\documentclass[12pt,a4paper]{scrbook}

\usepackage{amsfonts}
\usepackage{amsmath}
\usepackage{amssymb}
\usepackage{amsthm}

\usepackage{a4wide}
\usepackage{geometry}
\usepackage{graphicx}
\usepackage{longtable}
\usepackage{multirow}
\usepackage{xcolor}

\newtheorem{theorem}{Theorem}
\newtheorem{corollary}{Corollary}

\theoremstyle{definition}
\newtheorem{definition}{Definition}
\newtheorem{example}{Example}
\newtheorem{experiment}{Experimental Result}
\newtheorem{remark}{Remark}

\newtheorem{conjecture}{Conjecture}

\newtheorem{algorithm}{Algorithm}
\newtheorem{warning}{Warning}

\numberwithin{equation}{section}

\begin{document}

\subject{Four Plenary Lectures and Exercises}

\title{Theoretical and Experimental Approach to \\ p-Class Field Towers of Cyclic Cubic Number Fields}

\author{by Daniel C. Mayer}

\thanks{Research supported by the Austrian Science Fund (FWF): projects J0497-PHY, P26008-N25, and the Research Executive Agency of the European Union (EUREA)}

\begin{dedication}
\qquad \qquad \qquad Four Plenary Lectures and Exercises, Respectfully Dedicated to \\
\phantom{.} \qquad \qquad \qquad \qquad \qquad \textbf{Professor Mohammed Ayadi}.
\end{dedication}

\date{Sixi\`emes Journ\'ees d'Alg\`ebre, Th\'eorie des Nombres et leurs Applications (JATNA), Oujda, 23.---26. November 2022}

\maketitle

\frontmatter

\tableofcontents


\chapter{Preface}
\label{ch:Preface}

\noindent
This \textit{course in algebraic number theory} with title
\lq\lq Theoretical and Experimental Approach to \(p\)-Class Field Towers of Cyclic Cubic Number Fields\rq\rq\
is intended for an audience consisting of graduate students, doctorands, postdocs, and professional scientists.
It is divided into four plenary lectures
concerning the theory of \textit{cyclic number fields \(F/\mathbb{Q}\) of odd prime degree}
\(\ell:=\lbrack F:\mathbb{Q}\rbrack\)
and their unramified \(p\)-extension fields \(E/F\)
for various exemplary cases of \(p\)-class groups \(\mathrm{Cl}_p{F}=\mathrm{Syl}_p(\mathrm{Cl}(F))\)
with an assigned prime number \(p\in\mathbb{P}\).
The plenary lectures are separated by three experimental sessions with exercises in the laboratory,
where the theory is illuminated by means of the \textit{computational algebra system Magma}.
Experiments cover both, algebra and arithmetic as well as group theory.

Particular emphasis lies on the \textit{cubic situation} \(\ell=3\).
For a prime number \(p\),
denote by \(G=\mathrm{Gal}(F_p^{(\infty)}/F)\) the relative pro-\(p\) Galois group of
the unramified Hilbert \(p\)-class field tower \(F_p^{(\infty)}\)
of a cyclic cubic field \(F\) with \(p\)-class group \(\mathrm{Cl}_p{F}\).
For \(p\in\lbrace 2,5\rbrace\) and elementary bicyclic \(\mathrm{Cl}_p{F}\simeq C_p\times C_p\),
it is shown (Thm.
\ref{thm:SigmaGroupsDegree3})
that the \textit{action of the absolute group} \(\mathrm{Gal}(F/\mathbb{Q})\) on \(G\)
severely restricts the possibilities for the metabelianization \(\mathfrak{M}=G/G^{\prime\prime}\) of \(G\).
The length \(\ell_p{F}\) of the \(p\)-tower,
which coincides with the soluble length \(\mathrm{sl}(G)\) of \(G\),
is given by \(\ell_2{F}\in\lbrace 1,2\rbrace\) for \(p=2\) generally,
and by \(\ell_5{F}=2\) for \(p=5\) in all known examples.
Let \(t\) denote the number of prime divisors of the conductor \(c\) of \(F/\mathbb{Q}\),
viewed as a subfield of the \(3\)-\textit{ray class field modulo} \(c\) of \(\mathbb{Q}\).
For \(p=3\),
the statistical distribution of the two well-known (Thm.
\ref{thm:TwoPrimeCond})
possibilities \(\ell_3{F}\in\lbrace 1,2\rbrace\) is determined
in the case \(t=2\) and \(\mathrm{Cl}_3{F}\simeq C_3\times C_3\),
and a \textit{broad variety of new scenarios} (Thms.
\ref{thm:Cat3Gph679} --
\ref{thm:Cat1Gph12})
is presented in the case \(t=3\)
with elementary bicyclic or tricyclic \(\mathrm{Cl}_3{F}\),
in dependence on the mutual cubic residue conditions between the prime divisors of \(c\).
The coronation of this course is
the first rigorous proof of cyclic cubic fields \(F\) with
\(3\)-class field towers of exact length \(\ell_3{F}=3\),
established by means of Artin patterns and relation ranks (Thms.
\ref{thm:ImgQdr} --
\ref{thm:d19}).
For elementary tricyclic \(\mathrm{Cl}_3{F}\simeq C_3\times C_3\times C_3\),
finite \textit{three-stage towers} (Thm.
\ref{thm:333})
were completely unknown up to now,
for any kind of algebraic number fields.


\mainmatter

\part{Foundations}
\label{pt:Foundations}

\chapter{Construction of Cyclic Fields}
\label{ch:Construction}


\section{Multiplicity of Conductors and Discriminants}
\label{s:Multiplicity}

\noindent
For a fixed odd prime number \(\ell\ge 3\),
let \(F\) be a \textit{cyclic number field} of degree \(\ell\),
that is,
\(F/\mathbb{Q}\) is a Galois extension of degree
\(\lbrack F:\mathbb{Q}\rbrack=\ell\)
with absolute automorphism group 
\(\mathrm{Gal}(F/\mathbb{Q})=\langle\sigma\mid\sigma^\ell=1\rangle\).
According to the \textbf{Theorem of Kronecker, Weber and Hilbert}
on abelian extensions of the rational number field \(\mathbb{Q}\),
the \textit{conductor} \(c\) of \(F\)
is the smallest positive integer such that
\(F=:F_c\) is contained in the cyclotomic field
\(\mathbb{Q}(\zeta_c)\),
where \(\zeta_c\) denotes a primitive \(c\)-th root of unity,
more precisely, in the \(\ell\)-\textit{ray class field modulo} \(c\) of \(\mathbb{Q}\),
which lies in the maximal real subfield \(\mathbb{Q}(\zeta_c+\zeta_c^{-1})\).
It is well known that
\(c=\ell^e\cdot q_1\cdots q_\tau\),
where
\(e\in\lbrace 0,2\rbrace\)
and the \(q_i\) are pairwise distinct prime numbers
\(q_i\equiv +1\,(\mathrm{mod}\,\ell)\),
for \(1\le i\le\tau\).
The \textit{discriminant} of \(F=F_c\) is the perfect \((\ell-1)\)-th power
\(d_F=c^{\ell-1}\),
and the number of rational primes which are (totally) ramified in \(F\)
is given by

\begin{equation}
\label{eqn:Ramification}
t:=
\begin{cases}
\tau   & \text{ if } e=0\ (\ell \text{ is unramified in } F),\\
\tau+1 & \text{ if } e=2\ (\ell \text{ is ramified in } F).
\end{cases}
\end{equation}

\noindent
In the last case, we formally put \(q_{\tau+1}:=\ell^2\).
The number of non-isomorphic cyclic number fields
\(F_{c,1},\ldots,F_{c,m}\)
of degree \(\ell\), sharing the common conductor \(c\),
is given by the \textit{multiplicity} formula
\cite[p. 831]{Ma1992}

\begin{equation}
\label{eqn:Multiplicity}
m=m(c)=(\ell-1)^{t-1}.
\end{equation}


\section{Construction as Ray Class Fields}
\label{s:RayClassFields}

\noindent
For the construction of all cyclic number fields \(F=F_c\) of degree \(\ell\)
with ascending conductors \(b\le c\le B\)
between an assigned lower bound \(b\) and upper bound \(B\)
by means of the computational algebra system Magma
\cite{MAGMA2022},
the class field theoretic routines by Fieker
\cite{Fi2001}
can be used without the need of preparing a list of
suitable generating polynomials of \(\ell\)-th degree.
The big advantage of this technique is that
the cyclic number fields of degree \(\ell\) are produced as \textit{multiplets}
\((F_{c,1},\ldots,F_{c,m})\)
of pairwise non-isomorphic fields sharing the common conductor \(c\)
with \textit{multiplicities} \(m\in\lbrace 1,\ell-1,(\ell-1)^2,(\ell-1)^3,\ldots\rbrace\)
in dependence on the number \(t\in\lbrace 1,2,3,4,\ldots\rbrace\)
of primes dividing the conductor \(c\),
according to Formula
\eqref{eqn:Multiplicity}.


Throughout this course, we design algorithms
by the \textit{principle of successive refinement}.
We begin with a loop which sieves \(\ell\)-admissible conductors,
for a given odd prime number \(\ell\ge 3\).
If \(c=q_1^{n_1}\cdots q_t^{n_t}\) is the prime factorization of \(c\),
the function \texttt{Factorization()} returns a collection of pairs
\(((q_1,n_1),\ldots,(q_t,n_t))\).
In order to obtain \textbf{cyclic cubic fields}, we must put \(\ell=3\).

\begin{algorithm}
\label{alg:AdmissibleConductors}
(Filtering \(\ell\)-admissible conductors.) \\
\textbf{Input:}
prime \texttt{l}, lower bound \texttt{b}, upper bound \texttt{B}. \\
\textbf{Code:}
uses the subroutine \texttt{Process()}.
{\scriptsize
\texttt{
\begin{tabbing}
for \= for \= for \= for \= for \= \kill
fld := 0; // counter of fields, inclusively multiplicity\\
for c in [b..B] do\+\\
   cPD := Factorization(c);\\
   adm := true; // admissible\\
   for i in [1..\#cPD] do\+\\
      if not ( ((1 eq cPD[i][1] mod l) and (1 eq cPD[i][2])) or\+\+\\
               ((l eq cPD[i][1]) and (2 eq cPD[i][2])) ) then\-\\
         adm := false;\\
         break; // save CPU time\-\\
      end if;\-\\
   end for; // i\\
   if adm then\+\\
      Process(l,c,fld); // subroutine\-\\
   end if; // adm\-\\
end for; // c
\end{tabbing}
}
}
\noindent
\textbf{Output:}
managed by the subroutine \texttt{Process()}.
\end{algorithm}


\noindent
Now we implement the subroutine \texttt{Process()}.
Although the conductor \(c=1\) is not admissible
in the present context,
we suggest to avoid the ray class group modulo \(1\),
and to replace it by the ordinary class group,
if necessary.
Initially,
the wider multiplet of all abelian \(\ell\)-extensions
of the rational number field \(\mathbb{Q}\),
with conductors \textit{dividing} \(c\), is constructed and optimized.
In terms of the \(\ell\)-ray class rank modulo \(c\), \(\varrho=\varrho_\ell(c)\),
it consists of \(\frac{\ell^{\varrho}-1}{\ell-1}\) members.
This is the number of subgroups of index \(\ell\)
in the ray class group modulo \(c\).
By testing the correct discriminant \(d=c^{\ell-1}\),
the narrower multiplet \((F_{c,1},\ldots,F_{c,m})\)
of all cyclic extensions of degree \(\ell\)
with \textit{precise} conductor \(c\) is extracted,
and the \textit{abelian type invariants} of
the complete class group \(\mathrm{Cl}(F_{c,\mu})\)
and of the \(\ell\)-class group \(\mathrm{Cl}_\ell{F_{c,\mu}}\)
are determined for \(1\le\mu\le m\).

\newpage

\begin{algorithm}
\label{alg:CyclicMultiplets}
(Construction of multiplets.) \\
\textbf{Input:}
prime \texttt{l}, \texttt{l}-admissible conductor \texttt{c}, counter of fields \texttt{fld}. \\
\textbf{Code:}
{\scriptsize
\texttt{
\begin{tabbing}
for \= for \= for \= for \= for \= \kill
      SetClassGroupBounds("GRH");\\
      Q := RationalsAsNumberField();\\
      OQ := MaximalOrder(Q); // this is Z\\
      CQ,mCQ := ClassGroup(OQ);\\
      RQ,mRQ := RayClassGroup(c\(\ast\)OQ); // modulo c\(\ast\)Z\\
      if (1 eq c) then // avoid the ray mod 1\+\\
         XQ := CQ;\\
         mXQ := mCQ;\-\\
      else\+\\
         XQ := RQ;\\
         mXQ := mRQ;\-\\
      end if; // c=1\\
      sS1 := Subgroups(XQ: Quot := [l]); // full multiplet of l-ray class fields\\
      sA1 := [AbelianExtension(Inverse(mQQ)\(\ast\)mXQ)\+\+\\
              where QQ,mQQ := quo<XQ|x\({}\,\grave{}\,{}\)subgroup>: x in sS1];\-\-\\
      sN1 := [NumberField(x): x in sA1]; // relative fields\\
      sR1 := [MaximalOrder(x): x in sA1];\\
      sF1 := [AbsoluteField(x): x in sN1]; // absolute fields\\
      sM1 := [MaximalOrder(x): x in sF1];\\
      sB1 := [OptimizedRepresentation(x): x in sF1]; // first optimization\\
      sK1 := [NumberField(DefiningPolynomial(x)): x in sB1];\\
      sO1 := [Simplify(LLL(MaximalOrder(x))): x in sK1]; // second optimization\\
      // initialization of multiplets\\
      totMult := \#sO1; // inclusively wrong partial conductors\\
      locMult := 0; // only with correct conductor\\
      Collection := []; // all l-class groups of the multiplet\\
      for j in [1..totMult] do\+\\
         ON := sO1[j];\\
         N := NumberField(ON);\\
         G,Z,D:= GaloisGroup(N);\\
         if ((l eq \#G) and IsAbelian(G)) then // cyclic Galois group\+\\
            C,mC := ClassGroup(ON);\\
            CN := AbelianInvariants(C);\\
            SN := pPrimaryInvariants(C,l); // the Sylow l-subgroup\\
            Disc := c\({}\,\hat{}\,{}\)2; // root discriminant\\
            if (Disc\({}\,\hat{}\,{}\)((l-1) div 2) eq AbsoluteDiscriminant(ON)) then\+\\
               fld := fld + 1; // one field more\\
               locMult := locMult + 1; // local counter within multiplet\\
               printf "Nr \%4o",fld; // 4 digits\\
               printf ": c=\%6o",c; // 6 digits\\
               printf "=\%o",Factorization(c);\\
               printf ", j/m=\%o/\%o",j,totMult; // component nr\\
               printf ", SN=\%o, CN=\%o\(\backslash\)n",SN,CN;\\
               Append(\(\sim\)Collection,SN); // save l-class group\-\\
            end if; // correct discriminant\-\\
         end if; // cyclic\-\\
      end for; // j (end of multiplet)
\end{tabbing}
}
}
\noindent
\textbf{Output:}
counter \texttt{fld}, conductor \texttt{c}, factors, class group \texttt{CN}, \texttt{l}-class group \texttt{SN}.
\end{algorithm}


\noindent
When a \textit{statistic evaluation} of all multiplets is desired,
the loop in Algorithm
\ref{alg:AdmissibleConductors}
can be surrounded by the following initialization and finalization,
which is designed for the special case of cyclic cubic fields, \(\ell=3\).
The middle part must be inserted at the end of Algorithm
\ref{alg:CyclicMultiplets}.

\begin{algorithm}
\label{alg:StatisticMultiplets}
(Statistics of multiplets.) \\
\textbf{Input:}
local counter \texttt{locMult} and \texttt{Collection}. \\
\textbf{Code:}
{\scriptsize
\texttt{
\begin{tabbing}
for \= for \= for \= for \= for \= \kill
// distribution of multiplicities from 1 to 20\\
Mult := [0,0,0,0,0,0,0,0,0,0,0,0,0,0,0,0,0,0,0,0];\\
// collections of cyclic cubic multiplets\\
Singlets := []; // m=1\\
Doublets := []; // m=2\\
Quartets := []; // m=4\\
Octets := []; // m=8\\
Hexadecuplets := []; // m=16\\
\phantom{for}\\
         // after (end of multiplet)\\
         // finalization of multiplets\\
         if (20 ge locMult) then\+\\
            Mult[locMult] := Mult[locMult] + 1;\-\\
         end if;\\
         if (1 le locMult) then\+\\
            if (1 eq locMult) then\+\\
               Append(\(\sim\)Singlets,Collection);\\
               printf "m=1, no \%o\(\backslash\)n",Mult[1];\-\\
            elif (2 eq locMult) then\+\\
               Append(\(\sim\)Doublets,Collection);\\
               printf "m=2, no \%o\(\backslash\)n",Mult[2];\-\\
            elif (4 eq locMult) then\+\\
               Append(\(\sim\)Quartets,Collection);\\
               printf "m=4, no \%o\(\backslash\)n",Mult[4];\-\\
            elif (8 eq locMult) then\+\\
               Append(\(\sim\)Octets,Collection);\\
               printf "m=8, no \%o\(\backslash\)n",Mult[8];\-\\
            elif (16 eq locMult) then\+\\
               Append(\(\sim\)Hexadecuplets,Collection);\\
               printf "m=16, no \%o\(\backslash\)n",Mult[16];\-\\
            end if;\-\\
         end if; // locMult\\
\phantom{for}\\
printf "multiplicity:\(\backslash\)n";\\
for i in [1..20] do\+\\
   printf "\%2o:\%o\(\backslash\)n",i,Mult[i];\-\\
end for; // i\\
printf "\(\backslash\)nSinglets:\(\backslash\)n";\\
for i in [1..\#Singlets] do\+\\
   printf "\%o: ( ",i;\\
   for j in [1..\#Singlets[i]] do\+\\
      printf "\%o ",Singlets[i][j];\-\\
   end for;\\
   printf ")\(\backslash\)n";\-\\
end for;\\
printf "\(\backslash\)nDoublets:\(\backslash\)n";\\
for i in [1..\#Doublets] do\+\\
   printf "\%o: ( ",i;\\
   for j in [1..\#Doublets[i]] do\+\\
      printf "\%o ",Doublets[i][j];\-\\
   end for;\\
   printf ")\(\backslash\)n";\-\\
end for;\\
// and so on until\\
printf "\(\backslash\)nHexadecuplets:\(\backslash\)n";\\
for i in [1..\#Hexadecuplets] do\+\\
   printf "\%o: ( ",i;\\
   for j in [1..\#Hexadecuplets[i]] do\+\\
      printf "\%o ",Hexadecuplets[i][j];\-\\
   end for;\\
   printf ")\(\backslash\)n";\-\\
end for;\\
printf "\(\backslash\)n\(\backslash\)n\(\backslash\)n";
\end{tabbing}
}
}
\noindent
\textbf{Output:}
singlets, doublets, quartets, octets, hexadecuplets, with \texttt{l}-class groups.
\end{algorithm}


\begin{warning}
\label{wrn:RayClassFields}
It should be pointed out that
the totally complex cyclotomic field
\(\mathbb{Q}(\zeta_c)\) with \(c\ge 3\)
is the \textit{ray class field modulo} \(c\cdot w\) of \(\mathbb{Q}\),
where \(\zeta_c\) denotes a primitive \(c\)-th root of unity,
e.g. \(\zeta_c=\exp(2\pi\sqrt{-1}/c)\),
and \(w\) is the unique real archimedean place of \(\mathbb{Q}\),
which splits into \(\frac{\phi(c)}{2}\) pairs of conjugate complex archimedean places in the cyclotomic field.
However,
the multiplet of cyclic fields \(F_{c,\mu}\), \(1\le\mu\le m\), of degree \(\ell\)
with \(\ell\)-admissible conductor \(c\) is contained in 
the \(\ell\)-\textit{ray class field modulo} \(c\) of \(\mathbb{Q}\),
which is contained in the maximal real subfield \(\mathbb{Q}(\zeta_c+\zeta_c^{-1})\)
of \(\mathbb{Q}(\zeta_c)\).
\end{warning}


\section{Census of Multiplets and Class Groups}
\label{s:ClassGroups}

\noindent
With the aid of the computational algebra system Magma
\cite{BCP1997,BCFS2022,MAGMA2022},
we have implemented the Algorithms
\ref{alg:AdmissibleConductors},
\ref{alg:CyclicMultiplets},
\ref{alg:StatisticMultiplets}
and constructed the multiplets \((F_{c,1},\ldots,F_{c,m})\)
of all cyclic \textit{quintic}, respectively cyclic \textit{cubic}, number fields \(F\)
with conductors \(0<c<100\,000\)
as subfields of ray class fields
\cite{Fi2001}
modulo \(5\)-admissible, respectively \(3\)-admissible, conductors \(c\)
over the rational number field \(\mathbb{Q}\).


\subsection{Statistics of Multiplets}
\label{ss:StatMult}


\noindent
The complete census of multiplets for cyclic \textit{quintic} fields is given in Table
\ref{tbl:StatMult5}.
Obviously, the fields \(F\) in \textit{quartets} are dominating, but not the corresponding conductors \(c\).
We always emphasize the \textit{minimal conductor} \(c_{\text{min}}\) as a prototype.

\renewcommand{\arraystretch}{1.1}

\begin{table}[ht]
\caption{Statistics of Cyclic Quintic Multiplets}
\label{tbl:StatMult5}
\begin{center}
\begin{tabular}{|l|r||r||r||r|}
\hline
 Multiplets    &  \(m\) & \(\#\) Conductors \(c\) & \(\#\) Fields \(F\) & \(c_{\text{min}}\) \\
\hline
 \(c<\)        &   \(\) & \(100000\) & \(100000\) & \(\) \\
\hline
 Singlets      &  \(1\) & \(2388\) & \(2388\) &   \(11\) \\
 Quartets      &  \(4\) &  \(845\) & \(3380\) &  \(275\) \\
 Hexadecuplets & \(16\) &   \(49\) &  \(784\) & \(8525\) \\
\hline
 Total         &        & \(3282\) & \(6552\) &          \\
\hline
\end{tabular}
\end{center}
\end{table}


\noindent
For cyclic \textit{cubic} fields, we split the information into four ranges of conductors,
with lengths increasing in steps of \(25000\).
The complete census of multiplets is given in Table
\ref{tbl:StatMult3}.
Here, fields \(F\) in \textit{doublets} are dominating, but not conductors \(c\).

\renewcommand{\arraystretch}{1.1}

\begin{table}[ht]
\caption{Statistics of Cyclic Cubic Multiplets}
\label{tbl:StatMult3}
\begin{center}

{\small

\begin{tabular}{|l|r||r|r|r|r||r|r|r|r||r|}
\hline
 Multiplets & \(m\) & \multicolumn{4}{|c||}{\(\#\) Conductors \(c\)} & \multicolumn{4}{|c||}{\(\#\) Fields \(F\)} & \(c_{\text{min}}\) \\
\hline
 \(c<\)     &  \(\) &\(25000\) &\(50000\) &\(75000\) & \(10^5\) &\(25000\) &\(50000\) & \(75000\) &  \(10^5\) & \(\) \\
\hline
 Singlets   & \(1\) & \(1372\) & \(2557\) & \(3682\) & \(4785\) & \(1372\) & \(2557\) &  \(3682\) &  \(4785\) &     \(7\) \\
 Doublets   & \(2\) &  \(993\) & \(1953\) & \(2921\) & \(3863\) & \(1986\) & \(3906\) &  \(5842\) &  \(7726\) &    \(63\) \\
 Quartets   & \(4\) &  \(149\) &  \(351\) &  \(556\) &  \(783\) &  \(596\) & \(1404\) &  \(2224\) &  \(3132\) &   \(819\) \\
 Octets     & \(8\) &    \(1\) &    \(7\) &   \(17\) &   \(26\) &    \(8\) &   \(56\) &   \(136\) &   \(208\) & \(15561\) \\
\hline
 Total      &       & \(2515\) & \(4868\) & \(7176\) & \(9457\) & \(3962\) & \(7923\) & \(11884\) & \(15851\) &           \\
\hline
\end{tabular}

}

\end{center}
\end{table}


\subsection{Statistics of Class Group Structures}
\label{ss:StatStruc}


\renewcommand{\arraystretch}{1.1}

\begin{table}[hb]
\caption{Statistics of Cyclic Quintic \(5\)-Class Group Structures}
\label{tbl:StatStruc5}
\begin{center}
\begin{tabular}{|r||r||r|}
\hline
 \(\mathrm{Cl}_5{F}\) & \(\#\) Fields \(F\) & \(c_{\text{min}}\) \\
\hline
 \(c<\)      & \(100000\) & \(\) \\
\hline
 \(1\)       & \(2388\) &    \(11\) \\
 \((5)\)     & \(3260\) &   \(275\) \\
 \((5,5)\)   &  \(869\) &  \(2651\) \\
 \((5,5,5)\) &   \(35\) & \(13981\) \\
\hline
 Total       & \(6552\) &           \\
\hline
\end{tabular}
\end{center}
\end{table}

\noindent
Statistics of the structures of
\(5\)-class groups \(\mathrm{Cl}_5{F}\), respectively
\(3\)-class groups \(\mathrm{Cl}_3{F}\),
is shown in Table
\ref{tbl:StatStruc5}
for cyclic \textit{cubic} fields, respectively
\ref{tbl:StatStruc3}
for cyclic \textit{quintic} fields.
The \(3\)-class groups reveal a wealth of various structures
up to order \(3^5\).


\renewcommand{\arraystretch}{1.1}

\begin{table}[ht]
\caption{Statistics of Cyclic Cubic \(3\)-Class Group Structures}
\label{tbl:StatStruc3}
\begin{center}

{\small

\begin{tabular}{|r||r||r|}
\hline
 \(\mathrm{Cl}_3{F}\) & \(\#\) Fields \(F\) & \(c_{\text{min}}\) \\
\hline
 \(c<\)        & \(100000\) & \(\) \\
\hline
 \(1\)         &   \(4785\) &     \(7\) \\
 \((3)\)       &   \(6910\) &    \(63\) \\
 \((3,3)\)     &   \(3498\) &   \(657\) \\
 \((3,3,3)\)   &    \(481\) &  \(3913\) \\
 \((3,3,3,3)\) &     \(13\) & \(25389\) \\
\hline
 \((3,9)\)     &    \(105\) &  \(4711\) \\
 \((3,3,9)\)   &     \(43\) &  \(7657\) \\
 \((3,3,3,9)\) &      \(6\) & \(15561\) \\
\hline
 \((9,9)\)     &      \(5\) & \(41977\) \\
 \((3,9,9)\)   &      \(3\) & \(66157\) \\
\hline
 \((9,27)\)    &      \(2\) & \(36667\) \\
\hline
 Total         &  \(15851\) &           \\
\hline
\end{tabular}

}

\end{center}
\end{table}

\newpage

\subsection{Class Group Structures Grouped by Multiplets}
\label{ss:StatMultStruc}

\noindent
The complete census of the structures of
\(5\)-class groups \(\mathrm{Cl}_5{F}\), respectively
\(3\)-class groups \(\mathrm{Cl}_3{F}\),
\textit{grouped by multiplets},
is shown in Table
\ref{tbl:StatMultStruc5},
respectively
\ref{tbl:StatMultStruc3}.
This point of view admits
considerably deeper insight into the connections
between the number \(t\) of primes dividing the conductor \(c\)
and the structure of \(3\)-class groups \(\mathrm{Cl}_3{F_{c,\mu}}\)
of all components of the multiplet
\((F_{c,\mu})_{1\le\mu\le m}\), \(m=2^{t-1}\),
of cyclic cubic fields sharing the common conductor \(c\).


\renewcommand{\arraystretch}{1.1}

\begin{table}[ht]
\caption{Statistics of \(5\)-Class Group Structures of Cyclic Quintic Multiplets}
\label{tbl:StatMultStruc5}
\begin{center}
\begin{tabular}{|r|r||r||r||r|}
\hline
 \((\mathrm{Cl}_5{F_{c,\mu}})_{\mu=1}^m\) & \(m\) & \(\#\) Conductors \(c\) & \(\#\) Fields \(F\) & \(c_{\text{min}}\) \\
\hline
 \(c<\)                   &   \(\) & \(100000\) & \(100000\) & \(\) \\
\hline
 \(1\)                    &  \(1\) & \(2388\) &  \(2388\) &    \(11\) \\
\hline
 \((5)^4\)                &  \(4\) &  \(815\) &  \(3260\) &   \(275\) \\
 \((5,5)^4\)              &  \(4\) &   \(28\) &   \(112\) &  \(2651\) \\
 \((5,5,5)^4\)            &  \(4\) &    \(2\) &     \(8\) & \(23411\) \\
\hline
 \((5,5)^{16}\)           & \(16\) &   \(41\) &   \(656\) &  \(8525\) \\
 \((5,5,5)^3,(5,5)^{13}\) & \(16\) &    \(5\) & \(15+65\) & \(13981\) \\
 \((5,5,5)^4,(5,5)^{12}\) & \(16\) &    \(3\) & \(12+36\) & \(47275\) \\
\hline
 Total                    &        & \(3282\) &  \(6552\) &           \\
\hline
\end{tabular}
\end{center}
\end{table}

\newpage

\begin{experiment}
\label{exp:Plan}
Table
\ref{tbl:StatMultStruc3}
will be the leading principle
for the entire layout of this course.
\begin{itemize}
\item
\(3\)-class groups \(\mathrm{Cl}_3{F_{c,\mu}}\) of rank \(4\)
appear in \(14\) among the \(26\) \textbf{octets} only.
One, two or three components of an octet may have
\(3\)-class groups of type \((3,3,3,3)\) or \((3,3,3,9)\).
Octets will be reserved for future investigations.
\item
Since they are trivial, we shall not deal with
\(4785\) \textit{singlets} having the trivial \(3\)-class group \(1\),
and with \(3455\) \textit{doublets} having cyclic \(3\)-class group \((3)\).
\item
In this course,
we pay attention to \(352\), respectively \(56\), \textbf{doublets}
with \(3\)-class number \(h_3=9\) and elementary bicyclic \(3\)-class group \((3,3)\),
respectively \(h_3\ge 27\) and non-elementary bicyclic \(3\)-class groups
\((3,9)\), \((9,9)\), \((9,27)\),
in Section \S\
\ref{ch:Classical},
respectively Section \S\
\ref{ch:Singular},
and to all \(783\) \textbf{quartets}
in the Sections \S\S\
\ref{ss:Cat3Gph1To4},
\ref{ss:Cat3Gph5To9},
\ref{ss:Cat2Gph1},
\ref{ss:Cat2Gph2},
\ref{ss:Cat1Gph1},
\ref{ss:Cat1Gph2}.
\end{itemize}
\end{experiment}


\renewcommand{\arraystretch}{1.1}

\begin{table}[ht]
\caption{Statistics of \(3\)-Class Group Structures of Cyclic Cubic Multiplets}
\label{tbl:StatMultStruc3}
\begin{center}
\begin{tabular}{|r|r||r||r||r|}
\hline
 \((\mathrm{Cl}_3{F_{c,\mu}})_{\mu=1}^m\) & \(m\) & \(\#\) Conductors \(c\) & \(\#\) Fields \(F\) & \(c_{\text{min}}\) \\
\hline
 \(c<\)                              &   \(\) & \(100000\) & \(100000\) & \(\) \\
\hline
 \(1\)                               &  \(1\) & \(4785\) &    \(4785\) &     \(7\) \\
\hline
 \((3)^2\)                           &  \(2\) & \(3455\) &    \(6910\) &    \(63\) \\
 \((3,3)^2\)                         &  \(2\) &  \(352\) &     \(704\) &   \(657\) \\
 \((3,9)^2\)                         &  \(2\) &   \(50\) &     \(100\) &  \(4711\) \\
 \((9,9),(3,9)\)                     &  \(2\) &    \(4\) &     \(4+4\) & \(41977\) \\
 \((9,27),(3,9)\)                    &  \(2\) &    \(1\) &     \(1+1\) & \(36667\) \\
 \((9,27),(9,9)\)                    &  \(2\) &    \(1\) &     \(1+1\) & \(42127\) \\
\hline
 \((3,3)^4\)                         &  \(4\) &  \(579\) &    \(2316\) &   \(819\) \\
 \((3,3,3),(3,3)^3\)                 &  \(4\) &   \(79\) &  \(79+237\) &  \(4977\) \\
 \((3,3,3)^2,(3,3)^2\)               &  \(4\) &   \(80\) & \(160+160\) &  \(3913\) \\
 \((3,3,9),(3,3)^3\)                 &  \(4\) &   \(18\) &   \(18+54\) &  \(7657\) \\
 \((3,3,9)^2,(3,3)^2\)               &  \(4\) &   \(10\) &   \(20+20\) & \(27873\) \\
 \((3,9,9),(3,3)^3\)                 &  \(4\) &    \(1\) &     \(1+3\) & \(67347\) \\
 \((3,9,9),(3,3,9),(3,3)^2\)         &  \(4\) &    \(2\) &   \(2+2+4\) & \(66157\) \\
 \((3,3,3)^4\)                       &  \(4\) &   \(12\) &      \(48\) & \(38311\) \\
 \((3,3,9),(3,3,3)^3\)               &  \(4\) &    \(1\) &     \(1+3\) & \(91819\) \\
 \((3,3,9)^2,(3,3,3)^2\)             &  \(4\) &    \(1\) &     \(2+2\) & \(97747\) \\
\hline
 \((3,3,3)^8\)                       &  \(8\) &   \(12\) &      \(96\) & \(30303\) \\
 \((3,3,3,3),(3,3,3)^7\)             &  \(8\) &    \(9\) &    \(9+63\) & \(25389\) \\
 \((3,3,3,9),(3,3,3)^7\)             &  \(8\) &    \(2\) &    \(2+14\) & \(15561\) \\
 \((3,3,3,9)^2,(3,3,3)^6\)           &  \(8\) &    \(1\) &     \(2+6\) & \(49959\) \\
 \((3,3,3,9),(3,3,3,3)^2,(3,3,3)^5\) &  \(8\) &    \(2\) &  \(2+4+10\) & \(54873\) \\
\hline
 Total                               &        & \(9457\) &   \(15851\) &           \\
\hline
\end{tabular}
\end{center}
\end{table}


\chapter{Arithmetic of Cyclic Cubic Fields}
\label{ch:Arithmetic}

\section{Rank of 3-Class Groups of Cyclic Cubic Fields}
\label{s:3ClassRank}

\noindent
Since the rank \(\varrho_3{F}\) of the \(3\)-class group \(\mathrm{Cl}_3{F}\) of a cyclic cubic field \(F\)
depends on the mutual cubic residue conditions between the prime divisors \(q_1,\ldots,q_t\) of the conductor \(c\),
Gras
\cite[pp. 21--22]{Gr1973}
has introduced directed graphs with \(t\) vertices \(q_1,\ldots,q_t\)
whose directed edges \(q_i\to q_j\) describe values of cubic residue symbols.
We use a simplified notation of these graphs, fitting in a single line,
but occasionally requiring the repetition of a vertex.


\begin{definition}
\label{dfn:CubicResidues}
Let \(\zeta_3\) be a fixed primitive third root of unity.
For each pair \((q_i,q_j)\) with \(1\le i\ne j\le t\),
the value of the \textit{cubic residue symbol}
\(\left(\frac{q_i}{q_j}\right)_3=\zeta_3^{a_{i,j}}\)
is determined uniquely by the integer \(a_{i,j}\in\lbrace 0,1,2\rbrace\).
Let a \textit{directed edge} \(q_i\to q_j\) be defined if and only if \(\left(\frac{q_i}{q_j}\right)_3=1\),
that is, \(q_i\) is a cubic residue modulo \(q_j\) (and thus \(a_{i,j}=0\)).
The \textbf{combined cubic residue symbol} \(\lbrack q_1,\ldots,q_t\rbrack_3:=\)

\begin{equation}
\label{eqn:CubicResidues}
\Biggl\lbrace q_i\to q_j\Biggm\vert i\ne j,\left(\frac{q_i}{q_j}\right)_3=1\Biggr\rbrace
\bigcup
\Biggl\lbrace q_i\Biggm\vert (\forall j\ne i)\,\left(\frac{q_i}{q_j}\right)_3\ne 1,\left(\frac{q_j}{q_i}\right)_3\ne 1 \Biggr\rbrace
\end{equation}

\noindent
where the subscripts \(i\) and \(j\) run from \(1\) to \(t\),
is defined as the union of the set of all directed edges
which occur in the graph associated with \(q_1,\ldots,q_t\) in the sense of Gras,
and the set of all isolated vertices.
For \(t=3\),
we additionally need the invariant \(\delta:=a_{1,2}a_{2,3}a_{3,1}-a_{1,3}a_{3,2}a_{2,1}\)
in order to distinguish two subcases of the case with three isolated vertices.
\end{definition}


\begin{theorem}
\label{thm:3ClassRank}
(Rank Distribution, \textbf{G. Gras},
\cite{Gr1973}.) \\
Let \(F\) be a cyclic cubic field of conductor \(c=q_1\cdots q_t\) with \(1\le t\le 3\).
We indicate \textbf{mutual cubic residues} simply by writing \(q_1\leftrightarrow q_2\) instead of \(q_1\to q_2\to q_1\).

\begin{itemize}

\item
If \(t=1\), then \(m=1\), \(F\) forms a singlet, \(\lbrack q_1\rbrack_3=\lbrace q_1\rbrace\), and \(\varrho_3{F}=0\).

\item
If \(t=2\), then \(m=2\), \(F\) is member of a doublet \((F_{1},F_{2})\), and there arise two possibilities.

\begin{enumerate}
\item
\((\varrho_3{F_{1}},\varrho_3{F_{2}})=(1,1)\), if
\begin{equation}
\label{eqn:Rank1}
\lbrack q_1,q_2\rbrack_3=
\begin{cases}
\lbrace q_1,q_2\rbrace & \text{ or} \\
\lbrace q_i\to q_j\rbrace & \text{ with } i\ne j.
\end{cases}
\end{equation}
\item
\((\varrho_3{F_{1}},\varrho_3{F_{2}})=(2,2)\), if
\begin{equation}
\label{eqn:Rank2}
\lbrack q_1,q_2\rbrack_3=\lbrace q_1\leftrightarrow q_2\rbrace.
\end{equation}
\end{enumerate}

\item
If \(t=3\), then \(m=4\), \(F\) is member of a quartet \((F_{1},\ldots,F_{4})\), and there arise five cases.

\begin{enumerate}
\item
\((\varrho_3{F_{1}},\varrho_3{F_{2}},\varrho_3{F_{3}},\varrho_3{F_{4}})=(2,2,2,2)\), if
\begin{equation}
\label{eqn:Category3}
\lbrack q_1,q_2,q_3\rbrack_3=
\begin{cases}
\lbrace q_1,q_2,q_3;\delta\ne 0\rbrace & \text{ or} \\
\lbrace q_i\rightarrow q_j;q_k\rbrace & \text{ or} \\
\lbrace q_i\rightarrow q_j\rightarrow q_k\rbrace & \text{ or} \\
\lbrace q_i\rightarrow q_j\rightarrow q_k\rightarrow q_i\rbrace & \text{ or} \\
\lbrace q_i\leftrightarrow q_j;q_k\rbrace & \text{ or} \\
\lbrace q_i\leftrightarrow q_j\rightarrow q_k\rbrace & \text{ or} \\
\lbrace q_i\leftrightarrow q_j\leftarrow q_k\rbrace & \text{ or} \\
\lbrace q_k\rightarrow q_i\leftrightarrow q_j\leftarrow q_k\rbrace & \text{ or} \\
\lbrace q_k\rightarrow q_i\leftrightarrow q_j\rightarrow q_k\rbrace & \\
\end{cases}
\end{equation}
with \(i,j,k\) pairwise distinct.
\item
\((\varrho_3{F_{1}},\varrho_3{F_{2}},\varrho_3{F_{3}},\varrho_3{F_{4}})=(3,2,2,2)\), if
\begin{equation}
\label{eqn:Category1}
\lbrack q_1,q_2,q_3\rbrack_3=
\begin{cases}
\lbrace q_1,q_2,q_3;\delta=0\rbrace & \text{ or} \\
\lbrace q_i\leftarrow q_j\rightarrow q_k\rbrace & \text{ with } i,j,k \text{ pairwise distinct}.
\end{cases}
\end{equation}
\item
\((\varrho_3{F_{1}},\varrho_3{F_{2}},\varrho_3{F_{3}},\varrho_3{F_{4}})=(3,3,2,2)\), if
\begin{equation}
\label{eqn:Category2}
\lbrack q_1,q_2,q_3\rbrack_3=
\begin{cases}
\lbrace q_i\rightarrow q_j\leftarrow q_k\rbrace & \text{ or} \\
\lbrace q_i\rightarrow q_j\leftarrow q_k\rightarrow q_i\rbrace &
\end{cases}
\end{equation}
with \(i,j,k\) pairwise distinct.
\item
\((\varrho_3{F_{1}},\varrho_3{F_{2}},\varrho_3{F_{3}},\varrho_3{F_{4}})=(3,3,3,3)\), if
\begin{equation}
\label{eqn:Rank3}
\lbrack q_1,q_2,q_3\rbrack_3=
\begin{cases}
\lbrace q_i\leftarrow q_j\leftrightarrow q_k\rightarrow q_i\rbrace & \text{ or} \\
\lbrace q_i\leftrightarrow q_j\leftrightarrow q_k\rbrace & \text{ or} \\
\lbrace q_i\leftrightarrow q_j\leftrightarrow q_k\rightarrow q_i\rbrace &
\end{cases}
\end{equation}
with \(i,j,k\) pairwise distinct.
\item
\((\varrho_3{F_{1}},\varrho_3{F_{2}},\varrho_3{F_{3}},\varrho_3{F_{4}})=(4,4,4,4)\), if
\begin{equation}
\label{eqn:Rank4}
 \lbrack q_1,q_2,q_3\rbrack_3=\lbrace q_1\leftrightarrow q_2\leftrightarrow q_3\leftrightarrow q_1\rbrace.
\end{equation}
\end{enumerate}

\end{itemize}

\end{theorem}

\begin{proof}
See
\cite[Prp. VI.5, pp. 21--22]{Gr1973}.
\end{proof}


\begin{remark}
\label{rmk:3ClassRank}
Ayadi
\cite[pp. 45--47]{Ay1995}
investigated the cases
\(t=2\), Formula
\eqref{eqn:Rank2},
and
\(t=3\), Formulas
\eqref{eqn:Category3},
\eqref{eqn:Category1},
\eqref{eqn:Category2},
in Theorem
\ref{thm:3ClassRank}.
For \(t=3\), he denoted
the nine subcases of Formula
\eqref{eqn:Category3}
by Graph 1,2,3,4,5,6,7,8,9 of Category III,
the two subcases of Formula
\eqref{eqn:Category1}
by Graph 1 and 2 of Category I,
and the two subcases of Formula
\eqref{eqn:Category2}
by Graph 1 and 2 of Category II.
For the Categories I and II,
Ayadi did not consider the fields with \(3\)-class rank \(\varrho_3{F_{\mu}}=3\).
\end{remark}


\section{Power Residues, Categories and Graphs}
\label{s:PowResCatGph}

\noindent
The following variant of
the Euler criterion for \(\ell\)-th power residues,
with an odd prime \(\ell\),
is the foundation of all further algorithms in this section.

\begin{algorithm}
\label{alg:PowerResidues}
(Euler criterion for \(\ell\)-th power residues.) \\
\textbf{Input:}
prime \texttt{l}, prime(-power) module \texttt{m}, residue \texttt{r}. \\
\textbf{Code:}
{\tiny
\texttt{
\begin{tabbing}
for \= for \= for \= for \= for \= \kill
chr := 0; // residue character\\
g := PrimitiveRoot(m);\\
if (0 eq g) then\+\\
   printf "no primitive root g modulo \%o\(\backslash\)n",m;\-\\
else\+\\
   i := 0;\\
   bool := false;\\
   for e in [1..m-1] do\+\\
      pot := g\({}\,\hat{}\,{}\)e;\\
      dif := pot - r;\\
      if (0 eq dif mod m) then\+\\
         bool := true;\\
         i := e; // storage of exponent e\\
         break;\-\\
      end if;\-\\
   end for; // e\\
   // EULER criterion\\
   if (true eq bool) and (0 lt Valuation(i,l)) then\+\\
      chr := +1; // residue\-\\
   elif not (0 eq r mod m) then\+\\
      chr := -1; // (coarse) non-residue\-\\
   end if;\\
   printf "(\%o/\%o)\_\%o=\%o, g=\%o, i=\%o\(\backslash\)n",r,m,l,chr,g,i;\-\\
end if;
\end{tabbing}
}
}
\noindent
\textbf{Output:}
\texttt{l}-th residue character \texttt{chr} of \texttt{r} mod \texttt{m}, primitive root \texttt{g}, exponent \texttt{i}.
\end{algorithm}


\noindent
The preceding Algorithm
\ref{alg:PowerResidues}
can be implemented as a function \texttt{PowerResidue()},
which returns the \(\ell\)-th power residue character, 
rather than printing it as output.
Now we design a procedure for \(\ell=3\)
which determines the category and graph
of any \(3\)-admissible conductor \(c\) with \(t=3\) prime divisors.
The graph consists of three prime divisors, \(q_1,q_2,q_3\), as vertices,
and a directed edge \(q_i\to q_k\) whenever
\(i\ne k\) and \((q_i/q_k)_3=1\),
for instance \(13\to 7\), since \(3^3-13=14=2\cdot 7\)
and thus \(13\) is cubic residue mod \(7\), i.e. \((13/7)_3=1\).
The procedure determines the number of bidirectional edges
(\textit{mutual cubic residues}), all directed edges,
\textit{universally attractive} vertices, and \textit{universally repulsive} vertices,
which is sufficient in order to determine the category and graph,
according to the Gras Theorem
\ref{thm:3ClassRank},
except for one special configuration without any edges:
since we do not classify non-residues by the character values \(\zeta_3\) and \(\zeta_3^2\),
we cannot determine the invariant \(\delta\).
Thus the decision if Category III, graph 1,
has to be replaced by Category I, graph 1,
will be supplemented later by means of the rank distribution.


\begin{algorithm}
\label{alg:CategoryGraph}
(Automatic determination of category and graph.) \\
\textbf{Input:}
critical prime \texttt{iCrit}, three prime divisors \texttt{iPrm1}, \texttt{iPrm2}, \texttt{iPrm3} of the conductor \texttt{c}. \\
\textbf{Code:}
uses the function \texttt{PowerResidue()}.
{\scriptsize
\texttt{
\begin{tabbing}
for \= for \= for \= for \= for \= \kill
P := iCrit; // critical prime for p-th power residues\\
Q := iPrm1; // first prime divisor q\\
R := iPrm2; // second prime divisor r\\
S := iPrm3; // third prime divisor s\\
// six tests (in three pairs)\\
Char := [];\\
C12 := PowerResidue(P,Q,R); // Q -> R\\
C21 := PowerResidue(P,R,Q); // R -> Q\\
C13 := PowerResidue(P,Q,S); // Q -> S\\
C31 := PowerResidue(P,S,Q); // S -> Q\\
C23 := PowerResidue(P,R,S); // R -> S\\
C32 := PowerResidue(P,S,R); // S -> R\\
Append(\(\sim\)Char,C12);\\
Append(\(\sim\)Char,C21);\\
Append(\(\sim\)Char,C13);\\
Append(\(\sim\)Char,C31);\\
Append(\(\sim\)Char,C23);\\
Append(\(\sim\)Char,C32);\\
// number of bidirectional edges\\
AnzBid := 0;\\
if (1 eq C12) and (1 eq C21) then\+\\
   AnzBid := AnzBid + 1;\-\\
end if;\\
if (1 eq C13) and (1 eq C31) then\+\\
   AnzBid := AnzBid + 1;\-\\
end if;\\
if (1 eq C23) and (1 eq C32) then\+\\
   AnzBid := AnzBid + 1;\-\\
end if;\\
// three bidirectional edges\\
if (3 eq AnzBid) then\+\\
   kat := 5;\\
   grp := 1;\-\\
else\+\\
   // number of edges\\
   AnzKnt := 0;\\
   for j in [1..\#Char] do\+\\
      if (1 eq Char[j]) then\+\\
         AnzKnt := AnzKnt + 1;\-\\
      end if;\-\\
   end for; // j\\
   // number of attractive edges\\
   Att := 0;\\
   if (1 eq C21) and (1 eq C31) then\+\\
      Att := Att + 1;\-\\
   end if;
   if (1 eq C12) and (1 eq C32) then\+\\
      Att := Att + 1;\-\\
   end if;
   if (1 eq C13) and (1 eq C23) then\+\\
      Att := Att + 1;\-\\
   end if;\\
   // number of repulsive edges\\
   Rep := 0;\\
   if (1 eq C12) and (1 eq C13) then\+\\
      Rep := Rep + 1;\-\\
   end if;
   if (1 eq C21) and (1 eq C23) then\+\\
      Rep := Rep + 1;\-\\
   end if;
   if (1 eq C31) and (1 eq C32) then\+\\
      Rep := Rep + 1;\-\\
   end if;\\
   if (2 eq AnzBid) then\+\\
      if (4 eq AnzKnt) then\+\\
         kat := 4;\\
         grp := 2;\-\\
      else // then AnzKnt=5\+\\
         kat := 4;\\
         grp := 3;\-\\
      end if;\-\\
   elif (1 eq AnzBid) then\+\\
      if (2 eq AnzKnt) then\+\\
         kat := 3;\\
         grp := 5;\-\\
      elif (3 eq AnzKnt) then\+\\
         if (1 eq Rep) then\+\\
            kat := 3;\\
            grp := 6;\-\\
         elif (1 eq Att) then\+\\
            kat := 3;\\
            grp := 7;\-\\
         end if;\-\\
      elif (4 eq AnzKnt) then\+\\
         if (1 eq Att) and (1 eq Rep) then\+\\
            kat := 3;\\
            grp := 9;\-\\
         elif (1 eq Att) and (2 eq Rep) then\+\\
            kat := 4;\\
            grp := 1;\-\\
         elif (2 eq Att) and (1 eq Rep) then\+\\
            kat := 3;\\
            grp := 8;\-\\
         end if;\-\\
      end if;\-\\
   elif (0 eq AnzBid) then\+\\
      if (0 eq AnzKnt) then\+\\
         kat := 3; // or 1 (final decision later)\\
         grp := 1;\-\\
      elif (1 eq AnzKnt) then\+\\
         kat := 3;\\
         grp := 2;\-\\
      elif (2 eq AnzKnt) then\+\\
         if (1 eq Rep) then\+\\
            kat := 1;\\
            grp := 2;\-\\
         elif (1 eq Att) then\+\\
            kat := 2;\\
            grp := 1;\-\\
         else\+\\
            kat := 3;\\
            grp := 3;\-\\
         end if;\-\\
      elif (3 eq AnzKnt) then\+\\
         if (1 eq Att) and (1 eq Rep) then\+\\
            kat := 2;\\
            grp := 2;\-\\
         else\+\\
            kat := 3;\\
            grp := 4;\-\\
         end if;\-\\
      end if;\-\\  
   end if;\-\\
end if;\\
printf "Category=\%o, Graph=\%o\(\backslash\)n",kat,grp;\\
printf "\%o-symbols of (\%o,\%o,\%o): ",P,Q,R,S;\\
for j in [1..\#Char] do\+\\
   printf "\%o,",Char[j];\-\\
end for;\\
printf "\(\backslash\)n";
\end{tabbing}
}
}
\noindent
\textbf{Output:}
category, graph, and power residue symbols.
\end{algorithm}

\newpage

\section{Detailed Statistics of Categories and Graphs}
\label{s:StatMultCatGph3}

\noindent
Statistics of categories and graphs for the cyclic cubic quartets 
(with \(m=4\)) is shown in Table
\ref{tbl:StatCatGph3}.


\renewcommand{\arraystretch}{1.1}

\begin{table}[hb]
\caption{Statistics of Categories and Graphs of Cyclic Cubic Quartets}
\label{tbl:StatCatGph3}
\begin{center}

{\footnotesize

\begin{tabular}{|c|c||r|r|r|r||r|r|r|r||r|}

\hline
 Category & Gph & \multicolumn{4}{|c||}{\(\#\) Conductors \(c\)} & \multicolumn{4}{|c||}{\(\#\) Fields \(F\)} & \(c_{\text{min}}\) \\
\hline
 \(c<\) &  \(\) & \(25000\) & \(50000\) & \(75000\) & \(10^5\) & \(25000\) & \(50000\) & \(75000\) & \(10^5\) & \(\) \\
\hline
 I      & \(1\) &   \(7\) &  \(19\) &  \(27\) &  \(38\) &  \(28\) &   \(76\) &  \(108\) &  \(152\) &  \(4977\) \\
 I      & \(2\) &  \(15\) &  \(29\) &  \(44\) &  \(60\) &  \(60\) &  \(116\) &  \(176\) &  \(240\) &  \(7657\) \\
\hline
 Subtotal &     &  \(22\) &  \(48\) &  \(71\) &  \(98\) &  \(88\) &  \(192\) &  \(284\) &  \(392\) &           \\
\hline
 II     & \(1\) &  \(10\) &  \(25\) &  \(36\) &  \(47\) &  \(40\) &  \(100\) &  \(144\) &  \(188\) &  \(3913\) \\
 II     & \(2\) &   \(7\) &  \(19\) &  \(34\) &  \(45\) &  \(28\) &   \(76\) &  \(136\) &  \(180\) &  \(6327\) \\
\hline
 Subtotal &     &  \(17\) &  \(44\) &  \(70\) &  \(92\) &  \(68\) &  \(176\) &  \(280\) &  \(368\) &           \\
\hline
 III    & \(1\) &  \(11\) &  \(22\) &  \(41\) &  \(52\) &  \(44\) &   \(88\) &  \(164\) &  \(208\) &  \(1953\) \\
 III    & \(2\) &  \(57\) & \(125\) & \(181\) & \(262\) & \(228\) &  \(500\) &  \(724\) & \(1048\) &   \(819\) \\
 III    & \(3\) &  \(20\) &  \(50\) &  \(81\) & \(124\) &  \(80\) &  \(200\) &  \(324\) &  \(496\) &  \(1197\) \\
 III    & \(4\) &   \(5\) &   \(8\) &  \(16\) &  \(17\) &  \(20\) &   \(32\) &   \(64\) &   \(68\) &  \(6643\) \\
 III    & \(5\) &   \(4\) &  \(17\) &  \(27\) &  \(37\) &  \(16\) &   \(68\) &  \(108\) &  \(148\) & \(14049\) \\
 III    & \(6\) &   \(5\) &  \(16\) &  \(27\) &  \(31\) &  \(20\) &   \(64\) &  \(108\) &  \(124\) &  \(8541\) \\
 III    & \(7\) &   \(4\) &   \(8\) &  \(20\) &  \(34\) &  \(16\) &   \(32\) &   \(80\) &  \(136\) &  \(4599\) \\
 III    & \(8\) &   \(1\) &   \(4\) &   \(6\) &   \(7\) &   \(4\) &   \(16\) &   \(24\) &   \(28\) & \(20293\) \\
 III    & \(9\) &   \(3\) &   \(7\) &  \(11\) &  \(15\) &  \(12\) &   \(28\) &   \(44\) &   \(60\) & \(16471\) \\
\hline
 Subtotal &     & \(110\) & \(257\) & \(410\) & \(579\) & \(440\) & \(1028\) & \(1640\) & \(2316\) &           \\
\hline
 IV     & \(1\) &   \(0\) &   \(0\) &   \(2\) &   \(7\) &   \(0\) &    \(0\) &    \(8\) &   \(28\) & \(61579\) \\
 IV     & \(2\) &   \(0\) &   \(1\) &   \(1\) &   \(2\) &   \(0\) &    \(4\) &    \(4\) &    \(8\) & \(49543\) \\
 IV     & \(3\) &   \(0\) &   \(1\) &   \(2\) &   \(5\) &   \(0\) &    \(4\) &    \(8\) &   \(20\) & \(38311\) \\
\hline
 Subtotal &     &   \(0\) &   \(2\) &   \(5\) &  \(14\) &   \(0\) &    \(8\) &   \(20\) &   \(56\) &           \\
\hline
 Total    &     & \(149\) & \(351\) & \(556\) & \(783\) & \(596\) & \(1404\) & \(2224\) & \(3132\) &           \\
\hline

\end{tabular}

}

\end{center}
\end{table}

\noindent
In the range \(1<c<10^5\) of \(3\)-admissible conductors,
with a total of \(9457\) hits,
there are \(783\) with exactly \(t=3\) prime divisors,
which give rise to quartets,
according to the multiplicity formula
\eqref{eqn:Multiplicity}.
With \(262\) occurrences (\(33.5\%\)), category III, graph 2, \((q_i\to q_j;q_k)\), is most frequent.
On the second place with \(124\) (\(15.8\%\)) we have category III, graph 3, \((q_i\to q_j\to q_k)\).

\newpage

\section{Ambiguous Principal Ideals}
\label{s:AmbiguousIdeals}

\noindent
The number of \textit{primitive ambiguous ideals} of a cyclic cubic field \(F\)
increases with the number \(t\) of prime factors of the conductor \(c\).
According to Hilbert's Theorem \(93\), we have:
\begin{equation}
\label{eqn:AmbiguousIdeals}
\#\left(\mathcal{I}_F^{\langle\sigma\rangle}/\mathcal{I}_{\mathbb{Q}}\right)=3^t.
\end{equation}
However, the number of \textit{ambiguous principal ideals} of \(F\)
is a fixed invariant of all cyclic cubic fields, regardless of the number \(t\).
It is given by the Theorem on the \textbf{Herbrand quotient} of the unit group \(U_F\) of \(F\)
as a Galois module over the group \(\mathrm{Gal}(F/\mathbb{Q})=\langle\sigma\rangle\),
which can be expressed by abstract cohomology groups
\(\#\mathrm{H}^{-1}(\langle\sigma\rangle,U_F)/\#\hat{\mathrm{H}}^{0}(\langle\sigma\rangle,U_F)=\lbrack F:\mathbb{Q}\rbrack\)
or more ostensively by:
\begin{equation}
\label{eqn:AmbiguousPrincipalIdeals}
\#\left(\mathcal{P}_F^{\langle\sigma\rangle}/\mathcal{P}_{\mathbb{Q}}\right)
=\#\left(E_{F/\mathbb{Q}}/U_F^{1-\sigma}\right)
=\lbrack F:\mathbb{Q}\rbrack\cdot\#\left(U_{\mathbb{Q}}/\mathrm{N}_{F/\mathbb{Q}}{U_F}\right)=3,
\end{equation}
since the unit norm index is given by \(\left(U_{\mathbb{Q}}:\mathrm{Norm}_{F/\mathbb{Q}}{U_F}\right)=1\).
Consequently, if we speak about a \textit{non-trivial primitive ambiguous principal ideal} of \(F\),
then we either mean \((\alpha)=\alpha\mathcal{O}_F\) or \((\alpha^2/b)=(\alpha^2/b)\mathcal{O}_F\),
where \(\mathcal{P}_F^{\langle\sigma\rangle}/\mathcal{P}_{\mathbb{Q}}=\lbrace 1,(\alpha),(\alpha^2/b)\rbrace\).
The norms of these two elements are divisors of the square \(c^2=q_1^2\cdots q_t^2\)
of the conductor \(c\) of \(F\),
where \(q_t\) must be replaced by \(3\) if \(q_t=9\).
When \(\mathrm{N}_{F/\mathbb{Q}}{\alpha}=a\cdot b^2\)
with square free coprime integers \(a,b\), then
\(\mathrm{N}_{F/\mathbb{Q}}(\alpha^2/b)=a^2\cdot b^4/b^3=a^2\cdot b\).
It follows that both norms are cube free integers.


\begin{definition}
\label{dfn:PrincipalFactor}
The minimum of the two norms
of non-trivial primitive ambiguous principal ideals
\((\alpha),(\alpha^2/b)\)
of a cyclic cubic field \(F\) is called the \textbf{principal factor}
(of the discriminant \(d_F=c^2\) of the field)
\(B(F):=\min\lbrace a\cdot b^2,a^2\cdot b\rbrace\) of \(F\), that is
\begin{equation}
\label{eqn:PrincipalFactor}
B(F)=
\begin{cases}
a\cdot b^2 & \text{ if } b<a, \\
a^2\cdot b & \text{ if } a<b.
\end{cases}
\end{equation}
\end{definition}

\noindent
Ayadi 
\cite{Ay1995,AAI2001}
speaks about the \textit{Parry constant} or \textit{Parry invariant} \(B(F)\) of \(F\),
and Derhem
\cite{Dh2002}
calls \(B(F)\) the \textit{Kummer resolvent} of \(F\).


\begin{algorithm}
\label{alg:PrincipalFactor}
(Principal Factor.) \\
\textbf{Input:}
critical prime \texttt{p}, conductor \texttt{f}, cyclic number field \texttt{N}. \\
\textbf{Code:}
{\scriptsize
\texttt{
\begin{tabbing}
for \= for \= for \= for \= for \= for \= \kill
    ON := MaximalOrder(N); // algebraic integers of the base field N\\
    // trivial initialization of prime ideals\\
    Pid := Decomposition(ON,p)[1][1];\\
    Qid := Decomposition(ON,p)[1][1];\\
    Rid := Decomposition(ON,p)[1][1];\\
    Sid := Decomposition(ON,p)[1][1];\\
    // number of prime divisors of f totally ramified in N\\
    nPrimes := 0;\\
    q := 1;\\
    while (f ge q) do\+\\
       q := NextPrime(q);\\
       if (0 eq f mod q) then\+\\
          nPrimes := nPrimes + 1;\\
          // full decomposition of q in N: necessarily a singlet\\
          sPrimesN := Decomposition(ON,q);\\
          if (1 eq nPrimes) then\+\\
             Pid := sPrimesN[1][1];\-\\
          end if;\\
          if (2 eq nPrimes) then\+\\
             Qid := sPrimesN[1][1];\-\\
          end if;\\
          if (3 eq nPrimes) then\+\\
             Rid := sPrimesN[1][1];\-\\
          end if;\\
          // and so on for bigger number nPrimes\-\\
       end if; // q divides f\-\\
    end while; // q\\
    // possible absolute principal factorizations in N\\
    nPrincipalN := 0;\\
    if (1 eq nPrimes) then\+\\
       for e in [0..p-1] do\+\\
          I := (Pid\({}\,\hat{}\,{}\)e);\\
          if IsPrincipal(I) then\+\\
             nPrincipalN := nPrincipalN + 1;\\
             printf "   \%o: [\%o]\(\backslash\)n",nPrincipalN,e;\-\\
          end if;\-\\
       end for;\\
       printf "   \%o power products of primes among \%o\(\backslash\)n",nPrincipalN,p;\-\\
    end if; // 1\\
    if (2 eq nPrimes) then\+\\
       for e in [0..p-1] do\+\\
          for f in [0..p-1] do\+\\
             I := (Pid\({}\,\hat{}\,{}\)e)\(\ast\)(Qid\({}\,\hat{}\,{}\)f);\\
             if IsPrincipal(I) then\+\\
                nPrincipalN := nPrincipalN + 1;\\
                printf "   \%o: [\%o,\%o]\(\backslash\)n",nPrincipalN,e,f;\-\\
             end if;\-\\
          end for;\-\\
       end for;\\
       printf "   \%o power products of primes among \%o\(\backslash\)n",nPrincipalN,p\({}\,\hat{}\,{}\)2;\-\\
    end if; // 2\\
    if (3 eq nPrimes) then\+\\
       for e in [0..p-1] do\+\\
          for f in [0..p-1] do\+\\
             for g in [0..p-1] do\+\\
                I := (Pid\({}\,\hat{}\,{}\)e)\(\ast\)(Qid\({}\,\hat{}\,{}\)f)\(\ast\)(Rid\({}\,\hat{}\,{}\)g);\\
                if IsPrincipal(I) then\+\\
                   nPrincipalN := nPrincipalN + 1;\\
                   printf "   \%o: [\%o,\%o,\%o]\(\backslash\)n",nPrincipalN,e,f,g;\-\\
                end if;\-\\
             end for;\-\\
          end for;\-\\
       end for;\\
       printf "   \%o power products of primes among \%o\(\backslash\)n",nPrincipalN,p\({}\,\hat{}\,{}\)3;\-\\
    end if; // 3\\
    // and so on for bigger number nPrimes
\end{tabbing}
}
}
\noindent
\textbf{Output:}
prime factorization of norms of ambiguous principal ideals.
\end{algorithm}


\chapter{Unramified Extensions}
\label{ch:UnramifiedExtensions}

\section{The Absolute 3-Genus Field}
\label{s:GenusField}

\noindent
The \textit{absolute \(3\)-genus field} \(F^\ast=(F/\mathbb{Q})^\ast\)
of a cyclic cubic field \(F\) is the maximal unramified \(3\)-extension \(F^\ast/F\)
with abelian absolute Galois group \(\mathrm{Gal}(F^\ast/\mathbb{Q})\).
If the conductor \(c=q_1\cdots q_t\) of \(F=F_c\) has \(t\) prime divisors,
then \(F^\ast\) is the compositum of the multiplet \((F_{c,1},\ldots,F_{c,m})\)
of all cyclic cubic fields sharing the common conductor \(c\),
where \(m=m(c)=2^{t-1}\), according to the multiplicity formula
\eqref{eqn:Multiplicity}.
The absolute Galois group \(\mathrm{Gal}(F^\ast/\mathbb{Q})\)
is the elementary abelian \(3\)-group \(C_3^t\).
In particular, if \(t=1\) then \(F^\ast=F\),
and if \(t=2\), \(c=q_1q_2\), then \(F^\ast=F_{c,1}\cdot F_{c,2}\) is a bicyclic bicubic field
with conductor \(c\) and discriminant
\begin{equation}
\label{eqn:Discriminant}
d(F^\ast)=d(F_{q_1})\cdot d(F_{q_2})\cdot d(F_{c,1})\cdot d(F_{c,2})=q_1^2\cdot q_2^2\cdot(q_1q_2)^2\cdot(q_1q_2)^2=c^6.
\end{equation}
Parry
\cite{Pa1990}
has investigated the arithmetic of a general \textit{bicyclic bicubic field} \(K\)
with four cyclic cubic subfields \(K_1,\ldots,K_4\).
In particular, he determined the \textit{class number relation}
\begin{equation}
\label{eqn:ParryFormula}
h=\frac{I}{3^5}\cdot\prod_{i=1}^4\,h_i,
\end{equation}
where \(I:=(U:V)\) with \(V:=\prod_{i=1}^4\,U_i\) denotes the \textit{index of subfield units} of \(K\)
\cite[Prp. 7, p. 496]{Pa1990}.
Generally, \(I\) is a divisor of \(27=3^3\)
\cite[Lem. 11, p. 500]{Pa1990}.


For a cyclic cubic field \(F\) with \(t=2\),
the \(3\)-class numbers of the \(3\)-genus field \(F^\ast\) and its subfields
can be summarized as follows.

\begin{theorem}
\label{thm:TwoPrimeGenus}
Let \(F^\ast=F_{q_1}\cdot F_{q_2}\cdot F_{c,1}\cdot F_{c,2}\) be the genus field
of the two cyclic cubic fields \(F_{c,1}\) and \(F_{c,2}\) with conductor \(c=q_1q_2\).
Denote the \(3\)-valuations of the class numbers \(h\), \(h_1\), \(h_2\), \(h_3\), \(h_4\)
of \(F^\ast\), \(F_{q_1}\), \(F_{q_2}\), \(F_{c,1}\), \(F_{c,2}\), respectively,
by \(v\), \(v_1\), \(v_2\), \(v_3\), \(v_4\). Then \(v_1=v_2=0\), and
\begin{equation}
\label{eqn:TwoPrimeGenus}
v
\begin{cases}
=0,\ v_3=v_4=1,\ I=27, & \text{ if not } q_1\leftrightarrow q_2, \\
=1, & \text{ if } q_1\leftrightarrow q_2,\ v_3=v_4=2,\ I=9, \\
=2, & \text{ if } q_1\leftrightarrow q_2,\ v_3=v_4=2,\ I=27, \\
\ge 3, & \text{ if } q_1\leftrightarrow q_2,\ v_3\ge 3,\ v_4\ge 3,\ I\ge 9.
\end{cases}
\end{equation}
\end{theorem}

\begin{proof}
According to Theorem
\ref{thm:3ClassRank},
we generally have \(v_1=v_2=0\),
\(v_3\ge 1\), \(v_4\ge 1\) if not \(q_1\leftrightarrow q_2\),
and \(v_3\ge 2\), \(v_4\ge 2\) if \(q_1\leftrightarrow q_2\).
Now, the claim is an immediate consequence of Formula
\eqref{eqn:ParryFormula},
which yields
\[v=v_3h=v_3I-5+\sum_{i=1}^4\,v_3h_i=v_3I-5+v_1+v_2+v_3+v_4=v_3I-5+v_3+v_4.\]
The combination of
\cite[Thm. 9, p. 497]{Pa1990}
and
\cite[Cor. 1, p. 498]{Pa1990}
shows that \(v=0\) if and only if not \(q_1\leftrightarrow q_2\),
and \(v=0\) implies \(v_3I=3\),
whence necessarily \(v_3=v_4=1\).
However, if \(q_1\leftrightarrow q_2\), then
\(v_3=2\) is equivalent with \(v_4=2\), according to
\cite[Thm. 4.1, p. 472]{AAI2001}.
\end{proof}

\begin{remark}
\label{rmk:TwoPrimeGenus}
For \(v_3=v_4=2\), we have \(\mathrm{Cl}_3{F_{q_1q_2,\mu}}\simeq (3,3)\).
The smallest occurrences of \(v_3=v_4=3\) are the conductors \(7\cdot 673=4\,711\)
(Eau de Cologne) and \(7\cdot 769=5\,383\)
with \(\mathrm{Cl}_3{F_{q_1q_2,\mu}}\simeq (9,3)\), \(\mu=1,2\).
They will be considered in Section \S\
\ref{ch:Singular}.
\end{remark}


\noindent
For a cyclic cubic field \(F\) with \(t=3\) and conductor \(c=q_1q_2q_3\),
the \(3\)-genus field \(F^\ast\) contains \(13\) bicyclic bicubic subfields.
Three of them are the \textit{sub genus fields} \(F_i^\ast\), \(1\le i\le 3\),
of the cyclic cubic fields with conductors \(q_1q_2\), \(q_2q_3\), \(q_3q_1\),
respectively.
In the sequel, we always start with the leading three sub genus fields \(F_i^\ast\), \(1\le i\le 3\),
separated by a semicolon from the trailing ten remaining bicyclic bicubic subfields,
when we give a family of invariants for these \(13\) subfields \(S_1,\ldots,S_{13}\),
\begin{equation}
\label{eqn:SubGenusFields}
\text{in particular, }
\left\lbrack\mathrm{Cl}_3{S_i}\right\rbrack_{1\le i\le 13}:=\lbrack \mathrm{Cl}_3{F_1^\ast},\ldots,\mathrm{Cl}_3{F_3^\ast};\mathrm{Cl}_3{S_4},\ldots,\mathrm{Cl}_3{S_{13}}\rbrack.
\end{equation}


\begin{algorithm}
\label{alg:BicyclicBicubic}
(Class groups of \(13\) bicyclic bicubic fields.) \\
\textbf{Input:}
critical prime \texttt{p}, prime power factorization \texttt{cPD} of the conductor \texttt{c}. \\
\textbf{Code:}
{\scriptsize
\texttt{
\begin{tabbing}
for \= for \= for \= for \= for \= \kill
      numPF := \#cPD;\\
      if (1 le numPF) then\+\\
         prime1 := cPD[1][1];\\
         if (p eq prime1) then\+\\
            prime1 := p\({}\,\hat{}\,{}\)2;\-\\
         end if;\\
         if (2 le numPF) then\+\\
            prime2 := cPD[2][1];\\
            if (3 le numPF) then\+\\
               prime3 := cPD[3][1];\\
               if (4 le numPF) then\+\\
                  prime4 := cPD[4][1];\-\\
               end if; // 4\-\\
            end if; // 3\-\\
         end if; // 2\-\\
      end if; // 1\\
            f1    := prime1\(\ast\)prime2;\\
            Disc1 := f1\({}\,\hat{}\,{}\)2;\\
            f2    := prime1\(\ast\)prime3;\\
            Disc2 := f2\({}\,\hat{}\,{}\)2;\\     
            // and so on until\\
            Disc  := c\({}\,\hat{}\,{}\)2; // root discriminant\\
            if (Disc\({}\,\hat{}\,{}\)((p-1) div 2) eq AbsoluteDiscriminant(ON)) then\+\\
               fld := fld + 1;\\
               locMult := locMult + 1;\\
               if (1 eq locMult) then\+\\
                  Component1 := N;\-\\
               elif (2 eq locMult) then\+\\
                  Component2 := N;\-\\
               elif (3 eq locMult) then\+\\
                  Component3 := N;\-\\
               else\+\\
                  Component4 := N;\-\\
               end if;\-\\
            elif (Disc1\({}\,\hat{}\,{}\)((p-1) div 2) eq AbsoluteDiscriminant(ON)) then\+\\
               locMult1 := locMult1 + 1;\\
               if (1 eq locMult1) then\+\\
                  Sub11 := N;\-\\
               elif (2 eq locMult1) then\+\\
                  Sub12 := N;\-\\
               end if;\-\\
            elif (Disc2\({}\,\hat{}\,{}\)((p-1) div 2) eq AbsoluteDiscriminant(ON)) then\+\\
               locMult2 := locMult2 + 1;\\
               if (1 eq locMult2) then\+\\
                  Sub21 := N;\-\\
               elif (2 eq locMult2) then\+\\
                  Sub22 := N;\-\\
               end if;\-\\
            end if;
            // and so on until\\
            elif (2 eq locMult) then\+\\
               Append(\(\sim\)Doublets,Collection);\\
               printf "\(\backslash\)nm=2, no. \%o\(\backslash\)n",Mult[2];\\
               GenusField := Compositum(Component1,Component2);\\
               OGK := MaximalOrder(GenusField);\\
               CGK,mCGK := ClassGroup(OGK);\\
               AGK := AbelianInvariants(CGK);\\
               PGK := pPrimaryInvariants(CGK,p);\\
               VGK := Valuation(\#CGK,p);\\
               printf "Genus Field 2: SGK=\%o, CGK=\%o\(\backslash\)n",PGK,AGK;\-\\
            elif (4 eq locMult) then\+\\
               Append(\(\sim\)Quartets,Collection);\\
               printf "\(\backslash\)nm=4, no. \%o\(\backslash\)n",Mult[4];\\
               // ten bicyclic bicubic fields\\
               B1  := Compositum(Sub11,Sub21);\\
               B2  := Compositum(Sub11,Sub22);\\
               B3  := Compositum(Sub12,Sub22);\\
               B4  := Compositum(Sub12,Sub21);\\
               B5  := Compositum(Component1,Component3);\\
               B6  := Compositum(Component1,Component4);\\
               B7  := Compositum(Component1,Component2);\\
               B8  := Compositum(Component2,Component4);\\
               B9  := Compositum(Component2,Component3);\\
               B10 := Compositum(Component3,Component4);\\
               // their class groups\\
               OB := MaximalOrder(B1);\\
               CB,mCB := ClassGroup(OB);\\
               AB := AbelianInvariants(CB);\\
               PB := pPrimaryInvariants(CB,p);\\
               VB := Valuation(\#CB,p);\\
               printf "Bicubic 1:  SB=\%o, CB=\%o\(\backslash\)n",PB,AB;\\
               // and so on until\\
               OB := MaximalOrder(B10);\\
               CB,mCB := ClassGroup(OB);\\
               AB := AbelianInvariants(CB);\\
               PB := pPrimaryInvariants(CB,p);\\
               VB := Valuation(\#CB,p);\\
               printf "Bicubic 10: SB=\%o, CB=\%o\(\backslash\)n",PB,AB;\-\\
            end if;
\end{tabbing}
}
}
\noindent
\textbf{Output:}
class groups of ten bicyclic bicubic subfields of the genus field.
\end{algorithm}


\section{Capitulation Kernels of Cyclic Cubic Fields}
\label{s:Capitulation}

\noindent
Finally, we recall the connection between the size of the capitulation kernel \(\ker(T_{E/F})\)
and the unit norm index \((U_F:\mathrm{Norm}_{E/F}U_E)\)
of an unramified cyclic cubic extension \(E/F\) of a cyclic cubic field \(F\).
Here, \(T_{E/F}:\,\mathrm{Cl}_3{F}\to\mathrm{Cl}_3{E}\),
\(\mathfrak{a}\mathcal{P}_F\mapsto(\mathfrak{a}\mathcal{O}_E)\mathcal{P}_E\),
denotes the transfer of \(3\)-classes from \(F\) to \(E\).
\begin{theorem}
\label{thm:Capitulation}
The order of the \(3\)-capitulation kernel of \(E/F\) is given by
\begin{equation}
\label{eqn:Capitulation}
\#\ker(T_{E/F})=
\begin{cases}
3, \\
9, \\
27,
\end{cases}
\text{ if and only if }\quad
(U_F:\mathrm{Norm}_{E/F}U_E)=
\begin{cases}
1, \\
3, \\
9.
\end{cases}
\end{equation}
\end{theorem}
\begin{proof}
According to the Herbrand Theorem on the cohomology of the unit group \(U_E\)
as a Galois module with respect to \(G=\mathrm{Gal}(E/F)\),
we have the relation
\(\#\ker(T_{E/F})=\lbrack E:F\rbrack\cdot(U_F:\mathrm{Norm}_{E/F}U_E)\),
since \(\ker(T_{E/F})\simeq H^1(G,U_E)\)
when \(E/F\) is unramified of odd prime degree \(\lbrack E:F\rbrack=3\)
and \(U_F/\mathrm{Norm}_{E/F}U_E\simeq\hat{H}^0(G,U_E)\).
The cyclic cubic base field \(F\) has signature \((r_1,r_2)=(3,0)\)
and torsionfree Dirichlet unit rank \(r=r_1+r_2-1=3+0-1=2\).
Thus, there are three possibilities for the unit norm index
\((U_F:\mathrm{Norm}_{E/F}U_E)\in\lbrace 1,3,9\rbrace\).
\end{proof}

\begin{remark}
\label{rmk:Capitulation}
When \(F\) is a cyclic cubic field
with \(3\)-class group \(O:=\mathrm{Cl}_3{F}\) of elementary tricyclic type \((3,3,3)\),
viewed as a vector space of dimension \(3\) over the finite field \(\mathbb{F}_3\), then
\(\#\ker(T_{E/F})=3\) if and only if \((\exists\,1\le i\le 13)\) \(\ker(T_{E/F})=L_i\)
is a line in Table
\ref{tbl:Lines},
\(\#\ker(T_{E/F})=9\) if and only if \((\exists\,1\le i\le 13)\) \(\ker(T_{E/F})=P_i\)
is a plane in Table
\ref{tbl:Planes}, and
\(\#\ker(T_{E/F})=27\) if and only if \(\ker(T_{E/F})=O\)
is the entire vector space.
See Section \S\
\ref{s:ElmAblRk3}.
\end{remark}


\begin{algorithm}
\label{alg:ArtinPatternPP}
(Artin pattern for type \((3,3)\).) \\
\textbf{Input:}
critical prime \texttt{p}, algebraic number field \texttt{F} of type \((p,p)\). \\
\textbf{Code:}
{\scriptsize
\texttt{
\begin{tabbing}
for \= for \= for \= for \= for \= for \= \kill
O := MaximalOrder(F); // algebraic integers of the base field F\\
C,mC := ClassGroup(O);\\
epsilon := 0;       // counter of (3,3,3)\\
polarization1 := 3;\\
polarization2 := 3; // co-polarization\\
fixedpoints := 0;\\
capitulations := 0; // occupation number of 0\\
occupation := 0;    // cardinality of occupation support (including 0)\\
repetitions := 0;   // maximal occupation number (except for 0)\\
intersection := 0;  // meet of repetitions and fixed points\\
sS := Subgroups(C: Quot := [p]); // subgroups of index p\\
sI := []; // fixed ordering of subgroups\\
for j in [1..p+1] do\+\\
   Append(\(\sim\)sI,0);\-\\
end for; // j\\
n := Ngens(C);\\
H := (Order(C.(n-1)) div p)\(\ast\)C.(n-1); // 1st p-generator\\
K := (Order(C.n) div p)\(\ast\)C.n;         // 2nd p-generator\\
ct := 0; // local counter\\
for x in sS do\+\\
   ct := ct+1;\\
   if H in x\({}\,\grave{}\,{}\)subgroup then\+\\
      sI[1] := ct;\-\\
   end if; // n-1\\
   if K in x\({}\,\grave{}\,{}\)subgroup then\+\\
      sI[2] := ct;\-\\
   end if; // n\\
   for e in [1..p-1] do\+\\
      if (H+(e\(\ast\)K)) in x\({}\,\grave{}\,{}\)subgroup then\+\\
         sI[e+2] := ct;\-\\
      end if; // product\-\\
   end for; // e\-\\
end for; // x\\
// p+1 unramified extensions of degree p\\
sA := [AbelianExtension(Inverse(mQ)\(\ast\)mC)\+\+\\
       where Q,mQ := quo<C|x\({}\,\grave{}\,{}\)subgroup>: x in sS];\-\-\\
sN := [NumberField(x): x in sA];   // relative extensions\\
sR := [MaximalOrder(x): x in sA];\\
sF := [AbsoluteField(x): x in sN]; // absolute extensions\\
sM := [MaximalOrder(x): x in sF];\\
sM := [OptimizedRepresentation(x): x in sF];     // 1st optimization\\
sA := [NumberField(DefiningPolynomial(x)): x in sM];\\
sO := [Simplify(LLL(MaximalOrder(x))): x in sA]; // 2nd optimization\\
// Artin pattern with two components TTT and TKT\\
TTT := []; // transfer target type\\
for j in [1..\#sO] do\+\\
   CO := ClassGroup(sO[j]);\\
   Append(\(\sim\)TTT,pPrimaryInvariants(CO,p));\\
   if (3 eq \#pPrimaryInvariants(CO,p)) then\+\\
      epsilon := epsilon + 1;\-\\
   end if;\\
   val := Valuation(Order(CO),p);\\
   if (2 eq val) then\+\\
      polarization2 := val;\-\\
   elif (4 le val) then\+\\
      if (3 eq polarization1) then\+\\
         polarization1 := val;\-\\
      else\+\\
         polarization2 := val;\-\\
      end if;\-\\
   end if; // val\-\\
end for; // j\\
// using the fixed ordering of subgroups\\
TKT := []; // transfer kernel type\\
for j in [1..\#sR] do // use sR rather than sO\+\\
   Collector := [];\\
   I := sR[j]!!mC(H);\\
   if IsPrincipal(I) then\+\\
      Append(\(\sim\)Collector,sI[1]);\-\\
   end if; // I\\
   I := sR[j]!!mC(K);\\
   if IsPrincipal(I) then\+\\
     Append(\(\sim\)Collector,sI[2]);\-\\
   end if; // I\\
   for e in [1..p-1] do\+\\
      I := sR[j]!!mC(H+(e\(\ast\)K));\\
      if IsPrincipal(I) then\+\\
         Append(\(\sim\)Collector,sI[e+2]);\-\\
      end if; // I\-\\
   end for; // e\\
   if (2 le \#Collector) then // total capitulation\+\\
      Append(\(\sim\)TKT,0);\-\\
   else                       // partial capitulation\+\\
      Append(\(\sim\)TKT,Collector[1]);\-\\
   end if; // \#Collector\-\\
end for; // j\\
TAB := []; // Taussky conditions A and B\\
image := []; // local image collection\\
for j in [1..\#TKT] do\+\\
   if (j eq TKT[j]) then // fixed point\+\\
      Append(\(\sim\)TAB,"A");\\
      fixedpoints := fixedpoints + 1;\-\\
   elif (0 eq TKT[j]) then // total\+\\
      Append(\(\sim\)TAB,"A");\\
      capitulations := capitulations + 1;\-\\
   else // non-fixed point partial\+\\
      Append(\(\sim\)TAB,"B");\-\\
   end if; // fixed point or total or partial\\
   if not (TKT[j] in image) then\+\\
      Append(\(\sim\)image,TKT[j]);\-\\
   end if; // image\-\\
end for; // j\\
occupation := \#image;\\
doublet := 0; // local memory\\
for digit in [1..p+1] do\+\\
   counter := 0; // local counter\\
   for j in [1..\#TKT] do\+\\
      if (digit eq TKT[j]) then\+\\
         counter := counter + 1;\-\\
      end if;\-\\
   end for; // j\\
   if (counter ge 2) then\+\\
      doublet := digit; // last assignment persists\-\\
   end if;\\
   if (counter gt repetitions) then\+\\
      repetitions := counter;\-\\
   end if;\-\\
end for; // digit\\
if (doublet ge 1) then\+\\
   if (doublet eq TKT[doublet]) then\+\\
      intersection := 1;\-\\
   end if;\-\\
end if;\\
printf "TKT \%o; TAB \%o; INV \%o; \%o; \%o; TTT ( ",TKT,TAB,epsilon,polarization1,polarization2;\\
for j in [1..\#TTT] do\+\\
   printf "\%o ",TTT[j];\-\\
end for; // j\\
printf ")\(\backslash\)n";
\end{tabbing}
}
}
\noindent
\textbf{Output:}
transfer kernel type \texttt{TKT}, transfer target type \texttt{TTT}, other invariants.
\end{algorithm}


\begin{algorithm}
\label{alg:ArtinPatternPPP}
(Artin pattern for type \((3,3,3)\).) \\
\textbf{Input:}
critical prime \texttt{p}, algebraic number field \texttt{N} of type \((p,p,p)\). \\
\textbf{Code:}
{\tiny
\texttt{
\begin{tabbing}
for \= for \= for \= for \= for \= for \= \kill
ON := MaximalOrder(N); // algebraic integers of the base field N\\
CN,mCN := ClassGroup(ON);\\
PN := pPrimaryInvariants(CN,p);\\
if ([p,p,p] eq PN) then\+\\
   sS2 := Subgroups(C: Quot := [p]);\\
   sA2 := [AbelianExtension(Inverse(mQ1)\(\ast\)mC)\+\+\\
           where Q1,mQ1 := quo<C|x\({}\,\grave{}\,{}\)subgroup>: x in sS2];\-\-\\
   sN2 := [NumberField(x): x in sA2];   // relative extensions\\
   sO  := [MaximalOrder(x): x in sA2];\\
   sF2 := [AbsoluteField(x): x in sN2]; // absolute extensions\\
   sM2 := [MaximalOrder(x): x in sF2];\\
   sB2 := [OptimizedRepresentation(x): x in sF2];     // 1st optimization\\
   sK2 := [NumberField(DefiningPolynomial(x)): x in sB2];\\
   sO2 := [Simplify(LLL(MaximalOrder(x))): x in sK2]; // 2nd optimization\\
   n := Ngens(C);\\
   u := (Order(C.n) div p)\(\ast\)C.n;         // 1st p-generator\\
   v := (Order(C.(n-1)) div p)\(\ast\)C.(n-1); // 2nd p-generator\\
   w := (Order(C.(n-2)) div p)\(\ast\)C.(n-2); // 3rd p-generator\\
   ArtinMap := []; // fixed ordering of subgroups\\
   for x in sS2 do\+\\
      Collector := [];\\
      if u in x\({}\,\grave{}\,{}\)subgroup then\+\\
         Append(\(\sim\)Collector,1);\-\\
      end if; // n\\
      if v in x\({}\,\grave{}\,{}\)subgroup then\+\\
         Append(\(\sim\)Collector,2);\-\\
      end if; // n-1\\
      if w in x\({}\,\grave{}\,{}\)subgroup then\+\\
         Append(\(\sim\)Collector,3);\-\\
      end if; // n-2\\
      for e in [1..p-1] do\+\\
         if u+e\(\ast\)v in x\({}\,\grave{}\,{}\)subgroup then\+\\
            Append(\(\sim\)Collector,3+e);\-\\
         end if;\-\\
      end for; // e\\
      for e in [1..p-1] do\+\\
         if v+e\(\ast\)w in x\({}\,\grave{}\,{}\)subgroup then\+\\
            Append(\(\sim\)Collector,2+p+e);\-\\
         end if;\-\\
      end for; // e\\
      for e in [1..p-1] do\+\\
         if w+e\(\ast\)u in x\({}\,\grave{}\,{}\)subgroup then\+\\
            Append(\(\sim\)Collector,1+2\(\ast\)p+e);\-\\
         end if;\-\\
      end for; // e\\
      for e in [1..p-1] do\+\\
         for f in [1..p-1] do\+\\
            if u+e\(\ast\)v+f\(\ast\)w in x\({}\,\grave{}\,{}\)subgroup then\+\\
               Append(\(\sim\)Collector,3\(\ast\)p+(e-1)\(\ast\)(p-1)+f);\-\\
            end if;\-\\
         end for; // f\-\\
      end for; // e\\
      Append(\(\sim\)ArtinMap,Collector);\-\\
   end for; // x\\
   printf "Artin Map: ";\\
   for i in [1..p\({}\,\hat{}\,{}\)2+p+1] do\+\\
      printf "(";
      for k in [1..\#ArtinMap[i]] do\+\\
         printf "\%o,",ArtinMap[i][k];\-\\
      end for;
      printf ")";\-\\
   end for;\\
   printf "\(\backslash\)n";\\
   ATI := []; // transfer target type\\
   for j in [1..\#sO2] do\+\\
      C2,mC2 := ClassGroup(sO2[j]);\\
      Append(\(\sim\)ATI,pPrimaryInvariants(C2,p));\-\\
   end for; // j
   printf "TTT = (";\\
   for i in [1..p\({}\,\hat{}\,{}\)2+p+1] do\+\\
      printf "\%o,",ATI[i];\-\\
   end for;\\
   printf ")\(\backslash\)n";\\
   // using the fixed ordering of subgroups\\
   TKT := []; // transfer kernel type\\
   for j in [1..\#sO] do\+\\
      Collector := [];\\
      I := sO[j]!!mC(u);\\
      if IsPrincipal(I) then\+\\
         Append(\(\sim\)Collector,1);\-\\
      end if;\\
      I := sO[j]!!mC(v);\\
      if IsPrincipal(I) then\+\\
         Append(\(\sim\)Collector,2);\-\\
      end if;\\
      I := sO[j]!!mC(w);\\
      if IsPrincipal(I) then\+\\
         Append(\(\sim\)Collector,3);\-\\
      end if;\\
      for e in [1..p-1] do\+\\
         I := sO[j]!!mC(u+e\(\ast\)v);\\
         if IsPrincipal(I) then\+\\
            Append(\(\sim\)Collector,3+e);\-\\
         end if;\-\\
      end for; // e\\
      for e in [1..p-1] do\+\\
         I := sO[j]!!mC(v+e\(\ast\)w);\\
         if IsPrincipal(I) then\+\\
            Append(\(\sim\)Collector,2+p+e);\-\\
         end if;\-\\
      end for; // e\\
      for e in [1..p-1] do\+\\
         I := sO[j]!!mC(w+e\(\ast\)u);\\
         if IsPrincipal(I) then\+\\
            Append(\(\sim\)Collector,1+2\(\ast\)p+e);\-\\
         end if;\-\\
      end for; // e\\
      for e in [1..p-1] do\+\\
         for f in [1..p-1] do\+\\
            I := sO[j]!!mC(u+e\(\ast\)v+f\(\ast\)w);\\
            if IsPrincipal(I) then\+\\
               Append(\(\sim\)Collector,3\(\ast\)p+(e-1)\(\ast\)(p-1)+f);\-\\
            end if;\-\\
         end for; // f\-\\
      end for; // e\\
      if (p\({}\,\hat{}\,{}\)2+p+1 eq \#Collector) then // total capitulation\+\\
         Append(\(\sim\)TKT,[0]);\-\\
      else\+\\
         Append(\(\sim\)TKT,Collector);\-\\
      end if;\-\\
   end for; // j\\
   printf "TKT = (";\\
   for i in [1..p\({}\,\hat{}\,{}\)2+p+1] do\+\\
      printf "\%o,",TKT[i];\-\\
   end for;\\
   printf ")\(\backslash\)n";\-\\
end if; // (p,p,p)
\end{tabbing}
}
}
\noindent
\textbf{Output:}
\texttt{ArtinMap}, transfer target type \texttt{TTT}, transfer kernel type \texttt{TKT}.
\end{algorithm}


\chapter{Classical Results}
\label{ch:Classical}

\section{Cyclic Cubic Doublets of Type (3,3)}
\label{s:CycCub3x3Classical}


\subsection{Conductors with Two Prime Divisors}
\label{ss:TwoPrimeCond}

\noindent
First we focus on conductors \(c\) with \textit{two} prime divisors, \(t=2\).
According to the multiplicity formula \(m(c)=(3-1)^{t-1}\),
there are \(2\) cyclic cubic fields \(F_{c,1},F_{c,2}\) sharing the common conductor \(c\).
If one of them has \(3\)-class group \((3,3)\),
then the same is true for the other
\cite[Thm. 4.1, p. 472]{AAI2001}.
Necessarily, the prime divisors of the conductor \(c\)
are mutual cubic residues with respect to each other,
\(\left(\frac{q_1}{q_2}\right)_3=\left(\frac{q_2}{q_1}\right)_3=+1\),
with graph \(q_1\leftrightarrow q_2\), by Theorem
\ref{thm:3ClassRank}.
Further, it turns out that both members \(F=F_{c,\mu}\), \(\mu\in\lbrace 1,2\rbrace\),
of the \textit{doublet} have the same \(3\)-capitulation type \(\varkappa(F)\)
\cite[Prp. 3.3, p. 27]{Ay1995}
and the same \(3\)-class field tower group
\(G:=\mathrm{G}_3^{(\infty)}{F}:=\mathrm{Gal}(F_3^{(\infty)}/F)\),
which can be determined with the aid of the following Theorem
\ref{thm:TwoPrimeCond}.
Thus, we have only one graph, \(q_1\leftrightarrow q_2\), and do not need categories in the case \(t=2\).
In the sequel, we always use \textit{identifiers} \(\langle o,i\rangle\) of the SmallGroups Library
\cite{BEO2002,BEO2005}
in order to characterize a group \(G\) by its order \(o=\#G\) and a counting number \(i\), enclosed in angle brackets.


\begin{theorem}
\label{thm:TwoPrimeCond}
(Principal factor criterion, \textbf{Ayadi, 1995},
\cite{Ay1995,AAI2001}.) \\
Let \(c\) be a conductor divisible by two primes, \(t=2\), such that
\(\mathrm{Cl}_3{F_{c,\mu}}\simeq (3,3)\)
for both cyclic cubic fields
\(F_{c,\mu}\), \(1\le\mu\le 2\),
with conductor \(c\).
Denote by \(n\) the number of prime divisors of the norm
\(B(F)=\mathrm{N}_{F/\mathbb{Q}}(\alpha)\)
of a non-trivial primitive ambiguous principal ideal \((\alpha)\),
i.e. a \textbf{principal factor},
of any of the two fields \(F_{c,\mu}\). \\
Then \(n\in\lbrace 1,2\rbrace\), and the second \(3\)-class group
\(\mathfrak{M}:=\mathrm{G}_3^{(2)}{F}:=\mathrm{Gal}(F_3^{(2)}/F)\)
of both fields \(F=F_{c,\mu}\) is given by
\begin{equation}
\label{eqn:TwoPrimeCondAmbiguous}
\mathrm{G}_3^{(2)}{F}\simeq
\begin{cases}
\langle 9,2\rangle \text{ with capitulation type } \mathrm{a}.1,\ \varkappa{F}=(0000),  & \text{ if } n=2, \\
\langle 27,4\rangle \text{ with capitulation type } \mathrm{A}.1,\ \varkappa{F}=(1111), & \text{ if } n=1.
\end{cases}
\end{equation}
The length of the Hilbert \(3\)-class field tower is
\(\ell_3{F}=1\) with \(F_3^{(\infty)}=F_3^{(1)}\) if \(n=2\), and
\(\ell_3{F}=2\) with \(F_3^{(\infty)}=F_3^{(2)}\) if \(n=1\).
In both cases, \(G:=\mathrm{G}_3^{(\infty)}{F}=\mathrm{G}_3^{(2)}{F}\).
\end{theorem}

\begin{proof}
See
\cite[Prp. 3.6, p. 32, Thm. 3.1, p. 34, Thm. 3.3, p. 37]{Ay1995}.
\end{proof}

Concerning the \(3\)-capitulation types \(\mathrm{a}.1\) and \(\mathrm{A}.1\),
viewed as transfer kernel types (TKT),
and the related concept of transfer target types (TTT), see
\cite{Ma2012}.


The criteria in Theorem
\ref{thm:TwoPrimeCond}
can also be expressed with the aid of the \(3\)-class number
\(h_3{F^\ast}=\#\mathrm{Cl}_3{F^\ast}\)
of the \(3\)-\textbf{genus field} \(F^\ast\)
instead of principal factors.

\begin{corollary}
\label{cor:TwoPrimeCond}
(Genus field criterion, \textbf{Ayadi, 2001},
\cite{Ay1995,AAI2001}.) \\
Under the assumptions of Theorem
\ref{thm:TwoPrimeCond},
the second \(3\)-class group of both fields \(F=F_{c,\mu}\), \(1\le\mu\le 2\),
is characterized by the \(3\)-valuation \(v\) of the class number of \(F^\ast\).
\begin{equation}
\label{eqn:TwoPrimeCondGenus}
\mathrm{G}_3^{(2)}{F}\simeq
\begin{cases}
\langle 9,2\rangle  \text{ with transfer target type } \tau{F}=\lbrack 1,1,1,1\rbrack,  & \text{ if } v=1, \\
\langle 27,4\rangle \text{ with transfer target type } \tau{F}=\lbrack 11,2,2,2\rbrack, & \text{ if } v=2.
\end{cases}
\end{equation}
\end{corollary}

\begin{proof}
See
\cite[Thm. 4.2 and Thm. 4.3, p. 473]{AAI2001}.
\end{proof}


\begin{remark}
\label{rmk:TwoPrimeCond}
Note that the criteria in Theorem
\ref{thm:TwoPrimeCond}
are not rational (in terms of the prime divisors of the conductor).
A non-trivial primitive ambiguous principal ideal \((\alpha)\) must be determined
either directly or by means of a fundamental system of units of \(F\),
for instance with Magma
\cite{MAGMA2022}.
\end{remark}

\begin{example}
\label{exm:TwoPrimeCond}
The first occurrences for both situations, \(9 \mid c\) and \(\gcd(c,3) = 1\), are: \\
\(B(F)=\mathrm{N}_{F/\mathbb{Q}}(\alpha)=3^2\cdot 73\), \(n=2\), and thus \(\mathrm{G}_3^{(2)}{F_{c,\mu}}\simeq\langle 9,2\rangle\), for \(c=657=3^2\cdot 73\), \\
\(B(F)=\mathrm{N}_{F/\mathbb{Q}}(\alpha)=7\cdot 181\), \(n=2\), and thus \(\mathrm{G}_3^{(2)}{F_{c,\mu}}\simeq\langle 9,2\rangle\), for \(c=1\,267=7\cdot 181\), \\
\(B(F)=\mathrm{N}_{F/\mathbb{Q}}(\alpha)=271\), \(n=1\), and thus \(\mathrm{G}_3^{(2)}{F_{c,\mu}}\simeq\langle 27,4\rangle\), for \(c=2\,439=3^2\cdot 271\), \\
\(B(F)=\mathrm{N}_{F/\mathbb{Q}}(\alpha)=853\), \(n=1\), and thus \(\mathrm{G}_3^{(2)}{F_{c,\mu}}\simeq\langle 27,4\rangle\), for \(c=5\,971=7\cdot 853\), \\
according to Ayadi, Azizi, and Ismaili
\cite[p. 474]{AAI2001},
and confirmed by our own computations.
\end{example}


\begin{remark}
\label{rmk:Probability}
The statistical distribution of the three possible graphs for \(t=2\),
according to Theorem
\ref{thm:3ClassRank},
in the range \(c<10^5\) of conductors is the following. From Table
\ref{tbl:StatMult3}
we know there are \(3863\) doublets with conductors \(c=q_1q_2\) in this range.
They can be partitioned as follows.
Graph \(1\) with symbol
\(\lbrack q_1,q_2\rbrack_3=\lbrace q_1,q_2\rbrace\)
occurs \(1740\) times, graph \(2\) with symbol
\(\lbrack q_1,q_2\rbrack_3=\lbrace q_1\rightarrow q_2\rbrace\) or \(\lbrace q_2\rightarrow q_1\rbrace\)
occurs \(1715\) times, and graph \(3\) with symbol
\(\lbrack q_1,q_2\rbrack_3=\lbrace q_1\leftrightarrow q_2\rbrace\)
occurs \(408\) times.
As required, the sum is \(1740+1715+408=3863\).

This statistical result is in good accordance
with the \textit{probability tree} in Table
\ref{tbl:Probability},
because the proportion \(4:4:1\) for the three graphs
is nearly satisfied by the census
\(1740:1715:408\),
since \(4\cdot 408=1632\approx 1715\approx 1740\).

Another partition according to the \(3\)-class number
\(h=h_3{F^\ast}=\#\mathrm{Cl}_3{F^\ast}\)
of the \(3\)-\textit{genus field} \(F^\ast=F_{c,1}\cdot F_{c,2}\)
is also very illuminating.
We have \(h=1\) for \(3455=1740+1715\) hits,
which correspond to the union of graphs \(1\) and \(2\).
On the other hand, the \(408\) occurrences of graph \(3\)
split into \(210\) cases with \(h=3\),
corresponding to the group \(\mathrm{G}_3^{(2)}{F_{c,\mu}}\simeq\langle 9,2\rangle\),
\(142\) cases with \(h=9\),
corresponding to the group \(\mathrm{G}_3^{(2)}{F_{c,\mu}}\simeq\langle 27,4\rangle\),
and \(56\) cases with \(h\ge 27\),
corresponding to elevated \(3\)-class groups
\(\mathrm{Cl}_3{F_{c,\mu}}\simeq (9,3)\) and bigger.
The sum is \(210+142+56=408\), as required.
So the \textbf{overall proportion} of the \textit{abelian group} \(\langle 9,2\rangle\) is
\(\frac{210}{352}\approx 59.7\%\),
as opposed to the \textit{extra special group} \(\langle 27,4\rangle\) with
\(\frac{142}{352}\approx 40.3\%\).
See, however, section
\ref{ss:Statistics}
in the Conclusion
for an \textit{inverted population} in the special case \(9\mid c\).

\end{remark}
 

\renewcommand{\arraystretch}{1.2}

\begin{table}[ht]
\caption{Probability tree for the symbol \(\lbrack q_1,q_2\rbrack_3\)}
\label{tbl:Probability}
\begin{center}
\begin{tabular}{cc|cc||cl}
 \(P_1\) & \(\left(\frac{q_1}{q_2}\right)_3\) & \(P_2\) & \(\left(\frac{q_2}{q_1}\right)_3\) & \(P_1\cdot P_2\) & for \(\lbrack q_1,q_2\rbrack_3\) \\
              &                   &              &                   &               & \\
\hline
\hline
              &                   & \(2/3\)      & \(\zeta,\zeta^2\) & \quad \(4/9\) & for \(\lbrace q_1,q_2\rbrace\)\\
              &                   & \(\nearrow\) &                   &               & \\
\hline
 \(2/3\)      & \(\zeta,\zeta^2\) & \(\searrow\) &                   &               & \\
 \(\nearrow\) &                   & \(1/3\)      & \(1\)             & \quad \(2/9\) & for \(\lbrace q_1\leftarrow q_2\rbrace\) \\
\hline
 \(\searrow\) &                   & \(2/3\)      & \(\zeta,\zeta^2\) & \quad \(2/9\) & for \(\lbrace q_1\rightarrow q_2\rbrace\) \\
 \(1/3\)      & \(1\)             & \(\nearrow\) &                   &               & \\
\hline
              &                   & \(\searrow\) &                   &               & \\
              &                   & \(1/3\)      & \(1\)             & \quad \(1/9\) & for \(\lbrace q_1\leftrightarrow q_2\rbrace\) \\
\end{tabular}
\end{center}
\end{table}


\part{Current Research}
\label{pt:Current}

\chapter{Singular and Super-Singular Doublets}
\label{ch:Singular}

In his thesis
\cite{Ay1995},
Ayadi investigated doublets of cyclic cubic fields \((F_{c,\mu})_{1\le\mu\le 2}\)
with conductor \(c=q_1q_2\) divisible by two primes,
which are mutual cubic residues with respect to each other \(q_1\leftrightarrow q_2\),
and \(3\)-class group \(\mathrm{Cl}_3{F_{c,\mu}}=(3,3)\).
However, since the constitution of the \(3\)-genus fields \(F^\ast\)
of cyclic cubic fields \(F\) is recursive,
we shall see that the sub genus fields \(F_c^\ast\)
with partial conductors \(c\in\lbrace q_1q_2,q_2q_3,q_3q_1\rbrace\)
exert a considerable impact on the quartets of cyclic cubic fields
\((F_{f,\mu})_{1\le\mu\le 4}\)
with conductor \(f=q_1q_2q_3\) divisible by three primes
and \(3\)-genus field \(F_f^\ast\).
Consequently, the behavior of quartets cannot be analyzed satisfactorily,
if one restricts to doublets with elementary bicyclic \(3\)-class group.
For doublets \((F_{c,\mu})_{1\le\mu\le 2}\)
with conductor \(c=q_1q_2\)
and non-elementary bicyclic \(3\)-class group,
a further distinction arises from the \(3\)-valuation \(v\) of the
class number \(\#\mathrm{Cl}{F_c^\ast}\) of the \(3\)-genus field \(F_c^\ast\):


\begin{definition}
\label{dfn:Singular}
A doublet \((F_{c,\mu})_{1\le\mu\le 2}\) of cyclic cubic fields
is called
\begin{equation}
\label{eqn:Singular}
\begin{cases}
\textbf{regular}        & \text{ if } v\in\lbrace 0,1,2\rbrace, \\
\textbf{singular}       & \text{ if } v=3,                      \\
\textbf{super-singular} & \text{ if } v\in\lbrace 4,5,6\rbrace.
\end{cases}
\end{equation}
\end{definition}


\noindent
Table
\ref{tbl:Singular}
shows all \textit{singular} doublets \((F_{c,\mu})_{1\le\mu\le 2}\)
of cyclic cubic fields with conductors \(c<10^5\).
Both components have \(3\)-class group
\(\mathrm{Cl}_3{F_{c,\mu}}=(3,9)\),
the \(3\)-genus field \(F^\ast\) has \(3\)-class group
\(\mathrm{Cl}_3{F^\ast}=(3,3,3)\),
and thus valuation \(v=3\).
The second \(3\)-class group is always
\(G_3^{(2)}{F_{c,\mu}}\simeq\langle 81,3\rangle\),
for both fields \(1\le\mu\le 2\).
This group has punctured Artin pattern \(\mathrm{AP}=(\varkappa,\tau)\)
with transfer kernel type \(\varkappa=(000;0)\), \(\mathrm{a}.1\), and
transfer target type \(\tau=\lbrack (21)^3;1^3\rbrack\). See
\cite{Ma2022}.

Within the frame \(c<10^5\) of our present computations,
there are \(13\) singular doublets 
but no super-singular doublets
with conductor \(c\) divisible by \(9\).


\renewcommand{\arraystretch}{1.1}

\begin{table}[ht]
\caption{Thirty-Seven Singular Doublets with \(v=3\)}
\label{tbl:Singular}
\begin{center}
{\small
\begin{tabular}{|rrc||rrc|}
\hline
 No.    & \(c\)       & Factors                    & No.    & \(c\)       & Factors                    \\
\hline
  \(1\) &  \(4\,711\) & \(7\leftrightarrow 673\)   & \(20\) & \(58\,329\) & \(9\leftrightarrow 6481\)  \\
  \(2\) & \(11\,167\) & \(13\leftrightarrow 859\)  & \(21\) & \(62\,257\) & \(13\leftrightarrow 4789\) \\
  \(3\) & \(12\,439\) & \(7\leftrightarrow 1777\)  & \(22\) & \(64\,971\) & \(9\leftrightarrow 7219\)  \\
  \(4\) & \(16\,177\) & \(7\leftrightarrow 2311\)  & \(23\) & \(65\,383\) & \(151\leftrightarrow 433\) \\
  \(5\) & \(17\,593\) & \(73\leftrightarrow 241\)  & \(24\) & \(66\,829\) & \(7\leftrightarrow 9547\)  \\
  \(6\) & \(20\,421\) & \(9\leftrightarrow 2269\)  & \(25\) & \(69\,183\) & \(9\leftrightarrow 7687\)  \\
  \(7\) & \(25\,963\) & \(7\leftrightarrow 3709\)  & \(26\) & \(71\,611\) & \(19\leftrightarrow 3769\) \\
  \(8\) & \(27\,571\) & \(79\leftrightarrow 349\)  & \(27\) & \(72\,099\) & \(9\leftrightarrow 8011\)  \\
  \(9\) & \(32\,689\) & \(97\leftrightarrow 337\)  & \(28\) & \(73\,873\) & \(31\leftrightarrow 2383\) \\
 \(10\) & \(35\,163\) & \(9\leftrightarrow 3907\)  & \(29\) & \(77\,281\) & \(109\leftrightarrow 709\) \\
 \(11\) & \(37\,933\) & \(7\leftrightarrow 5419\)  & \(30\) & \(78\,093\) & \(9\leftrightarrow 8677\)  \\
 \(12\) & \(40\,573\) & \(13\leftrightarrow 3121\) & \(31\) & \(81\,009\) & \(9\leftrightarrow 9001\)  \\
 \(13\) & \(40\,873\) & \(7\leftrightarrow 5839\)  & \(32\) & \(89\,109\) & \(9\leftrightarrow 9901\)  \\
 \(14\) & \(43\,081\) & \(67\leftrightarrow 643\)  & \(33\) & \(89\,863\) & \(73\leftrightarrow 1231\) \\
 \(15\) & \(44\,397\) & \(9\leftrightarrow 4933\)  & \(34\) & \(94\,357\) & \(157\leftrightarrow 601\) \\
 \(16\) & \(49\,743\) & \(9\leftrightarrow 5527\)  & \(35\) & \(95\,913\) & \(9\leftrightarrow 10657\) \\
 \(17\) & \(51\,847\) & \(139\leftrightarrow 373\) & \(36\) & \(96\,709\) & \(97\leftrightarrow 997\)  \\
 \(18\) & \(55\,951\) & \(7\leftrightarrow 7993\)  & \(37\) & \(96\,817\) & \(7\leftrightarrow 13831\) \\
 \(19\) & \(56\,223\) & \(9\leftrightarrow 6247\)  & \(\) & \(\) & \(\) \\
\hline
\end{tabular}
}
\end{center}
\end{table}

\newpage

\noindent
Table
\ref{tbl:SuperSingular}
shows all \textit{super-singular} doublets \((F_{c,\mu})_{1\le\mu\le 2}\)
of cyclic cubic fields with conductors \(c<10^5\).
On 23 December \(2001\), Aissa Derhem
\cite{Dh2002}
asked me,
if there exist such fields with \(3\)-class numbers \(h_3{F_{c,\mu}}\in\lbrace 81,243\rbrace\).
On Christmas Day \(2001\), I answered that in the \(1982\) tables of Ennola and Turunen,
there only appears \(c=5\,383\) with \(h_3{F_{c,\mu}}=27\),
and I announced that I am going to launch an extensive computation
of regulators and class numbers of cyclic cubic fields
by means of Voronoi's algorithm
\cite{Vo1896}
and the Euler product method for the analytical class number formula.
I needed three months until I arrived at \(h_3{F_{c,\mu}}\in\lbrace 81,243\rbrace\)
for \(c\in\lbrace 36\,667,41\,977,42\,127,42\,991\rbrace\) on 2 April \(2002\).
See the details in \\
\textbf{(1)} \texttt{http://www.algebra.at/AissaDan17.htm} \quad and\\
\textbf{(2)} \texttt{http://www.algebra.at/KarimAissaDan23.htm}. \\
Later I saw that \(c=36\,667\) was known to Georges Gras in \(1973\) already
\cite[Exm. VI.7, pp. 36--38]{Gr1973}.


\renewcommand{\arraystretch}{1.1}

\begin{table}[ht]
\caption{Nineteen Super-Singular Doublets with \(v\in\lbrace 4,5,6\rbrace\)}
\label{tbl:SuperSingular}
\begin{center}
{\footnotesize
\begin{tabular}{|rrc||rc|c|c|}
\hline
 No.    & \(c\) & Factors & \(\mathrm{Cl}_3{F^\ast}\) & \(v\) & \(\mathrm{Cl}_3{F_{c,\mu}}\) & \(G_3^{(2)}{F_{c,\mu}}\)    \\
\hline
  \(1\) &  \(5\,383\) & \(7\leftrightarrow 769\)   &   \((3,3,9)\) & \(4\) & \((3,9),(3,9)\) & \(\langle 243,14\rangle^2\) \\
  \(2\) & \(12\,403\) & \(79\leftrightarrow 157\)  &   \((3,3,9)\) & \(4\) & \((3,9),(3,9)\) & \(\langle 243,14\rangle^2\) \\ 
  \(3\) & \(21\,763\) & \(7\leftrightarrow 3109\)  & \((3,3,3,3)\) & \(4\) & \((3,9),(3,9)\) & \(\langle 243,13\rangle,\langle 729,12\rangle\) \\
  \(4\) & \(28\,177\) & \(19\leftrightarrow 1483\) &   \((3,3,9)\) & \(4\) & \((3,9),(3,9)\) & \(\langle 243,14\rangle^2\) \\
  \(5\) & \(32\,311\) & \(79\leftrightarrow 409\)  &   \((3,3,9)\) & \(4\) & \((3,9),(3,9)\) & \(\langle 243,14\rangle^2\) \\
  \(6\) & \(36\,667\) & \(37\leftrightarrow 991\)  &  \((3,9,27)\) & \(6\) & \((3,9),(9,27)\)& \(\langle 243,14\rangle,\ast\ast\) \\
  \(7\) & \(38\,503\) & \(139\leftrightarrow 277\) &   \((3,3,9)\) & \(4\) & \((3,9),(3,9)\) & \(\langle 243,14\rangle,\langle 729,17\rangle\) \\
  \(8\) & \(41\,977\) & \(13\leftrightarrow 3229\) &   \((3,3,9)\) & \(4\) & \((3,9),(9,9)\) & \(\langle 243,15\rangle,\ast\) \\
  \(9\) & \(42\,127\) & \(103\leftrightarrow 409\) & \((3,3,9,9)\) & \(6\) & \((9,9),(9,27)\) & \(\ast,\ast\ast,\) \\
 \(10\) & \(42\,991\) & \(13\leftrightarrow 3307\) &   \((3,3,9)\) & \(4\) & \((3,9),(9,9)\) & \(\langle 243,15\rangle,\ast\) \\
 \(11\) & \(49\,849\) & \(79\leftrightarrow 631\)  &   \((3,3,9)\) & \(4\) & \((3,9),(3,9)\) & \(\langle 243,14\rangle,\langle 729,17\rangle\) \\
 \(12\) & \(55\,657\) & \(7\leftrightarrow 7951\)  &   \((3,3,9)\) & \(4\) & \((3,9),(3,9)\) & \(\langle 243,14\rangle,\langle 729,17\rangle\) \\
 \(13\) & \(57\,811\) & \(13\leftrightarrow 4447\) & \((3,3,3,3)\) & \(4\) & \((3,9),(3,9)\) & \(\langle 243,13\rangle^2\) \\
 \(14\) & \(59\,803\) & \(79\leftrightarrow 757\)  &   \((3,3,9)\) & \(4\) & \((3,9),(3,9)\) & \(\langle 243,14\rangle^2\) \\
 \(15\) & \(59\,911\) & \(181\leftrightarrow 331\) & \((3,3,3,3)\) & \(4\) & \((3,9),(3,9)\) & \(\langle 243,13\rangle^2\) \\
 \(16\) & \(68\,857\) & \(37\leftrightarrow 1861\) &   \((3,9,9)\) & \(5\) & \((3,9),(9,9)\) & \(\langle 243,14\rangle,\ast\) \\
 \(17\) & \(75\,859\) & \(7\leftrightarrow 10837\) & \((3,3,3,3)\) & \(4\) & \((3,9),(3,9)\) & \(\langle 729,12\rangle^2\) \\
 \(18\) & \(84\,103\) & \(31\leftrightarrow 2713\) &   \((3,3,9)\) & \(4\) & \((3,9),(3,9)\) & \(\langle 243,14\rangle^2\) \\
 \(19\) & \(97\,249\) & \(79\leftrightarrow 1231\) &   \((3,9,9)\) & \(5\) & \((3,9),(9,9)\) & \(\langle 243,14\rangle,\ast\) \\
\hline
\end{tabular}
}
\end{center}
\end{table}


\chapter{Recent Results}
\label{ch:Recent}

\section{Finite 3-Groups of Type (3,3)}
\label{s:Groups3x3}

\noindent
In the following tables,
we list those invariants of finite \(3\)-groups
with elementary bicyclic commutator quotient \((3,3)\)
which qualify metabelian groups \(\mathfrak{M}\) as second \(3\)-class groups \(\mathrm{Gal}(F_3^{(2)}/F)\)
and non-metabelian groups \(G\) as \(3\)-class tower groups \(\mathrm{Gal}(F_3^{(\infty)}/F)\)
of cyclic cubic number fields \(F\).
The process of searching for suitable groups in descendant trees
by means of the strategy of pattern recognition
is governed by the \textit{Artin transfer pattern}
\(\mathrm{AP}=(\tau,\varkappa)\),
where \(\tau=\tau_1\), resp. \(\varkappa=\varkappa_1\), denotes the first layer
of the transfer target type (TTT), resp. transfer kernel type (TKT).
Additionally, we give the top layer of the TTT, which consists of the abelian quotient invariants
of the commutator subgroup \(G^\prime\),
corresponding to the \(3\)-class group of the first Hilbert \(3\)-class field \(F_3^{(1)}\).
The \textit{nuclear rank} \(\nu\) is responsible for the search complexity.
The \(p\)-multiplicator rank \(\mu\) of a group \(G\) is precisely its \textit{relation rank}
\(d_2(G)=\dim_{\mathbb{F}_3}\mathrm{H}^2(G,\mathbb{F}_3)\),
which decides whether \(G\) is admissible as \(\mathrm{Gal}(F_3^{(\infty)}/F)\),
according to the Shafarevich Theorem
\cite{Sh1964}.
In the case of cyclic cubic fields \(F\), it is limited by the \textit{Shafarevich bound}
\(\mu\le\varrho+r+\theta\),
where
\(\varrho=d_1(G)=\dim_{\mathbb{F}_3}\mathrm{H}^1(G,\mathbb{F}_3)\)
denotes the \textit{generator rank} of \(G\),
which coincides with the \(3\)-class rank \(\varrho\) of \(F\),
\(r=r_1+r_2-1=2\) is the torsion free Dirichlet unit rank of the field \(F\) with signature \((r_1,r_2)=(3,0)\),
and \(\theta=0\) indicates the absence of a (complex) primitive third root of unity in the totally real field \(F\).
Finally, \(\pi(\mathfrak{M})=\mathfrak{M}/\gamma_c\mathfrak{M}\) denotes the parent of \(\mathfrak{M}\).


In Table
\ref{tbl:Metabelian33},
we begin with metabelian groups \(\mathfrak{M}\) of generator rank \(d_1(\mathfrak{M})=2\).
The Shafarevich bound
\cite{Ma2016}
is given by \(\mu\le\varrho+r+\theta=2+2+0=4\).

\renewcommand{\arraystretch}{1.1}

\begin{table}[ht]
\caption{Invariants of Metabelian \(3\)-Groups \(\mathfrak{M}\) with \(\mathfrak{M}/\mathfrak{M}^\prime\simeq (3,3)\)}
\label{tbl:Metabelian33}
\begin{center}
{\scriptsize
\begin{tabular}{|c|r|l|c|c|c|r|r|c|}
\hline
\(\mathfrak{M}\)                     & cc    & Type & \(\varkappa\) & \(\tau\)           & \(\tau_2\) & \(\nu\) & \(\mu\) & \(\pi(\mathfrak{M})\) \\
\hline
\(\langle 9,2\rangle\)               & \(1\) & a.1  & \((0000)\)    & \((1)^4\)          & \(0\)      &   \(3\) &   \(3\) & \\
\(\langle 27,4\rangle\)              & \(1\) & A.1  & \((1111)\)    & \(1^2,(2)^3\)      & \(1\)      &   \(0\) &   \(2\) & \(\langle 9,2\rangle\) \\
\(\langle 81,7\rangle\)              & \(1\) & a.3  & \((2000)\)    & \(1^3,(1^2)^3\)    & \(1^2\)    &   \(0\) &   \(3\) & \(\langle 27,3\rangle\) \\
\(\langle 81,8\rangle\)              & \(1\) & a.3  & \((2000)\)    & \(21,(1^2)^3\)     & \(1^2\)    &   \(0\) &   \(3\) & \(\langle 27,3\rangle\) \\
\(\langle 81,10\rangle\)             & \(1\) & a.2  & \((1000)\)    & \(21,(1^2)^3\)     & \(1^2\)    &   \(0\) &   \(3\) & \(\langle 27,3\rangle\) \\
\(\langle 81,9\rangle\)              & \(1\) & a.1  & \((0000)\)    & \(21,(1^2)^3\)     & \(1^2\)    &   \(1\) &   \(4\) & \(\langle 27,3\rangle\) \\
\(\langle 243,28\ldots 30\rangle\)   & \(1\) & a.1  & \((0000)\)    & \(21,(1^2)^3\)     & \(21\)     &   \(0\) &   \(3\) & \(\langle 81,9\rangle\) \\
\(\langle 243,25\rangle\)            & \(1\) & a.3  & \((2000)\)    & \(2^2,(1^2)^3\)    & \(21\)     &   \(0\) &   \(3\) & \(\langle 81,9\rangle\) \\
\(\langle 243,27\rangle\)            & \(1\) & a.2  & \((1000)\)    & \(2^2,(1^2)^3\)    & \(21\)     &   \(0\) &   \(3\) & \(\langle 81,9\rangle\) \\
\hline
\(\langle 243,8\rangle\)             & \(2\) & c.21 & \((0231)\)    & \((21)^4\)         & \(1^3\)    &   \(1\) &   \(3\) & \(\langle 27,3\rangle\) \\
\(\langle 729,54\rangle=U\)          & \(2\) & c.21 & \((0231)\)    & \(22,(21)^3\)      & \(21^2\)   &   \(2\) &   \(4\) & \(\langle 243,8\rangle\) \\
\(\langle 243,3\rangle\)             & \(2\) & b.10 & \((0043)\)    & \((21)^2,(1^3)^2\) & \(1^3\)    &   \(2\) &   \(4\) & \(\langle 27,3\rangle\) \\
\(\langle 729,40\rangle=B\)          & \(2\) & b.10 & \((0043)\)    & \(2^2,21,(1^3)^2\) & \(21^2\)   &   \(2\) &   \(5\) & \(\langle 243,3\rangle\) \\
\(\langle 729,34\rangle=H\)          & \(2\) & b.10 & \((0043)\)    & \((21)^2,(1^3)^2\) & \(1^4\)    &   \(2\) &   \(5\) & \(\langle 243,3\rangle\) \\
\(\langle 729,35\rangle=I\)          & \(2\) & b.10 & \((0043)\)    & \((21)^2,(1^3)^2\) & \(1^4\)    &   \(1\) &   \(4\) & \(\langle 243,3\rangle\) \\
\(\langle 729,37\rangle=A\)          & \(2\) & b.10 & \((0043)\)    & \((21)^2,(1^3)^2\) & \(21^2\)   &   \(2\) &   \(5\) & \(\langle 243,3\rangle\) \\
\(\langle 729,38\rangle=C\)          & \(2\) & b.10 & \((0043)\)    & \((21)^2,(1^3)^2\) & \(21^2\)   &   \(1\) &   \(4\) & \(\langle 243,3\rangle\) \\
\(\langle 729,41\rangle=D\)          & \(2\) & d.19 & \((4043)\)    & \(32,21,(1^3)^2\)  & \(21^2\)   &   \(1\) &   \(4\) & \(\langle 243,3\rangle\) \\
\(\langle 2187,301\vert 305\rangle\) & \(2\) & G.16 & \((4231)\)    & \(32,(21)^3\)      & \(2^21\)   &   \(1\) &   \(4\) & \(\langle 729,54\rangle\) \\
\hline
\end{tabular}
}
\end{center}
\end{table}


For the metabelian groups \(\mathfrak{M}\) with non-trivial cover \(\mathrm{cov}(\mathfrak{M})\),
we need non-metabelian groups \(G\) in the cover, which are given in Table
\ref{tbl:NonMetabelian33},
where we begin with groups \(G\) of generator rank \(d_1(G)=2\).
For \(d_1(G)=3\), see \S\
\ref{s:Groups3x3x3}.

\renewcommand{\arraystretch}{1.1}

\begin{table}[ht]
\caption{Invariants of Non-Metabelian \(3\)-Groups \(G\) with \(G/G^\prime\simeq (3,3)\)}
\label{tbl:NonMetabelian33}
\begin{center}
{\scriptsize
\begin{tabular}{|c|r|l|c|c|c|r|r|c|}
\hline
\(G\)                                 & cc    & Type & \(\varkappa\) & \(\tau\)           & \(\tau_2\) & \(\nu\) & \(\mu\) & \(G/G^{\prime\prime}\) \\
\hline
\(\langle 2187,263\ldots 265\rangle\) & \(2\) & d.19 & \((4043)\)    & \(32,21,(1^3)^2\)  & \(21^2\)   &   \(0\) &   \(3\) & \(\langle 729,41\rangle\) \\
\(\langle 2187,307\vert 308\rangle\)  & \(2\) & c.21 & \((0231)\)    & \(32,(21)^3\)      & \(21^2\)   &   \(0\) &   \(3\) & \(\langle 729,54\rangle\) \\
\hline
\(\langle 6561,619\vert 623\rangle\)  & \(3\) & G.16 & \((4231)\)    & \(32,(21)^3\)      & \(2^21\)   &   \(1\) &   \(3\) & \(\langle 2187,301\vert 305\rangle\) \\
\hline
\end{tabular}
}
\end{center}
\end{table}

\newpage

\begin{theorem}
\label{thm:ClassGroup33}
Let \(F\) be a cyclic cubic number field
with elementary bicyclic \(3\)-class group \(\mathrm{Cl}_3{F}\simeq (3,3)\).
Denote by \(\mathfrak{M}=\mathrm{Gal}(F_3^{(2)}/F)\) the second \(3\)-class group of \(F\),
and by \(G=\mathrm{Gal}(F_3^{(\infty)}/F)\) the \(3\)-class field tower group of \(F\).
Then, the Artin pattern \((\tau,\varkappa)\) of \(F\)
identifies the groups \(\mathfrak{M}\) and \(G\),
and determines the length \(\ell_3{F}\) of the \(3\)-class field tower of \(F\),
according to the following \textbf{deterministic laws}.
\begin{enumerate}
\item
If \(\tau=\lbrack (1)^4\rbrack\), \(\varkappa=(0000)\) (type \(\mathrm{a}.1\)), then \(G\simeq\langle 9,2\rangle\) and \(\ell_3{F}=1\).
\item
If \(\tau\sim\lbrack 1^2,(2)^3\rbrack\), \(\varkappa\sim (1111)\) (type \(\mathrm{A}.1\)), then \(G\simeq\langle 27,4\rangle\).
\item
If \(\tau\sim\lbrack 1^3,(1^2)^3\rbrack\), \(\varkappa\sim (2000)\) (type \(\mathrm{a}.3\)), then \(G\simeq\langle 81,7\rangle\).
\item
If \(\tau\sim\lbrack 21,(1^2)^3\rbrack\), \(\varkappa\sim (2000)\) (type \(\mathrm{a}.3\)), then \(G\simeq\langle 81,8\rangle\).
\item
If \(\tau\sim\lbrack 21,(1^2)^3\rbrack\), \(\varkappa\sim (1000)\) (type \(\mathrm{a}.2\)), then \(G\simeq\langle 81,10\rangle\).
\item
If \(\tau\sim\lbrack 2^2,(1^2)^3\rbrack\), \(\varkappa\sim (2000)\) (type \(\mathrm{a}.3\)), then \(G\simeq\langle 243,25\rangle\).
\item
If \(\tau\sim\lbrack 2^2,(1^2)^3\rbrack\), \(\varkappa\sim (1000)\) (type \(\mathrm{a}.2\)), then \(G\simeq\langle 243,27\rangle\).
\end{enumerate} 
Except for the abelian tower in item 1, the tower is metabelian with \(\ell_3{F}=2\).
\end{theorem}

\begin{proof}
Generally, a cyclic cubic field \(F\)
has signature \((r_1,r_2)=(3,0)\), torsion free unit rank \(r=r_1+r_2-1=2\),
does not contain primitive third roots of unity,
and thus possesses the maximal admissible relation rank
\(d_2\le d_1+r=4\) for the group \(G\),
when its \(3\)-class rank, i.e. the generator rank of \(G\), is \(d_1=\varrho=2\).
Consequently, \(\ell_3{F}\ge 3\) in the case of \(d_2{\mathfrak{M}}\ge 5\).

For item 1, we have
\(\mathfrak{M}=\mathrm{Gal}(F_3^{(2)}/F)\simeq\langle 9,2\rangle\simeq (3,3)\simeq\mathrm{Cl}_3{F}\simeq\mathrm{Gal}(F_3^{(1)}/F)\),
whence \(\ell_3{F}=1\).

For item 2 to item 7,
the group \(\mathfrak{M}\) is of maximal class (coclass \(\mathrm{cc}(\mathfrak{M})=1\)),
and thus coincides with \(G\), whence \(\ell_3{F}=2\).

In each case, the Artin pattern \((\tau,\varkappa)\) identifies \(\mathfrak{M}=G\) uniquely,
and the relation ranks are
\(d_2\langle 9,2\rangle=3\), 
\(d_2\langle 27,4\rangle=2\), 
\(d_2\langle 81,7\rangle=3\), 
\(d_2\langle 81,8\rangle=3\), 
\(d_2\langle 81,10\rangle=3\), 
\(d_2\langle 243,25\rangle=3\),
\(d_2\langle 243,27\rangle=3\),
each of them less than \(4\).
\end{proof}


\begin{corollary}
\label{cor:ClassGroup33}
Under the assumptions of Theorem
\ref{thm:ClassGroup33},
the abelian type invariants \(\tau_2\) of the \(3\)-class group
\(\mathrm{Cl}_3{F_3^{(1)}}\)
of the first Hilbert \(3\)-class field of \(F\)
are required for the unambiguous identification
of the groups \(G\) respectively \(\mathfrak{M}\). \\
If \(\tau\sim\lbrack 21,(1^2)^3\rbrack\), \(\varkappa\sim (0000)\), \(\mathrm{a}.1\), then
\(G\simeq
\begin{cases}
\langle 81,9\rangle & \text{ for } \tau_2=\lbrack 1^2\rbrack, \\
\langle 243,28\ldots 30\rangle & \text{ for } \tau_2\sim\lbrack 21\rbrack.
\end{cases}\) \\
If \(\tau\sim\lbrack (21)^2,(1^3)^2\rbrack\), \(\varkappa\sim (0043)\), \(\mathrm{b}.10\), then
\(\mathfrak{M}\simeq
\begin{cases}
\langle 729,34\ldots 36\rangle & \text{ for } \tau_2=\lbrack 1^4\rbrack, \\
\langle 729,37\ldots 39\rangle & \text{ for } \tau_2\sim\lbrack 21^2\rbrack.
\end{cases}\) 
\end{corollary}

\begin{proof}
The Artin pattern \((\tau,\varkappa)\) of \(F\) alone is not able
to identify the groups \(\mathfrak{M}\) and \(G\) unambiguously.
Ascione uses the notation
\(\langle 729,34\rangle=H\),
\(\langle 729,35\rangle=I\),
\(\langle 729,37\rangle=A\),
\(\langle 729,38\rangle=C\).
\end{proof}

\noindent
Since many new groups will arise for \(t=3\), we insert a section on \(p\)-group theory.


\subsection{Descendant Trees of Finite 3-Groups}
\label{ss:Trees}

\noindent
Basic definitions, facts, and notation concerning \textit{descendant trees} of finite \(p\)-groups
are summarized briefly in
\cite[\S\ 2, pp. 410--411]{Ma2013}.
They are discussed thoroughly in the broadest detail in the initial sections of
\cite{Ma2015a,Ma2018}.
Trees are crucial for recent progress in the theory of \(p\)-class field towers
\cite{Ma2016a,Ma2016c,Ma2016d},
in particular in order to describe the mutual location of
\(\mathrm{G}_3^{(2)}{K}\) and \(\mathrm{G}_3^{(\infty)}{K}\)
for number fields \(K\).

\newpage

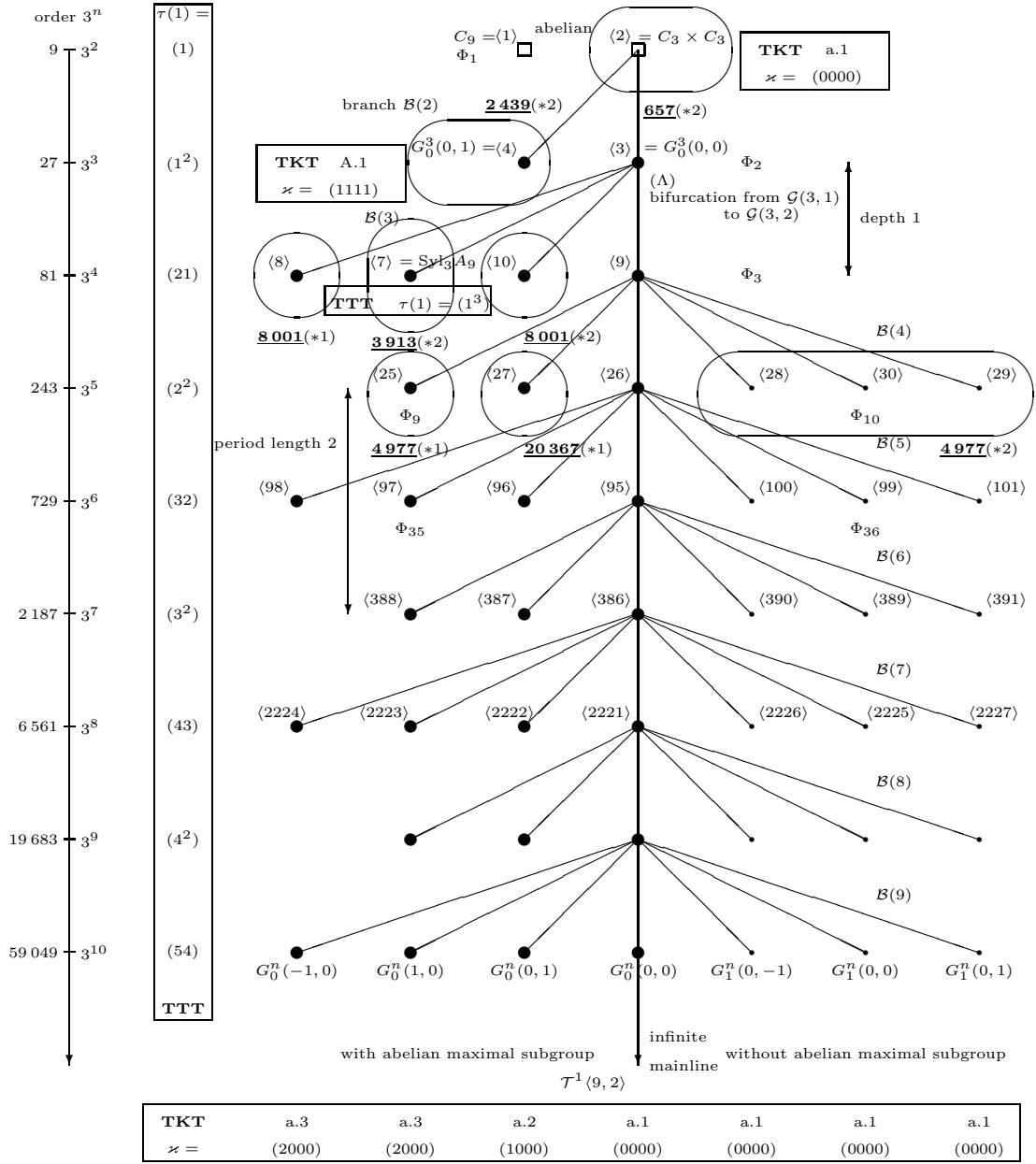
\begin{figure}[ht]
\caption{Distribution of Conductors for \(\mathrm{G}_3^{(2)}{F}\) on the Coclass Tree \(\mathcal{T}^1\langle 9,2\rangle\)}
\label{fig:MinDiscTyp33TreeCc1}

{\tiny


\setlength{\unitlength}{0.8cm}
\begin{picture}(15,18)(-11,-17)

\put(-10,0.5){\makebox(0,0)[cb]{order \(3^n\)}}

\put(-10,0){\line(0,-1){16}}
\multiput(-10.1,0)(0,-2){9}{\line(1,0){0.2}}

\put(-10.2,0){\makebox(0,0)[rc]{\(9\)}}
\put(-9.8,0){\makebox(0,0)[lc]{\(3^2\)}}
\put(-10.2,-2){\makebox(0,0)[rc]{\(27\)}}
\put(-9.8,-2){\makebox(0,0)[lc]{\(3^3\)}}
\put(-10.2,-4){\makebox(0,0)[rc]{\(81\)}}
\put(-9.8,-4){\makebox(0,0)[lc]{\(3^4\)}}
\put(-10.2,-6){\makebox(0,0)[rc]{\(243\)}}
\put(-9.8,-6){\makebox(0,0)[lc]{\(3^5\)}}
\put(-10.2,-8){\makebox(0,0)[rc]{\(729\)}}
\put(-9.8,-8){\makebox(0,0)[lc]{\(3^6\)}}
\put(-10.2,-10){\makebox(0,0)[rc]{\(2\,187\)}}
\put(-9.8,-10){\makebox(0,0)[lc]{\(3^7\)}}
\put(-10.2,-12){\makebox(0,0)[rc]{\(6\,561\)}}
\put(-9.8,-12){\makebox(0,0)[lc]{\(3^8\)}}
\put(-10.2,-14){\makebox(0,0)[rc]{\(19\,683\)}}
\put(-9.8,-14){\makebox(0,0)[lc]{\(3^9\)}}
\put(-10.2,-16){\makebox(0,0)[rc]{\(59\,049\)}}
\put(-9.8,-16){\makebox(0,0)[lc]{\(3^{10}\)}}

\put(-10,-16){\vector(0,-1){2}}

\put(-8,0.5){\makebox(0,0)[cb]{\(\tau(1)=\)}}

\put(-8,0){\makebox(0,0)[cc]{\((1)\)}}
\put(-8,-2){\makebox(0,0)[cc]{\((1^2)\)}}
\put(-8,-4){\makebox(0,0)[cc]{\((21)\)}}
\put(-8,-6){\makebox(0,0)[cc]{\((2^2)\)}}
\put(-8,-8){\makebox(0,0)[cc]{\((32)\)}}
\put(-8,-10){\makebox(0,0)[cc]{\((3^2)\)}}
\put(-8,-12){\makebox(0,0)[cc]{\((43)\)}}
\put(-8,-14){\makebox(0,0)[cc]{\((4^2)\)}}
\put(-8,-16){\makebox(0,0)[cc]{\((54)\)}}

\put(-8,-17){\makebox(0,0)[cc]{\textbf{TTT}}}
\put(-8.5,-17.2){\framebox(1,18){}}

\put(3.7,-3){\vector(0,1){1}}
\put(3.9,-3){\makebox(0,0)[lc]{depth \(1\)}}
\put(3.7,-3){\vector(0,-1){1}}

\put(-5.1,-8){\vector(0,1){2}}
\put(-5.3,-7){\makebox(0,0)[rc]{period length \(2\)}}
\put(-5.1,-8){\vector(0,-1){2}}

\put(-0.1,-0.1){\framebox(0.2,0.2){}}
\put(-2.1,-0.1){\framebox(0.2,0.2){}}

\multiput(0,-2)(0,-2){8}{\circle*{0.2}}
\multiput(-2,-2)(0,-2){8}{\circle*{0.2}}
\multiput(-4,-4)(0,-2){7}{\circle*{0.2}}
\multiput(-6,-4)(0,-4){4}{\circle*{0.2}}

\multiput(2,-6)(0,-2){6}{\circle*{0.1}}
\multiput(4,-6)(0,-2){6}{\circle*{0.1}}
\multiput(6,-6)(0,-2){6}{\circle*{0.1}}

\multiput(0,0)(0,-2){8}{\line(0,-1){2}}
\multiput(0,0)(0,-2){8}{\line(-1,-1){2}}
\multiput(0,-2)(0,-2){7}{\line(-2,-1){4}}
\multiput(0,-2)(0,-4){4}{\line(-3,-1){6}}
\multiput(0,-4)(0,-2){6}{\line(1,-1){2}}
\multiput(0,-4)(0,-2){6}{\line(2,-1){4}}
\multiput(0,-4)(0,-2){6}{\line(3,-1){6}}

\put(0,-16){\vector(0,-1){2}}
\put(0.2,-17.5){\makebox(0,0)[lc]{infinite}}
\put(0.2,-18){\makebox(0,0)[lc]{mainline}}
\put(-0.2,-18.3){\makebox(0,0)[rc]{\(\mathcal{T}^1\langle 9,2\rangle\)}}

\put(0,0.1){\makebox(0,0)[lb]{\(=C_3\times C_3\)}}
\put(-2.5,0.1){\makebox(0,0)[rb]{\(C_9=\)}}
\put(-0.8,0.3){\makebox(0,0)[rb]{abelian}}
\put(0.1,-1.9){\makebox(0,0)[lb]{\(=G^3_0(0,0)\)}}
\put(-2.5,-1.9){\makebox(0,0)[rb]{\(G^3_0(0,1)=\)}}
\put(-4.2,-3.9){\makebox(0,0)[lb]{\(=\mathrm{Syl}_3A_9\)}}
\put(-4,-4.5){\makebox(0,0)[cc]{\textbf{TTT} \quad \(\tau(1)=(1^3)\)}}
\put(-5.5,-4.7){\framebox(2.9,0.5){}}

\put(0.2,-2.2){\makebox(0,0)[lt]{\((\Lambda)\)}}
\put(0.2,-2.5){\makebox(0,0)[lt]{bifurcation from \(\mathcal{G}(3,1)\)}}
\put(1.5,-2.8){\makebox(0,0)[lt]{to \(\mathcal{G}(3,2)\)}}

\put(-3,0){\makebox(0,0)[ct]{\(\Phi_1\)}}
\put(-2.1,0.1){\makebox(0,0)[rb]{\(\langle 1\rangle\)}}
\put(-0.1,0.1){\makebox(0,0)[rb]{\(\langle 2\rangle\)}}

\put(-3.5,-1){\makebox(0,0)[rc]{branch \(\mathcal{B}(2)\)}}
\put(2,-2){\makebox(0,0)[cc]{\(\Phi_2\)}}
\put(-2.1,-1.9){\makebox(0,0)[rb]{\(\langle 4\rangle\)}}
\put(-0.1,-1.9){\makebox(0,0)[rb]{\(\langle 3\rangle\)}}

\put(-4.5,-3){\makebox(0,0)[cc]{\(\mathcal{B}(3)\)}}
\put(2,-4){\makebox(0,0)[cc]{\(\Phi_3\)}}
\put(-6.1,-3.9){\makebox(0,0)[rb]{\(\langle 8\rangle\)}}
\put(-4.3,-3.9){\makebox(0,0)[rb]{\(\langle 7\rangle\)}}
\put(-2.1,-3.9){\makebox(0,0)[rb]{\(\langle 10\rangle\)}}
\put(-0.1,-3.9){\makebox(0,0)[rb]{\(\langle 9\rangle\)}}

\put(-4,-6.5){\makebox(0,0)[cc]{\(\Phi_9\)}}
\put(-4.1,-5.9){\makebox(0,0)[rb]{\(\langle 25\rangle\)}}
\put(-2.1,-5.9){\makebox(0,0)[rb]{\(\langle 27\rangle\)}}
\put(-0.1,-5.9){\makebox(0,0)[rb]{\(\langle 26\rangle\)}}

\put(4.5,-5){\makebox(0,0)[cc]{\(\mathcal{B}(4)\)}}
\put(4,-6.5){\makebox(0,0)[cc]{\(\Phi_{10}\)}}
\put(2.1,-5.9){\makebox(0,0)[lb]{\(\langle 28\rangle\)}}
\put(4.1,-5.9){\makebox(0,0)[lb]{\(\langle 30\rangle\)}}
\put(6.1,-5.9){\makebox(0,0)[lb]{\(\langle 29\rangle\)}}

\put(-4,-8.5){\makebox(0,0)[cc]{\(\Phi_{35}\)}}
\put(-6.1,-7.9){\makebox(0,0)[rb]{\(\langle 98\rangle\)}}
\put(-4.1,-7.9){\makebox(0,0)[rb]{\(\langle 97\rangle\)}}
\put(-2.1,-7.9){\makebox(0,0)[rb]{\(\langle 96\rangle\)}}
\put(-0.1,-7.9){\makebox(0,0)[rb]{\(\langle 95\rangle\)}}

\put(4.5,-7){\makebox(0,0)[cc]{\(\mathcal{B}(5)\)}}
\put(4,-8.5){\makebox(0,0)[cc]{\(\Phi_{36}\)}}
\put(2.1,-7.9){\makebox(0,0)[lb]{\(\langle 100\rangle\)}}
\put(4.1,-7.9){\makebox(0,0)[lb]{\(\langle 99\rangle\)}}
\put(6.1,-7.9){\makebox(0,0)[lb]{\(\langle 101\rangle\)}}

\put(-4.1,-9.9){\makebox(0,0)[rb]{\(\langle 388\rangle\)}}
\put(-2.1,-9.9){\makebox(0,0)[rb]{\(\langle 387\rangle\)}}
\put(-0.1,-9.9){\makebox(0,0)[rb]{\(\langle 386\rangle\)}}

\put(4.5,-9){\makebox(0,0)[cc]{\(\mathcal{B}(6)\)}}
\put(2.1,-9.9){\makebox(0,0)[lb]{\(\langle 390\rangle\)}}
\put(4.1,-9.9){\makebox(0,0)[lb]{\(\langle 389\rangle\)}}
\put(6.1,-9.9){\makebox(0,0)[lb]{\(\langle 391\rangle\)}}

\put(-5.8,-11.9){\makebox(0,0)[rb]{\(\langle 2224\rangle\)}}
\put(-4,-11.9){\makebox(0,0)[rb]{\(\langle 2223\rangle\)}}
\put(-1.8,-11.9){\makebox(0,0)[rb]{\(\langle 2222\rangle\)}}
\put(-0.1,-11.9){\makebox(0,0)[rb]{\(\langle 2221\rangle\)}}

\put(4.5,-11){\makebox(0,0)[cc]{\(\mathcal{B}(7)\)}}
\put(2.1,-11.9){\makebox(0,0)[lb]{\(\langle 2226\rangle\)}}
\put(4,-11.9){\makebox(0,0)[lb]{\(\langle 2225\rangle\)}}
\put(5.8,-11.9){\makebox(0,0)[lb]{\(\langle 2227\rangle\)}}

\put(4.5,-13){\makebox(0,0)[cc]{\(\mathcal{B}(8)\)}}
\put(4.5,-15){\makebox(0,0)[cc]{\(\mathcal{B}(9)\)}}

\put(0.1,-16.2){\makebox(0,0)[ct]{\(G^n_0(0,0)\)}}
\put(-2,-16.2){\makebox(0,0)[ct]{\(G^n_0(0,1)\)}}
\put(-4,-16.2){\makebox(0,0)[ct]{\(G^n_0(1,0)\)}}
\put(-6,-16.2){\makebox(0,0)[ct]{\(G^n_0(-1,0)\)}}
\put(2,-16.2){\makebox(0,0)[ct]{\(G^n_1(0,-1)\)}}
\put(4,-16.2){\makebox(0,0)[ct]{\(G^n_1(0,0)\)}}
\put(6,-16.2){\makebox(0,0)[ct]{\(G^n_1(0,1)\)}}

\put(-3,-17.7){\makebox(0,0)[ct]{with abelian maximal subgroup}}
\put(4,-17.7){\makebox(0,0)[ct]{without abelian maximal subgroup}}

\put(2.5,0){\makebox(0,0)[cc]{\textbf{TKT}}}
\put(3.5,0){\makebox(0,0)[cc]{a.1}}
\put(2.5,-0.5){\makebox(0,0)[cc]{\(\varkappa=\)}}
\put(3.5,-0.5){\makebox(0,0)[cc]{\((0000)\)}}
\put(1.8,-0.7){\framebox(2.6,1){}}
\put(-6,-2){\makebox(0,0)[cc]{\textbf{TKT}}}
\put(-5,-2){\makebox(0,0)[cc]{A.1}}
\put(-6,-2.5){\makebox(0,0)[cc]{\(\varkappa=\)}}
\put(-5,-2.5){\makebox(0,0)[cc]{\((1111)\)}}
\put(-6.7,-2.7){\framebox(2.6,1){}}

\put(-8,-19){\makebox(0,0)[cc]{\textbf{TKT}}}
\put(0,-19){\makebox(0,0)[cc]{a.1}}
\put(-2,-19){\makebox(0,0)[cc]{a.2}}
\put(-4,-19){\makebox(0,0)[cc]{a.3}}
\put(-6,-19){\makebox(0,0)[cc]{a.3}}
\put(2,-19){\makebox(0,0)[cc]{a.1}}
\put(4,-19){\makebox(0,0)[cc]{a.1}}
\put(6,-19){\makebox(0,0)[cc]{a.1}}
\put(-8,-19.5){\makebox(0,0)[cc]{\(\varkappa=\)}}
\put(0,-19.5){\makebox(0,0)[cc]{\((0000)\)}}
\put(-2,-19.5){\makebox(0,0)[cc]{\((1000)\)}}
\put(-4,-19.5){\makebox(0,0)[cc]{\((2000)\)}}
\put(-6,-19.5){\makebox(0,0)[cc]{\((2000)\)}}
\put(2,-19.5){\makebox(0,0)[cc]{\((0000)\)}}
\put(4,-19.5){\makebox(0,0)[cc]{\((0000)\)}}
\put(6,-19.5){\makebox(0,0)[cc]{\((0000)\)}}
\put(-8.7,-19.7){\framebox(15.4,1){}}

\put(0.4,0){\oval(2.5,1.5)}
\put(-2.8,-2){\oval(2.5,1.5)}
\put(0.1,-1.1){\makebox(0,0)[lc]{\underbar{\textbf{657}}\((\ast 2)\)}}
\put(-2,-1){\makebox(0,0)[cc]{\underbar{\textbf{2\,439}}\((\ast 2)\)}}

\put(-4,-4){\oval(1.5,2)}
\multiput(-6,-4)(4,0){2}{\oval(1.5,1.5)}
\multiput(-4,-6.1)(0,-4){1}{\oval(1.5,1.5)}
\multiput(-2,-6.1)(0,-4){1}{\oval(1.5,1.5)}
\multiput(4,-6.1)(0,-4){1}{\oval(5.9,1.5)}
\put(-6,-5.1){\makebox(0,0)[cc]{\underbar{\textbf{8\,001}}\((\ast 1)\)}}
\put(-4,-5.2){\makebox(0,0)[cc]{\underbar{\textbf{3\,913}}\((\ast 2)\)}}
\put(-2,-5.1){\makebox(0,0)[lc]{\underbar{\textbf{8\,001}}\((\ast 2)\)}}
\put(-4,-7.1){\makebox(0,0)[cc]{\underbar{\textbf{4\,977}}\((\ast 1)\)}}
\put(-2,-7.1){\makebox(0,0)[lc]{\underbar{\textbf{20\,367}}\((\ast 1)\)}}
\put(6,-7.1){\makebox(0,0)[cc]{\underbar{\textbf{4\,977}}\((\ast 2)\)}}

\end{picture}

}

\end{figure}

\newpage

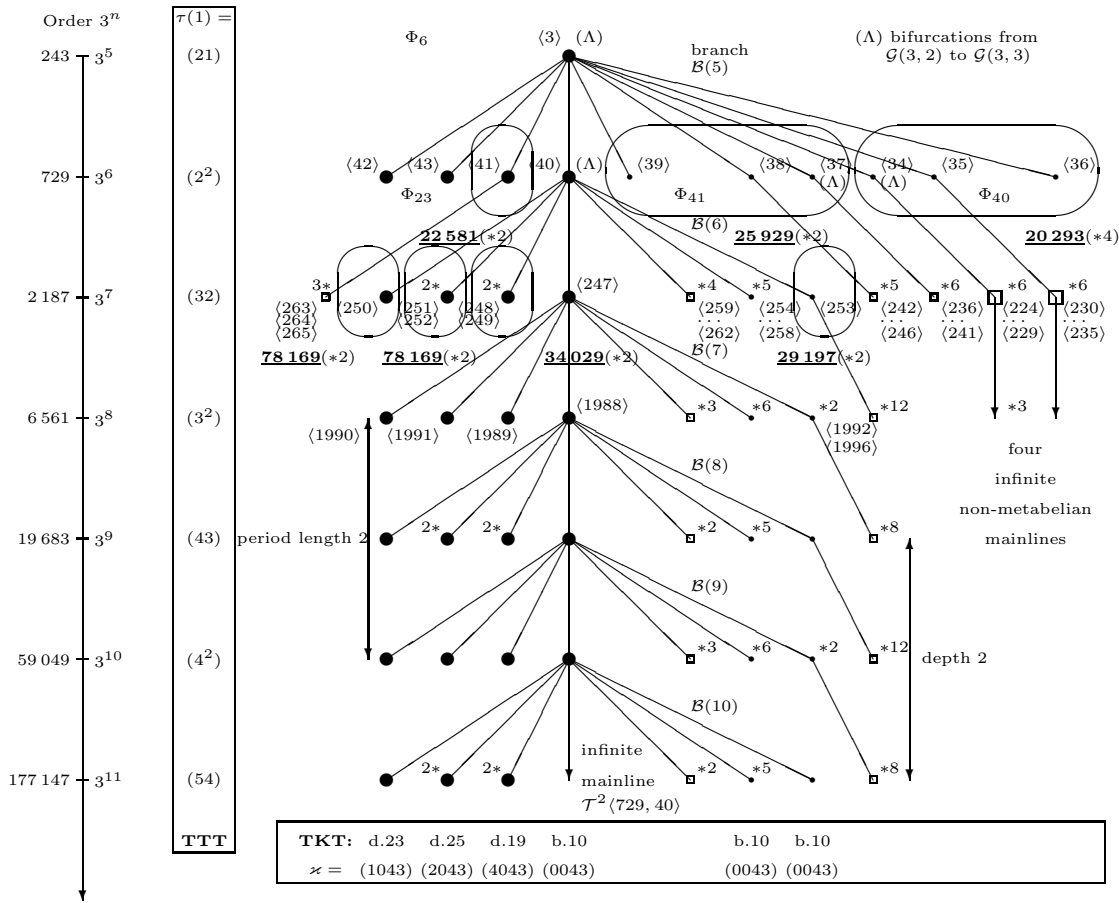
\begin{figure}[ht]
\caption{Distribution of Conductors for \(\mathrm{G}_3^{(2)}{F}\) on the Coclass Tree \(\mathcal{T}^2\langle 729,40\rangle\)}
\label{fig:MinDiscTyp33TreeBCc2}

{\tiny


\setlength{\unitlength}{0.8cm}
\begin{picture}(16,15)(-9,-14)

\put(-8,0.5){\makebox(0,0)[cb]{Order \(3^n\)}}
\put(-8,0){\line(0,-1){12}}
\multiput(-8.1,0)(0,-2){7}{\line(1,0){0.2}}
\put(-8.2,0){\makebox(0,0)[rc]{\(243\)}}
\put(-7.8,0){\makebox(0,0)[lc]{\(3^5\)}}
\put(-8.2,-2){\makebox(0,0)[rc]{\(729\)}}
\put(-7.8,-2){\makebox(0,0)[lc]{\(3^6\)}}
\put(-8.2,-4){\makebox(0,0)[rc]{\(2\,187\)}}
\put(-7.8,-4){\makebox(0,0)[lc]{\(3^7\)}}
\put(-8.2,-6){\makebox(0,0)[rc]{\(6\,561\)}}
\put(-7.8,-6){\makebox(0,0)[lc]{\(3^8\)}}
\put(-8.2,-8){\makebox(0,0)[rc]{\(19\,683\)}}
\put(-7.8,-8){\makebox(0,0)[lc]{\(3^9\)}}
\put(-8.2,-10){\makebox(0,0)[rc]{\(59\,049\)}}
\put(-7.8,-10){\makebox(0,0)[lc]{\(3^{10}\)}}
\put(-8.2,-12){\makebox(0,0)[rc]{\(177\,147\)}}
\put(-7.8,-12){\makebox(0,0)[lc]{\(3^{11}\)}}
\put(-8,-12){\vector(0,-1){2}}

\put(-6,0.5){\makebox(0,0)[cb]{\(\tau(1)=\)}}
\put(-6,0){\makebox(0,0)[cc]{\((21)\)}}
\put(-6,-2){\makebox(0,0)[cc]{\((2^2)\)}}
\put(-6,-4){\makebox(0,0)[cc]{\((32)\)}}
\put(-6,-6){\makebox(0,0)[cc]{\((3^2)\)}}
\put(-6,-8){\makebox(0,0)[cc]{\((43)\)}}
\put(-6,-10){\makebox(0,0)[cc]{\((4^2)\)}}
\put(-6,-12){\makebox(0,0)[cc]{\((54)\)}}
\put(-6,-13){\makebox(0,0)[cc]{\textbf{TTT}}}
\put(-6.5,-13.2){\framebox(1,14){}}

\put(5.6,-10){\vector(0,1){2}}
\put(5.8,-10){\makebox(0,0)[lc]{depth \(2\)}}
\put(5.6,-10){\vector(0,-1){2}}

\put(-3.3,-8){\vector(0,1){2}}
\put(-3.3,-8){\makebox(0,0)[rc]{period length \(2\)}}
\put(-3.3,-8){\vector(0,-1){2}}

\put(4.7,0.3){\makebox(0,0)[lc]{\((\Lambda)\) bifurcations from}}
\put(5.2,0){\makebox(0,0)[lc]{\(\mathcal{G}(3,2)\) to \(\mathcal{G}(3,3)\)}}

\multiput(0,0)(0,-2){6}{\circle*{0.2}}
\multiput(0,0)(0,-2){5}{\line(0,-1){2}}

\put(0,-10){\vector(0,-1){2}}
\put(0.2,-11.5){\makebox(0,0)[lc]{infinite}}
\put(0.2,-12){\makebox(0,0)[lc]{mainline}}
\put(0.2,-12.4){\makebox(0,0)[lc]{\(\mathcal{T}^2\langle 729,40\rangle\)}}

\put(7.2,-5.8){\makebox(0,0)[lc]{\(*3\)}}
\multiput(7,-4)(1,0){2}{\vector(0,-1){2}}
\put(7.5,-6.5){\makebox(0,0)[cc]{four}}
\put(7.5,-7){\makebox(0,0)[cc]{infinite}}
\put(7.5,-7.5){\makebox(0,0)[cc]{non-metabelian}}
\put(7.5,-8){\makebox(0,0)[cc]{mainlines}}

\multiput(-1,-2)(0,-2){6}{\circle*{0.2}}
\multiput(-2,-2)(0,-2){6}{\circle*{0.2}}
\multiput(-3,-2)(0,-2){6}{\circle*{0.2}}
\multiput(1.95,-4.05)(0,-2){5}{\framebox(0.1,0.1){}}
\multiput(3,-2)(0,-2){6}{\circle*{0.1}}
\multiput(4,-2)(0,-2){6}{\circle*{0.1}}
\multiput(5,-2)(1,0){2}{\circle*{0.1}}
\multiput(1,-2)(7,0){2}{\circle*{0.1}}

\multiput(0,0)(0,-2){6}{\line(-1,-2){1}}
\multiput(0,0)(0,-2){6}{\line(-1,-1){2}}
\multiput(0,0)(0,-2){6}{\line(-3,-2){3}}
\multiput(0,-2)(0,-2){5}{\line(1,-1){2}}
\multiput(0,0)(0,-2){6}{\line(3,-2){3}}
\multiput(0,0)(0,-2){6}{\line(2,-1){4}}
\put(0,0){\line(1,-2){1}}
\put(0,0){\line(5,-2){5}}
\put(0,0){\line(3,-1){6}}
\put(0,0){\line(4,-1){8}}

\put(-4.05,-4.05){\framebox(0.1,0.1){}}
\put(5.95,-4.05){\framebox(0.1,0.1){}}
\multiput(6.9,-4.1)(1,0){2}{\framebox(0.2,0.2){}}
\multiput(4.95,-4.05)(0,-2){5}{\framebox(0.1,0.1){}}

\put(-1,-2){\line(-3,-2){3}}
\multiput(3,-2)(1,0){4}{\line(1,-1){2}}
\multiput(4,-4)(0,-2){4}{\line(1,-2){1}}

\put(2,0.1){\makebox(0,0)[lc]{branch}}
\put(2,-0.2){\makebox(0,0)[lc]{\(\mathcal{B}(5)\)}}
\put(2,-2.8){\makebox(0,0)[lc]{\(\mathcal{B}(6)\)}}
\put(2,-4.9){\makebox(0,0)[lc]{\(\mathcal{B}(7)\)}}
\put(2,-6.8){\makebox(0,0)[lc]{\(\mathcal{B}(8)\)}}
\put(2,-8.8){\makebox(0,0)[lc]{\(\mathcal{B}(9)\)}}
\put(2,-10.8){\makebox(0,0)[lc]{\(\mathcal{B}(10)\)}}

\put(-2.5,0.3){\makebox(0,0)[cc]{\(\Phi_{6}\)}}
\put(-2.5,-2.3){\makebox(0,0)[cc]{\(\Phi_{23}\)}}
\put(2,-2.3){\makebox(0,0)[cc]{\(\Phi_{41}\)}}
\put(7,-2.3){\makebox(0,0)[cc]{\(\Phi_{40}\)}}

\put(-0.1,0.3){\makebox(0,0)[rc]{\(\langle 3\rangle\)}}
\put(0.1,0.3){\makebox(0,0)[lc]{\((\Lambda)\)}}

\put(-3.1,-1.8){\makebox(0,0)[rc]{\(\langle 42\rangle\)}}
\put(-2.1,-1.8){\makebox(0,0)[rc]{\(\langle 43\rangle\)}}
\put(-1.1,-1.8){\makebox(0,0)[rc]{\(\langle 41\rangle\)}}

\put(-0.1,-1.8){\makebox(0,0)[rc]{\(\langle 40\rangle\)}}
\put(0.1,-1.8){\makebox(0,0)[lc]{\((\Lambda)\)}}

\put(1.1,-1.8){\makebox(0,0)[lc]{\(\langle 39\rangle\)}}
\put(3.1,-1.8){\makebox(0,0)[lc]{\(\langle 38\rangle\)}}
\put(4.1,-1.8){\makebox(0,0)[lc]{\(\langle 37\rangle\)}}
\put(4.1,-2.1){\makebox(0,0)[lc]{\((\Lambda)\)}}
\put(5.1,-1.8){\makebox(0,0)[lc]{\(\langle 34\rangle\)}}
\put(5.1,-2.1){\makebox(0,0)[lc]{\((\Lambda)\)}}
\put(6.1,-1.8){\makebox(0,0)[lc]{\(\langle 35\rangle\)}}
\put(8.1,-1.8){\makebox(0,0)[lc]{\(\langle 36\rangle\)}}
\put(-4.1,-4.2){\makebox(0,0)[rc]{\(\langle 263\rangle\)}}
\put(-4.1,-4.4){\makebox(0,0)[rc]{\(\langle 264\rangle\)}}
\put(-4.1,-4.6){\makebox(0,0)[rc]{\(\langle 265\rangle\)}}
\put(-3.1,-4.2){\makebox(0,0)[rc]{\(\langle 250\rangle\)}}
\put(-2.1,-4.2){\makebox(0,0)[rc]{\(\langle 251\rangle\)}}
\put(-2.1,-4.4){\makebox(0,0)[rc]{\(\langle 252\rangle\)}}
\put(-1.1,-4.2){\makebox(0,0)[rc]{\(\langle 248\rangle\)}}
\put(-1.1,-4.4){\makebox(0,0)[rc]{\(\langle 249\rangle\)}}

\put(0.1,-3.8){\makebox(0,0)[lc]{\(\langle 247\rangle\)}}

\put(2.1,-4.2){\makebox(0,0)[lc]{\(\langle 259\rangle\)}}
\put(2.1,-4.4){\makebox(0,0)[lc]{\(\cdots\)}}
\put(2.1,-4.6){\makebox(0,0)[lc]{\(\langle 262\rangle\)}}
\put(3.1,-4.2){\makebox(0,0)[lc]{\(\langle 254\rangle\)}}
\put(3.1,-4.4){\makebox(0,0)[lc]{\(\cdots\)}}
\put(3.1,-4.6){\makebox(0,0)[lc]{\(\langle 258\rangle\)}}
\put(4.1,-4.2){\makebox(0,0)[lc]{\(\langle 253\rangle\)}}
\put(5.1,-4.2){\makebox(0,0)[lc]{\(\langle 242\rangle\)}}
\put(5.1,-4.4){\makebox(0,0)[lc]{\(\cdots\)}}
\put(5.1,-4.6){\makebox(0,0)[lc]{\(\langle 246\rangle\)}}
\put(6.1,-4.2){\makebox(0,0)[lc]{\(\langle 236\rangle\)}}
\put(6.1,-4.4){\makebox(0,0)[lc]{\(\cdots\)}}
\put(6.1,-4.6){\makebox(0,0)[lc]{\(\langle 241\rangle\)}}
\put(7.1,-4.2){\makebox(0,0)[lc]{\(\langle 224\rangle\)}}
\put(7.1,-4.4){\makebox(0,0)[lc]{\(\cdots\)}}
\put(7.1,-4.6){\makebox(0,0)[lc]{\(\langle 229\rangle\)}}
\put(8.1,-4.2){\makebox(0,0)[lc]{\(\langle 230\rangle\)}}
\put(8.1,-4.4){\makebox(0,0)[lc]{\(\cdots\)}}
\put(8.1,-4.6){\makebox(0,0)[lc]{\(\langle 235\rangle\)}}

\put(0.1,-5.8){\makebox(0,0)[lc]{\(\langle 1988\rangle\)}}

\put(-0.8,-6.3){\makebox(0,0)[rc]{\(\langle 1989\rangle\)}}
\put(-2.1,-6.3){\makebox(0,0)[rc]{\(\langle 1991\rangle\)}}
\put(-3.4,-6.3){\makebox(0,0)[rc]{\(\langle 1990\rangle\)}}

\put(4.2,-6.2){\makebox(0,0)[lc]{\(\langle 1992\rangle\)}}
\put(4.2,-6.5){\makebox(0,0)[lc]{\(\langle 1996\rangle\)}}

\put(-3.9,-3.8){\makebox(0,0)[rc]{\(3*\)}}
\multiput(-2.1,-3.8)(0,-4){3}{\makebox(0,0)[rc]{\(2*\)}}
\multiput(-1.1,-3.8)(0,-4){3}{\makebox(0,0)[rc]{\(2*\)}}
\put(2.1,-3.8){\makebox(0,0)[lc]{\(*4\)}}
\put(5.1,-3.8){\makebox(0,0)[lc]{\(*5\)}}
\put(6.1,-3.8){\makebox(0,0)[lc]{\(*6\)}}
\multiput(7.2,-3.8)(1,0){2}{\makebox(0,0)[lc]{\(*6\)}}
\multiput(2.1,-5.8)(0,-4){2}{\makebox(0,0)[lc]{\(*3\)}}
\multiput(2.1,-7.8)(0,-4){2}{\makebox(0,0)[lc]{\(*2\)}}
\multiput(3,-3.8)(0,-4){3}{\makebox(0,0)[lc]{\(*5\)}}
\multiput(3,-5.8)(0,-4){2}{\makebox(0,0)[lc]{\(*6\)}}
\multiput(4.1,-5.8)(0,-4){2}{\makebox(0,0)[lc]{\(*2\)}}
\multiput(5.1,-5.8)(0,-4){2}{\makebox(0,0)[lc]{\(*12\)}}
\multiput(5.1,-7.8)(0,-4){2}{\makebox(0,0)[lc]{\(*8\)}}

\put(-4,-13){\makebox(0,0)[cc]{\textbf{TKT:}}}
\put(-3,-13){\makebox(0,0)[cc]{d.23}}
\put(-2,-13){\makebox(0,0)[cc]{d.25}}
\put(-1,-13){\makebox(0,0)[cc]{d.19}}
\put(0,-13){\makebox(0,0)[cc]{b.10}}
\put(3,-13){\makebox(0,0)[cc]{b.10}}
\put(4,-13){\makebox(0,0)[cc]{b.10}}
\put(-4,-13.5){\makebox(0,0)[cc]{\(\varkappa=\)}}
\put(-3,-13.5){\makebox(0,0)[cc]{\((1043)\)}}
\put(-2,-13.5){\makebox(0,0)[cc]{\((2043)\)}}
\put(-1,-13.5){\makebox(0,0)[cc]{\((4043)\)}}
\put(0,-13.5){\makebox(0,0)[cc]{\((0043)\)}}
\put(3,-13.5){\makebox(0,0)[cc]{\((0043)\)}}
\put(4,-13.5){\makebox(0,0)[cc]{\((0043)\)}}
\put(-4.8,-13.7){\framebox(13.6,1){}}

\put(-1.1,-1.9){\oval(1,1.5)}
\put(2.6,-1.9){\oval(4,1.5)}
\put(6.7,-1.9){\oval(4,1.5)}
\put(-0.9,-3){\makebox(0,0)[rc]{\underbar{\textbf{22\,581}}\((\ast 2)\)}}
\put(3.5,-3){\makebox(0,0)[cc]{\underbar{\textbf{25\,929}}\((\ast 2)\)}}
\put(7.5,-3){\makebox(0,0)[lc]{\underbar{\textbf{20\,293}}\((\ast 4)\)}}
\put(-3.3,-3.9){\oval(1,1.5)}
\put(-3.5,-5){\makebox(0,0)[rc]{\underbar{\textbf{78\,169}}\((\ast 2)\)}}
\put(-2.2,-3.9){\oval(1,1.5)}
\put(-1.5,-5){\makebox(0,0)[rc]{\underbar{\textbf{78\,169}}\((\ast 2)\)}}
\put(-1.1,-3.9){\oval(1,1.5)}
\put(-0.4,-5){\makebox(0,0)[lc]{\underbar{\textbf{34\,029}}\((\ast 2)\)}}
\put(4.2,-3.9){\oval(1,1.5)}
\put(4.2,-5){\makebox(0,0)[cc]{\underbar{\textbf{29\,197}}\((\ast 2)\)}}

\end{picture}

}

\end{figure}

Generally, the \textit{vertices} of \textit{coclass trees} in the Figures
\ref{fig:MinDiscTyp33TreeCc1},
\ref{fig:MinDiscTyp33TreeBCc2}
and
\ref{fig:MinDiscTyp33TreeUCc2}
represent isomorphism classes of finite \(3\)-groups.
Two vertices are connected by a \textit{directed edge} \(G\to H\) if
\(H\) is isomorphic to the last lower central quotient \(G/\gamma_c(G)\),
where \(c=\mathrm{cl}(G)\) denotes the nilpotency class of \(G\),
and either \(\lvert G\rvert=3\lvert H\rvert\), that is,
the last lower central \(\gamma_c(G)\) is cyclic of order \(3\),
or \(\lvert G\rvert=9\lvert H\rvert\), that is,
\(\gamma_c(G)\) is bicyclic of type \((3,3)\).
See also
\cite[\S\ 2.2, pp. 410--411]{Ma2013}
and
\cite[\S\ 4, pp. 163--164]{Ma2015a}.

Vertices of the tree diagrams in Figure
\ref{fig:MinDiscTyp33TreeCc1}
are classified with various symbols:
\begin{enumerate}
\item
big full discs {\Large \(\bullet\)} represent metabelian groups \(\mathfrak{M}\)
with \textit{defect} \(k(\mathfrak{M})=0\)
\cite{Ma2012},
\item
small full discs {\footnotesize \(\bullet\)} represent metabelian groups \(\mathfrak{M}\)
with defect \(k(\mathfrak{M})=1\).
\end{enumerate}


\newpage

\noindent
In the Figures
\ref{fig:MinDiscTyp33TreeBCc2}
and
\ref{fig:MinDiscTyp33TreeUCc2},
\begin{enumerate}
\item
big full discs {\Large \(\bullet\)} represent \textit{metabelian} groups \(\mathfrak{M}\)
with bicyclic centre of type \((3,3)\) and defect \(k(\mathfrak{M})=0\)
\cite[\S\ 3.3.2, p. 429]{Ma2013},
\item
small full discs {\footnotesize \(\bullet\)} represent metabelian groups \(\mathfrak{M}\)
with cyclic centre of order \(3\) and defect \(k(\mathfrak{M})=1\),
\item
small contour squares {\tiny \(\square\)} represent \textit{non-metabelian} groups \(\mathfrak{G}\).
\end{enumerate}

\noindent
A symbol \(n\ast\) adjacent to a vertex
denotes the multiplicity of a \textit{batch} of \(n\) siblings,
that is, immediate descendants sharing a common parent.
The groups of particular importance are labelled by a number in angles,
which is the \textit{identifier} in the SmallGroups Library
\cite{BEO2002,BEO2005}
of GAP
\cite{GAP2021}
and Magma
\cite{MAGMA2022}.
We omit the orders, which are given on the left hand scale.
The transfer kernel type (TKT) \(\varkappa\)
\cite[Thm. 2.5, Tbl. 6--7]{Ma2012},
in the bottom rectangle concerns all
vertices located vertically above.
The first component \(\tau(1)\) of the transfer target type (TTT)
\cite[Dfn. 3.3, p. 288]{Ma2015b}
in the left rectangle
concerns vertices \(G\) on the same horizontal level with defect \(k(G)=0\).
The periodicity with length \(2\) of branches,
\(\mathcal{B}(j)\simeq\mathcal{B}(j+2)\) for \(j\ge 4\), respectively \(j\ge 7\),
sets in with branch \(\mathcal{B}(4)\), respectively \(\mathcal{B}(7)\),
having a root of order \(3^4\), respectively \(3^7\),
in Figure
\ref{fig:MinDiscTyp33TreeCc1},
resp. Figures
\ref{fig:MinDiscTyp33TreeBCc2},
\ref{fig:MinDiscTyp33TreeUCc2}.

\newpage

\section{Cyclic Cubic Quartets of Type (3,3)}
\label{s:CycCub3x3Recent}

\subsection{Conductors with Three Prime Divisors}
\label{ss:ThreePrimeCond3x3}

\noindent
Now we come to conductors \(c\) with \textit{three} prime divisors, \(t=3\).
According to the multiplicity formula \(m(c)=(3-1)^{t-1}\),
there are \(4\) cyclic cubic fields \(F_{c,1},\ldots,F_{c,4}\) sharing the common conductor \(c\).
The members of a \textit{quartet} \((F_{c,\mu})_{1\le\mu\le 4}\) do not necessarily have equal \(3\)-class ranks \(\varrho_3{F_{c,\mu}}\),
according to Theorem
\ref{thm:3ClassRank}.


\begin{definition}
\label{dfn:Categories}
This forces us to introduce \textit{three categories} of quartets \((F_{c,\mu})_{1\le\mu\le 4}\). \\
\textbf{Category} \(\mathrm{I}\): \\
one member \(F_{c,1}\) has \(3\)-class rank \(3\), the remaining three \(F_{c,2},\ldots,F_{c,4}\) have rank \(2\). \\
\textbf{Category} \(\mathrm{II}\): \\
two members \(F_{c,1},F_{c,2}\) have \(3\)-class rank \(3\), the other two \(F_{c,3},F_{c,4}\) have rank \(2\). \\
\textbf{Category} \(\mathrm{III}\): all four members \(F_{c,1},\ldots,F_{c,4}\) have \(3\)-class rank \(2\).

Additionally to these categories, defined in Ayadi's Thesis
\cite{Ay1995},
we introduce \textit{two further categories}. \\
\textbf{Category} \(\mathrm{IV}\): all four members \(F_{c,1},\ldots,F_{c,4}\) have \(3\)-class rank \(3\). \\
\textbf{Category} \(\mathrm{V}\): all four members \(F_{c,1},\ldots,F_{c,4}\) have \(3\)-class rank \(4\).
\end{definition}

\noindent
However, there arises an additional complication:
In general, the members of rank \(2\) neither have the same \(3\)-capitulation type \(\varkappa(F_{c,\mu})\)
nor the same \(3\)-class tower group \(\mathrm{G}_3^{(\infty)}{F_{c,\mu}}\).
Therefore we occasionally use formal powers with exponents denoting iteration.
There is still one classical result, Theorem
\ref{thm:Cat3Gph1To4},
for \(t=3\).


\subsection{Graphs 1,2,3,4 of Category III}
\label{ss:Cat3Gph1To4}

\noindent
Rational criteria for \(\mathrm{G}_3^{(\infty)}{F_{c,\mu}}\) to be abelian
are given in the following theorem which refers to the first four subcases of Formula
\eqref{eqn:Category3}.


\begin{theorem}
\label{thm:Cat3Gph1To4}
(Cubic residue criterion, \textbf{Ayadi, 1995},
\cite{Ay1995,AAI2001}.) \\
Let \(c\) be a conductor divisible by exactly three primes, \(t=3\), such that
\(\mathrm{Cl}_3{F_{c,\mu}}\simeq (3,3)\) for all four cyclic cubic fields \(F_{c,\mu}\), \(1\le\mu\le 4\), with conductor \(c\).
If there are no mutual cubic residues among the prime divisors of \(c\),
that is, if \(c=q_1q_2q_3\) belongs to the Graphs \(1,2,3,4\) of Category \(\mathrm{III}\), i.e.,
\begin{equation}
\label{eqn:Cat3Gph1To4}
\lbrack q_1,q_2,q_3\rbrack_3=
\begin{cases}
\lbrace q_1,q_2,q_3;\delta\ne 0\rbrace & \text{ or} \\
\lbrace q_i\rightarrow q_j;q_k\rbrace & \text{ or} \\
\lbrace q_i\rightarrow q_j\rightarrow q_k\rbrace & \text{ or} \\
\lbrace q_i\rightarrow q_j\rightarrow q_k\rightarrow q_i\rbrace &
\end{cases}
\end{equation}
with \(i,j,k\) pairwise distinct,
then the second \(3\)-class group \(\mathrm{G}_3^{(2)}{F_{c,\mu}}\) of all four fields \(F_{c,\mu}\)
is isomorphic to the abelian group \(\langle 9,2\rangle \simeq (3,3)\)
with capitulation type \(\mathrm{a}.1\), \(\varkappa(F_{c,\mu}) = (0000)\),
and transfer target type \(\tau(F_{c,\mu}) = \lbrack 1,1,1,1\rbrack\),
the \(3\)-class tower has length \(\ell_3{F_{c,\mu}}=1\),
and the \(3\)-class groups of the \(13\) bicyclic bicubic subfields \(S_i\)
of the \(3\)-genus field \(F^\ast\) are given by
the logarithmic abelian type invariants (ATI)
\begin{equation}
\label{eqn:Cat3Gph1To4Genus}
\left\lbrack\mathrm{Cl}_3{S_i}\right\rbrack_{1\le i\le 13}=\lbrack (0)^3;(1)^{10}\rbrack.
\end{equation}
\end{theorem}

\begin{conjecture}
\label{cnj:Cat3Gph1To4}
We conjecture that the converse of Theorem
\ref{thm:Cat3Gph1To4}
is also true, i.e.,
that \(\mathrm{G}_3^{(2)}{F}\) is never abelian for conductors \(c\) in the Categories \(\mathrm{I}\) and \(\mathrm{II}\),
and in the remaining Graphs \(5,6,7,8,9\) of Category \(\mathrm{III}\).
\end{conjecture}


\begin{example}
\label{exm:Cat3Gph1To4}
The decisive advantage of the criteria in Theorem
\ref{thm:Cat3Gph1To4}
is that they are \textit{rational}, i.e., expressible in terms of the prime divisors of the conductor.
In Table
\ref{tbl:Cat3Gph1To4},
the first occurrence of each of the Graphs \(1,2,3,4\) of Category \(\mathrm{III}\) is presented.
Again, we give the smallest examples for both, \(9\mid c\) and \(\gcd(c,3)=1\), separately,
indicated by the \(3\)-valuation \(v_3{c}\) of \(c=q_1q_2q_3\).
\end{example}


\renewcommand{\arraystretch}{1.1}

\begin{table}[ht]
\caption{The First Examples for the Graphs \(1,2,3,4\) of Category \(\mathrm{III}\)}
\label{tbl:Cat3Gph1To4}
\begin{center}
{\small
\begin{tabular}{|c|crc|c|}
\hline
 Graph & \(v_3{c}\) &       \(c\) &                                \(\lbrack q_1,q_2,q_3\rbrack_3\) & \(\mathrm{G}_3^{(2)}{F_{c,\mu}}\) \\
\hline
 \(1\) &      \(2\) &  \(1\,953\) &                         \(\lbrace 9,7,31;\ \delta\ne 0\rbrace\) &          \(\langle 9,2\rangle^4\) \\
 \(1\) &      \(0\) & \(14\,539\) &                        \(\lbrace 7,31,67;\ \delta\ne 0\rbrace\) &          \(\langle 9,2\rangle^4\) \\
\hline
 \(2\) &      \(2\) &     \(819\) &                            \(\lbrace 13\rightarrow 7;9\rbrace\) &          \(\langle 9,2\rangle^4\) \\
 \(2\) &      \(0\) &  \(3\,367\) &                           \(\lbrace 13\rightarrow 7;37\rbrace\) &          \(\langle 9,2\rangle^4\) \\
\hline
 \(3\) &      \(2\) &  \(1\,197\) &                 \(\lbrace 7\rightarrow 19\rightarrow 9\rbrace\) &          \(\langle 9,2\rangle^4\) \\
 \(3\) &      \(0\) &  \(1\,729\) &                \(\lbrace 13\rightarrow 7\rightarrow 19\rbrace\) &          \(\langle 9,2\rangle^4\) \\
\hline
 \(4\) &      \(0\) &  \(6\,643\) &  \(\lbrace 13\rightarrow 7\rightarrow 73\rightarrow 13\rbrace\) &          \(\langle 9,2\rangle^4\) \\
 \(4\) &      \(2\) & \(17\,613\) & \(\lbrace 19\rightarrow 9\rightarrow 103\rightarrow 19\rbrace\) &          \(\langle 9,2\rangle^4\) \\
\hline
\end{tabular}
}
\end{center}
\end{table}


\subsection{Graphs 5,6,7,8,9 of Category III}
\label{ss:Cat3Gph5To9}
\noindent
Now we begin with interesting \textbf{new quartets} of cyclic cubic fields \(F\),
having \(t=3\) and \(\mathrm{Cl}_3{F_{c,\mu}}\simeq (3,3)\) for \(1\le\mu\le 4\),
where the second \(3\)-class groups \(\mathrm{G}_3^{(2)}{F_{c,\mu}}\) were completely \textbf{unknown} up to now.
According to Theorems
\ref{thm:3ClassRank}
and
\ref{thm:Cat3Gph1To4},
there must occur mutual cubic residues among the prime divisors \(q_1,q_2,q_3\) of \(c\),
say \(q_1,q_2\),
and thus it is also interesting to consider the coinciding group \(\mathrm{G}_3^{(2)}{F_{f,\nu}}\)
of the two members \(F_{f,\nu}\), \(1\le\nu\le 2\),
of the doublet with \textit{partial conductor} \(f=q_1\cdot q_2\), \(f\mid c\),
as given by Theorem
\ref{thm:TwoPrimeCond},
and listed in the tables of \\
\((\ast)\) \texttt{http://www.algebra.at/ResearchFrontier2013ThreeByThree.htm}.

Within the frame of our computations,
conductors \(c\) belonging to the Graphs \(6\), \(7\) and \(9\) of Category \(\mathrm{III}\)
reveal a completely \textit{uniform behaviour} of the associated second \(3\)-class groups \(\mathrm{G}_3^{(2)}{F_{c,\mu}}\)
of the four members \(F_{c,\mu}\), \(1\le\mu\le 4\), of the quartet with conductor \(c=q_1\cdot q_2\cdot q_3\),
provided that \(\mathrm{G}_3^{(2)}{F_{f,\nu}}\simeq\langle 9,2\rangle\) is \textit{abelian}, black color in \((\ast)\),
for the two members \(F_{f,\nu}\), \(1\le\nu\le 2\), of the doublet with conductor \(f=q_1\cdot q_2\).
Therefore, we can summarize their properties in Theorem
\ref{thm:Cat3Gph679}.
However, exceptions arise if either \(\mathrm{G}_3^{(2)}{F_{f,\nu}}\simeq\langle 27,4\rangle\) is \textit{extra-special}, red color in \((\ast)\),
or \(\mathrm{Cl}_3{F_{f,\nu}}\simeq (9,3)\) is \textit{non-elementary},
listed in the tables of \\
\((\ast\ast)\) \texttt{http://www.algebra.at/ResearchFrontier2013NineByThree.htm}.

Exceptions are indicated by \textbf{boldface} font in the Tables
\ref{tbl:Cat3Gph6},
\ref{tbl:Cat3Gph7},
\ref{tbl:Cat3Gph9}.


\renewcommand{\arraystretch}{1.1}

\begin{table}[ht]
\caption{Thirty-One Examples for Graph \(6\) of Category \(\mathrm{III}\)}
\label{tbl:Cat3Gph6}
\begin{center}
{\tiny
\begin{tabular}{|crc|c|rc|c|c|}
\hline
    No. &       \(c\) &                       \(\lbrack q_1,q_2,q_3\rbrack_3\) & \(\mathrm{G}_3^{(2)}{F_{c,\mu}}\) & \(f\) &          \(\lbrack q_1,q_2\rbrack_3\) & \(\mathrm{G}_3^{(2)}{F_{f,\nu}}\) & \(\left\lbrack\mathrm{Cl}_3{S_i}\right\rbrack_{1\le i\le 13}\) \\
\hline
  \(1\) &  \(8\,541\) &   \(\lbrace 9\leftrightarrow 73\rightarrow 13\rbrace\) & \(\langle 81,7\rangle^4\) &    \(657\) &   \(\lbrace 9\leftrightarrow 73\rbrace\) & \(\langle 9,2\rangle^2\) & \(\lbrack (0)^2,1;(1^2)^8,(1^3)^2\rbrack\) \\
  \(2\) &  \(9\,373\) &  \(\lbrace 103\leftrightarrow 13\rightarrow 7\rbrace\) & \(\langle 81,7\rangle^4\) & \(1\,339\) & \(\lbrace 13\leftrightarrow 103\rbrace\) & \(\langle 9,2\rangle^2\) & \(\lbrack (0)^2,1;(1^2)^8,(1^3)^2\rbrack\) \\
  \(3\) & \(11\,403\) &   \(\lbrace 7\leftrightarrow 181\rightarrow 9\rbrace\) & \(\langle 81,7\rangle^4\) & \(1\,267\) &  \(\lbrace 7\leftrightarrow 181\rbrace\) & \(\langle 9,2\rangle^2\) & \(\lbrack (0)^2,1;(1^2)^8,(1^3)^2\rbrack\) \\
  \(4\) & \(19\,341\) &   \(\lbrace 9\leftrightarrow 307\rightarrow 7\rbrace\) & \(\langle 81,7\rangle^4\) & \(2\,763\) &  \(\lbrace 9\leftrightarrow 307\rbrace\) & \(\langle 9,2\rangle^2\) & \(\lbrack (0)^2,1;(1^2)^8,(1^3)^2\rbrack\) \\
  \(5\) & \(20\,839\) &  \(\lbrace 229\leftrightarrow 13\rightarrow 7\rbrace\) & \(\langle 81,7\rangle^4\) & \(2\,977\) & \(\lbrace 13\leftrightarrow 229\rbrace\) & \(\langle 9,2\rangle^2\) & \(\lbrack (0)^2,1;(1^2)^8,(1^3)^2\rbrack\) \\
  \(6\) & \(25\,441\) & \(\lbrace 13\leftrightarrow 103\rightarrow 19\rbrace\) & \(\langle 81,7\rangle^4\) & \(1\,339\) & \(\lbrace 13\leftrightarrow 103\rbrace\) & \(\langle 9,2\rangle^2\) & \(\lbrack (0)^2,1;(1^2)^8,(1^3)^2\rbrack\) \\
  \(7\) & \(29\,659\) &  \(\lbrace 223\leftrightarrow 7\rightarrow 19\rbrace\) & \(\langle 81,7\rangle^4\) & \(1\,561\) &  \(\lbrace 7\leftrightarrow 223\rbrace\) & \(\langle 9,2\rangle^2\) & \(\lbrack (0)^2,1;(1^2)^8,(1^3)^2\rbrack\) \\
  \(8\) & \(35\,919\) &  \(\lbrace 9\leftrightarrow 307\rightarrow 13\rbrace\) & \(\langle 81,7\rangle^4\) & \(2\,763\) &  \(\lbrace 9\leftrightarrow 307\rbrace\) & \(\langle 9,2\rangle^2\) & \(\lbrack (0)^2,1;(1^2)^8,(1^3)^2\rbrack\) \\
  \(9\) & \(37\,297\) & \(\lbrace 19\leftrightarrow 151\rightarrow 13\rbrace\) & \(\langle 81,7\rangle^4\) & \(2\,869\) & \(\lbrace 19\leftrightarrow 151\rbrace\) & \(\langle 9,2\rangle^2\) & \(\lbrack (0)^2,1;(1^2)^8,(1^3)^2\rbrack\) \\
 \(10\) & \(40\,077\) &   \(\lbrace 73\leftrightarrow 9\rightarrow 61\rbrace\) & \(\langle 81,7\rangle^4\) &    \(657\) &   \(\lbrace 9\leftrightarrow 73\rbrace\) & \(\langle 9,2\rangle^2\) & \(\lbrack (0)^2,1;(1^2)^8,(1^3)^2\rbrack\) \\
 \(11\) & \(44\,019\) &   \(\lbrace 73\leftrightarrow 9\rightarrow 67\rbrace\) & \(\langle 81,7\rangle^4\) &    \(657\) &   \(\lbrace 9\leftrightarrow 73\rbrace\) & \(\langle 9,2\rangle^2\) & \(\lbrack (0)^2,1;(1^2)^8,(1^3)^2\rbrack\) \\
 \(12\) & \(44\,821\) &  \(\lbrace 337\leftrightarrow 7\rightarrow 19\rbrace\) & \(\langle 81,7\rangle^4\) & \(2\,359\) &  \(\lbrace 7\leftrightarrow 337\rbrace\) & \(\langle 9,2\rangle^2\) & \(\lbrack (0)^2,1;(1^2)^8,(1^3)^2\rbrack\) \\
 \(13\) & \(45\,409\) &  \(\lbrace 499\leftrightarrow 13\rightarrow 7\rbrace\) & \(\langle 81,7\rangle^4\) & \(6\,487\) & \(\lbrace 13\leftrightarrow 499\rbrace\) & \(\langle 9,2\rangle^2\) & \(\lbrack (0)^2,1;(1^2)^8,(1^3)^2\rbrack\) \\
 \(14\) & \(45\,477\) &  \(\lbrace 31\leftrightarrow 163\rightarrow 9\rbrace\) & \(\langle 81,7\rangle^4\) & \(5\,053\) & \(\lbrace 31\leftrightarrow 169\rbrace\) & \(\langle 9,2\rangle^2\) & \(\lbrack (0)^2,1;(1^2)^8,(1^3)^2\rbrack\) \\
 \(15\) & \(47\,367\) &  \(\lbrace 277\leftrightarrow 19\rightarrow 9\rbrace\) & \(\langle 81,7\rangle^4\) & \(5\,263\) & \(\lbrace 19\leftrightarrow 277\rbrace\) & \(\langle 9,2\rangle^2\) & \(\lbrack (0)^2,1;(1^2)^8,(1^3)^2\rbrack\) \\
 \(16\) & \(47\,691\) &   \(\lbrace 9\leftrightarrow 757\rightarrow 7\rbrace\) & \(\langle 81,7\rangle^4\) & \(6\,813\) &  \(\lbrace 9\leftrightarrow 757\rbrace\) & \(\langle 9,2\rangle^2\) & \(\lbrack (0)^2,1;(1^2)^8,(1^3)^2\rbrack\) \\
 \(17\) & \(51\,093\) &   \(\lbrace 7\leftrightarrow 811\rightarrow 9\rbrace\) & \(\langle 81,7\rangle^4\) & \(5\,677\) &  \(\lbrace 7\leftrightarrow 811\rbrace\) & \(\langle 9,2\rangle^2\) & \(\lbrack (0)^2,1;(1^2)^8,(1^3)^2\rbrack\) \\
 \(18\) & \(53\,641\) &   \(\lbrace 79\leftrightarrow 97\rightarrow 7\rbrace\) & \(\langle 81,7\rangle^4\) & \(7\,663\) &  \(\lbrace 79\leftrightarrow 97\rbrace\) & \(\langle 9,2\rangle^2\) & \(\lbrack (0)^2,1;(1^2)^8,(1^3)^2\rbrack\) \\
 \(19\) & \(55\,629\) &   \(\lbrace 7\leftrightarrow 883\rightarrow 9\rbrace\) & \(\langle 81,7\rangle^4\) & \(6\,181\) &  \(\lbrace 7\leftrightarrow 883\rbrace\) & \(\langle 9,2\rangle^2\) & \(\lbrack (0)^2,1;(1^2)^8,(1^3)^2\rbrack\) \\
 \(20\) & \(56\,329\) &  \(\lbrace 619\leftrightarrow 13\rightarrow 7\rbrace\) & \(\langle 81,7\rangle^4\) & \(8\,047\) & \(\lbrace 13\leftrightarrow 619\rbrace\) &\(\langle 27,4\rangle^2\) & \(\lbrack (0)^2,1^2;(1^2)^8,(1^3)^2\rbrack\) \\
 \(21\) & \(56\,563\) & \(\lbrace 13\leftrightarrow 229\rightarrow 19\rbrace\) & \(\langle 81,7\rangle^4\) & \(2\,977\) & \(\lbrace 13\leftrightarrow 229\rbrace\) & \(\langle 9,2\rangle^2\) & \(\lbrack (0)^2,1;(1^2)^8,(1^3)^2\rbrack\) \\
 \(22\) & \(57\,757\) &  \(\lbrace 7\leftrightarrow 223\rightarrow 37\rbrace\) & \(\langle 81,7\rangle^4\) & \(1\,561\) &  \(\lbrace 7\leftrightarrow 223\rbrace\) & \(\langle 9,2\rangle^2\) & \(\lbrack (0)^2,1;(1^2)^8,(1^3)^2\rbrack\) \\
 \(23\) & \(58\,093\) &  \(\lbrace 193\leftrightarrow 43\rightarrow 7\rbrace\) & \(\langle 81,7\rangle^4\) & \(8\,299\) & \(\lbrace 43\leftrightarrow 193\rbrace\) &\(\langle 27,4\rangle^2\) & \(\lbrack (0)^2,1^2;(1^2)^8,(1^3)^2\rbrack\) \\
 \(24\) & \(63\,567\) &  \(\lbrace 9\leftrightarrow 1009\rightarrow 7\rbrace\) & \(\langle 81,7\rangle^4\) & \(9\,081\) & \(\lbrace 9\leftrightarrow 1009\rbrace\) & \(\langle 9,2\rangle^2\) & \(\lbrack (0)^2,1;(1^2)^8,(1^3)^2\rbrack\) \\
 \(25\) & \(63\,783\) &  \(\lbrace 373\leftrightarrow 19\rightarrow 9\rbrace\) & \(\langle 81,7\rangle^4\) & \(7\,087\) & \(\lbrace 19\leftrightarrow 373\rbrace\) &\(\langle 27,4\rangle^2\) & \(\lbrack (0)^2,1^2;(1^2)^8,(1^3)^2\rbrack\) \\
 \(26\) & \(67\,509\) &  \(\lbrace 9\leftrightarrow 577\rightarrow 13\rbrace\) & \(\langle 81,7\rangle^4\) & \(5\,193\) &  \(\lbrace 9\leftrightarrow 577\rbrace\) & \(\langle 9,2\rangle^2\) & \(\lbrack (0)^2,1;(1^2)^8,(1^3)^2\rbrack\) \\
 \(27\) & \(73\,129\) &  \(\lbrace 7\leftrightarrow 337\rightarrow 31\rbrace\) & \(\langle 81,7\rangle^4\) & \(2\,359\) &  \(\lbrace 7\leftrightarrow 337\rbrace\) & \(\langle 9,2\rangle^2\) & \(\lbrack (0)^2,1;(1^2)^8,(1^3)^2\rbrack\) \\
 \(28\) & \(75\,733\) &  \(\lbrace 31\leftrightarrow 349\rightarrow 7\rbrace\) & \(\langle 81,7\rangle^4\) &\(10\,819\) & \(\lbrace 31\leftrightarrow 349\rbrace\) & \(\langle 9,2\rangle^2\) & \(\lbrack (0)^2,1;(1^2)^8,(1^3)^2\rbrack\) \\
 \(29\) & \(\mathbf{78\,169}\) & \(\lbrace 859\leftrightarrow 13\rightarrow 7\rbrace\) & \(\langle 3^7,250\rangle^2,\langle 3^7,251\rangle^2\) & \(\mathbf{11\,167}\) & \(\lbrace 13\leftrightarrow 859\rbrace\) & \(\langle 81,3\rangle^2\) & \(\lbrack (0)^2,1^3;(21)^3,(32)^3,(1^3)^4\rbrack\) \\
 \(30\) & \(92\,491\) &  \(\lbrace 881\leftrightarrow 7\rightarrow 73\rbrace\) & \(\langle 81,7\rangle^4\) &\(6\,181\) & \(\lbrace 7\leftrightarrow 883\rbrace\) & \(\langle 9,2\rangle^2\) & \(\lbrack (0)^2,1;(1^2)^8,(1^3)^2\rbrack\) \\
 \(31\) & \(99\,619\) &  \(\lbrace 97\leftrightarrow 79\rightarrow 13\rbrace\) & \(\langle 81,7\rangle^4\) &\(7\,663\) & \(\lbrace 79\leftrightarrow 97\rbrace\) & \(\langle 9,2\rangle^2\) & \(\lbrack (0)^2,1;(1^2)^8,(1^3)^2\rbrack\) \\
\hline
\end{tabular}
}
\end{center}
\end{table}

\noindent
Table
\ref{tbl:Cat3Gph6}
shows all \(31\) hits of quartets
belonging to Graph \(6\) in Category \(\mathrm{III}\),
as mentioned among the computational results for \(c<10^5\) in Table
\ref{tbl:StatCatGph3}.
There is only one irregular case with conductor \(c=78\,169\) and
\textit{singular} partial conductor \(f=11\,167\),
where \(\mathrm{Cl}_3{F_{f,\nu}}\simeq (9,3)\) and
\(\mathrm{Cl}_3{F_0^\ast}\simeq (3,3,3)\), i.e. \(v=3\).
Three cases with \(\mathrm{G}_3^{(2)}{F_{f,\nu}}\simeq\langle 27,4\rangle\)
fit in unceremoniously without causing irregular behavior,
except that \(\mathrm{Cl}_3{F_0^\ast}\simeq (3,3)\), i.e. \(v=2\).

Table
\ref{tbl:Cat3Gph7}
shows all \(34\) occurrences of quartets
belonging to Graph \(7\) in Category \(\mathrm{III}\),
as mentioned in Table
\ref{tbl:StatCatGph3}.
There are two irregular cases with conductors \(c=69\,979\), \(86\,821\) and
\textit{super-singular} partial conductors \(f=5\,383\), \(12\,403\),
where \(\mathrm{Cl}_3{F_{f,\nu}}\simeq (9,3)\) and
\(\mathrm{Cl}_3{F_0^\ast}\simeq (9,3,3)\), i.e. \(v=4\).
There is one irregular case with conductor \(c=61\,243\) and
\textit{singular} partial conductor \(f=4\,711\),
where \(\mathrm{Cl}_3{F_{f,\nu}}\simeq (9,3)\) and
\(\mathrm{Cl}_3{F_0^\ast}\simeq (3,3,3)\), i.e. \(v=3\).
Seven cases with \(\mathrm{G}_3^{(2)}{F_{f,\nu}}\simeq\langle 27,4\rangle\)
fit in unspectacularly without causing irregular behavior,
except that \(\mathrm{Cl}_3{F_0^\ast}\simeq (3,3)\), i.e. \(v=2\).
However, the behavior is not uniform, since four conductors
\(c=76\,741\), \(89\,433\), \(90\,243\), \(99\,801\)
with \(\mathrm{G}_3^{(2)}{F_{f,\nu}}\simeq\langle 27,4\rangle\)
cause exceptions.


\renewcommand{\arraystretch}{1.1}

\begin{table}[ht]
\caption{Thirty-Four Examples for Graph \(7\) of Category \(\mathrm{III}\)}
\label{tbl:Cat3Gph7}
\begin{center}
{\tiny
\begin{tabular}{|crc|c|rc|c|c|}
\hline
    No. &       \(c\) &                     \(\lbrack q_1,q_2,q_3\rbrack_3\) & \(\mathrm{G}_3^{(2)}{F_{c,\mu}}\) & \(f\) &          \(\lbrack q_1,q_2\rbrack_3\) & \(\mathrm{G}_3^{(2)}{F_{f,\nu}}\) & \(\left\lbrack\mathrm{Cl}_3{S_i}\right\rbrack_{1\le i\le 13}\) \\
\hline
  \(1\) &  \(4\,599\) &   \(\lbrace 9\leftrightarrow 73\leftarrow 7\rbrace\) & \(\langle 81,7\rangle^4\) &    \(657\) &   \(\lbrace 9\leftrightarrow 73\rbrace\) & \(\langle 9,2\rangle^2\) & \(\lbrack (0)^2,1;(1^2)^8,(1^3)^2\rbrack\) \\
  \(2\) & \(12\,051\) & \(\lbrace 13\leftrightarrow 103\leftarrow 9\rbrace\) & \(\langle 81,7\rangle^4\) & \(1\,339\) & \(\lbrace 13\leftrightarrow 103\rbrace\) & \(\langle 9,2\rangle^2\) & \(\lbrack (0)^2,1;(1^2)^8,(1^3)^2\rbrack\) \\
  \(3\) & \(12\,483\) &  \(\lbrace 73\leftrightarrow 9\leftarrow 19\rbrace\) & \(\langle 81,7\rangle^4\) &    \(657\) &   \(\lbrace 9\leftrightarrow 73\rbrace\) & \(\langle 9,2\rangle^2\) & \(\lbrack (0)^2,1;(1^2)^8,(1^3)^2\rbrack\) \\
  \(4\) & \(20\,083\) & \(\lbrace 151\leftrightarrow 19\leftarrow 7\rbrace\) & \(\langle 81,7\rangle^4\) & \(2\,869\) & \(\lbrace 19\leftrightarrow 151\rbrace\) & \(\langle 9,2\rangle^2\) & \(\lbrack (0)^2,1;(1^2)^8,(1^3)^2\rbrack\) \\
  \(5\) & \(28\,251\) &  \(\lbrace 9\leftrightarrow 73\leftarrow 43\rbrace\) & \(\langle 81,7\rangle^4\) &    \(657\) &   \(\lbrace 9\leftrightarrow 73\rbrace\) & \(\langle 9,2\rangle^2\) & \(\lbrack (0)^2,1;(1^2)^8,(1^3)^2\rbrack\) \\
  \(6\) & \(31\,707\) & \(\lbrace 9\leftrightarrow 271\leftarrow 13\rbrace\) & \(\langle 81,7\rangle^4\) & \(2\,439\) &  \(\lbrace 9\leftrightarrow 271\rbrace\) &\(\langle 27,4\rangle^2\) & \(\lbrack (0)^2,1^2;(1^2)^8,(1^3)^2\rbrack\) \\
  \(7\) & \(36\,841\) & \(\lbrace 277\leftrightarrow 19\leftarrow 7\rbrace\) & \(\langle 81,7\rangle^4\) & \(5\,263\) & \(\lbrace 19\leftrightarrow 277\rbrace\) & \(\langle 9,2\rangle^2\) & \(\lbrack (0)^2,1;(1^2)^8,(1^3)^2\rbrack\) \\
  \(8\) & \(39\,277\) & \(\lbrace 7\leftrightarrow 181\leftarrow 31\rbrace\) & \(\langle 81,7\rangle^4\) & \(1\,267\) &  \(\lbrace 7\leftrightarrow 181\rbrace\) & \(\langle 9,2\rangle^2\) & \(\lbrack (0)^2,1;(1^2)^8,(1^3)^2\rbrack\) \\
  \(9\) & \(52\,633\) & \(\lbrace 103\leftrightarrow 73\leftarrow 7\rbrace\) & \(\langle 81,7\rangle^4\) & \(7\,519\) & \(\lbrace 73\leftrightarrow 103\rbrace\) &\(\langle 27,4\rangle^2\) & \(\lbrack (0)^2,1^2;(1^2)^8,(1^3)^2\rbrack\) \\
 \(10\) & \(53\,739\) &  \(\lbrace 7\leftrightarrow 853\leftarrow 9\rbrace\) & \(\langle 81,7\rangle^4\) & \(5\,971\) &  \(\lbrace 7\leftrightarrow 853\rbrace\) &\(\langle 27,4\rangle^2\) & \(\lbrack (0)^2,1^2;(1^2)^8,(1^3)^2\rbrack\) \\
 \(11\) & \(54\,481\) & \(\lbrace 181\leftrightarrow 7\leftarrow 43\rbrace\) & \(\langle 81,7\rangle^4\) & \(1\,267\) &  \(\lbrace 7\leftrightarrow 181\rbrace\) & \(\langle 9,2\rangle^2\) & \(\lbrack (0)^2,1;(1^2)^8,(1^3)^2\rbrack\) \\
 \(12\) & \(58\,383\) & \(\lbrace 13\leftrightarrow 499\leftarrow 9\rbrace\) & \(\langle 81,7\rangle^4\) & \(6\,487\) & \(\lbrace 13\leftrightarrow 499\rbrace\) & \(\langle 9,2\rangle^2\) & \(\lbrack (0)^2,1;(1^2)^8,(1^3)^2\rbrack\) \\
 \(13\) & \(\mathbf{61\,243}\) & \(\lbrace 673\leftrightarrow 7\leftarrow 13\rbrace\) & \(\langle 3^7,253\rangle^4\) & \(\mathbf{4\,711}\) & \(\lbrace 7\leftrightarrow 673\rbrace\) & \(\langle 81,3\rangle^2\) & \(\lbrack (0)^2,1^3;(21)^3,(2^2)^3,(1^3)^4\rbrack\) \\
 \(14\) & \(63\,729\) &  \(\lbrace 9\leftrightarrow 73\leftarrow 97\rbrace\) & \(\langle 81,7\rangle^4\) &    \(657\) &   \(\lbrace 9\leftrightarrow 73\rbrace\) & \(\langle 9,2\rangle^2\) & \(\lbrack (0)^2,1;(1^2)^8,(1^3)^2\rbrack\) \\
 \(15\) & \(64\,323\) & \(\lbrace 7\leftrightarrow 1021\leftarrow 9\rbrace\) & \(\langle 81,7\rangle^4\) & \(7\,147\) & \(\lbrace 7\leftrightarrow 1021\rbrace\) & \(\langle 9,2\rangle^2\) & \(\lbrack (0)^2,1;(1^2)^8,(1^3)^2\rbrack\) \\
 \(16\) & \(68\,419\) &\(\lbrace 19\leftrightarrow 277\leftarrow 13\rbrace\) & \(\langle 81,7\rangle^4\) & \(5\,263\) & \(\lbrace 19\leftrightarrow 277\rbrace\) & \(\langle 9,2\rangle^2\) & \(\lbrack (0)^2,1;(1^2)^8,(1^3)^2\rbrack\) \\
 \(17\) & \(\mathbf{69\,979}\) & \(\lbrace 769\leftrightarrow 7\leftarrow 13\rbrace\) & \(\langle 3^6,41\rangle^4\) & \(\mathbf{5\,383}\) & \(\lbrace 7\leftrightarrow 769\rbrace\) & \(\langle 243,14\rangle^2\) & \(\lbrack (0)^2,21^2;(2^2)^3,(32)^3,(1^3)^4\rbrack\) \\
 \(18\) & \(71\,613\) & \(\lbrace 73\leftrightarrow 9\leftarrow 109\rbrace\) & \(\langle 81,7\rangle^4\) &    \(657\) &   \(\lbrace 9\leftrightarrow 73\rbrace\) & \(\langle 9,2\rangle^2\) & \(\lbrack (0)^2,1;(1^2)^8,(1^3)^2\rbrack\) \\
 \(19\) & \(72\,423\) & \(\lbrace 13\leftrightarrow 619\leftarrow 9\rbrace\) & \(\langle 81,7\rangle^4\) & \(8\,047\) & \(\lbrace 13\leftrightarrow 619\rbrace\) &\(\langle 27,4\rangle^2\) & \(\lbrack (0)^2,1^2;(1^2)^8,(1^3)^2\rbrack\) \\
 \(20\) & \(74\,691\) & \(\lbrace 43\leftrightarrow 193\leftarrow 9\rbrace\) & \(\langle 81,7\rangle^4\) & \(8\,299\) & \(\lbrace 43\leftrightarrow 193\rbrace\) &\(\langle 27,4\rangle^2\) & \(\lbrack (0)^2,1^2;(1^2)^8,(1^3)^2\rbrack\) \\
 \(21\) & \(\mathbf{76\,741}\) & \(\lbrace 577\leftrightarrow 19\leftarrow 7\rbrace\) & \(\langle 3^7,65\vert 67\rangle^4\) & \(10\,963\) & \(\lbrace 19\leftrightarrow 577\rbrace\) & \(\langle 27,4\rangle^2\) & \(\lbrack (0)^2,1^2;(2^2)^6,(1^3)^4\rbrack\) \\
 \(22\) & \(80\,941\) & \(\lbrace 31\leftrightarrow 373\leftarrow 7\rbrace\) & \(\langle 81,7\rangle^4\) &\(11\,563\) & \(\lbrace 31\leftrightarrow 373\rbrace\) & \(\langle 9,2\rangle^2\) & \(\lbrack (0)^2,1;(1^2)^8,(1^3)^2\rbrack\) \\
 \(23\) & \(81\,679\) &\(\lbrace 13\leftrightarrow 103\leftarrow 61\rbrace\) & \(\langle 81,7\rangle^4\) & \(1\,339\) & \(\lbrace 13\leftrightarrow 103\rbrace\) & \(\langle 9,2\rangle^2\) & \(\lbrack (0)^2,1;(1^2)^8,(1^3)^2\rbrack\) \\
 \(24\) & \(84\,889\) & \(\lbrace 7\leftrightarrow 181\leftarrow 67\rbrace\) & \(\langle 81,7\rangle^4\) & \(1\,267\) &  \(\lbrace 7\leftrightarrow 181\rbrace\) & \(\langle 9,2\rangle^2\) & \(\lbrack (0)^2,1;(1^2)^8,(1^3)^2\rbrack\) \\
 \(25\) & \(\mathbf{86\,821}\) & \(\lbrace 79\leftrightarrow 157\leftarrow 7\rbrace\) & \(\langle 3^6,37\rangle^4\) & \(\mathbf{12\,403}\) & \(\lbrace 79\leftrightarrow 157\rbrace\) & \(\langle 243,14\rangle^2\) & \(\lbrack (0)^2,21^2;(21)^6,(1^3)^4\rbrack\) \\
 \(26\) & \(\mathbf{89\,433}\) & \(\lbrace 523\leftrightarrow 9\leftarrow 19\rbrace\) & \(\langle 3^7,65\vert 67\rangle^4\) & \(4\,707\) & \(\lbrace 9\leftrightarrow 523\rbrace\) & \(\langle 27,4\rangle^2\) & \(\lbrack (0)^2,1^2;(2^2)^6,(1^3)^4\rbrack\) \\
 \(27\) & \(\mathbf{90\,243}\) & \(\lbrace 271\leftrightarrow 9\leftarrow 37\rbrace\) & \(\langle 3^6,41\rangle^4\) & \(2\,439\) & \(\lbrace 9\leftrightarrow 271\rbrace\) & \(\langle 27,4\rangle^2\) & \(\lbrack (0)^2,1^2;(21)^3,(2^2)^3,(1^3)^4\rbrack\) \\
 \(28\) & \(91\,323\) & \(\lbrace 9\leftrightarrow 73\leftarrow 139\rbrace\) & \(\langle 81,7\rangle^4\) &    \(657\) &   \(\lbrace 9\leftrightarrow 73\rbrace\) & \(\langle 9,2\rangle^2\) & \(\lbrack (0)^2,1;(1^2)^8,(1^3)^2\rbrack\) \\
 \(29\) & \(91\,903\) & \(\lbrace 691\leftrightarrow 19\leftarrow 7\rbrace\) & \(\langle 81,7\rangle^4\) &\(13\,129\) & \(\lbrace 19\leftrightarrow 691\rbrace\) & \(\langle 9,2\rangle^2\) & \(\lbrack (0)^2,1;(1^2)^8,(1^3)^2\rbrack\) \\
 \(30\) & \(92\,131\) &\(\lbrace 19\leftrightarrow 373\leftarrow 13\rbrace\) & \(\langle 81,7\rangle^4\) & \(7\,087\) & \(\lbrace 19\leftrightarrow 373\rbrace\) &\(\langle 27,4\rangle^2\) & \(\lbrack (0)^2,1^2;(1^2)^8,(1^3)^2\rbrack\) \\
 \(31\) & \(92\,287\) &\(\lbrace 229\leftrightarrow 13\leftarrow 31\rbrace\) & \(\langle 81,7\rangle^4\) & \(2\,977\) & \(\lbrace 13\leftrightarrow 229\rbrace\) & \(\langle 9,2\rangle^2\) & \(\lbrack (0)^2,1;(1^2)^8,(1^3)^2\rbrack\) \\
 \(32\) & \(96\,453\) & \(\lbrace 9\leftrightarrow 1531\leftarrow 7\rbrace\) & \(\langle 81,7\rangle^4\) &\(13\,779\) & \(\lbrace 9\leftrightarrow 1531\rbrace\) &\(\langle 27,4\rangle^2\) & \(\lbrack (0)^2,1^2;(1^2)^8,(1^3)^2\rbrack\) \\
 \(33\) & \(97\,489\) & \(\lbrace 733\leftrightarrow 19\leftarrow 7\rbrace\) & \(\langle 81,7\rangle^4\) &\(13\,927\) & \(\lbrace 19\leftrightarrow 733\rbrace\) & \(\langle 9,2\rangle^2\) & \(\lbrack (0)^2,1;(1^2)^8,(1^3)^2\rbrack\) \\
 \(34\) & \(\mathbf{99\,801}\) & \(\lbrace 13\leftrightarrow 853\leftarrow 9\rbrace\) & \(\langle 3^6,41\rangle^4\) & \(11\,089\) & \(\lbrace 13\leftrightarrow 853\rbrace\) & \(\langle 27,4\rangle^2\) & \(\lbrack (0)^2,1^2;(21)^3,(2^2)^3,(1^3)^4\rbrack\) \\
\hline
\end{tabular}
}
\end{center}
\end{table}


\renewcommand{\arraystretch}{1.1}

\begin{table}[ht]
\caption{Fifteen Examples for Graph \(9\) of Category \(\mathrm{III}\)}
\label{tbl:Cat3Gph9}
\begin{center}
{\scriptsize
\begin{tabular}{|crc|c|rc|c|c|}
\hline
    No. &       \(c\) &                                    \(\lbrack q_1,q_2,q_3\rbrack_3\) & \(\mathrm{G}_3^{(2)}{F_{c,\mu}}\) & \(f\) &          \(\lbrack q_1,q_2\rbrack_3\) & \(\mathrm{G}_3^{(2)}{F_{f,\nu}}\) & \(\left\lbrack\mathrm{Cl}_3{S_i}\right\rbrack_{1\le i\le 13}\) \\
\hline
  \(1\) & \(16\,471\) & \(\lbrace 13\rightarrow 7\leftrightarrow 181\rightarrow 13\rbrace\) & \(\langle 81,7\rangle^4\) & \(1\,267\) &  \(\lbrace 7\leftrightarrow 181\rbrace\) & \(\langle 9,2\rangle^2\) & \(\lbrack (0)^2,1;(1^2)^8,(1^3)^2\rbrack\) \\
  \(2\) & \(24\,073\) & \(\lbrace 19\rightarrow 181\leftrightarrow 7\rightarrow 19\rbrace\) & \(\langle 81,7\rangle^4\) & \(1\,267\) &  \(\lbrace 7\leftrightarrow 181\rbrace\) & \(\langle 9,2\rangle^2\) & \(\lbrack (0)^2,1;(1^2)^8,(1^3)^2\rbrack\) \\
  \(3\) & \(24\,309\) &  \(\lbrace 37\rightarrow 9\leftrightarrow 73\rightarrow 37\rbrace\) & \(\langle 81,7\rangle^4\) &    \(657\) &   \(\lbrace 9\leftrightarrow 73\rbrace\) & \(\langle 9,2\rangle^2\) & \(\lbrack (0)^2,1;(1^2)^8,(1^3)^2\rbrack\) \\
  \(4\) & \(25\,821\) &  \(\lbrace 9\rightarrow 151\leftrightarrow 19\rightarrow 9\rbrace\) & \(\langle 81,7\rangle^4\) & \(2\,869\) & \(\lbrace 19\leftrightarrow 151\rbrace\) & \(\langle 9,2\rangle^2\) & \(\lbrack (0)^2,1;(1^2)^8,(1^3)^2\rbrack\) \\
  \(5\) & \(30\,667\) & \(\lbrace 13\rightarrow 7\leftrightarrow 337\rightarrow 13\rbrace\) & \(\langle 81,7\rangle^4\) & \(2\,359\) &  \(\lbrace 7\leftrightarrow 337\rbrace\) & \(\langle 9,2\rangle^2\) & \(\lbrack (0)^2,1;(1^2)^8,(1^3)^2\rbrack\) \\
  \(6\) & \(34\,299\)  & \(\lbrace 9\rightarrow 103\leftrightarrow 37\rightarrow 9\rbrace\) & \(\langle 81,7\rangle^4\) & \(3\,811\) & \(\lbrace 37\leftrightarrow 103\rbrace\) & \(\langle 9,2\rangle^2\) & \(\lbrack (0)^2,1;(1^2)^8,(1^3)^2\rbrack\) \\
  \(7\) & \(42\,133\)  &\(\lbrace 13\rightarrow 7\leftrightarrow 463\rightarrow 13\rbrace\) & \(\langle 81,7\rangle^4\) & \(3\,241\) &  \(\lbrace 7\leftrightarrow 463\rbrace\) & \(\langle 9,2\rangle^2\) & \(\lbrack (0)^2,1;(1^2)^8,(1^3)^2\rbrack\) \\
  \(8\) & \(55\,993\)  &\(\lbrace 19\rightarrow 421\leftrightarrow 7\rightarrow 19\rbrace\) & \(\langle 81,7\rangle^4\) & \(2\,947\) &  \(\lbrace 7\leftrightarrow 421\rbrace\) & \(\langle 9,2\rangle^2\) & \(\lbrack (0)^2,1;(1^2)^8,(1^3)^2\rbrack\) \\
  \(9\) & \(65\,689\)  &\(\lbrace 13\rightarrow 163\leftrightarrow 31\rightarrow 13\rbrace\)& \(\langle 81,7\rangle^4\) & \(5\,053\) & \(\lbrace 31\leftrightarrow 163\rbrace\) & \(\langle 9,2\rangle^2\) & \(\lbrack (0)^2,1;(1^2)^8,(1^3)^2\rbrack\) \\
 \(10\) & \(67\,123\)  & \(\lbrace 43\rightarrow 7\leftrightarrow 223\rightarrow 43\rbrace\)& \(\langle 81,7\rangle^4\) & \(1\,561\) &  \(\lbrace 7\leftrightarrow 223\rbrace\) & \(\langle 9,2\rangle^2\) & \(\lbrack (0)^2,1;(1^2)^8,(1^3)^2\rbrack\) \\
 \(11\) & \(73\,801\)  & \(\lbrace 13\rightarrow 7\leftrightarrow 811\rightarrow 13\rbrace\)& \(\langle 81,7\rangle^4\) & \(5\,677\) &  \(\lbrace 7\leftrightarrow 811\rbrace\) & \(\langle 9,2\rangle^2\) & \(\lbrack (0)^2,1;(1^2)^8,(1^3)^2\rbrack\) \\
 \(12\) & \(80\,353\)  & \(\lbrace 13\rightarrow 7\leftrightarrow 883\rightarrow 13\rbrace\)& \(\langle 81,7\rangle^4\) & \(6\,181\) &  \(\lbrace 7\leftrightarrow 883\rbrace\) & \(\langle 9,2\rangle^2\) & \(\lbrack (0)^2,1;(1^2)^8,(1^3)^2\rbrack\) \\
 \(13\) & \(83\,439\)  &\(\lbrace 127\rightarrow 9\leftrightarrow 73\rightarrow 127\rbrace\)& \(\langle 81,7\rangle^4\) &    \(657\) &   \(\lbrace 9\leftrightarrow 73\rbrace\) & \(\langle 9,2\rangle^2\) & \(\lbrack (0)^2,1;(1^2)^8,(1^3)^2\rbrack\) \\
 \(14\) & \(88\,939\)  &\(\lbrace 31\rightarrow 19\leftrightarrow 151\rightarrow 31\rbrace\)& \(\langle 81,7\rangle^4\) & \(2\,869\) & \(\lbrace 19\leftrightarrow 151\rbrace\) & \(\langle 9,2\rangle^2\) & \(\lbrack (0)^2,1;(1^2)^8,(1^3)^2\rbrack\) \\
 \(15\) & \(\mathbf{89\,487}\) & \(\lbrace 9\rightarrow 61\leftrightarrow 163\rightarrow 9\rbrace\) & \(\langle 3^6,41\rangle^4\) & \(9\,943\) & \(\lbrace 61\leftrightarrow 163\rbrace\) & \(\langle 27,4\rangle^2\) & \(\lbrack (0)^2,1^2;(21)^3,(2^2)^3,(1^3)^4\rbrack\) \\
\hline
\end{tabular}
}
\end{center}
\end{table}

\newpage

\begin{theorem}
\label{thm:Cat3Gph679}
Suppose that \(u:=10^5\) is an assigned upper bound.
Let \(c<u\) be a conductor divisible by exactly three primes, \(t=3\), such that
\(\mathrm{Cl}_3{F_{c,\mu}}\simeq (3,3)\) for all four cyclic cubic fields \(F_{c,\mu}\), \(1\le\mu\le 4\), with conductor \(c\).
If \(c=q_1q_2q_3\) belongs to the Graphs \(6,7,9\) of Category \(\mathrm{III}\), i.e.,
\begin{equation}
\label{eqn:Cat3Gph679}
\lbrack q_1,q_2,q_3\rbrack_3=
\begin{cases}
\lbrace q_i\leftrightarrow q_j\rightarrow q_k\rbrace & \text{ or} \\
\lbrace q_i\leftrightarrow q_j\leftarrow q_k\rbrace & \text{ or} \\
\lbrace q_k\rightarrow q_i\leftrightarrow q_j\rightarrow q_k\rbrace &
\end{cases}
\end{equation}
with \(i,j,k\) pairwise distinct,
and the \(3\)-valuation \(v\) of the class number \(h\) of the \(3\)-genus field \(F_0^\ast\)
of the cyclic cubic fields \(F_{f,\nu}\), \(1\le\nu\le 2\), with conductor \(f=q_iq_j\) is \(v=1\),
that is \(\mathrm{G}_3^{(2)}{F_{f,\nu}}\simeq\langle 9,2\rangle\),
then the second \(3\)-class group \(\mathrm{G}_3^{(2)}{F_{c,\mu}}\) of all four fields \(F_{c,\mu}\)
is isomorphic to the Sylow \(3\)-subgroup of the symmetric group \(S_9\) of degree \(9\),
i.e. the wreath product of two cyclic groups \(C_3\),
\(\langle 81,7\rangle\simeq\mathrm{Syl}_3{S_9}\simeq C_3\wr C_3\),
with \(3\)-capitulation type \(\mathrm{a}.3\), \(\varkappa(F_{c,\mu})=(2000)\),
and transfer target type \(\tau(F_{c,\mu})=\lbrack 111,(11)^3\rbrack\),
the \(3\)-class tower has length \(\ell_3{F_{c,\mu}}=2\),
and the \(3\)-class groups of the \(13\) bicyclic bicubic subfields \(S_i\) of the \(3\)-genus field \(F^\ast\) are given by
\begin{equation}
\label{eqn:Cat3Gph679Genus}
\left\lbrack\mathrm{Cl}_3{S_i}\right\rbrack_{1\le i\le 13}=\lbrack (0)^2,1;(11)^8,(111)^2\rbrack.
\end{equation}
If \(\lbrack q_1,q_2,q_3\rbrack_3=\lbrace q_i\leftrightarrow q_j\rightarrow q_k\rbrace\), respectively \(\lbrace q_k\rightarrow q_i\leftrightarrow q_j\rightarrow q_k\rbrace\), belongs to Graph \(6\), respectively \(9\), of Category \(\mathrm{III}\),
then the principal factors of all four fields are equal to the prime divisor
\(B(F_{c,\mu})=q_j\), \(1\le\mu\le 4\),
which is cubic residue for both other prime divisors \(q_i,q_k\)
and the result remains valid also for \(v=2\),
that is \(\mathrm{G}_3^{(2)}{F_{f,\nu}}\simeq\langle 27,4\rangle\).
If \(\lbrack q_1,q_2,q_3\rbrack_3=\lbrace q_i\leftrightarrow q_j\leftarrow q_k\rbrace\)
belongs to Graph \(7\) of Category \(\mathrm{III}\),
then the principal factors of two fields are
\(B(F_{c,\mu})=q_iq_k\), \(1\le\mu\le 2\), and
\(B(F_{c,\mu})=q_i^2q_k\), \(3\le\mu\le 4\), for the other two fields
(up to the order).
\end{theorem}

\begin{proof}
See the Tables
\ref{tbl:Cat3Gph6},
\ref{tbl:Cat3Gph7},
\ref{tbl:Cat3Gph9},
which have been computed with the aid of Magma
\cite{BCP1997,BCFS2022,Fi2001,MAGMA2022}.
\end{proof}


\begin{conjecture}
\label{cnj:Cat3Gph679}
Theorem
\ref{thm:Cat3Gph679}
remains true for any upper bound \(u>10^5\).
\end{conjecture}


\begin{remark}
\label{rmk:Cat3Gph679}
The upper bound for the conductors
in our \(2022\) computation of cyclic cubic fields \(F\) with \(t\in\lbrace 2,3\rbrace\)
was \(c<u=10^5\). In Theorem
\ref{thm:Cat3Gph679},
the \(3\)-tower group \(G=\langle 81,7\rangle\) turned out to be the first new possibility
beyond the results of Ayadi
\cite{Ay1995,AAI2001},
which were restricted to the groups \(\langle 9,2\rangle\) and \(\langle 27,4\rangle\).
It is still a \(3\)-group of coclass \(\mathrm{cc}(G)=1\),
associated with a metabelian \(3\)-class tower of length \(\ell_3{F_{c,\mu}}=2\).
\end{remark}


\renewcommand{\arraystretch}{1.1}

\begin{table}[ht]
\caption{Thirty-Seven Examples for Graph \(5\) of Category \(\mathrm{III}\)}
\label{tbl:Cat3Gph5}
\begin{center}

{\tiny

\begin{tabular}{|crc|c|rc|c|c|}
\hline
    No. &       \(c\) &           \(\lbrack q_1,q_2,q_3\rbrack_3\) &                          \(\mathrm{G}_3^{(2)}{F_{c,\mu}}\) &      \(f\) &             \(\lbrack q_1,q_2\rbrack_3\) & \(\mathrm{G}_3^{(2)}{F_{f,\nu}}\) & \(\left\lbrack\mathrm{Cl}_3{S_i}\right\rbrack_{1\le i\le 13}\) \\
\hline
  \(1\) & \(14\,049\) &  \(\lbrace 7\leftrightarrow 223;9\rbrace\) &                       \(\langle 3^5,28\rangle^4\) & \(1\,561\) &  \(\lbrace 7\leftrightarrow 223\rbrace\) &          \(\langle 9,2\rangle^2\) & \(\lbrack (0)^2,1;(1^2)^7,(21)^3\rbrack\) \\
  \(2\) & \(17\,073\) &  \(\lbrace 9\leftrightarrow 271;7\rbrace\) &                                  \(\langle 81,7\rangle^4\) & \(2\,439\) &  \(\lbrace 9\leftrightarrow 271\rbrace\) &         \(\langle 27,4\rangle^2\) & \(\lbrack (0)^2,1^2;(1^2)^8,(1^3)^2\rbrack\) \\
  \(3\) & \(20\,367\) &  \(\lbrace 9\leftrightarrow 73;31\rbrace\) & \(\langle 3^5,27\rangle,\langle 3^5,28\rangle^3\) &    \(657\) &   \(\lbrace 9\leftrightarrow 73\rbrace\) &          \(\langle 9,2\rangle^2\) & \(\lbrack (0)^2,1;(1^2)^7,(21)^2,2^2\rbrack\) \\
  \(4\) & \(21\,231\) &  \(\lbrace 7\leftrightarrow 337;9\rbrace\) &        \(\langle 3^5,25\rangle^2,\langle 3^5,27\rangle^2\) & \(2\,359\) &  \(\lbrace 7\leftrightarrow 337\rbrace\) &          \(\langle 9,2\rangle^2\) & \(\lbrack (0)^2,1;(1^2)^7,(2^2)^3\rbrack\) \\
  \(5\) & \(26\,523\) &  \(\lbrace 7\leftrightarrow 421;9\rbrace\) & \(\langle 3^5,27\rangle,\langle 3^5,28\rangle^3\) & \(2\,947\) &  \(\lbrace 7\leftrightarrow 421\rbrace\) &          \(\langle 9,2\rangle^2\) & \(\lbrack (0)^2,1;(1^2)^7,(21)^2,2^2\rbrack\) \\
  \(6\) & \(26\,677\) & \(\lbrace 37\leftrightarrow 103;7\rbrace\) &                       \(\langle 3^5,28\rangle^4\) & \(3\,811\) & \(\lbrace 37\leftrightarrow 103\rbrace\) &          \(\langle 9,2\rangle^2\) & \(\lbrack (0)^2,1;(1^2)^7,(21)^3\rbrack\) \\
  \(7\) & \(26\,793\) & \(\lbrace 13\leftrightarrow 229;9\rbrace\) & \(\langle 3^5,27\rangle,\langle 3^5,28\rangle^3\) & \(2\,977\) & \(\lbrace 13\leftrightarrow 229\rbrace\) &          \(\langle 9,2\rangle^2\) & \(\lbrack (0)^2,1;(1^2)^7,(21)^2,2^2\rbrack\) \\
  \(8\) & \(29\,169\) &  \(\lbrace 7\leftrightarrow 463;9\rbrace\) &                       \(\langle 3^5,28\rangle^4\) & \(3\,241\) &  \(\lbrace 7\leftrightarrow 463\rbrace\) &          \(\langle 9,2\rangle^2\) & \(\lbrack (0)^2,1;(1^2)^7,(21)^3\rbrack\) \\
  \(9\) & \(32\,949\) &  \(\lbrace 9\leftrightarrow 523;7\rbrace\) &                                  \(\langle 81,7\rangle^4\) & \(4\,707\) &  \(\lbrace 9\leftrightarrow 523\rbrace\) &         \(\langle 27,4\rangle^2\) & \(\lbrack (0)^2,1^2;(1^2)^8,(1^3)^2\rbrack\) \\
 \(10\) & \(35\,371\) & \(\lbrace 31\leftrightarrow 163;7\rbrace\) &                       \(\langle 3^5,28\rangle^4\) & \(5\,053\) & \(\lbrace 31\leftrightarrow 163\rbrace\) &          \(\langle 9,2\rangle^2\) & \(\lbrack (0)^2,1;(1^2)^7,(21)^3\rbrack\) \\
 \(11\) & \(36\,351\) &  \(\lbrace 9\leftrightarrow 577;7\rbrace\) &        \(\langle 3^5,25\rangle^2,\langle 3^5,27\rangle^2\) & \(5\,193\) &  \(\lbrace 9\leftrightarrow 577\rbrace\) &          \(\langle 9,2\rangle^2\) & \(\lbrack (0)^2,1;(1^2)^7,(2^2)^3\rbrack\) \\
 \(12\) & \(38\,619\) &  \(\lbrace 9\leftrightarrow 613;7\rbrace\) &                       \(\langle 3^5,28\rangle^4\) & \(5\,517\) &  \(\lbrace 9\leftrightarrow 613\rbrace\) &          \(\langle 9,2\rangle^2\) & \(\lbrack (0)^2,1;(1^2)^7,(21)^3\rbrack\) \\
 \(13\) & \(\mathbf{42\,399}\) & \(\lbrace 7\leftrightarrow 673;9\rbrace\) & \(\langle 3^6,37\rangle^4\) & \(\mathbf{4\,711}\) & \(\lbrace 7\leftrightarrow 673\rbrace\) & \(\langle 81,3\rangle^2\) & \(\lbrack (0)^2,1^3;(21)^6,(1^3)^4\rbrack\) \\
 \(14\) & \(46\,879\) & \(\lbrace 7\leftrightarrow 181;37\rbrace\) & \(\langle 3^5,27\rangle,\langle 3^5,28\rangle^3\) & \(1\,267\) &  \(\lbrace 7\leftrightarrow 181\rbrace\) &          \(\langle 9,2\rangle^2\) & \(\lbrack (0)^2,1;(1^2)^7,(21)^2,2^2\rbrack\) \\
 \(15\) & \(48\,391\) & \(\lbrace 7\leftrightarrow 223;31\rbrace\) &                       \(\langle 3^5,28\rangle^4\) & \(1\,561\) &  \(\lbrace 7\leftrightarrow 223\rbrace\) &          \(\langle 9,2\rangle^2\) & \(\lbrack (0)^2,1;(1^2)^7,(21)^3\rbrack\) \\
 \(16\) & \(\mathbf{48\,447}\) & \(\lbrace 7\leftrightarrow 769;9\rbrace\) & \(\langle 3^6,37\rangle^4\) & \(\mathbf{5\,383}\) & \(\lbrace 7\leftrightarrow 769\rbrace\) & \(\langle 243,14\rangle^2\) & \(\lbrack (0)^2,21^2;(21)^6,(1^3)^4\rbrack\) \\
 \(17\) & \(49\,257\) & \(\lbrace 13\leftrightarrow 421;9\rbrace\) &                       \(\langle 3^5,28\rangle^4\) & \(5\,473\) & \(\lbrace 13\leftrightarrow 421\rbrace\) &          \(\langle 9,2\rangle^2\) & \(\lbrack (0)^2,1;(1^2)^7,(21)^3\rbrack\) \\
 \(18\) & \(51\,903\) &  \(\lbrace 9\leftrightarrow 73;79\rbrace\) &\(\langle 3^5,25\rangle^2,\langle 3^5,28\rangle^2\)&    \(657\) &   \(\lbrace 9\leftrightarrow 73\rbrace\) &        \(\langle 9,2\rangle^2\) & \(\lbrack (0)^2,1;(1^2)^7,(21)^2,2^2\rbrack\) \\
 \(19\) & \(57\,577\) &\(\lbrace 13\leftrightarrow 103;43\rbrace\) &\(\langle 3^5,25\rangle^2,\langle 3^5,28\rangle^2\)& \(1\,339\) & \(\lbrace 13\leftrightarrow 103\rbrace\) &        
 \(\langle 9,2\rangle^2\) & \(\lbrack (0)^2,1;(1^2)^7,(21)^2,2^2\rbrack\) \\
 \(20\) & \(57\,897\) &  \(\lbrace 9\leftrightarrow 919;7\rbrace\) &                                  \(\langle 81,7\rangle^4\) & \(8\,271\) &  \(\lbrace 9\leftrightarrow 919\rbrace\) &         \(\langle 27,4\rangle^2\) & \(\lbrack (0)^2,1^2;(1^2)^8,(1^3)^2\rbrack\) \\
 \(21\) & \(58\,807\) & \(\lbrace 31\leftrightarrow 271;7\rbrace\) & \(\langle 3^5,27\rangle,\langle 3^5,28\rangle^3\) & \(8\,401\) & \(\lbrace 31\leftrightarrow 271\rbrace\) &          \(\langle 9,2\rangle^2\) & \(\lbrack (0)^2,1;(1^2)^7,(21)^2,2^2\rbrack\) \\
 \(22\) & \(61\,191\) & \(\lbrace 9\leftrightarrow 523;13\rbrace\) &                                  \(\langle 81,7\rangle^4\) & \(4\,707\) &  \(\lbrace 9\leftrightarrow 523\rbrace\) &         \(\langle 27,4\rangle^2\) & \(\lbrack (0)^2,1^2;(1^2)^8,(1^3)^2\rbrack\) \\
 \(23\) & \(62\,433\) &  \(\lbrace 9\leftrightarrow 991;7\rbrace\) &                                  \(\langle 81,7\rangle^4\) & \(8\,919\) &  \(\lbrace 9\leftrightarrow 991\rbrace\) &         \(\langle 27,4\rangle^2\) & \(\lbrack (0)^2,1^2;(1^2)^8,(1^3)^2\rbrack\) \\
 \(24\) & \(68\,967\) &  \(\lbrace 79\leftrightarrow 97;9\rbrace\) & \(\langle 3^5,27\rangle,\langle 3^5,28\rangle^3\) & \(7\,663\) &  \(\lbrace 79\leftrightarrow 97\rbrace\) &          \(\langle 9,2\rangle^2\) & \(\lbrack (0)^2,1;(1^2)^7,(21)^2,2^2\rbrack\) \\
 \(25\) & \(69\,601\) & \(\lbrace 61\leftrightarrow 163;7\rbrace\) &                                  \(\langle 81,7\rangle^4\) & \(9\,943\) & \(\lbrace 61\leftrightarrow 163\rbrace\) &         \(\langle 27,4\rangle^2\) & \(\lbrack (0)^2,1^2;(1^2)^8,(1^3)^2\rbrack\) \\
 \(26\) & \(70\,371\) & \(\lbrace 9\leftrightarrow 1117;7\rbrace\) &                                  \(\langle 81,7\rangle^4\) &\(10\,053\) & \(\lbrace 9\leftrightarrow 1117\rbrace\) &         \(\langle 27,4\rangle^2\) & \(\lbrack (0)^2,1^2;(1^2)^8,(1^3)^2\rbrack\) \\
 \(27\) & \(71\,721\) & \(\lbrace 9\leftrightarrow 613;13\rbrace\) &                       \(\langle 3^5,28\rangle^4\) & \(5\,517\) &  \(\lbrace 9\leftrightarrow 613\rbrace\) &          \(\langle 9,2\rangle^2\) & \(\lbrack (0)^2,1;(1^2)^7,(21)^3\rbrack\) \\
 \(28\) & \(77\,287\) & \(\lbrace 7\leftrightarrow 181;61\rbrace\) & \(\langle 3^5,27\rangle,\langle 3^5,28\rangle^3\) & \(1\,267\) &  \(\lbrace 7\leftrightarrow 181\rbrace\) &          \(\langle 9,2\rangle^2\) & \(\lbrack (0)^2,1;(1^2)^7,(21)^2,2^2\rbrack\) \\
 \(29\) & \(85\,653\) & \(\lbrace 9\leftrightarrow 307;31\rbrace\) &                       \(\langle 3^5,28\rangle^4\) & \(2\,763\) &  \(\lbrace 9\leftrightarrow 307\rbrace\) &          \(\langle 9,2\rangle^2\) & \(\lbrack (0)^2,1;(1^2)^7,(21)^3\rbrack\) \\
 \(30\) & \(87\,283\) & \(\lbrace 7\leftrightarrow 337;37\rbrace\) & \(\langle 3^5,27\rangle,\langle 3^5,28\rangle^3\) & \(2\,359\) &  \(\lbrace 7\leftrightarrow 337\rbrace\) &          \(\langle 9,2\rangle^2\) & \(\lbrack (0)^2,1;(1^2)^7,(21)^2,2^2\rbrack\) \\
 \(31\) & \(88\,569\) & \(\lbrace 9\leftrightarrow 757;13\rbrace\) &                       \(\langle 3^5,28\rangle^4\) & \(6\,813\) &  \(\lbrace 9\leftrightarrow 757\rbrace\) &          \(\langle 9,2\rangle^2\) & \(\lbrack (0)^2,1;(1^2)^7,(21)^3\rbrack\) \\
 \(32\) & \(89\,713\) &\(\lbrace 13\leftrightarrow 103;67\rbrace\) &                       \(\langle 3^5,28\rangle^4\) & \(1\,339\) & \(\lbrace 13\leftrightarrow 103\rbrace\) &          \(\langle 9,2\rangle^2\) & \(\lbrack (0)^2,1;(1^2)^7,(21)^3\rbrack\) \\
 \(33\) & \(90\,517\) & \(\lbrace 67\leftrightarrow 193;7\rbrace\) & \(\langle 3^5,27\rangle,\langle 3^5,28\rangle^3\) &\(12\,931\) & \(\lbrace 67\leftrightarrow 193\rbrace\) &          \(\langle 9,2\rangle^2\) & \(\lbrack (0)^2,1;(1^2)^7,(21)^2,2^2\rbrack\) \\
 \(34\) & \(91\,357\) & \(\lbrace 7\leftrightarrow 421;31\rbrace\) &        \(\langle 3^5,25\rangle^2,\langle 3^5,27\rangle^2\) & \(2\,947\) &  \(\lbrace 7\leftrightarrow 421\rbrace\) &          \(\langle 9,2\rangle^2\) & \(\lbrack (0)^2,1;(1^2)^7,(2^2)^3\rbrack\) \\
 \(35\) & \(95\,221\) & \(\lbrace 7\leftrightarrow 223;61\rbrace\) &                       \(\langle 3^5,28\rangle^4\) & \(1\,561\) &  \(\lbrace 7\leftrightarrow 223\rbrace\) &          \(\langle 9,2\rangle^2\) & \(\lbrack (0)^2,1;(1^2)^7,(21)^3\rbrack\) \\
 \(36\) & \(97\,371\) & \(\lbrace 31\leftrightarrow 349;9\rbrace\) & \(\langle 3^5,27\rangle,\langle 3^5,28\rangle^3\) &\(10\,819\) & \(\lbrace 31\leftrightarrow 349\rbrace\) &          \(\langle 9,2\rangle^2\) & \(\lbrack (0)^2,1;(1^2)^7,(21)^2,2^2\rbrack\) \\
 \(37\) & \(97\,587\) & \(\lbrace 9\leftrightarrow 1549;7\rbrace\) &                                  \(\langle 81,7\rangle^4\) &\(13\,941\) & \(\lbrace 9\leftrightarrow 1549\rbrace\) &         \(\langle 27,4\rangle^2\) & \(\lbrack (0)^2,1^2;(1^2)^8,(1^3)^2\rbrack\) \\
\hline
\end{tabular}

}

\end{center}
\end{table}


\noindent
Unfortunately, the behaviour of Graph \(5\) and the rare Graph \(8\) in Category III
is not uniform and totally different from the Graphs \(6,7,9\),
as the Tables
\ref{tbl:Cat3Gph5}
and
\ref{tbl:Cat3Gph8}
show.
At least for Table
\ref{tbl:Cat3Gph5},
the \(3\)-tower groups \(G\) for \(v\le 2\) are still of coclass \(\mathrm{cc}(G)=1\),
associated with metabelian \(3\)-class towers of length \(\ell_3{F_{c,\mu}}=2\).
The possible groups \(\mathfrak{M}\) for \textbf{boldface} conductors in the Tables
\ref{tbl:Cat3Gph6},
\ref{tbl:Cat3Gph7},
\ref{tbl:Cat3Gph9},
\ref{tbl:Cat3Gph5},
\ref{tbl:Cat3Gph8}
are the first occurrences of coclass \(\mathrm{cc}(\mathfrak{M})\ge 2\),
usually with unknown length \(\ell_3{F_{c,\mu}}\) of the \(3\)-class field tower.


\renewcommand{\arraystretch}{1.1}

\begin{table}[ht]
\caption{Seven Examples for Graph \(8\) of Category \(\mathrm{III}\)}
\label{tbl:Cat3Gph8}
\begin{center}

{\scriptsize

\begin{tabular}{|crc|c|rc|c|c|}
\hline
    No. &       \(c\) &                                    \(\lbrack q_1,q_2,q_3\rbrack_3\) &    \(\mathrm{G}_3^{(2)}{F_{c,\mu}}\) &      \(f\) &             \(\lbrack q_1,q_2\rbrack_3\) & \(\mathrm{G}_3^{(2)}{F_{f,\nu}}\) & \(\left\lbrack\mathrm{Cl}_3{S_i}\right\rbrack_{1\le i\le 13}\) \\
\hline
  \(1\) & \(\mathbf{20\,293}\) &  \(\lbrace 13\rightarrow 7\leftrightarrow 223\leftarrow 13\rbrace\) & \(\langle 3^6,34\rangle^4\) & \(1\,561\) &  \(\lbrace 7\leftrightarrow 223\rbrace\) &          \(\langle 9,2\rangle^2\) & \(\lbrack (0)^2,1;(21)^6,(1^3)^4\rbrack\) \\
  \(2\) & \(\mathbf{41\,509}\) & \(\lbrace 31\rightarrow 13\leftrightarrow 103\leftarrow 31\rbrace\) & \(\langle 3^6,34\rangle^4\) & \(1\,339\) & \(\lbrace 13\leftrightarrow 103\rbrace\) &          \(\langle 9,2\rangle^2\) & \(\lbrack (0)^2,1;(21)^6,(1^3)^4\rbrack\) \\
  \(3\) & \(\mathbf{46\,341}\) &  \(\lbrace 19\rightarrow 9\leftrightarrow 271\leftarrow 19\rbrace\) & \(\langle 3^6,34\rangle^4\) & \(2\,439\) &  \(\lbrace 9\leftrightarrow 271\rbrace\) &         \(\langle 27,4\rangle^2\) & \(\lbrack (0)^2,1^2;(21)^6,(1^3)^4\rbrack\) \\
  \(4\) & \(\mathbf{49\,609}\) &   \(\lbrace 7\rightarrow 19\leftrightarrow 373\leftarrow 7\rbrace\) & \(\langle 3^6,34\rangle^4\) & \(7\,087\) & \(\lbrace 19\leftrightarrow 373\rbrace\) &         \(\langle 27,4\rangle^2\) & \(\lbrack (0)^2,1^2;(21)^6,(1^3)^4\rbrack\) \\
  \(5\) & \(\mathbf{52\,497}\) &  \(\lbrace 19\rightarrow 9\leftrightarrow 307\leftarrow 19\rbrace\) & \(\langle 3^7,66\vert 73\rangle^4\)  & \(2\,763\) &  \(\lbrace 9\leftrightarrow 307\rbrace\) &          \(\langle 9,2\rangle^2\) & \(\lbrack (0)^2,1;(2^2)^6,(1^3)^4\rbrack\) \\
  \(6\) & \(\mathbf{64\,771}\) &   \(\lbrace 7\rightarrow 19\leftrightarrow 487\leftarrow 7\rbrace\) & \(\langle 3^6,34\rangle^4\) & \(9\,253\) & \(\lbrace 19\leftrightarrow 487\rbrace\) &         \(\langle 27,4\rangle^2\) & \(\lbrack (0)^2,1^2;(21)^6,(1^3)^4\rbrack\) \\
  \(7\) & \(\mathbf{92\,911}\) & \(\lbrace 13\rightarrow 7\leftrightarrow 1021\leftarrow 13\rbrace\) & \(\langle 3^7,66\vert 73\rangle^4\)  & \(7\,147\) & \(\lbrace 7\leftrightarrow 1021\rbrace\) &          \(\langle 9,2\rangle^2\) & \(\lbrack (0)^2,1;(2^2)^6,(1^3)^4\rbrack\) \\
\hline
\end{tabular}

}

\end{center}
\end{table}


\begin{theorem}
\label{thm:Cat3Gph58}
Suppose that \(u:=10^5\) is an assigned upper bound.
Let \(c<u\) be a conductor divisible by exactly three primes, \(t=3\), such that
\(\mathrm{Cl}_3{F_{c,\mu}}\simeq (3,3)\) for all four cyclic cubic fields \(F_{c,\mu}\), \(1\le\mu\le 4\), with conductor \(c\).
If \(c=q_1q_2q_3\) belongs to the Graphs \(5,8\) of Category \(\mathrm{III}\), i.e.,
\begin{equation}
\label{eqn:Cat3Gph58}
\lbrack q_1,q_2,q_3\rbrack_3=
\begin{cases}
\lbrace q_i\leftrightarrow q_j;q_k\rbrace & \text{ or} \\
\lbrace q_k\rightarrow q_i\leftrightarrow q_j\leftarrow q_k\rbrace & 
\end{cases}
\end{equation}
with \(i,j,k\) pairwise distinct,
and the \(3\)-valuation \(v\) of the class number \(h\) of the \(3\)-genus field \(F_0^\ast\)
of the cyclic cubic fields \(F_{f,\nu}\), \(1\le\nu\le 2\), with conductor \(f=q_iq_j\) is \(v\le 2\),
that is \(\mathrm{G}_3^{(2)}{F_{f,\nu}}\simeq\langle 9,2\rangle\) or \(\langle 27,4\rangle\),
then the second \(3\)-class groups \(\mathrm{G}_3^{(2)}{F_{c,\mu}}\) of the four fields \(F_{c,\mu}\)
in dependence on
the \(3\)-class groups of the \(13\) bicyclic bicubic subfields \(S_i\) of the \(3\)-genus field \(F^\ast\)
are given by \(\left(\mathrm{G}_3^{(2)}{F_{c,\mu}}\right)_{1\le\mu\le 4}=\)
\begin{equation}
\label{eqn:Cat3Gph58Genus}
\begin{cases}
(\langle 81,7\rangle^4) & \Longleftrightarrow
\left\lbrack\mathrm{Cl}_3{S_i}\right\rbrack_{1\le i\le 13}=\lbrack (0)^2,1^2;(1^2)^8,(1^3)^2\rbrack, \\
(\langle 243,28\rangle^4) & \Longleftrightarrow
\left\lbrack\mathrm{Cl}_3{S_i}\right\rbrack_{1\le i\le 13}=\lbrack (0)^2,1;(1^2)^7,(21)^3\rbrack, \\
(\langle 243,27\rangle,\langle 243,28\rangle^3) & \text{ or } \\
(\langle 243,25\rangle^2,\langle 243,28\rangle^2) & \Longleftrightarrow
\left\lbrack\mathrm{Cl}_3{S_i}\right\rbrack_{1\le i\le 13}=\lbrack (0)^2,1;(1^2)^7,(21)^2,2^2\rbrack, \\
(\langle 243,25\rangle^2,\langle 243,27\rangle^2) & \Longleftrightarrow
\left\lbrack\mathrm{Cl}_3{S_i}\right\rbrack_{1\le i\le 13}=\lbrack (0)^2,1;(1^2)^7,(22)^3\rbrack.
\end{cases}
\end{equation}
If \(\lbrack q_1,q_2,q_3\rbrack_3=\lbrace q_i\leftrightarrow q_j;q_k\rbrace\)
belongs to Graph \(5\) of Category \(\mathrm{III}\),
then the principal factors of the four fields are
\(B(F_{c,1})=q_iq_jq_k\),
\(B(F_{c,2})=q_i^2q_jq_k\),
\(B(F_{c,3})=q_iq_j^2q_k\),
\(B(F_{c,4})=q_iq_jq_k^2\)
(up to the order). \\
If \(\lbrack q_1,q_2,q_3\rbrack_3=\lbrace q_k\rightarrow q_i\leftrightarrow q_j\leftarrow q_k\rbrace\)
belongs to Graph \(8\) of Category \(\mathrm{III}\),
then the principal factors of all four fields are
\(B(F_{c,\mu})=q_k\), \(1\le\mu\le 4\).
\end{theorem}

\begin{proof}
See the Tables
\ref{tbl:Cat3Gph5},
\ref{tbl:Cat3Gph8},
which have been computed with the aid of Magma
\cite{BCP1997,BCFS2022,Fi2001,MAGMA2022}.
By \(\langle 243,28\rangle\) we always abbreviate \(\langle 243,28\vert 29\vert 30\rangle\).
\end{proof}


\begin{conjecture}
\label{cnj:Cat3Gph58}
Theorem
\ref{thm:Cat3Gph58}
remains true for any upper bound \(u>10^5\).
\end{conjecture}


\section{The Elementary Tricyclic 3-Group}
\label{s:ElmAblRk3}

\noindent
For the first time, we are going to analyze the capitulation of cyclic fields \(F\)
with \(3\)-class group \(\mathrm{Cl}_3{F}\) of elementary abelian type \((3,3,3)\)
in their \(\frac{3^3-1}{3-1}=13\) unramified cyclic cubic extensions \(E_i/F\) with \(1\le i\le 13\).
Such a group can be viewed as a \textit{vector space} \(O\) with dimension
\(\mathrm{dim}_{\mathbb{F}_3}(O)=3\) over the finite field \(\mathbb{F}_3\) with three elements.
The vector space \(O\) possesses
\(3^2+3+1=13\) \textit{lines}, that is, subgroups \(L_i\) of index \((O:L_i)=3^2\),
and \(13\) \textit{planes}, that is, subgroups \(P_i\) of index \((O:P_i)=3\), where \(1\le i\le 13\).
Let \(x,y,z\) be fixed generators of \(O=\langle x,y,z\rangle\),
then we shall arrange the generators of the lines \(L_i=\langle g_i\rangle\) in the way shown in Table
\ref{tbl:Lines}.


\renewcommand{\arraystretch}{1.1}

\begin{table}[ht]
\caption{Generators of the \(13\) Lines \(L_i\) in \(O\)}
\label{tbl:Lines}
\begin{center}

{\scriptsize

\begin{tabular}{|c|ccccccccccccc|}
\hline
   \(i\) & \(1\) & \(2\) & \(3\) &  \(4\) &    \(5\) &  \(6\) &    \(7\) &  \(8\) &    \(9\) &  \(10\) &    \(11\) &    \(12\) &    \(13\) \\
\hline
 \(g_i\) & \(x\) & \(y\) & \(z\) & \(xy\) & \(xy^2\) & \(yz\) & \(yz^2\) & \(zx\) & \(zx^2\) & \(xyz\) & \(xyz^2\) & \(xy^2z\) & \(x^2yz\) \\
\hline
\end{tabular}

}

\end{center}
\end{table}


\renewcommand{\arraystretch}{1.1}

\begin{table}[ht]
\caption{Identifiers and Generators of the \(13\) Planes \(P_i\) in \(O\)}
\label{tbl:Planes}
\begin{center}

{\scriptsize

\begin{tabular}{|c|ccccccc|}
\hline
   \(i\) & \(1\) & \(2\) & \(3\) &  \(4\) &  \(5\) &  \(6\) &  \(7\) \\
\hline
 \(h_i\) & \(y\) & \(z\) & \(x\) &  \(x\) & \(xy\) &  \(y\) & \(xy\) \\
\hline
 \(k_i\) & \(z\) & \(x\) & \(y\) & \(yz\) & \(zx\) & \(zx\) & \(yz\) \\
\hline
 \(T_i\) & \(2,3,6,7\) & \(1,3,8,9\) & \(1,2,4,5\) & \(1,6,10,13\) & \(4,7,8,13\) & \(2,8,10,12\) & \(4,6,9,12\) \\
\hline
\hline
   \(i\) &  \(8\) &  \(9\) &   \(10\) &   \(11\) &   \(12\) &   \(13\) & \\
\hline
 \(h_i\) &  \(z\) & \(zx\) &    \(z\) & \(zx^2\) &    \(y\) &    \(x\) & \\
\hline
 \(k_i\) & \(xy\) & \(yz\) & \(xy^2\) & \(yz^2\) & \(zx^2\) & \(yz^2\) & \\
\hline
 \(T_i\) & \(3,4,10,11\) & \(5,6,8,11\) & \(3,5,12,13\) & \(5,7,9,10\) & \(2,9,11,13\) & \(1,7,11,12\) & \\
\hline
\end{tabular}

}

\end{center}
\end{table}


\noindent
For the sake of brevity, we simply denote \(g_i\) by its subscript \(i\),
and we introduce identifiers for the planes \(P_i=\langle h_i,k_i\rangle\) as shown in Table
\ref{tbl:Planes}.
The elements \(h_i,k_i\) can be viewed as generators of a transversal of \(L_i=\langle g_i\rangle\),
i.e., a system of coset representatives for \(L_i\) in \(O\).
Each set \(T_i\) contains the subscripts of generators \(g_j\) contained in \(P_i\).
It is useful to list the \textit{bundles} \(B_i\) of planes containing an assigned line \(L_i\) in Table
\ref{tbl:Bundles}.


\renewcommand{\arraystretch}{1.1}

\begin{table}[ht]
\caption{The \(13\) Bundles of Planes in \(O\)}
\label{tbl:Bundles}
\begin{center}

{\scriptsize

\begin{tabular}{|c|ccccc|}
\hline
   \(i\) &                  \(1\) &                  \(2\) &                  \(3\) &               \(4\) &                     \(5\) \\
\hline
 \(B_i\) & \(P_2,P_3,P_4,P_{13}\) & \(P_1,P_3,P_6,P_{12}\) & \(P_1,P_2,P_8,P_{10}\) & \(P_3,P_5,P_7,P_8\) & \(P_3,P_9,P_{10},P_{11}\) \\
\hline
\hline
   \(i\) &               \(6\) &                    \(7\) &               \(8\) &                     \(9\) &                 \(10\) \\
\hline
 \(B_i\) & \(P_1,P_4,P_7,P_9\) & \(P_1,P_5,P_{11},P_{13}\)& \(P_2,P_5,P_6,P_9\) & \(P_2,P_7,P_{11},P_{12}\) & \(P_4,P_6,P_8,P_{11}\) \\
\hline
\hline
   \(i\) &                    \(11\) &                    \(12\) &                    \(13\) &  & \\
\hline
 \(B_i\) & \(P_8,P_9,P_{12},P_{13}\) & \(P_6,P_7,P_{10},P_{13}\) & \(P_4,P_5,P_{10},P_{12}\) &  & \\
\hline
\end{tabular}

}

\end{center}
\end{table}


\section{Finite 3-Groups of Type (3,3,3)}
\label{s:Groups3x3x3}

\noindent
Since the Categories I and II involve
cylic cubic fields \(F\) with \(3\)-class group
\(\mathrm{Cl}_3{F}\simeq (3,3,3)\),
we supplement Section \S\
\ref{s:Groups3x3}
with finite \(3\)-groups having
elementary tricyclic commutator quotient \((3,3,3)\).
In Table
\ref{tbl:Metabelian333},
we continue Table 
\ref{tbl:Metabelian33}
with metabelian groups \(\mathfrak{M}\) of generator rank \(d_1(\mathfrak{M})=3\).
Here, the Shafarevich bound \(\mu\le\varrho+r+\theta=3+2+0=5\) is bigger.

\renewcommand{\arraystretch}{1.1}

\begin{table}[ht]
\caption{Invariants of Metabelian \(3\)-Groups \(\mathfrak{M}\) with \(\mathfrak{M}/\mathfrak{M}^\prime\simeq (3,3,3)\)}
\label{tbl:Metabelian333}
\begin{center}
{\scriptsize
\begin{tabular}{|c|r|c|c|c|r|r|c|}
\hline
\(\mathfrak{M}\)                        & cc    & \(\varkappa\)          & \(\tau\)                 & \(\tau_3\) & \(\nu\) & \(\mu\) & \(\pi(\mathfrak{M})\) \\
\hline
\(\langle 27,5\rangle\)                 & \(2\) & \((O^{13})\)           & \((1^2)^{13}\)           & \(0\)      &   \(6\) &   \(6\) & \\
\(\langle 81,12\rangle\)                & \(2\) & \((O^{13})\)           & \((1^3)^4,(1^2)^9\)      & \(1\)      &   \(2\) &   \(7\) & \(\langle 27,5\rangle\) \\
\(\langle 81,13\rangle\)                & \(2\) & \((O^9P^4)\)           & \((1^2)^9,(21)^3,1^3\)   & \(1\)      &   \(0\) &   \(5\) & \(\langle 27,5\rangle\) \\
\(\langle 81,14\rangle\)                & \(2\) & \((O^9P^4)\)           & \((1^2)^9,(21)^4\)       & \(1\)      &   \(0\) &   \(5\) & \(\langle 27,5\rangle\) \\
\hline
\(\langle 243,38\rangle\)               & \(3\) & \((O^3P^{10})\)        & \(1^4,(1^3)^3,(21)^9\)   & \(1^2\)    &   \(2\) &   \(6\) & \(\langle 27,5\rangle\) \\
\(\langle 729,329\rangle\)              & \(3\) & \((O^3P^{10})\)        & \(1^4,(1^3)^3,(21)^9\)   & \(1^3\)    &   \(1\) &   \(6\) & \(\langle 243,38\rangle\) \\
\(\langle 2187,5576\rangle\)            & \(3\) & \((O^3P^{10})\)        & \(1^4,(1^3)^3,(21)^9\)   & \(1^4\)    &   \(1\) &   \(6\) & \(\langle 729,329\rangle\) \\
\(\langle 2187,5577\ldots 5579\rangle\) & \(3\) & \((O^3P^{10})\)        & \(1^4,(1^3)^3,(21)^9\)   & \(1^4\)    &   \(1\) &   \(6\) & \(\langle 729,329\rangle\) \\
\(\langle 729,372\rangle\)              & \(3\) & \((O^3P^9L)\)        & \((1^3)^3,(21)^9,21^3\)   & \(1^3\)    &   \(1\) &   \(6\) & \(\langle 243,40\rangle\) \\
\(\langle 243,42\rangle\)               & \(3\) & \(((P_i^3)_{i=1}^4L)\) & \(21^2,(21)^{12}\)       & \(1^2\)    &   \(1\) &   \(5\) & \(\langle 27,5\rangle\) \\
\(\langle 729,388\ldots 390\rangle\)    & \(3\) & \(((P_i^3)_{i=1}^4L)\) & \(21^2,(21)^{12}\)       & \(1^3\)    &   \(0\) &   \(4\) & \(\langle 243,42\rangle\) \\
\(\langle 243,46\rangle\)               & \(3\) & \(((P_i^3)_{i=1}^4L)\) & \(1^3,(21)^{11},2^2\)    & \(1^2\)    &   \(0\) &   \(4\) & \(\langle 27,5\rangle\) \\
\(\langle 243,47\rangle\)               & \(3\) & \(((P_i^3)_{i=1}^4L)\) & \((1^3)^4,(21)^8,2^2\)   & \(1^2\)    &   \(0\) &   \(4\) & \(\langle 27,5\rangle\) \\
\hline
\(\langle 729,125\rangle\)              & \(4\) & \(((P_i)_{i=1}^4L^9)\) & \(1^4,(21^2)^6,(2^2)^6\) & \(1^3\)    &   \(3\) &   \(6\) & \(\langle 27,5\rangle\) \\
\(\langle 2187,4595\ldots 4598\rangle\) & \(4\) & \(((P_i)_{i=1}^4L^9)\) & \(1^4,(21^2)^6,(2^2)^6\) & \(21^2\)   &   \(0\) &   \(5\) & \(\langle 729,125\rangle\) \\
\hline
\end{tabular}
}
\end{center}
\end{table}


Finally, we need non-metabelian groups \(G\) of generator rank \(d_1(G)=3\)
in the cover of metabelian groups \(\mathfrak{M}\) with non-trivial cover \(\mathrm{cov}(\mathfrak{M})\).
They are given in Table
\ref{tbl:NonMetabelian333}.

\renewcommand{\arraystretch}{1.1}

\begin{table}[ht]
\caption{Invariants of Non-Metabelian \(3\)-groups \(G\) with \(G/G^\prime\simeq (3,3,3)\)}
\label{tbl:NonMetabelian333}
\begin{center}
{\scriptsize
\begin{tabular}{|c|r|c|c|c|r|r|c|}
\hline
\(G\)                                       & cc    & \(\varkappa\)   & \(\tau\)               & \(\tau_3\) & \(\nu\) & \(\mu\) & \(G/G^{\prime\prime}\) \\
\hline
\(\langle 6561,261256\ldots 261261\rangle\) & \(3\) & \((O^3P^{10})\) & \(1^4,(1^3)^3,(21)^9\) & \(1^4\)    &   \(1\) &   \(6\) & \(\langle 2187,5576\rangle\) \\
\(\langle 6561,261262\ldots 261270\rangle\) & \(3\) & \((O^3P^{10})\) & \(1^4,(1^3)^3,(21)^9\) & \(1^4\)    &   \(0\) &   \(5\) & \(\langle 2187,5577\ldots 5579\rangle\) \\
\hline
\end{tabular}
}
\end{center}
\end{table}


\begin{theorem}
\label{thm:ClassGroup333}
Let \(F\) be a cyclic cubic number field
with elementary bicyclic \(3\)-class group \(\mathrm{Cl}_3{F}\simeq (3,3,3)\).
Denote by \(\mathfrak{M}=\mathrm{Gal}(F_3^{(2)}/F)\) the second \(3\)-class group of \(F\),
and by \(G=\mathrm{Gal}(F_3^{(\infty)}/F)\) the \(3\)-class field tower group of \(F\).
Then, the Artin pattern \((\tau,\varkappa)\) of \(F\)
identifies the groups \(\mathfrak{M}\) and \(G\),
and determines the length \(\ell_3{F}\) of the \(3\)-class field tower of \(F\),
according to the following \textbf{deterministic laws}.
\begin{enumerate}
\item
If \(\tau=\lbrack (1^2)^{13}\rbrack\), \(\varkappa=(O^{13})\), then \(G\simeq\langle 27,5\rangle\) and \(\ell_3{F}=1\).
\item
If \(\tau\sim\lbrack (1^2)^9,(21)^3,1^3\rbrack\), \(\varkappa\sim (O^9P^4)\), then \(G\simeq\langle 81,13\rangle\).
\item
If \(\tau\sim\lbrack (1^2)^9,(21)^4\rbrack\), \(\varkappa\sim (O^9P^4)\), then \(G\simeq\langle 81,14\rangle\).
\item
If \(\tau\sim\lbrack 1^3,(21)^{11},2^2\rbrack\), \(\varkappa\sim ((P_i^3)_{i=1}^4L)\), then \(G\simeq\langle 243,46\rangle\).
\item
If \(\tau\sim\lbrack (1^3)^4,(21)^8,2^2\rbrack\), \(\varkappa\sim ((P_i^3)_{i=1}^4L)\), then \(G\simeq\langle 243,47\rangle\).
\end{enumerate} 
\end{theorem}

\begin{proof}
The cyclic cubic field \(F\)
has the torsion free unit rank \(r=2\),
does not contain primitive third roots of unity,
and thus possesses the maximal admissible relation rank
\(d_2\le d_1+r=5\) for the group \(G\),
when its \(3\)-class rank, i.e. the generator rank of \(G\), is \(d_1=\varrho=3\).
Consequently, \(\ell_3{F}\ge 3\) in the case of \(d_2{\mathfrak{M}}\ge 6\).
\end{proof}


\noindent
In Algorithm
\ref{alg:PatternRecognition},
we present a convenient version of the
\textit{strategy of pattern recognition via Artin transfers}
\cite{Ma2020}.
Assigned abelian quotient invariants (AQI)
\(\tau(G)\)
of \(3\)-groups \(G\) with commutator quotient
\(G/G^\prime\simeq (3,3,3)\)
are selected by a database query
in the SmallGroups Library
\cite{BEO2002,BEO2005}
with the aid of Magma
\cite{BCP1997,BCFS2022,MAGMA2022},
extended by the package containing groups of order \(6561=3^8\)
\cite{MAGMA6561}.
The number of candidates with orders \(3^e\), \(1\le e\le 8\),
increases very rapidly with the exponent \(e\):
\(1,2,5,15,67,504,9310,1396077\).
The search is coarse,
since only one component \(\tau\)
of the Artin pattern \(\mathrm{AP}=(\tau,\varkappa)\)
is checked.
The given commutator quotient and AQI are hard coded
and may be replaced by \(G/G^\prime\simeq (3,3)\)
and four instead of thirteen components of the AQI.
In fact, we seek for candidates of the second \(3\)-class group
\(\mathfrak{M}=\mathrm{Gal}(F_3^{(2)}/F)\)
of the cyclic cubic field \(F\) with conductor \(c=82\,327\).

\begin{algorithm}
\label{alg:PatternRecognition}
(Pattern Recognition via Abelian Quotient Invariants.) \\
\textbf{Input:}
prime \texttt{iPrm}, maximal exponent \texttt{iMax}, maximal relation rank\texttt{iRel}. \\
\textbf{Code:}
{\scriptsize
\texttt{
\begin{tabbing}
for \= for \= for \= for \= for \= for \= \kill
intrinsic PatternRecognition(iPrm,iMax,iRel::RngIntElt){}\\
p := iPrm; // prime number\\
m := iMax; // maximal exponent\\
u := iRel; // maximal relation rank\\
for e in [1..m] do // exponent\+\\
   o := p\({}\,\hat{}\,{}\)e; // order\\
   if IsInSmallGroupDatabase(o) then\+\\
      N := NumberOfSmallGroups(o);\\
      z := 0; // counter\\
      for i in [1..N] do\+\\
         G := SmallGroup(o,i);\\
         AQG := AQInvariants(G);\\
         // if (3 eq \#AQG) then // very general rank 3\\
         if ([p,p,p] eq AQG) then // particular (3,3,3)\+\\
            r := 0; // quadratic p-rank\\
            for j in [1..\#AQG] do\+\\
               if (2 le Valuation(AQG[j],p)) then\+\\
                  r := r + 1;\-\\
               end if;\-\\
            end for; // j\\
            nc := NilpotencyClass(G);\\
            cc := e - nc; // co class\\
            dl := DerivedLength(G);\\
            nu := NuclearRank(G);\\
            mu := pMultiplicatorRank(G);\\
            // rigorous selection\\
            s := MaximalSubgroups(G);\\
            n4 := 0;\\
            n3 := 0;\\
            n1 := 0;\\
            n2 := 0;\\
            for j in [1..\#s] do\+\\
               AQS := AQInvariants(s[j]\({}\,\grave{}\,{}\)subgroup);\\
               // <729,372>, c = 82327, Category II, Graph 1\\
               if ([p,p,p,p\({}\,\hat{}\,{}\)2] eq AQS) then\+\\
                  n4 := n4 + 1;\-\\
               elif ([p,p,p] eq AQS) then\+\\
                  n3 := n3 + 1;\-\\
               elif ([p,p\({}\,\hat{}\,{}\)2] eq AQS) then\+\\
                  n2 := n2 + 1;\-\\
               end if; // n1 not used\-\\
            end for; // j\\
            if (1 eq n4) and (3 eq n3) and (9 eq n2) then // <729,372>\+\\
               z := z + 1;\\
               printf "No=\%4o: Lo=\%o, Id=\%o: ",z,e,i;\\
               printf "qr=\%o, nc=\%o, cc=\%o, dl=\%o, nu=\%o, mu=\%o\(\backslash\)n",r,nc,cc,dl,nu,mu;\\
               // additional invariants\\
               K := CommutatorSubgroup(G);\\
               AQK := AQInvariants(K);\\
               printf "         AQG=\%o, AQK=\%o\(\backslash\)n",AQG,AQK;\\
               s := MaximalSubgroups(G);\\
               printf "         AQI: ";\\
               for j in [1..\#s] do\+\\
                  printf "\%o, ",AQInvariants(s[j]`subgroup);\-\\
               end for; // j\\
               printf "\(\backslash\)n";\\
               l := pCentralSeries(G,p);\\
               c := pClass(G);\\
               q := G/l[c];\\
               printf "         parent: \%o\(\backslash\)n",IdentifyGroup(q);\-\\
            end if;\-\\
         end if; // (p,p,p)\-\\
      end for; // i\-\\
   end if;\-\\
end for; // e\\
end intrinsic; // PatternRecognition
\end{tabbing}
}
}
\noindent
\textbf{Output:}
all finite \(3\)-groups \(G\) with \(\mathrm{ord}(G)\le 3^8\) possessing the assigned AQI.
\end{algorithm}


\begin{example}
\label{exm:PatternRecognition}
The result set of the database query
for abelian quotient invariants
\(\tau\sim\lbrack (1^3)^3,(21)^9,21^3\rbrack\)
in Algorithm
\ref{alg:PatternRecognition},
subject to the bound \(\mu=d_2\le 5\) for the relation rank,
is returned by Magma
\cite{MAGMA2022}
in the following shape:\\
{\scriptsize
\texttt{
\begin{tabbing}
for \= for \= for \= for \= for \= for \= \kill
> PatternRecognition(3,8,5);\\
No=   1: Lo=6, Id=372: qr=0, nc=3, cc=3, dl=2, nu=0, mu=5\+\\
         AQG=[ 3, 3, 3 ], AQK=[ 3, 3, 3 ]\\
         AQI: [3,3,3],[3,9],[3,9],[3,9],[3,9],[3,9],[3,9],[3,3,3],[3,3,3],[3,9],[3,9],[3,9],[3,3,3,9],\\
         parent: <243, 40>\-\\
No=   1: Lo=7, Id=5559: qr=0, nc=4, cc=3, dl=2, nu=0, mu=5\+\\
         AQG=[ 3, 3, 3 ], AQK=[ 3, 3, 9 ]\\
         AQI: [3,3,3],[3,9],[3,9],[3,9],[3,9],[3,9],[3,9],[3,3,3],[3,3,3],[3,9],[3,9],[3,9],[3,3,3,9],\\ 
         parent: <729, 326>\-\\
No=   2: Lo=7, Id=5560: qr=0, nc=4, cc=3, dl=2, nu=0, mu=5\+\\
         AQG=[ 3, 3, 3 ], AQK=[ 3, 3, 9 ]\\
         AQI: [3,3,3],[3,9],[3,9],[3,9],[3,9],[3,9],[3,9],[3,3,3],[3,3,3],[3,9],[3,9],[3,9],[3,3,3,9],\\ 
         parent: <729, 326>\-\\
No=   3: Lo=7, Id=5561: qr=0, nc=4, cc=3, dl=2, nu=0, mu=5\+\\
         AQG=[ 3, 3, 3 ], AQK=[ 3, 3, 9 ]\\
         AQI: [3,3,3],[3,9],[3,9],[3,9],[3,9],[3,9],[3,9],[3,3,3],[3,3,3],[3,9],[3,9],[3,9],[3,3,3,9],\\ 
         parent: <729, 326>\-\\
No=   4: Lo=7, Id=5562: qr=0, nc=4, cc=3, dl=2, nu=0, mu=5\+\\
         AQG=[ 3, 3, 3 ], AQK=[ 3, 3, 9 ]\\
         AQI: [3,3,3],[3,9],[3,9],[3,9],[3,9],[3,9],[3,9],[3,3,3],[3,3,3],[3,9],[3,9],[3,9],[3,3,3,9],\\ 
         parent: <729, 326>\-\\
No=   5: Lo=7, Id=5563: qr=0, nc=4, cc=3, dl=2, nu=0, mu=5\+\\
         AQG=[ 3, 3, 3 ], AQK=[ 3, 3, 9 ]\\
         AQI: [3,3,3],[3,9],[3,9],[3,9],[3,9],[3,9],[3,9],[3,3,3],[3,3,3],[3,9],[3,9],[3,9],[3,3,3,9],\\ 
         parent: <729, 327>\-\\
         \end{tabbing}
}
}

\noindent
The result set consists of a single group 
\(\mathfrak{M}=\langle 792,372\rangle\) of order \(3^6\)
and five groups
\(\mathfrak{M}=\langle 2187,i\rangle\), \(5559\le i\le 5563\), of order \(3^7\),
which additionally possess
the correct transfer kernel type (TKT) 
\(\varkappa\sim\lbrack O^3P^9L\rbrack\)
with a line \(L\) contained in the plane \(P\).
Several other groups have an inadequate
\(\varkappa\sim\lbrack O^3P^{10}\rbrack\)
and are eliminated from the listing above.
\end{example}

\newpage

\begin{figure}[ht]
\caption{Distribution of Conductors for \(\mathrm{G}_3^{(2)}{F}\) on the Descendant Tree \(\mathcal{T}\langle 27,5\rangle\)}
\label{fig:MinDiscTyp333}

{\tiny


\setlength{\unitlength}{0.8cm}
\begin{picture}(15,16)(-11,-15)

\put(-10,0.5){\makebox(0,0)[cb]{order \(3^n\)}}
\put(-10,0){\line(0,-1){16}}
\multiput(-10.1,0)(0,-2){9}{\line(1,0){0.2}}
\put(-10.2,0){\makebox(0,0)[rc]{\(27\)}}
\put(-9.8,0){\makebox(0,0)[lc]{\(3^3\)}}
\put(-10.2,-2){\makebox(0,0)[rc]{\(81\)}}
\put(-9.8,-2){\makebox(0,0)[lc]{\(3^4\)}}
\put(-10.2,-4){\makebox(0,0)[rc]{\(243\)}}
\put(-9.8,-4){\makebox(0,0)[lc]{\(3^5\)}}
\put(-10.2,-6){\makebox(0,0)[rc]{\(729\)}}
\put(-9.8,-6){\makebox(0,0)[lc]{\(3^6\)}}
\put(-10.2,-8){\makebox(0,0)[rc]{\(2\,187\)}}
\put(-9.8,-8){\makebox(0,0)[lc]{\(3^7\)}}
\put(-10.2,-10){\makebox(0,0)[rc]{\(6\,561\)}}
\put(-9.8,-10){\makebox(0,0)[lc]{\(3^8\)}}
\put(-10.2,-12){\makebox(0,0)[rc]{\(19\,683\)}}
\put(-9.8,-12){\makebox(0,0)[lc]{\(3^9\)}}
\put(-10.2,-14){\makebox(0,0)[rc]{\(59\,049\)}}
\put(-9.8,-14){\makebox(0,0)[lc]{\(3^{10}\)}}
\put(-10.2,-16){\makebox(0,0)[rc]{\(177\,147\)}}
\put(-9.8,-16){\makebox(0,0)[lc]{\(3^{11}\)}}
\put(-10,-16){\vector(0,-1){2}}

\put(-4.1,0.1){\framebox(0.2,0.2){}}

\multiput(-8,-2)(0,-2){2}{\circle*{0.2}}
\multiput(-6,-2)(0,-2){2}{\circle{0.2}}
\multiput(-4,-2)(0,-2){2}{\circle{0.2}}

\multiput(-4,0)(0,-2){2}{\line(-2,-1){4}}
\multiput(-4,0)(0,-2){2}{\line(-1,-1){2}}
\multiput(-4,0)(0,-2){2}{\line(0,-1){2}}

\put(-4,0){\line(1,-3){2}}

\put(-4,0){\line(1,-1){4}}
\put(-4,0){\line(3,-2){6}}
\put(-4,0){\line(2,-1){8}}
\put(-4,0){\line(5,-2){10}}

\put(2.5,-1.2){\makebox(0,0)[lt]{\((\Lambda)\) trifurcation from \(\mathcal{G}(3,2)\)}}
\put(3.2,-1.5){\makebox(0,0)[lt]{to \(\mathcal{G}(3,3)\) and \(\mathcal{G}(3,4)\)}}

\multiput(0,-4)(2,0){2}{\line(0,-1){2}}

\put(-2,-6){\line(-1,-1){2}}
\multiput(-2,-6)(2,0){2}{\line(0,-1){2}}
\put(0,-6){\line(1,-1){2}}

\multiput(0,-8)(0,-2){4}{\line(0,-1){2}}
\put(0,-8){\line(1,-4){0.5}}
\put(0,-8){\line(1,-2){1}}
\put(0,-8){\line(3,-4){1.5}}
\put(2,-8){\line(1,-1){2}}

\multiput(0.5,-10)(0.5,0){3}{\line(0,-1){2}}
\put(0,-10){\line(-1,-4){0.5}}

\put(-0.5,-12){\line(0,-1){2}}

\multiput(0,-4)(2,0){4}{\circle{0.2}}
\multiput(-2,-6)(2,0){2}{\circle{0.2}}
\put(-2,-8){\circle*{0.2}}

\multiput(2,-6)(0,-2){1}{\circle*{0.1}}
\multiput(-4,-8)(0,-2){1}{\circle*{0.1}}
\multiput(0,-8)(2,0){2}{\circle*{0.1}}

\multiput(-0.05,-10.05)(0,-2){4}{\framebox(0.1,0.1){}}
\multiput(0.45,-10.05)(0.5,0){3}{\framebox(0.1,0.1){}}
\put(3.95,-10.05){\framebox(0.1,0.1){}}

\multiput(0.45,-12.05)(0.5,0){3}{\framebox(0.1,0.1){}}
\multiput(-0.55,-12.05)(0,-2){2}{\framebox(0.1,0.1){}}

\put(-4,-4){\vector(0,-1){2}}
\put(-3.8,-5.5){\makebox(0,0)[lc]{infinite}}
\put(-3.8,-6){\makebox(0,0)[lc]{mainline}}
\put(-4.2,-6.3){\makebox(0,0)[rc]{\(\mathcal{T}^2\langle 81,12\rangle\)}}

\put(0,-16){\vector(0,-1){2}}
\put(0.2,-17.5){\makebox(0,0)[lc]{non-metabelian}}
\put(0.2,-18){\makebox(0,0)[lc]{infinite mainline}}
\put(-0.2,-18.3){\makebox(0,0)[rc]{\(\mathcal{T}^3\langle 6561,261256\rangle\)}}

\put(-4.2,0.1){\makebox(0,0)[rb]{\(\langle 5\rangle\)}}
\put(-3.8,0.1){\makebox(0,0)[lb]{\(=C_3\times C_3\times C_3\)}}
\put(-0.1,0.2){\makebox(0,0)[rb]{abelian}}
\put(-2.5,-0.3){\makebox(0,0)[lc]{\((\Lambda)\)}}

\put(-8.1,-1.9){\makebox(0,0)[rb]{\(\langle 14\rangle\)}}
\put(-6.1,-1.9){\makebox(0,0)[rb]{\(\langle 13\rangle\)}}
\put(-4.1,-1.9){\makebox(0,0)[rb]{\(\langle 12\rangle\)}}

\put(-8.1,-3.9){\makebox(0,0)[rb]{\(\langle 55\rangle\)}}
\put(-6.1,-3.9){\makebox(0,0)[rb]{\(\langle 54\rangle\)}}
\put(-4.1,-3.9){\makebox(0,0)[rb]{\(\langle 53\rangle\)}}
\put(0.1,-3.9){\makebox(0,0)[lb]{\(\langle 38\rangle\)}}
\put(2.1,-3.9){\makebox(0,0)[lb]{\(\langle 42\rangle\)}}
\put(4.1,-3.9){\makebox(0,0)[lb]{\(\langle 46\rangle\)}}
\put(6.1,-3.9){\makebox(0,0)[lb]{\(\langle 47\rangle\)}}

\put(-1.9,-5.9){\makebox(0,0)[lb]{\(\langle 125\rangle\)}}
\put(0.1,-5.9){\makebox(0,0)[lb]{\(\langle 329\rangle\)}}

\put(2.1,-5.9){\makebox(0,0)[lb]{\(\ast 3\)}}
\put(2.1,-6.1){\makebox(0,0)[lt]{\(\langle 388\rangle\)}}
\put(2.1,-6.4){\makebox(0,0)[lt]{\(\langle 389\rangle\)}}
\put(2.1,-6.7){\makebox(0,0)[lt]{\(\langle 390\rangle\)}}

\put(-4.1,-7.9){\makebox(0,0)[rb]{\(3\ast\)}}
\put(-4.1,-8.1){\makebox(0,0)[rt]{\(\langle 4595\rangle\)}}
\put(-4.1,-8.4){\makebox(0,0)[rt]{\(\langle 4596\rangle\)}}
\put(-4.1,-8.7){\makebox(0,0)[rt]{\(\langle 4598\rangle\)}}

\put(-2.1,-8.1){\makebox(0,0)[rt]{\(\langle 4597\rangle\)}}

\put(-0.1,-8.1){\makebox(0,0)[rt]{\(\langle 5576\rangle\)}}

\put(2.1,-8){\makebox(0,0)[lb]{\(\ast 3\)}}
\put(1.9,-8.1){\makebox(0,0)[rt]{\(\langle 5577\rangle\)}}
\put(1.9,-8.4){\makebox(0,0)[rt]{\(\langle 5578\rangle\)}}
\put(1.9,-8.7){\makebox(0,0)[rt]{\(\langle 5579\rangle\)}}

\put(-0.5,-10.1){\makebox(0,0)[rt]{\(\langle 261256\rangle\)}}
\put(-0.9,-10.5){\makebox(0,0)[rc]{\(\cdots\)}}
\put(-0.5,-10.7){\makebox(0,0)[rt]{\(\langle 261261\rangle\)}}

\put(4.1,-10){\makebox(0,0)[lb]{\(\ast 3\)}}
\put(2.2,-10.1){\makebox(0,0)[lt]{\(\langle 261262\rangle\)}}
\put(2.2,-10.4){\makebox(0,0)[lt]{\(\langle 261263\rangle\)}}
\put(2.2,-10.7){\makebox(0,0)[lt]{\(\langle 261264\rangle\)}}
\put(4,-10.1){\makebox(0,0)[lt]{\(\langle 261265\rangle\)}}
\put(4,-10.4){\makebox(0,0)[lt]{\(\langle 261266\rangle\)}}
\put(4,-10.7){\makebox(0,0)[lt]{\(\langle 261267\rangle\)}}
\put(5.8,-10.1){\makebox(0,0)[lt]{\(\langle 261268\rangle\)}}
\put(5.8,-10.4){\makebox(0,0)[lt]{\(\langle 261269\rangle\)}}
\put(5.8,-10.7){\makebox(0,0)[lt]{\(\langle 261270\rangle\)}}

\put(1.1,-10.1){\makebox(0,0)[rt]{\(\ast 2\)}}
\put(1.6,-10.1){\makebox(0,0)[rt]{\(\ast 2\)}}

\put(0.6,-12.1){\makebox(0,0)[rt]{\(\ast 42\)}}
\put(1.1,-12.1){\makebox(0,0)[rt]{\(\ast 4\)}}
\put(1.6,-12.1){\makebox(0,0)[rt]{\(\ast 6\)}}

\put(-0.6,-14.1){\makebox(0,0)[rt]{\(\ast 27\)}}

\put(1,0.4){\makebox(0,0)[cc]{\textbf{TKT}}}
\put(1,-0.1){\makebox(0,0)[cc]{\(\varkappa=\)}}
\put(2,-0.1){\makebox(0,0)[cc]{\((O^{13})\)}}
\put(0.3,-0.3){\framebox(2.6,1){}}

\put(-8,0){\makebox(0,0)[cc]{\textbf{TKT}}}
\put(-8,-0.5){\makebox(0,0)[cc]{\(\varkappa=\)}}
\put(-7,-0.5){\makebox(0,0)[cc]{\((O^9P^4)\)}}
\put(-8.7,-0.7){\framebox(2.6,1){}}

\put(-10,-19){\makebox(0,0)[cc]{\textbf{TKT}}}
\put(-10,-19.5){\makebox(0,0)[cc]{\(\varkappa=\)}}
\put(-8,-19.5){\makebox(0,0)[cc]{\((O^{12}P)\)}}
\put(-6,-19.5){\makebox(0,0)[cc]{\((O^{12}P)\)}}
\put(-4,-19.5){\makebox(0,0)[cc]{\((O^{13})\)}}
\put(-2,-19.5){\makebox(0,0)[cc]{\(((P_i)_{i=1}^4L^9)\)}}
\put(0,-19.5){\makebox(0,0)[cc]{\((O^3P^{10})\)}}
\put(2,-19.5){\makebox(0,0)[cc]{\(((P_i^3)_{i=1}^4L)\)}}
\put(4,-19.5){\makebox(0,0)[cc]{\(((P_i^3)_{i=1}^4L)\)}}
\put(6,-19.5){\makebox(0,0)[cc]{\(((P_i^3)_{i=1}^4L)\)}}
\put(-10.7,-19.7){\framebox(17.6,1){}}

\multiput(-8,-2)(2,0){2}{\oval(1.5,1.5)}
\put(-8,-3.1){\makebox(0,0)[cc]{\underbar{\textbf{8\,001}}\((\ast 1)\)}}
\put(-6,-3.1){\makebox(0,0)[cc]{\underbar{\textbf{3\,913}}\((\ast 2)\)}}

\multiput(2.2,-5.3)(0,-4){1}{\oval(1.5,3.9)}
\put(3.7,-7){\makebox(0,0)[cc]{\underbar{\textbf{21\,049}}\((\ast 1)\)}}

\multiput(4.1,-4.1)(2,0){2}{\oval(1.5,1.5)}
\put(4,-5.1){\makebox(0,0)[cc]{\underbar{\textbf{10\,621}}\((\ast 1)\)}}
\put(6,-5.1){\makebox(0,0)[cc]{\underbar{\textbf{4\,977}}\((\ast 1)\)}}

\put(-3.4,-8.3){\oval(3.6,1.5)}
\put(-3,-9.3){\makebox(0,0)[cc]{\underbar{\textbf{22\,581}}\((\ast 1)\)}}

\put(4.5,-11.2){\makebox(0,0)[ct]{with high probability}}

\put(0.8,-8.3){\oval(4.0,1.5)}
\put(2,-9.3){\makebox(0,0)[cc]{\underbar{\textbf{22\,581}}\((\ast 1)\)}}
\put(1,-12.5){\makebox(0,0)[lt]{with lower probability}}

\end{picture}

}

\end{figure}
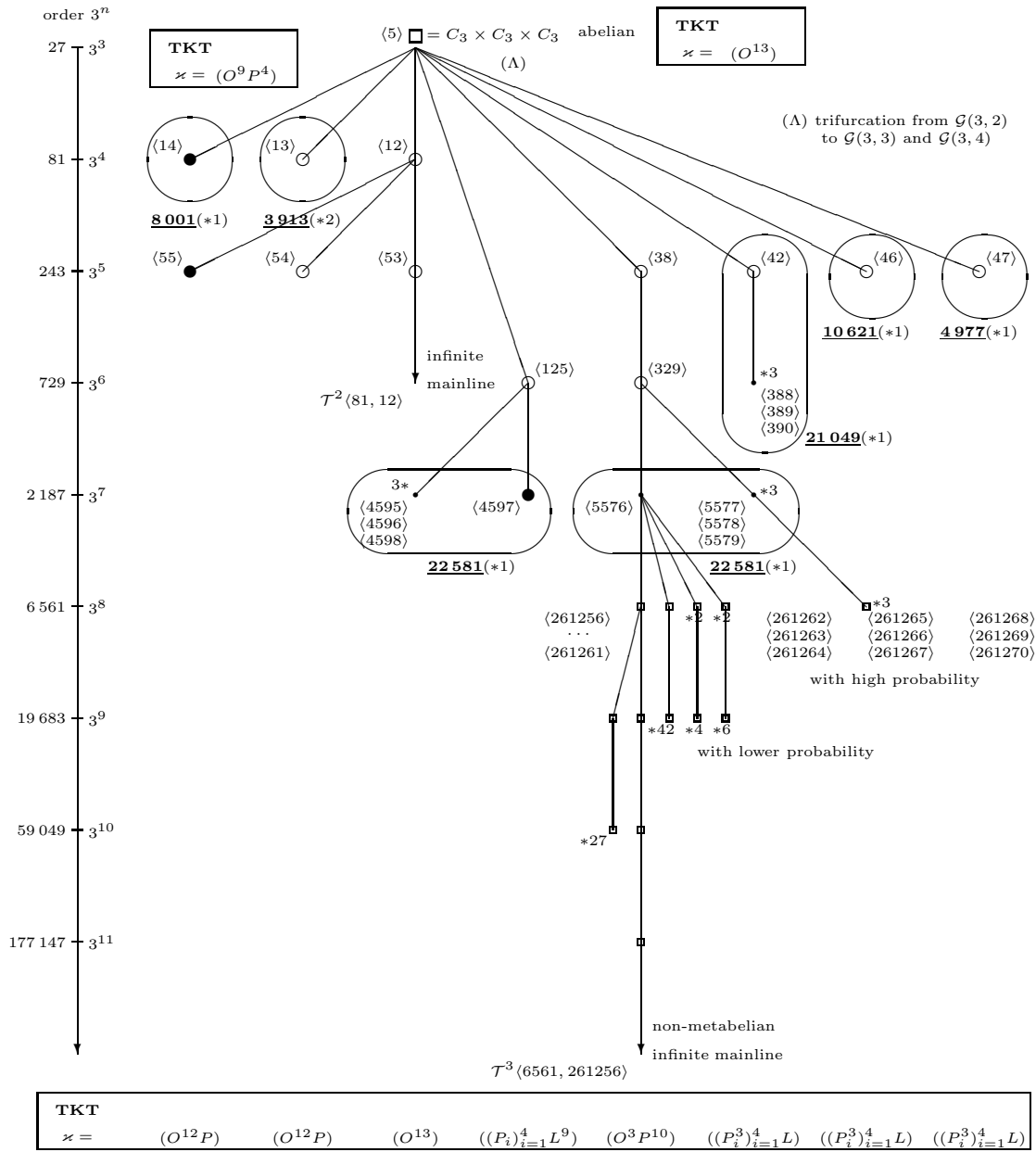

\newpage

\section{Cyclic Cubic Fields of Type (3,3,3)}
\label{s:CycCub3x3x3}


\subsection{Graph 1 of Category II}
\label{ss:Cat2Gph1}
\noindent
Since we have completed the discussion of conductors \(c\) in the Category \(\mathrm{III}\),
where all four fields of a quartet share a common \(3\)-class rank \(\varrho_3{F_{c,\mu}}=2\),
in the sections
\ref{ss:Cat3Gph1To4}
and
\ref{ss:Cat3Gph5To9},
we must now enter the realm of multiplets \((F_{c,1},\ldots,F_{c,m})\) with inhomogeneous \(3\)-class rank.
We begin with Category \(\mathrm{II}\),
where two fields have \(\varrho_3{F_{c,\mu}}=3\) and the remaining two fields have \(\varrho_3{F_{c,\mu}}=2\).
In the Tables
\ref{tbl:Cat2Gph1a},
\ref{tbl:Cat2Gph1b},
\ref{tbl:Cat2Gph2a}
and
\ref{tbl:Cat2Gph2b},
the second \(3\)-class groups of the two fields with \(\varrho_3{F_{c,\mu}}=3\) are given first
and separated by a semicolon.
Exceptional \(3\)-class groups of type \((9,3,3)\), resp. \((9,9,3)\),
are indicated by asterisks \(\ast\), resp. \(\ast\ast\).


\renewcommand{\arraystretch}{1.1}

\begin{table}[ht]
\caption{Fourty-Seven Examples for Graph \(1\) of Category \(\mathrm{II}\)}
\label{tbl:Cat2Gph1a}
\begin{center}

{\tiny

\begin{tabular}{|rrc|c|c|}
\hline
    No. &       \(c\) &                 \(\lbrack q_1,q_2,q_3\rbrack_3\) &                \(\mathrm{G}_3^{(2)}{F_{c,\mu}}\) & \(\left\lbrack\mathrm{Cl}_3{S_i}\right\rbrack_{1\le i\le 13}\) \\
\hline
  \(1\) &  \(3\,913\) &  \(\lbrace 13\rightarrow 7\leftarrow 43\rbrace\) & \(\langle 81,13\rangle^2;\langle 81,7\rangle^2\) &                        \(\lbrack(0)^3;(1^2)^8,(1^3)^2\rbrack\) \\
  \(2\) &  \(4\,123\) &  \(\lbrace 7\rightarrow 19\leftarrow 31\rbrace\) & \(\langle 81,13\rangle^2;\langle 81,7\rangle^2\) &                        \(\lbrack(0)^3;(1^2)^8,(1^3)^2\rbrack\) \\
  \(3\) &  \(4\,921\) &  \(\lbrace 7\rightarrow 19\leftarrow 37\rbrace\) & \(\langle 81,13\rangle^2;\langle 81,7\rangle^2\) &                        \(\lbrack(0)^3;(1^2)^8,(1^3)^2\rbrack\) \\
  \(4\) &  \(8\,827\) &  \(\lbrace 13\rightarrow 7\leftarrow 97\rbrace\) & \(\langle 81,13\rangle^2;\langle 81,7\rangle^2\) &                        \(\lbrack(0)^3;(1^2)^8,(1^3)^2\rbrack\) \\
  \(5\) & \(11\,557\) & \(\lbrace 13\rightarrow 7\leftarrow 127\rbrace\) & \(\langle 81,13\rangle^2;\langle 81,7\rangle^2\) &                        \(\lbrack(0)^3;(1^2)^8,(1^3)^2\rbrack\) \\
  \(6\) & \(12\,649\) & \(\lbrace 13\rightarrow 7\leftarrow 139\rbrace\) & \(\langle 81,13\rangle^2;\langle 81,7\rangle^2\) &                        \(\lbrack(0)^3;(1^2)^8,(1^3)^2\rbrack\) \\
  \(7\) & \(13\,699\) & \(\lbrace 7\rightarrow 19\leftarrow 103\rbrace\) & \(\langle 81,13\rangle^2;\langle 81,7\rangle^2\) &                        \(\lbrack(0)^3;(1^2)^8,(1^3)^2\rbrack\) \\
  \(8\) & \(21\,679\) & \(\lbrace 7\rightarrow 19\leftarrow 163\rbrace\) & \(\langle 81,13\rangle^2;\langle 81,7\rangle^2\) &                        \(\lbrack(0)^3;(1^2)^8,(1^3)^2\rbrack\) \\
  \(9\) & \(\mathbf{22\,581}\) & \(\lbrace 9\rightarrow 193\leftarrow 13\rbrace\) & \(\langle 3^7,5577\rangle,\langle 3^7,4595\rangle;\langle 3^6,41\rangle^2\) & \(\lbrack(0)^3;(21)^3,(2^2)^3,(1^3)^3,1^4\rbrack\) \\
 \(10\) & \(23\,121\) &  \(\lbrace 7\rightarrow 367\leftarrow 9\rbrace\) & \(\langle 81,13\rangle^2;\langle 81,7\rangle^2\) &                        \(\lbrack(0)^3;(1^2)^8,(1^3)^2\rbrack\) \\
 \(11\) & \(\mathbf{25\,929}\) &  \(\lbrace 9\rightarrow 67\leftarrow 43\rbrace\) & \(\langle 3^8,249242\rangle,\langle 3^8,249293\rangle;\langle 3^6,37\ldots 39\rangle^2\) & \(\lbrack(0)^3;(21)^6,(1^3)^3,1^4\rbrack\) \\
 \(12\) & \(27\,657\) &  \(\lbrace 7\rightarrow 439\leftarrow 9\rbrace\) & \(\langle 81,13\rangle^2;\langle 81,7\rangle^2\) &                        \(\lbrack(0)^3;(1^2)^8,(1^3)^2\rbrack\) \\
 \(13\) & \(28\,737\) & \(\lbrace 9\rightarrow 103\leftarrow 31\rbrace\) & \(\langle 81,13\rangle^2;\langle 81,7\rangle^2\) &                        \(\lbrack(0)^3;(1^2)^8,(1^3)^2\rbrack\) \\
 \(14\) & \(29\,419\) & \(\lbrace 31\rightarrow 13\leftarrow 73\rbrace\) & \(\langle 81,13\rangle^2;\langle 81,7\rangle^2\) &                        \(\lbrack(0)^3;(1^2)^8,(1^3)^2\rbrack\) \\
 \(15\) & \(\mathbf{30\,457}\) & \(\lbrace 7\rightarrow 19\leftarrow 229\rbrace\) & \(\ast^2;\langle 3^6,37\ldots 39\rangle^2\) & \(\lbrack(0)^3;(21)^6,(1^3)^3,21^3\rbrack\) \\
 \(16\) & \(31\,759\) & \(\lbrace 13\rightarrow 7\leftarrow 349\rbrace\) & \(\langle 81,13\rangle^2;\langle 81,7\rangle^2\) &                        \(\lbrack(0)^3;(1^2)^8,(1^3)^2\rbrack\) \\
 \(17\) & \(31\,837\) & \(\lbrace 31\rightarrow 13\leftarrow 79\rbrace\) & \(\langle 81,13\rangle^2;\langle 81,7\rangle^2\) &                        \(\lbrack(0)^3;(1^2)^8,(1^3)^2\rbrack\) \\
 \(18\) & \(\mathbf{34\,029}\) & \(\lbrace 19\rightarrow 9\leftarrow 199\rbrace\) & \(\ast^2;\langle 3^7,248\rangle^2\) & \(\lbrack(0)^3;(21)^3,(32)^3,(1^3)^3,21^3\rbrack\) \\
 \(19\) & \(\mathbf{34\,489}\) & \(\lbrace 13\rightarrow 7\leftarrow 379\rbrace\) & \(\langle 3^7,5577\rangle,\langle 3^7,4595\rangle;\langle 3^6,41\rangle^2\) & \(\lbrack(0)^3;(21)^3,(2^2)^3,(1^3)^3,1^4\rbrack\) \\
 \(20\) & \(36\,297\) & \(\lbrace 37\rightarrow 9\leftarrow 109\rbrace\) & \(\langle 81,13\rangle^2;\langle 81,7\rangle^2\) &                        \(\lbrack(0)^3;(1^2)^8,(1^3)^2\rbrack\) \\
 \(21\) & \(39\,403\) & \(\lbrace 13\rightarrow 7\leftarrow 433\rbrace\) & \(\langle 81,13\rangle^2;\langle 81,7\rangle^2\) &                        \(\lbrack(0)^3;(1^2)^8,(1^3)^2\rbrack\) \\
 \(22\) & \(41\,643\) &  \(\lbrace 7\rightarrow 661\leftarrow 9\rbrace\) & \(\langle 81,13\rangle^2;\langle 81,7\rangle^2\) &                        \(\lbrack(0)^3;(1^2)^8,(1^3)^2\rbrack\) \\
 \(23\) & \(\mathbf{41\,839}\) & \(\lbrace 43\rightarrow 7\leftarrow 139\rbrace\) & \(\ast^2;\langle 6561,673\rangle^2\) & \(\lbrack(0)^3;(2^2)^6,(1^3)^3,2^21^3\rbrack\) \\
 \(24\) & \(42\,291\) & \(\lbrace 37\rightarrow 9\leftarrow 127\rbrace\) & \(\langle 81,13\rangle^2;\langle 81,7\rangle^2\) &                        \(\lbrack(0)^3;(1^2)^8,(1^3)^2\rbrack\) \\
 \(25\) & \(49\,321\) & \(\lbrace 31\rightarrow 37\leftarrow 43\rbrace\) & \(\langle 81,13\rangle^2;\langle 81,7\rangle^2\) &                        \(\lbrack(0)^3;(1^2)^8,(1^3)^2\rbrack\) \\
 \(26\) & \(57\,421\) & \(\lbrace 13\rightarrow 7\leftarrow 631\rbrace\) & \(\langle 81,13\rangle^2;\langle 81,7\rangle^2\) &                        \(\lbrack(0)^3;(1^2)^8,(1^3)^2\rbrack\) \\
 \(27\) & \(58\,513\) & \(\lbrace 13\rightarrow 7\leftarrow 643\rbrace\) & \(\langle 81,13\rangle^2;\langle 81,7\rangle^2\) &                        \(\lbrack(0)^3;(1^2)^8,(1^3)^2\rbrack\) \\
 \(28\) & \(60\,273\) & \(\lbrace 37\rightarrow 9\leftarrow 181\rbrace\) & \(\langle 81,13\rangle^2;\langle 81,7\rangle^2\) &                        \(\lbrack(0)^3;(1^2)^8,(1^3)^2\rbrack\) \\
 \(29\) & \(60\,781\) & \(\lbrace 7\rightarrow 19\leftarrow 457\rbrace\) & \(\langle 81,13\rangle^2;\langle 81,7\rangle^2\) &                        \(\lbrack(0)^3;(1^2)^8,(1^3)^2\rbrack\) \\
 \(30\) & \(66\,937\) & \(\lbrace 13\rightarrow 271\leftarrow 19\rbrace\)& \(\langle 81,13\rangle^2;\langle 81,7\rangle^2\) &                        \(\lbrack(0)^3;(1^2)^8,(1^3)^2\rbrack\) \\
 \(31\) & \(67\,887\) & \(\lbrace 19\rightarrow 9\leftarrow 397\rbrace\) & \(\langle 81,13\rangle^2;\langle 81,7\rangle^2\) &                        \(\lbrack(0)^3;(1^2)^8,(1^3)^2\rbrack\) \\
 \(32\) & \(68\,887\) & \(\lbrace 13\rightarrow 7\leftarrow 757\rbrace\) & \(\langle 81,13\rangle^2;\langle 81,7\rangle^2\) &                        \(\lbrack(0)^3;(1^2)^8,(1^3)^2\rbrack\) \\
 \(33\) & \(73\,233\) & \(\lbrace 9\rightarrow 103\leftarrow 79\rbrace\) & \(\langle 81,13\rangle^2;\langle 81,7\rangle^2\) &                        \(\lbrack(0)^3;(1^2)^8,(1^3)^2\rbrack\) \\
 \(34\) & \(73\,633\) & \(\lbrace 7\rightarrow 157\leftarrow 67\rbrace\) & \(\langle 81,13\rangle^2;\langle 81,7\rangle^2\) &                        \(\lbrack(0)^3;(1^2)^8,(1^3)^2\rbrack\) \\
 \(35\) & \(\mathbf{74\,043}\) & \(\lbrace 19\rightarrow 9\leftarrow 433\rbrace\) & \(\langle 3^7,4606\rangle^2;\langle 3^7,65\vert 67\rangle^2\) & \(\lbrack(0)^3;(2^2)^6,(1^3)^3,21^3\rbrack\) \\
 \(36\) & \(74\,971\) & \(\lbrace 73\rightarrow 13\leftarrow 79\rbrace\) & \(\langle 81,13\rangle^2;\langle 81,7\rangle^2\) &                        \(\lbrack(0)^3;(1^2)^8,(1^3)^2\rbrack\) \\
 \(37\) & \(77\,337\) & \(\lbrace 9\rightarrow 661\leftarrow 13\rbrace\) & \(\langle 81,13\rangle^2;\langle 81,7\rangle^2\) &                        \(\lbrack(0)^3;(1^2)^8,(1^3)^2\rbrack\) \\
 \(38\) & \(78\,403\) & \(\lbrace 13\rightarrow 163\leftarrow 37\rbrace\)& \(\langle 81,13\rangle^2;\langle 81,7\rangle^2\) &                        \(\lbrack(0)^3;(1^2)^8,(1^3)^2\rbrack\) \\
 \(39\) & \(\mathbf{82\,327}\) & \(\lbrace 7\rightarrow 19\leftarrow 619\rbrace\) & \(\langle 3^6,372\rangle^2;\langle 3^6,37\ldots 39\rangle^2\) & \(\lbrack(0)^3;(21)^6,(1^3)^3,21^3\rbrack\) \\
 \(40\) & \(83\,327\) & \(\lbrace 7\rightarrow 73\leftarrow 163\rbrace\) & \(\langle 81,13\rangle^2;\langle 81,7\rangle^2\) &                        \(\lbrack(0)^3;(1^2)^8,(1^3)^2\rbrack\) \\
\hline
\end{tabular}

}

\end{center}
\end{table}


\renewcommand{\arraystretch}{1.1}

\begin{table}[ht]
\caption{Graph \(1\) of Category \(\mathrm{II}\) Continued}
\label{tbl:Cat2Gph1b}
\begin{center}

{\tiny

\begin{tabular}{|rrc|c|c|}
\hline
    No. &       \(c\) &                 \(\lbrack q_1,q_2,q_3\rbrack_3\) &                \(\mathrm{G}_3^{(2)}{F_{c,\mu}}\) & \(\left\lbrack\mathrm{Cl}_3{S_i}\right\rbrack_{1\le i\le 13}\) \\
\hline
 \(41\) & \(83\,731\) & \(\lbrace 31\rightarrow 37\leftarrow 73\rbrace\) & \(\langle 81,13\rangle^2;\langle 81,7\rangle^2\) &                        \(\lbrack(0)^3;(1^2)^8,(1^3)^2\rbrack\) \\
 \(42\) & \(\mathbf{83\,817}\) & \(\lbrace 9\rightarrow 67\leftarrow 139\rbrace\) & \(\ast\ast,\ast;\langle 6561,673\rangle^2\) & \(\lbrack(0)^3;(2^2)^6,(1^3)^3,21^4\rbrack\) \\
 \(43\) & \(89\,053\) & \(\lbrace 19\rightarrow 109\leftarrow 43\rbrace\)& \(\langle 81,13\rangle^2;\langle 81,7\rangle^2\) &                        \(\lbrack(0)^3;(1^2)^8,(1^3)^2\rbrack\) \\
 \(44\) & \(92\,407\) & \(\lbrace 43\rightarrow 7\leftarrow 307\rbrace\) & \(\langle 81,13\rangle^2;\langle 81,7\rangle^2\) &                        \(\lbrack(0)^3;(1^2)^8,(1^3)^2\rbrack\) \\
 \(45\) & \(92\,511\) & \(\lbrace 19\rightarrow 9\leftarrow 541\rbrace\) & \(\langle 81,13\rangle^2;\langle 81,7\rangle^2\) &                        \(\lbrack(0)^3;(1^2)^8,(1^3)^2\rbrack\) \\
 \(46\) & \(95\,641\) & \(\lbrace 13\rightarrow 7\leftarrow 1051\rbrace\)& \(\langle 81,13\rangle^2;\langle 81,7\rangle^2\) &                        \(\lbrack(0)^3;(1^2)^8,(1^3)^2\rbrack\) \\
 \(47\) & \(97\,209\) & \(\lbrace 7\rightarrow 1543\leftarrow 9\rbrace\) & \(\langle 81,13\rangle^2;\langle 81,7\rangle^2\) &                        \(\lbrack(0)^3;(1^2)^8,(1^3)^2\rbrack\) \\
\hline
\end{tabular}

}

\end{center}
\end{table}

\newpage

\phantom{strut}
\vskip 1.2in

\subsection{Graph 2 of Category II}
\label{ss:Cat2Gph2}
\noindent
Similarly as for Graph \(1\),
the behaviour of the dominating part of conductors \(c\) belonging to Graph \(2\)
seems to be uniform and can be summarized in Theorem
\ref{thm:Cat2Gph2}.


\renewcommand{\arraystretch}{1.1}

\begin{table}[ht]
\caption{Fourty-Five Examples for Graph \(2\) of Category \(\mathrm{II}\)}
\label{tbl:Cat2Gph2a}
\begin{center}

{\tiny

\begin{tabular}{|rrc|c|c|}
\hline
    No. &       \(c\) &                                \(\lbrack q_1,q_2,q_3\rbrack_3\) &                \(\mathrm{G}_3^{(2)}{F_{c,\mu}}\) & \(\left\lbrack\mathrm{Cl}_3{S_i}\right\rbrack_{1\le i\le 13}\) \\
\hline
  \(1\) &  \(6\,327\) &   \(\lbrace 19\rightarrow 9\leftarrow 37\rightarrow 19\rbrace\) & \(\langle 81,13\rangle^2;\langle 81,7\rangle^2\) & \(\lbrack(0)^3;(1^2)^8,(1^3)^2\rbrack\) \\
  \(2\) & \(18\,639\) & \(\lbrace 109\rightarrow 9\leftarrow 19\rightarrow 109\rbrace\) & \(\langle 81,13\rangle^2;\langle 81,7\rangle^2\) & \(\lbrack(0)^3;(1^2)^8,(1^3)^2\rbrack\) \\
  \(3\) & \(19\,201\) & \(\lbrace 211\rightarrow 7\leftarrow 13\rightarrow 211\rbrace\) & \(\langle 81,13\rangle^2;\langle 81,7\rangle^2\) & \(\lbrack(0)^3;(1^2)^8,(1^3)^2\rbrack\) \\
  \(4\) & \(20\,313\) &   \(\lbrace 61\rightarrow 9\leftarrow 37\rightarrow 61\rbrace\) & \(\langle 81,13\rangle^2;\langle 81,7\rangle^2\) & \(\lbrack(0)^3;(1^2)^8,(1^3)^2\rbrack\) \\
  \(5\) & \(21\,717\) &   \(\lbrace 9\rightarrow 127\leftarrow 19\rightarrow 9\rbrace\) & \(\langle 81,13\rangle^2;\langle 81,7\rangle^2\) & \(\lbrack(0)^3;(1^2)^8,(1^3)^2\rbrack\) \\
  \(6\) & \(21\,793\) &  \(\lbrace 37\rightarrow 19\leftarrow 31\rightarrow 37\rbrace\) & \(\langle 81,13\rangle^2;\langle 81,7\rangle^2\) & \(\lbrack(0)^3;(1^2)^8,(1^3)^2\rbrack\) \\
  \(7\) & \(21\,973\) &    \(\lbrace 7\rightarrow 73\leftarrow 43\rightarrow 7\rbrace\) & \(\langle 81,13\rangle^2;\langle 81,7\rangle^2\) & \(\lbrack(0)^3;(1^2)^8,(1^3)^2\rbrack\) \\
  \(8\) & \(\mathbf{27\,873}\) & \(\lbrace 19\rightarrow 9\leftarrow 163\rightarrow 19\rbrace\) & \(\ast^2;\langle 3^6,37\ldots 39\rangle^2\) & \(\lbrack(0)^3;(21)^6,(1^3)^3,1^5\rbrack\) \\
  \(9\) & \(27\,937\) &  \(\lbrace 13\rightarrow 7\leftarrow 307\rightarrow 13\rbrace\) & \(\langle 81,13\rangle^2;\langle 81,7\rangle^2\) & \(\lbrack(0)^3;(1^2)^8,(1^3)^2\rbrack\) \\
 \(10\) & \(\mathbf{29\,197}\) & \(\lbrace 43\rightarrow 7\leftarrow 97\rightarrow 43\rbrace\) & \(\ast^2;\langle 3^7,253\rangle^2\) & \(\lbrack(0)^3;(21)^3,(2^2)^3,(1^3)^3,21^3\rbrack\) \\
 \(11\) & \(30\,951\) & \(\lbrace 181\rightarrow 9\leftarrow 19\rightarrow 181\rbrace\) & \(\langle 81,13\rangle^2;\langle 81,7\rangle^2\) & \(\lbrack(0)^3;(1^2)^8,(1^3)^2\rbrack\) \\
 \(12\) & \(\mathbf{33\,943}\) & \(\lbrace 7\rightarrow 373\leftarrow 13\rightarrow 7\rbrace\) & \(\ast^2;\langle 3^6,37\ldots 39\rangle^2\) & \(\lbrack(0)^3;(21)^6,(1^3)^3,1^5\rbrack\) \\
 \(13\) & \(38\,227\) &  \(\lbrace 43\rightarrow 7\leftarrow 127\rightarrow 43\rbrace\) & \(\langle 81,13\rangle^2;\langle 81,7\rangle^2\) & \(\lbrack(0)^3;(1^2)^8,(1^3)^2\rbrack\) \\
 \(14\) & \(\mathbf{41\,629}\) & \(\lbrace 19\rightarrow 313\leftarrow 7\rightarrow 19\rbrace\) & \(\langle 3^8,249232\rangle^2;\langle 3^6,37\ldots 39\rangle^2\) & \(\lbrack(0)^3;(21)^6,(1^3)^3,1^4\rbrack\) \\
 \(15\) & \(43\,927\) & \(\lbrace 31\rightarrow 13\leftarrow 109\rightarrow 31\rbrace\) & \(\langle 81,13\rangle^2;\langle 81,7\rangle^2\) & \(\lbrack(0)^3;(1^2)^8,(1^3)^2\rbrack\) \\
 \(16\) & \(44\,023\) & \(\lbrace 331\rightarrow 19\leftarrow 7\rightarrow 331\rbrace\) & \(\langle 81,13\rangle^2;\langle 81,7\rangle^2\) & \(\lbrack(0)^3;(1^2)^8,(1^3)^2\rbrack\) \\
 \(17\) & \(\mathbf{46\,417}\) & \(\lbrace 7\rightarrow 19\leftarrow 349\rightarrow 7\rbrace\) & \(\ast^2;\langle 3^6,37\ldots 39\rangle^2\) & \(\lbrack(0)^3;(21)^6,(1^3)^3,1^5\rbrack\) \\
 \(18\) & \(49\,567\) &    \(\lbrace 7\rightarrow 73\leftarrow 97\rightarrow 7\rbrace\) & \(\langle 81,13\rangle^2;\langle 81,7\rangle^2\) & \(\lbrack(0)^3;(1^2)^8,(1^3)^2\rbrack\) \\
 \(19\) & \(49\,777\) &  \(\lbrace 13\rightarrow 7\leftarrow 547\rightarrow 13\rbrace\) & \(\langle 81,13\rangle^2;\langle 81,7\rangle^2\) & \(\lbrack(0)^3;(1^2)^8,(1^3)^2\rbrack\) \\
 \(20\) & \(50\,407\) &   \(\lbrace 7\rightarrow 19\leftarrow 379\rightarrow 7\rbrace\) & \(\langle 81,13\rangle^2;\langle 81,7\rangle^2\) & \(\lbrack(0)^3;(1^2)^8,(1^3)^2\rbrack\) \\
 \(21\) & \(54\,279\) & \(\lbrace 163\rightarrow 9\leftarrow 37\rightarrow 163\rbrace\) & \(\langle 81,13\rangle^2;\langle 81,7\rangle^2\) & \(\lbrack(0)^3;(1^2)^8,(1^3)^2\rbrack\) \\
 \(22\) & \(54\,691\) & \(\lbrace 601\rightarrow 7\leftarrow 13\rightarrow 601\rbrace\) & \(\langle 81,13\rangle^2;\langle 81,7\rangle^2\) & \(\lbrack(0)^3;(1^2)^8,(1^3)^2\rbrack\) \\
 \(23\) & \(\mathbf{56\,547}\) & \(\lbrace 61\rightarrow 103\leftarrow 9\rightarrow 61\rbrace\) & \(\langle 3^7,5577\rangle,\langle 3^7,4595\rangle;\langle 3^6,41\rangle^2\) & \(\lbrack(0)^3;(21)^3,(2^2)^3,(1^3)^3,1^4\rbrack\) \\
 \(24\) & \(60\,151\) &   \(\lbrace 7\rightarrow 661\leftarrow 13\rightarrow 7\rbrace\) & \(\langle 81,13\rangle^2;\langle 81,7\rangle^2\) & \(\lbrack(0)^3;(1^2)^8,(1^3)^2\rbrack\) \\
 \(25\) & \(60\,667\) & \(\lbrace 103\rightarrow 19\leftarrow 31\rightarrow 103\rbrace\)& \(\langle 81,13\rangle^2;\langle 81,7\rangle^2\) & \(\lbrack(0)^3;(1^2)^8,(1^3)^2\rbrack\) \\
\hline
\end{tabular}

}

\end{center}
\end{table}


\renewcommand{\arraystretch}{1.1}

\begin{table}[ht]
\caption{Graph \(2\) of Category \(\mathrm{II}\) Continued}
\label{tbl:Cat2Gph2b}
\begin{center}

{\tiny

\begin{tabular}{|rrc|c|c|}
\hline
    No. &       \(c\) &                                \(\lbrack q_1,q_2,q_3\rbrack_3\) &                \(\mathrm{G}_3^{(2)}{F_{c,\mu}}\) & \(\left\lbrack\mathrm{Cl}_3{S_i}\right\rbrack_{1\le i\le 13}\) \\
\hline
 \(26\) & \(60\,853\) & \(\lbrace 31\rightarrow 13\leftarrow 151\rightarrow 31\rbrace\) & \(\langle 81,13\rangle^2;\langle 81,7\rangle^2\) & \(\lbrack(0)^3;(1^2)^8,(1^3)^2\rbrack\) \\
 \(27\) & \(63\,271\) & \(\lbrace 31\rightarrow 13\leftarrow 157\rightarrow 31\rbrace\) & \(\langle 81,13\rangle^2;\langle 81,7\rangle^2\) & \(\lbrack(0)^3;(1^2)^8,(1^3)^2\rbrack\) \\
 \(28\) & \(\mathbf{63\,511}\) & \(\lbrace 43\rightarrow 7\leftarrow 211\rightarrow 43\rbrace\) & \(\ast^2;\langle 3^6,37\ldots 39\rangle^2\) & \(\lbrack(0)^3;(21)^6,(1^3)^3,21^3\rbrack\) \\
 \(29\) & \(64\,809\) &  \(\lbrace 19\rightarrow 9\leftarrow 379\rightarrow 19\rbrace\) & \(\langle 81,13\rangle^2;\langle 81,7\rangle^2\) & \(\lbrack(0)^3;(1^2)^8,(1^3)^2\rbrack\) \\
 \(30\) & \(\mathbf{65\,727}\) & \(\lbrace 9\rightarrow 67\leftarrow 109\rightarrow 9\rbrace\) & \(\ast^2;\langle 3^6,37\ldots 39\rangle^2\) & \(\lbrack(0)^3;(21)^6,(1^3)^3,21^3\rbrack\) \\
 \(31\) & \(\mathbf{66\,157}\) & \(\lbrace 13\rightarrow 7\leftarrow 727\rightarrow 13\rbrace\) & \(\ast\ast,\ast;\langle 6561,1989\rangle^2\) & \(\lbrack(0)^3;(21)^3,(3^2)^3,(1^3)^3,2^21^2\rbrack\) \\
 \(32\) & \(\mathbf{66\,267}\) & \(\lbrace 37\rightarrow 9\leftarrow 199\rightarrow 37\rbrace\) & \(\langle 3^8,249242\rangle,\langle 3^8,249293\rangle;\langle 3^6,37\ldots 39\rangle^2\) & \(\lbrack(0)^3;(21)^6,(1^3)^3,1^4\rbrack\) \\
 \(33\) & \(71\,029\) &   \(\lbrace 7\rightarrow 73\leftarrow 139\rightarrow 7\rbrace\) & \(\langle 81,13\rangle^2;\langle 81,7\rangle^2\) & \(\lbrack(0)^3;(1^2)^8,(1^3)^2\rbrack\) \\
 \(34\) & \(72\,943\) & \(\lbrace 181\rightarrow 13\leftarrow 31\rightarrow 181\rbrace\)& \(\langle 81,13\rangle^2;\langle 81,7\rangle^2\) & \(\lbrack(0)^3;(1^2)^8,(1^3)^2\rbrack\) \\
 \(35\) & \(75\,943\) & \(\lbrace 571\rightarrow 19\leftarrow 7\rightarrow 571\rbrace\) & \(\langle 81,13\rangle^2;\langle 81,7\rangle^2\) & \(\lbrack(0)^3;(1^2)^8,(1^3)^2\rbrack\) \\
 \(36\) & \(\mathbf{79\,933}\) & \(\lbrace 7\rightarrow 19\leftarrow 601\rightarrow 7\rbrace\) & \(\langle 3^7,5577\rangle,\langle 3^7,4595\rangle;\langle 3^6,41\rangle^2\) & \(\lbrack(0)^3;(21)^3,(2^2)^3,(1^3)^3,1^4\rbrack\) \\
 \(37\) & \(\mathbf{80\,731}\) & \(\lbrace 607\rightarrow 19\leftarrow 7\rightarrow 607\rbrace\) & \(\ast^2;\langle 3^6,37\ldots 39\rangle^2\) & \(\lbrack(0)^3;(21)^6,(1^3)^3,1^5\rbrack\) \\
 \(38\) & \(85\,267\) &  \(\lbrace 13\rightarrow 7\leftarrow 937\rightarrow 13\rbrace\) & \(\langle 81,13\rangle^2;\langle 81,7\rangle^2\) & \(\lbrack(0)^3;(1^2)^8,(1^3)^2\rbrack\) \\
 \(39\) & \(86\,233\) &  \(\lbrace 97\rightarrow 7\leftarrow 127\rightarrow 97\rbrace\) & \(\langle 81,13\rangle^2;\langle 81,7\rangle^2\) & \(\lbrack(0)^3;(1^2)^8,(1^3)^2\rbrack\) \\
 \(40\) & \(87\,913\) &  \(\lbrace 19\rightarrow 661\leftarrow 7\rightarrow 19\rbrace\) & \(\langle 81,13\rangle^2;\langle 81,7\rangle^2\) & \(\lbrack(0)^3;(1^2)^8,(1^3)^2\rbrack\) \\
 \(41\) & \(87\,997\) &  \(\lbrace 13\rightarrow 7\leftarrow 967\rightarrow 13\rbrace\) & \(\langle 81,13\rangle^2;\langle 81,7\rangle^2\) & \(\lbrack(0)^3;(1^2)^8,(1^3)^2\rbrack\) \\
 \(42\) & \(91\,053\) &  \(\lbrace 67\rightarrow 151\leftarrow 9\rightarrow 67\rbrace\) & \(\langle 81,13\rangle^2;\langle 81,7\rangle^2\) & \(\lbrack(0)^3;(1^2)^8,(1^3)^2\rbrack\) \\
 \(43\) & \(94\,381\) &  \(\lbrace 97\rightarrow 7\leftarrow 139\rightarrow 97\rbrace\) & \(\langle 81,13\rangle^2;\langle 81,7\rangle^2\) & \(\lbrack(0)^3;(1^2)^8,(1^3)^2\rbrack\) \\
 \(44\) & \(96\,733\) &\(\lbrace 1063\rightarrow 7\leftarrow 13\rightarrow 1063\rbrace\)& \(\langle 81,13\rangle^2;\langle 81,7\rangle^2\) & \(\lbrack(0)^3;(1^2)^8,(1^3)^2\rbrack\) \\
 \(45\) & \(99\,463\) & \(\lbrace 13\rightarrow 7\leftarrow 1093\rightarrow 13\rbrace\) & \(\langle 81,13\rangle^2;\langle 81,7\rangle^2\) & \(\lbrack(0)^3;(1^2)^8,(1^3)^2\rbrack\) \\
\hline
\end{tabular}

}

\end{center}
\end{table}

\newpage

\begin{theorem}
\label{thm:Cat2Gph2}
Suppose that \(u:=10^5\) is an assigned upper bound.
Let \(c<u\) be a conductor divisible by exactly three primes, \(t=3\), such that
\(\mathrm{Cl}_3{F_{c,\mu}}\simeq (3,3,3)\) for two cyclic cubic fields \(F_{c,1}\), \(F_{c,2}\) with conductor \(c\), and
\(\mathrm{Cl}_3{F_{c,\mu}}\simeq (3,3)\) for the other two cyclic cubic fields \(F_{c,3}\), \(F_{c,4}\) with conductor \(c\).
If \(c\) belongs either to Graph \(1\)
or to Graph \(2\) of Category \(\mathrm{II}\),
and if the \(3\)-class groups of the \(13\) bicyclic bicubic subfields \(S_i\) of the \(3\)-genus field \(F^\ast\) are given by
\begin{equation}
\label{eqn:Cat2Gph2Genus}
\left\lbrack\mathrm{Cl}_3{S_i}\right\rbrack_{1\le i\le 13}=\lbrack (0)^3;(11)^{8},(111)^2\rbrack,
\end{equation}
then the second \(3\)-class group \(\mathrm{G}_3^{(2)}{F_{c,\mu}}\) 
is isomorphic to \(\langle 81,13\rangle\) with \(3\)-capitulation type \(\varkappa(F_{c,\mu})=(O^9P^4)\)
for the fields \(F_{c,1}\), \(F_{c,2}\), and isomorphic to 
\(\langle 81,7\rangle\simeq\mathrm{Syl}_3{S_9}\) with \(3\)-capitulation type \(\mathrm{a}.3\), \(\varkappa(F_{c,\mu})=(2000)\),
for the fields \(F_{c,3}\), \(F_{c,4}\),
and the \(3\)-class tower has length \(\ell_3{F_{c,\mu}}=2\) for all fields.
\end{theorem}

\begin{proof}
See Tables
\ref{tbl:Cat2Gph1a},
\ref{tbl:Cat2Gph1b},
\ref{tbl:Cat2Gph2a}
and
\ref{tbl:Cat2Gph2b},
which have been computed with the aid of Magma
\cite{BCP1997,BCFS2022,Fi2001,MAGMA2022}.
\end{proof}


\begin{conjecture}
\label{cnj:Cat2Gph2}
Theorem
\ref{thm:Cat2Gph2}
remains true for any upper bound \(u>10^5\).
\end{conjecture}

\newpage

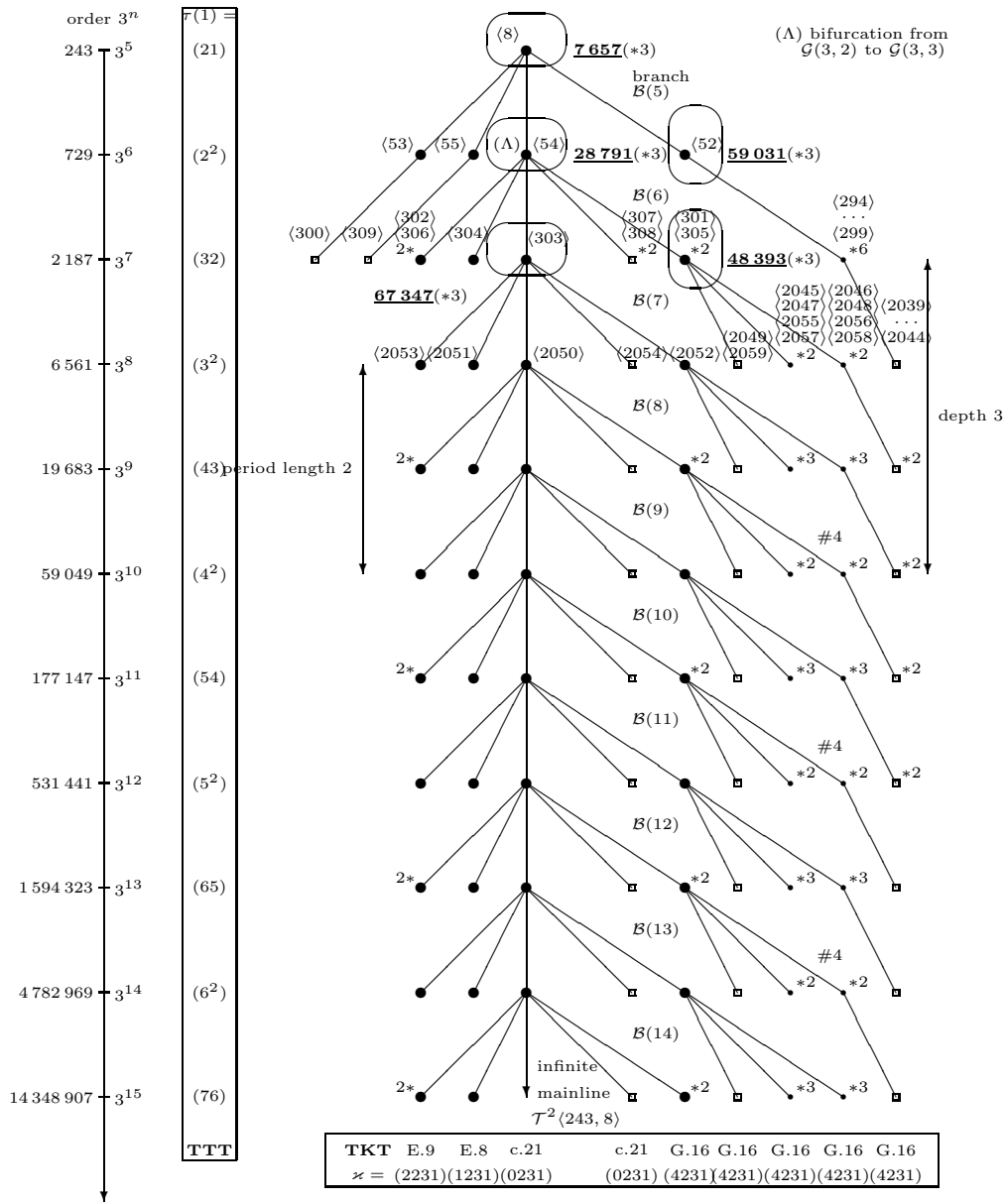
\begin{figure}[ht]
\caption{Distribution of Conductors for \(\mathrm{G}_3^{(2)}{F}\) on the Coclass Tree \(\mathcal{T}^2\langle 243,8\rangle\)}
\label{fig:MinDiscTyp33TreeUCc2}

{\tiny


\setlength{\unitlength}{0.7cm}
\begin{picture}(15,18)(-9,-17)

\put(-8,0.5){\makebox(0,0)[cb]{order \(3^n\)}}
\put(-8,0){\line(0,-1){20}}
\multiput(-8.1,0)(0,-2){11}{\line(1,0){0.2}}
\put(-8.2,0){\makebox(0,0)[rc]{\(243\)}}
\put(-7.8,0){\makebox(0,0)[lc]{\(3^5\)}}
\put(-8.2,-2){\makebox(0,0)[rc]{\(729\)}}
\put(-7.8,-2){\makebox(0,0)[lc]{\(3^6\)}}
\put(-8.2,-4){\makebox(0,0)[rc]{\(2\,187\)}}
\put(-7.8,-4){\makebox(0,0)[lc]{\(3^7\)}}
\put(-8.2,-6){\makebox(0,0)[rc]{\(6\,561\)}}
\put(-7.8,-6){\makebox(0,0)[lc]{\(3^8\)}}
\put(-8.2,-8){\makebox(0,0)[rc]{\(19\,683\)}}
\put(-7.8,-8){\makebox(0,0)[lc]{\(3^9\)}}
\put(-8.2,-10){\makebox(0,0)[rc]{\(59\,049\)}}
\put(-7.8,-10){\makebox(0,0)[lc]{\(3^{10}\)}}
\put(-8.2,-12){\makebox(0,0)[rc]{\(177\,147\)}}
\put(-7.8,-12){\makebox(0,0)[lc]{\(3^{11}\)}}
\put(-8.2,-14){\makebox(0,0)[rc]{\(531\,441\)}}
\put(-7.8,-14){\makebox(0,0)[lc]{\(3^{12}\)}}
\put(-8.2,-16){\makebox(0,0)[rc]{\(1\,594\,323\)}}
\put(-7.8,-16){\makebox(0,0)[lc]{\(3^{13}\)}}
\put(-8.2,-18){\makebox(0,0)[rc]{\(4\,782\,969\)}}
\put(-7.8,-18){\makebox(0,0)[lc]{\(3^{14}\)}}
\put(-8.2,-20){\makebox(0,0)[rc]{\(14\,348\,907\)}}
\put(-7.8,-20){\makebox(0,0)[lc]{\(3^{15}\)}}
\put(-8,-20){\vector(0,-1){2}}

\put(-6,0.5){\makebox(0,0)[cb]{\(\tau(1)=\)}}
\put(-6,0){\makebox(0,0)[cc]{\((21)\)}}
\put(-6,-2){\makebox(0,0)[cc]{\((2^2)\)}}
\put(-6,-4){\makebox(0,0)[cc]{\((32)\)}}
\put(-6,-6){\makebox(0,0)[cc]{\((3^2)\)}}
\put(-6,-8){\makebox(0,0)[cc]{\((43)\)}}
\put(-6,-10){\makebox(0,0)[cc]{\((4^2)\)}}
\put(-6,-12){\makebox(0,0)[cc]{\((54)\)}}
\put(-6,-14){\makebox(0,0)[cc]{\((5^2)\)}}
\put(-6,-16){\makebox(0,0)[cc]{\((65)\)}}
\put(-6,-18){\makebox(0,0)[cc]{\((6^2)\)}}
\put(-6,-20){\makebox(0,0)[cc]{\((76)\)}}
\put(-6,-21){\makebox(0,0)[cc]{\textbf{TTT}}}
\put(-6.5,-21.2){\framebox(1,22){}}

\put(7.6,-7){\vector(0,1){3}}
\put(7.8,-7){\makebox(0,0)[lc]{depth \(3\)}}
\put(7.6,-7){\vector(0,-1){3}}

\put(-3.1,-8){\vector(0,1){2}}
\put(-3.3,-8){\makebox(0,0)[rc]{period length \(2\)}}
\put(-3.1,-8){\vector(0,-1){2}}

\put(4.7,0.3){\makebox(0,0)[lc]{\((\Lambda)\) bifurcation from}}
\put(5.2,0){\makebox(0,0)[lc]{\(\mathcal{G}(3,2)\) to \(\mathcal{G}(3,3)\)}}


\multiput(0,0)(0,-2){10}{\circle*{0.2}}
\multiput(0,0)(0,-2){9}{\line(0,-1){2}}
\multiput(-1,-2)(0,-2){10}{\circle*{0.2}}
\multiput(-2,-2)(0,-2){10}{\circle*{0.2}}
\multiput(1.95,-4.05)(0,-2){9}{\framebox(0.1,0.1){}}
\multiput(3,-2)(0,-2){10}{\circle*{0.2}}
\multiput(0,0)(0,-2){10}{\line(-1,-2){1}}
\multiput(0,0)(0,-2){10}{\line(-1,-1){2}}
\multiput(0,-2)(0,-2){9}{\line(1,-1){2}}
\multiput(0,0)(0,-2){10}{\line(3,-2){3}}
\multiput(-3.05,-4.05)(-1,0){2}{\framebox(0.1,0.1){}}
\multiput(3.95,-6.05)(0,-2){8}{\framebox(0.1,0.1){}}
\multiput(5,-6)(0,-2){8}{\circle*{0.1}}
\multiput(6,-4)(0,-2){9}{\circle*{0.1}}
\multiput(-1,-2)(-1,0){2}{\line(-1,-1){2}}
\multiput(3,-4)(0,-2){8}{\line(1,-2){1}}
\multiput(3,-4)(0,-2){8}{\line(1,-1){2}}
\multiput(3,-2)(0,-2){9}{\line(3,-2){3}}
\multiput(6.95,-6.05)(0,-2){8}{\framebox(0.1,0.1){}}
\multiput(6,-4)(0,-2){8}{\line(1,-2){1}}

\put(2,-0.5){\makebox(0,0)[lc]{branch}}
\put(2,-0.8){\makebox(0,0)[lc]{\(\mathcal{B}(5)\)}}
\put(2,-2.8){\makebox(0,0)[lc]{\(\mathcal{B}(6)\)}}
\put(2,-4.8){\makebox(0,0)[lc]{\(\mathcal{B}(7)\)}}
\put(2,-6.8){\makebox(0,0)[lc]{\(\mathcal{B}(8)\)}}
\put(2,-8.8){\makebox(0,0)[lc]{\(\mathcal{B}(9)\)}}
\put(2,-10.8){\makebox(0,0)[lc]{\(\mathcal{B}(10)\)}}
\put(2,-12.8){\makebox(0,0)[lc]{\(\mathcal{B}(11)\)}}
\put(2,-14.8){\makebox(0,0)[lc]{\(\mathcal{B}(12)\)}}
\put(2,-16.8){\makebox(0,0)[lc]{\(\mathcal{B}(13)\)}}
\put(2,-18.8){\makebox(0,0)[lc]{\(\mathcal{B}(14)\)}}

\put(-0.1,0.3){\makebox(0,0)[rc]{\(\langle 8\rangle\)}}

\put(-2.1,-1.8){\makebox(0,0)[rc]{\(\langle 53\rangle\)}}
\put(-1.1,-1.8){\makebox(0,0)[rc]{\(\langle 55\rangle\)}}
\put(0.1,-1.8){\makebox(0,0)[lc]{\(\langle 54\rangle\)}}
\put(-0.1,-1.8){\makebox(0,0)[rc]{\((\Lambda)\)}}
\put(3.1,-1.8){\makebox(0,0)[lc]{\(\langle 52\rangle\)}}

\put(-4.1,-3.5){\makebox(0,0)[cc]{\(\langle 300\rangle\)}}
\put(-3.1,-3.5){\makebox(0,0)[cc]{\(\langle 309\rangle\)}}
\put(-2.1,-3.2){\makebox(0,0)[cc]{\(\langle 302\rangle\)}}
\put(-2.1,-3.5){\makebox(0,0)[cc]{\(\langle 306\rangle\)}}
\put(-1.1,-3.5){\makebox(0,0)[cc]{\(\langle 304\rangle\)}}
\put(0.0,-3.6){\makebox(0,0)[lc]{\(\langle 303\rangle\)}}
\put(2.2,-3.2){\makebox(0,0)[cc]{\(\langle 307\rangle\)}}
\put(2.2,-3.5){\makebox(0,0)[cc]{\(\langle 308\rangle\)}}
\put(3.2,-3.2){\makebox(0,0)[cc]{\(\langle 301\rangle\)}}
\put(3.2,-3.5){\makebox(0,0)[cc]{\(\langle 305\rangle\)}}
\put(6.2,-2.9){\makebox(0,0)[cc]{\(\langle 294\rangle\)}}
\put(6.2,-3.2){\makebox(0,0)[cc]{\(\cdots\)}}
\put(6.2,-3.5){\makebox(0,0)[cc]{\(\langle 299\rangle\)}}

\put(-2.4,-5.8){\makebox(0,0)[cc]{\(\langle 2053\rangle\)}}
\put(-1.4,-5.8){\makebox(0,0)[cc]{\(\langle 2051\rangle\)}}
\put(0.1,-5.8){\makebox(0,0)[lc]{\(\langle 2050\rangle\)}}
\put(2.2,-5.8){\makebox(0,0)[cc]{\(\langle 2054\rangle\)}}
\put(3.2,-5.8){\makebox(0,0)[cc]{\(\langle 2052\rangle\)}}
\put(4.2,-5.5){\makebox(0,0)[cc]{\(\langle 2049\rangle\)}}
\put(4.2,-5.8){\makebox(0,0)[cc]{\(\langle 2059\rangle\)}}
\put(5.2,-4.6){\makebox(0,0)[cc]{\(\langle 2045\rangle\)}}
\put(5.2,-4.9){\makebox(0,0)[cc]{\(\langle 2047\rangle\)}}
\put(5.2,-5.2){\makebox(0,0)[cc]{\(\langle 2055\rangle\)}}
\put(5.2,-5.5){\makebox(0,0)[cc]{\(\langle 2057\rangle\)}}
\put(6.2,-4.6){\makebox(0,0)[cc]{\(\langle 2046\rangle\)}}
\put(6.2,-4.9){\makebox(0,0)[cc]{\(\langle 2048\rangle\)}}
\put(6.2,-5.2){\makebox(0,0)[cc]{\(\langle 2056\rangle\)}}
\put(6.2,-5.5){\makebox(0,0)[cc]{\(\langle 2058\rangle\)}}
\put(7.2,-4.9){\makebox(0,0)[cc]{\(\langle 2039\rangle\)}}
\put(7.2,-5.2){\makebox(0,0)[cc]{\(\cdots\)}}
\put(7.2,-5.5){\makebox(0,0)[cc]{\(\langle 2044\rangle\)}}

\put(2.1,-3.8){\makebox(0,0)[lc]{\(*2\)}}
\multiput(-2.1,-3.8)(0,-4){5}{\makebox(0,0)[rc]{\(2*\)}}
\multiput(3.1,-3.8)(0,-4){5}{\makebox(0,0)[lc]{\(*2\)}}
\put(6.1,-3.8){\makebox(0,0)[lc]{\(*6\)}}
\multiput(5.1,-5.8)(0,-4){4}{\makebox(0,0)[lc]{\(*2\)}}
\multiput(5.5,-9.3)(0,-4){3}{\makebox(0,0)[lc]{\(\#4\)}}
\multiput(6.1,-5.8)(0,-4){4}{\makebox(0,0)[lc]{\(*2\)}}
\multiput(5.1,-7.8)(0,-4){4}{\makebox(0,0)[lc]{\(*3\)}}
\multiput(6.1,-7.8)(0,-4){4}{\makebox(0,0)[lc]{\(*3\)}}
\multiput(7.1,-7.8)(0,-2){4}{\makebox(0,0)[lc]{\(*2\)}}

\put(-3,-21){\makebox(0,0)[cc]{\textbf{TKT}}}
\put(-2,-21){\makebox(0,0)[cc]{E.9}}
\put(-1,-21){\makebox(0,0)[cc]{E.8}}
\put(0,-21){\makebox(0,0)[cc]{c.21}}
\put(2,-21){\makebox(0,0)[cc]{c.21}}
\put(3.1,-21){\makebox(0,0)[cc]{G.16}}
\put(4,-21){\makebox(0,0)[cc]{G.16}}
\put(5,-21){\makebox(0,0)[cc]{G.16}}
\put(6,-21){\makebox(0,0)[cc]{G.16}}
\put(7,-21){\makebox(0,0)[cc]{G.16}}
\put(-3,-21.5){\makebox(0,0)[cc]{\(\varkappa=\)}}
\put(-2,-21.5){\makebox(0,0)[cc]{\((2231)\)}}
\put(-1,-21.5){\makebox(0,0)[cc]{\((1231)\)}}
\put(0,-21.5){\makebox(0,0)[cc]{\((0231)\)}}
\put(2,-21.5){\makebox(0,0)[cc]{\((0231)\)}}
\put(3.1,-21.5){\makebox(0,0)[cc]{\((4231)\)}}
\put(4,-21.5){\makebox(0,0)[cc]{\((4231)\)}}
\put(5,-21.5){\makebox(0,0)[cc]{\((4231)\)}}
\put(6,-21.5){\makebox(0,0)[cc]{\((4231)\)}}
\put(7,-21.5){\makebox(0,0)[cc]{\((4231)\)}}
\put(-3.8,-21.7){\framebox(11.6,1){}}

\put(0,-18){\vector(0,-1){2}}
\put(0.2,-19.4){\makebox(0,0)[lc]{infinite}}
\put(0.2,-19.9){\makebox(0,0)[lc]{mainline}}
\put(1.8,-20.4){\makebox(0,0)[rc]{\(\mathcal{T}^2\langle 243,8\rangle\)}}


\multiput(0,0.2)(0,-2){3}{\oval(1.5,1)}
\put(0.9,0.0){\makebox(0,0)[lc]{\underbar{\textbf{7\,657}}\((\ast 3)\)}}
\put(0.9,-2.0){\makebox(0,0)[lc]{\underbar{\textbf{28\,791}}\((\ast 3)\)}}
\put(-1.1,-4.7){\makebox(0,0)[rc]{\underbar{\textbf{67\,347}}\((\ast 3)\)}}
\multiput(3.2,-1.8)(0,-2){2}{\oval(1,1.5)}
\put(4.7,-2.0){\makebox(0,0)[cc]{\underbar{\textbf{59\,031}}\((\ast 3)\)}}
\put(4.7,-4.0){\makebox(0,0)[cc]{\underbar{\textbf{48\,393}}\((\ast 3)\)}}


\end{picture}

}

\end{figure}

\newpage

\subsection{Graph 1 of Category I}
\label{ss:Cat1Gph1}
\noindent
We continue with Category \(\mathrm{I}\),
where three fields have \(\varrho_3{F_{c,\mu}}=2\) and the remaining single field has \(\varrho_3{F_{c,\mu}}=3\).
In the Tables
\ref{tbl:Cat1Gph1A},
\ref{tbl:Cat1Gph1B},
\ref{tbl:Cat1Gph2A},
and
\ref{tbl:Cat1Gph2B},
the second \(3\)-class group of the unique field with \(\varrho_3{F_{c,\mu}}=3\) is given first
and separated by a semicolon.
Again, we denote a \(3\)-class group
\(\mathrm{Cl}_3{F_{c,\mu}}=(9,3,3)\), respectively \((9,9,3)\),
by an asterisk \(\ast\), respectively two asterisks \(\ast\ast\).


\renewcommand{\arraystretch}{1.1}

\begin{table}[ht]
\caption{Thirty-Eight Examples for Graph \(1\) of Category \(\mathrm{I}\)}
\label{tbl:Cat1Gph1A}
\begin{center}

{\scriptsize

\begin{tabular}{|crc|c|c|}
\hline
    No. &       \(c\) &      \(\lbrack q_1,q_2,q_3\rbrack_3\) &                       \(\mathrm{G}_3^{(2)}{F_{c,\mu}}\) & \(\left\lbrack\mathrm{Cl}_3{S_i}\right\rbrack_{1\le i\le 13}\) \\
\hline
  \(1\) &  \(4\,977\) &    \(\lbrace 9,7,79;\delta=0\rbrace\) & \(\langle 243,47\rangle;\langle 243,25\rangle,\langle 243,28\rangle^2\) & \(\lbrack(0)^3;(1^2)^6,(21)^2,2^2,1^3\rbrack\) \\
  \(2\) & \(10\,621\) &  \(\lbrace 13,19,43;\delta=0\rbrace\) & \(\langle 243,46\rangle;\langle 243,25\rangle,\langle 243,28\rangle^2\) & \(\lbrack(0)^3;(1^2)^6,(21)^2,2^2,1^3\rbrack\) \\
  \(3\) & \(11\,349\) &   \(\lbrace 9,13,97;\delta=0\rbrace\) &     \(\langle 81,14\rangle;\langle 81,8\rangle,\langle 81,10\rangle^2\) &         \(\lbrack(0)^3;(1^2)^7,(21)^3\rbrack\) \\
  \(4\) & \(16\,263\) &  \(\lbrace 9,13,139;\delta=0\rbrace\) & \(\langle 243,46\rangle;\langle 243,25\rangle,\langle 243,28\rangle^2\) & \(\lbrack(0)^3;(1^2)^6,(21)^2,2^2,1^3\rbrack\) \\
  \(5\) & \(17\,353\) &   \(\lbrace 7,37,67;\delta=0\rbrace\) &     \(\langle 81,14\rangle;\langle 81,8\rangle,\langle 81,10\rangle^2\) &         \(\lbrack(0)^3;(1^2)^7,(21)^3\rbrack\) \\
  \(6\) & \(17\,829\) &   \(\lbrace 9,7,283;\delta=0\rbrace\) &     \(\langle 81,14\rangle;\langle 81,8\rangle,\langle 81,10\rangle^2\) &         \(\lbrack(0)^3;(1^2)^7,(21)^3\rbrack\) \\
  \(7\) & \(22\,041\) &   \(\lbrace 9,31,79;\delta=0\rbrace\) &     \(\langle 81,14\rangle;\langle 81,8\rangle,\langle 81,10\rangle^2\) &         \(\lbrack(0)^3;(1^2)^7,(21)^3\rbrack\) \\
  \(8\) & \(28\,197\) &  \(\lbrace 9,13,241;\delta=0\rbrace\) &     \(\langle 81,14\rangle;\langle 81,8\rangle,\langle 81,10\rangle^2\) &         \(\lbrack(0)^3;(1^2)^7,(21)^3\rbrack\) \\
  \(9\) & \(28\,609\) &   \(\lbrace 7,61,67;\delta=0\rbrace\) & \(\langle 243,46\rangle;\langle 243,25\rangle,\langle 243,28\rangle^2\) & \(\lbrack(0)^3;(1^2)^6,(21)^2,2^2,1^3\rbrack\) \\
 \(10\) & \(\mathbf{28\,791}\) &   \(\lbrace 9,7,457;\delta=0\rbrace\) &                               \(\ast;\langle 729,54\rangle^3\) &    \(\lbrack(0)^3;(21)^6,(2^2)^3,2^21\rbrack\) \\
 \(11\) & \(\mathbf{32\,227}\) &  \(\lbrace 13,37,67;\delta=0\rbrace\) &                               \(\ast;\langle 729,54\rangle^3\) &    \(\lbrack(0)^3;(21)^6,(2^2)^3,2^21\rbrack\) \\
 \(12\) & \(34\,099\) &  \(\lbrace 13,43,61;\delta=0\rbrace\) & \(\langle 243,47\rangle;\langle 243,25\rangle,\langle 243,28\rangle^2\) & \(\lbrack(0)^3;(1^2)^6,(21)^2,2^2,1^3\rbrack\) \\
 \(13\) & \(34\,333\) & \(\lbrace 13,19,139;\delta=0\rbrace\) & \(\langle 243,47\rangle;\langle 243,25\rangle,\langle 243,28\rangle^2\) & \(\lbrack(0)^3;(1^2)^6,(21)^2,2^2,1^3\rbrack\) \\
 \(14\) & \(38\,727\) &  \(\lbrace 9,13,331;\delta=0\rbrace\) &                                         \(\ast;\langle 243,8\rangle^3\) &            \(\lbrack(0)^3;(21)^9,21^2\rbrack\) \\
 \(15\) & \(40\,833\) &  \(\lbrace 9,13,349;\delta=0\rbrace\) &                                         \(\ast;\langle 243,8\rangle^3\) &            \(\lbrack(0)^3;(21)^9,21^2\rbrack\) \\
 \(16\) & \(43\,183\) &  \(\lbrace 7,31,199;\delta=0\rbrace\) &     \(\langle 81,14\rangle;\langle 81,8\rangle,\langle 81,10\rangle^2\) &         \(\lbrack(0)^3;(1^2)^7,(21)^3\rbrack\) \\
 \(17\) & \(43\,533\) &   \(\lbrace 9,7,691;\delta=0\rbrace\) &                                         \(\ast;\langle 243,8\rangle^3\) &            \(\lbrack(0)^3;(21)^9,21^2\rbrack\) \\
 \(18\) & \(46\,179\) &   \(\lbrace 9,7,733;\delta=0\rbrace\) &     \(\langle 81,14\rangle;\langle 81,8\rangle,\langle 81,10\rangle^2\) &         \(\lbrack(0)^3;(1^2)^7,(21)^3\rbrack\) \\
 \(19\) & \(49\,153\) & \(\lbrace 13,19,199;\delta=0\rbrace\) &     \(\langle 81,14\rangle;\langle 81,8\rangle,\langle 81,10\rangle^2\) &         \(\lbrack(0)^3;(1^2)^7,(21)^3\rbrack\) \\
 \(20\) & \(52\,297\) &  \(\lbrace 7,31,241;\delta=0\rbrace\) & \(\langle 243,46\rangle;\langle 243,25\rangle,\langle 243,28\rangle^2\) & \(\lbrack(0)^3;(1^2)^6,(21)^2,2^2,1^3\rbrack\) \\
 \(21\) & \(53\,793\) &  \(\lbrace 9,43,139;\delta=0\rbrace\) &     \(\langle 81,14\rangle;\langle 81,8\rangle,\langle 81,10\rangle^2\) &         \(\lbrack(0)^3;(1^2)^7,(21)^3\rbrack\) \\
 \(22\) & \(58\,869\) &  \(\lbrace 9,31,211;\delta=0\rbrace\) &     \(\langle 81,14\rangle;\langle 81,8\rangle,\langle 81,10\rangle^2\) &         \(\lbrack(0)^3;(1^2)^7,(21)^3\rbrack\) \\
 \(23\) & \(\mathbf{61\,087}\) & \(\lbrace 13,37,127;\delta=0\rbrace\) &                               \(\ast;\langle 729,54\rangle^3\) &    \(\lbrack(0)^3;(21)^6,(2^2)^3,2^21\rbrack\) \\
 \(24\) & \(64\,543\) &  \(\lbrace 19,43,79;\delta=0\rbrace\) &     \(\langle 81,14\rangle;\langle 81,8\rangle,\langle 81,10\rangle^2\) &         \(\lbrack(0)^3;(1^2)^7,(21)^3\rbrack\) \\
 \(25\) & \(65\,457\) &  \(\lbrace 9,7,1039;\delta=0\rbrace\) &     \(\langle 81,14\rangle;\langle 81,8\rangle,\langle 81,10\rangle^2\) &         \(\lbrack(0)^3;(1^2)^7,(21)^3\rbrack\) \\
 \(26\) & \(67\,239\) &  \(\lbrace 9,31,241;\delta=0\rbrace\) &                                         \(\ast;\langle 243,8\rangle^3\) &            \(\lbrack(0)^3;(21)^9,21^2\rbrack\) \\
 \(27\) & \(\mathbf{67\,347}\) &  \(\lbrace 9,7,1069;\delta=0\rbrace\) &                         \(\ast\ast;\langle 2187,303\rangle^3\) &      \(\lbrack(0)^3;(21)^6,(32)^3,2^3\rbrack\) \\
 \(28\) & \(77\,931\) &  \(\lbrace 9,7,1237;\delta=0\rbrace\) &     \(\langle 81,14\rangle;\langle 81,8\rangle,\langle 81,10\rangle^2\) &         \(\lbrack(0)^3;(1^2)^7,(21)^3\rbrack\) \\
 \(29\) & \(78\,589\) & \(\lbrace 7,103,109;\delta=0\rbrace\) &                                         \(\ast;\langle 243,8\rangle^3\) &            \(\lbrack(0)^3;(21)^9,21^2\rbrack\) \\
 \(30\) & \(80\,847\) &  \(\lbrace 9,13,691;\delta=0\rbrace\) & \(\langle 243,46\rangle;\langle 243,25\rangle,\langle 243,28\rangle^2\) & \(\lbrack(0)^3;(1^2)^6,(21)^2,2^2,1^3\rbrack\) \\
\hline
\end{tabular}

}

\end{center}
\end{table}


\renewcommand{\arraystretch}{1.1}

\begin{table}[ht]
\caption{Graph \(1\) of Category \(\mathrm{I}\) Continued}
\label{tbl:Cat1Gph1B}
\begin{center}

{\scriptsize

\begin{tabular}{|crc|c|c|}
\hline
    No. &       \(c\) &      \(\lbrack q_1,q_2,q_3\rbrack_3\) &                       \(\mathrm{G}_3^{(2)}{F_{c,\mu}}\) & \(\left\lbrack\mathrm{Cl}_3{S_i}\right\rbrack_{1\le i\le 13}\) \\
\hline
 \(31\) & \(81\,333\) &  \(\lbrace 9,7,1291;\delta=0\rbrace\) &                                         \(\ast;\langle 243,8\rangle^3\) &            \(\lbrack(0)^3;(21)^9,21^2\rbrack\) \\
 \(32\) & \(82\,411\) &  \(\lbrace 7,61,193;\delta=0\rbrace\) & \(\langle 243,47\rangle;\langle 243,25\rangle,\langle 243,28\rangle^2\) & \(\lbrack(0)^3;(1^2)^6,(21)^2,2^2,1^3\rbrack\) \\
 \(33\) & \(84\,607\) &  \(\lbrace 19,61,73;\delta=0\rbrace\) &     \(\langle 81,14\rangle;\langle 81,8\rangle,\langle 81,10\rangle^2\) &         \(\lbrack(0)^3;(1^2)^7,(21)^3\rbrack\) \\
 \(34\) & \(87\,327\) &  \(\lbrace 9,31,313;\delta=0\rbrace\) &     \(\langle 81,14\rangle;\langle 81,8\rangle,\langle 81,10\rangle^2\) &         \(\lbrack(0)^3;(1^2)^7,(21)^3\rbrack\) \\
 \(35\) & \(87\,867\) &  \(\lbrace 9,13,751;\delta=0\rbrace\) &     \(\langle 81,14\rangle;\langle 81,8\rangle,\langle 81,10\rangle^2\) &         \(\lbrack(0)^3;(1^2)^7,(21)^3\rbrack\) \\
 \(36\) & \(90\,649\) & \(\lbrace 13,19,367;\delta=0\rbrace\) &                                         \(\ast;\langle 243,8\rangle^3\) &            \(\lbrack(0)^3;(21)^9,21^2\rbrack\) \\
 \(37\) & \(92\,349\) &  \(\lbrace 9,31,331;\delta=0\rbrace\) &     \(\langle 81,14\rangle;\langle 81,8\rangle,\langle 81,10\rangle^2\) &         \(\lbrack(0)^3;(1^2)^7,(21)^3\rbrack\) \\
 \(38\) & \(96\,291\) &  \(\lbrace 9,13,823;\delta=0\rbrace\) &     \(\langle 81,14\rangle;\langle 81,8\rangle,\langle 81,10\rangle^2\) &         \(\lbrack(0)^3;(1^2)^7,(21)^3\rbrack\) \\
\hline
\end{tabular}

}

\end{center}
\end{table}

\newpage

\subsection{Graph 2 of Category I}
\label{ss:Cat1Gph2}


\begin{theorem}
\label{thm:Cat1Gph12}
Suppose that \(u:=10^5\) is an assigned upper bound.
Let \(c<u\) be a conductor divisible by exactly three primes, \(t=3\), such that
\(\mathrm{Cl}_3{F_{c,\mu}}\simeq (3,3,3)\) for the single cyclic cubic field \(\mu=1\), and
\(\mathrm{Cl}_3{F_{c,\mu}}\simeq (3,3)\) for the other three cyclic cubic fields \(2\le\mu\le 4\) with conductor \(c\).
If \(c=q_1q_2q_3\) belongs to Graph \(1\) or \(2\) of Category \(\mathrm{I}\), i.e.,
\begin{equation}
\label{eqn:Cat1Gph2}
\lbrack q_1,q_2,q_3\rbrack_3=
\begin{cases}
\lbrace q_1,q_2,q_3;\delta=0\rbrace & \text{ or} \\
\lbrace q_i\leftarrow q_j\rightarrow q_k\rbrace & \text{ with } i,j,k \text{ pairwise distinct},
\end{cases}
\end{equation}
then the second \(3\)-class groups \(\mathrm{G}_3^{(2)}{F_{c,\mu}}\) of the four fields \(1\le\mu\le 4\)
in dependence on
the \(3\)-class groups of the \(13\) bicyclic bicubic subfields \(S_i\), \(1\le i\le 13\), of the \(3\)-genus field \(F^\ast\)
are given by \(\left(\mathrm{G}_3^{(2)}{F_{c,\mu}}\right)_{1\le\mu\le 4}=\)
\begin{equation}
\label{eqn:Cat1Gph12Genus}
\begin{cases}
\langle 81,14\rangle;\langle 81,8\rangle,\langle 81,10\rangle^2 & \Longleftrightarrow
\left\lbrack\mathrm{Cl}_3{S_i}\right\rbrack=\lbrack(0)^3;(1^2)^7,(21)^3\rbrack, \\
\langle 243,46\rangle;\langle 243,25\rangle,\langle 243,28\rangle^2 & \Longleftrightarrow
\left\lbrack\mathrm{Cl}_3{S_i}\right\rbrack=\lbrack(0)^3;(1^2)^6,(21)^2,2^2,1^3\rbrack, \\
\langle 243,42\rangle;\langle 243,8\rangle^3 & \Longleftrightarrow
\left\lbrack\mathrm{Cl}_3{S_i}\right\rbrack=\lbrack(0)^3;(21)^9,21^2\rbrack.
\end{cases}
\end{equation}
\end{theorem}

\begin{proof}
Theorem
\ref{thm:Cat1Gph12}
is a consequence of the Tables
\ref{tbl:Cat1Gph1A},
\ref{tbl:Cat1Gph1B},
\ref{tbl:Cat1Gph2A}, and
\ref{tbl:Cat1Gph2B},
which have been computed with the aid of Magma
\cite{BCP1997,BCFS2022,Fi2001,MAGMA2022}.
\end{proof}


\begin{conjecture}
\label{cnj:Cat1Gph12}
Theorem
\ref{thm:Cat1Gph12}
remains true for any upper bound \(u>10^5\).
\end{conjecture}

\newpage

\renewcommand{\arraystretch}{1.1}

\begin{table}[ht]
\caption{Sixty Examples for Graph \(2\) of Category \(\mathrm{I}\)}
\label{tbl:Cat1Gph2A}
\begin{center}

{\scriptsize

\begin{tabular}{|rrc|c|c|}
\hline
    No. &       \(c\) &  \(\lbrack q_1,q_2,q_3\rbrack_3\) &                                   \(\mathrm{G}_3^{(2)}{F_{c,\mu}}\) & \(\left\lbrack\mathrm{Cl}_3{S_i}\right\rbrack_{1\le i\le 13}\) \\
\hline
  \(1\) &  \(7\,657\) & \(\lbrace 13\leftarrow 31\rightarrow 19\rbrace\) &                    \(\ast;\langle 243,8\rangle^3\) &                     \(\lbrack(0)^3;(21)^9,21^2\rbrack\) \\
  \(2\) &  \(8\,001\) &  \(\lbrace 9\leftarrow 127\rightarrow 7\rbrace\) & \(\langle 81,14\rangle;\langle 81,8\rangle,\langle 81,10\rangle^2\) & \(\lbrack(0)^3;(1^2)^7,(21)^3\rbrack\) \\
  \(3\) &  \(9\,709\) &  \(\lbrace 19\leftarrow 7\rightarrow 73\rbrace\) & \(\langle 81,14\rangle;\langle 81,8\rangle,\langle 81,10\rangle^2\) & \(\lbrack(0)^3;(1^2)^7,(21)^3\rbrack\) \\
  \(4\) & \(11\,137\) &  \(\lbrace 7\leftarrow 43\rightarrow 37\rbrace\) & \(\langle 81,14\rangle;\langle 81,8\rangle,\langle 81,10\rangle^2\) & \(\lbrack(0)^3;(1^2)^7,(21)^3\rbrack\) \\
  \(5\) & \(12\,753\) & \(\lbrace 9\leftarrow 109\rightarrow 13\rbrace\) & \(\langle 81,14\rangle;\langle 81,8\rangle,\langle 81,10\rangle^2\) & \(\lbrack(0)^3;(1^2)^7,(21)^3\rbrack\) \\
  \(6\) & \(14\,833\) & \(\lbrace 7\leftarrow 13\rightarrow 163\rbrace\) & \(\langle 81,14\rangle;\langle 81,8\rangle,\langle 81,10\rangle^2\) & \(\lbrack(0)^3;(1^2)^7,(21)^3\rbrack\) \\
  \(7\) & \(14\,911\) & \(\lbrace 13\leftarrow 31\rightarrow 37\rbrace\) & \(\langle 81,14\rangle;\langle 81,8\rangle,\langle 81,10\rangle^2\) & \(\lbrack(0)^3;(1^2)^7,(21)^3\rbrack\) \\
  \(8\) & \(16\,587\) &  \(\lbrace 9\leftarrow 19\rightarrow 97\rbrace\) & \(\langle 81,14\rangle;\langle 81,8\rangle,\langle 81,10\rangle^2\) & \(\lbrack(0)^3;(1^2)^7,(21)^3\rbrack\) \\
  \(9\) & \(17\,563\) & \(\lbrace 7\leftarrow 13\rightarrow 193\rbrace\) &                                     \(\ast;\langle 243,8\rangle^3\) &    \(\lbrack(0)^3;(21)^9,21^2\rbrack\) \\
 \(10\) & \(20\,167\) &  \(\lbrace 7\leftarrow 43\rightarrow 67\rbrace\) &                                     \(\ast;\langle 243,8\rangle^3\) &    \(\lbrack(0)^3;(21)^9,21^2\rbrack\) \\
 \(11\) & \(20\,881\) & \(\lbrace 19\leftarrow 7\rightarrow 157\rbrace\) &                                     \(\ast;\langle 243,8\rangle^3\) &    \(\lbrack(0)^3;(21)^9,21^2\rbrack\) \\
 \(12\) & \(21\,049\) &  \(\lbrace 7\leftarrow 97\rightarrow 31\rbrace\) &                    \(\langle 243,42\rangle;\langle 243,8\rangle^3\) &  \(\lbrack(0)^3;(21)^9,21^2\rbrack\) \\
 \(13\) & \(21\,177\) & \(\lbrace 9\leftarrow 181\rightarrow 13\rbrace\) &                                     \(\ast;\langle 243,8\rangle^3\) &    \(\lbrack(0)^3;(21)^9,21^2\rbrack\) \\
 \(14\) & \(23\,877\) &  \(\lbrace 7\leftarrow 379\rightarrow 9\rbrace\) &                                     \(\ast;\langle 243,8\rangle^3\) &    \(\lbrack(0)^3;(21)^9,21^2\rbrack\) \\
 \(15\) & \(24\,661\) & \(\lbrace 7\leftarrow 13\rightarrow 271\rbrace\) & \(\langle 81,14\rangle;\langle 81,8\rangle,\langle 81,10\rangle^2\) & \(\lbrack(0)^3;(1^2)^7,(21)^3\rbrack\) \\
 \(16\) & \(25\,123\) &  \(\lbrace 7\leftarrow 97\rightarrow 37\rbrace\) & \(\langle 81,14\rangle;\langle 81,8\rangle,\langle 81,10\rangle^2\) & \(\lbrack(0)^3;(1^2)^7,(21)^3\rbrack\) \\
 \(17\) & \(25\,207\) & \(\lbrace 7\leftarrow 13\rightarrow 277\rbrace\) & \(\langle 81,14\rangle;\langle 81,8\rangle,\langle 81,10\rangle^2\) & \(\lbrack(0)^3;(1^2)^7,(21)^3\rbrack\) \\
 \(18\) & \(27\,279\) &  \(\lbrace 7\leftarrow 433\rightarrow 9\rbrace\) & \(\langle 81,14\rangle;\langle 81,8\rangle,\langle 81,10\rangle^2\) & \(\lbrack(0)^3;(1^2)^7,(21)^3\rbrack\) \\
 \(19\) & \(30\,163\) & \(\lbrace 7\leftarrow 139\rightarrow 31\rbrace\) & \(\langle 81,14\rangle;\langle 81,8\rangle,\langle 81,10\rangle^2\) & \(\lbrack(0)^3;(1^2)^7,(21)^3\rbrack\) \\
 \(20\) & \(30\,411\) & \(\lbrace 9\leftarrow 109\rightarrow 31\rbrace\) & \(\langle 81,14\rangle;\langle 81,8\rangle,\langle 81,10\rangle^2\) & \(\lbrack(0)^3;(1^2)^7,(21)^3\rbrack\) \\
 \(21\) & \(32\,809\) & \(\lbrace 7\leftarrow 43\rightarrow 109\rbrace\) & \(\langle 81,14\rangle;\langle 81,8\rangle,\langle 81,10\rangle^2\) & \(\lbrack(0)^3;(1^2)^7,(21)^3\rbrack\) \\
 \(22\) & \(35\,113\) & \(\lbrace 13\leftarrow 73\rightarrow 37\rbrace\) & \(\langle 81,14\rangle;\langle 81,8\rangle,\langle 81,10\rangle^2\) & \(\lbrack(0)^3;(1^2)^7,(21)^3\rbrack\) \\
 \(23\) & \(36\,783\) &  \(\lbrace 61\leftarrow 9\rightarrow 67\rbrace\) & \(\langle 81,14\rangle;\langle 81,8\rangle,\langle 81,10\rangle^2\) & \(\lbrack(0)^3;(1^2)^7,(21)^3\rbrack\) \\
 \(24\) & \(37\,219\) & \(\lbrace 7\leftarrow 13\rightarrow 409\rbrace\) & \(\langle 81,14\rangle;\langle 81,8\rangle,\langle 81,10\rangle^2\) & \(\lbrack(0)^3;(1^2)^7,(21)^3\rbrack\) \\
 \(25\) & \(39\,753\) &  \(\lbrace 7\leftarrow 631\rightarrow 9\rbrace\) & \(\langle 81,14\rangle;\langle 81,8\rangle,\langle 81,10\rangle^2\) & \(\lbrack(0)^3;(1^2)^7,(21)^3\rbrack\) \\
 \(26\) & \(42\,883\) & \(\lbrace 19\leftarrow 37\rightarrow 61\rbrace\) & \(\langle 81,14\rangle;\langle 81,8\rangle,\langle 81,10\rangle^2\) & \(\lbrack(0)^3;(1^2)^7,(21)^3\rbrack\) \\
 \(27\) & \(\mathbf{48\,393}\) & \(\lbrace 9\leftarrow 19\rightarrow 283\rbrace\) & \(\ast;\langle 3^7,301\vert 305\rangle^3\) &             \(\lbrack(0)^3;(21)^6,(32)^3,2^21\rbrack\) \\
 \(28\) & \(48\,811\) & \(\lbrace 19\leftarrow 7\rightarrow 367\rbrace\) & \(\langle 81,14\rangle;\langle 81,8\rangle,\langle 81,10\rangle^2\) & \(\lbrack(0)^3;(1^2)^7,(21)^3\rbrack\) \\
 \(29\) & \(49\,149\) & \(\lbrace 9\leftarrow 127\rightarrow 43\rbrace\) & \(\langle 81,14\rangle;\langle 81,8\rangle,\langle 81,10\rangle^2\) & \(\lbrack(0)^3;(1^2)^7,(21)^3\rbrack\) \\
 \(30\) & \(53\,523\) & \(\lbrace 9\leftarrow 19\rightarrow 313\rbrace\) & \(\langle 81,14\rangle;\langle 81,8\rangle,\langle 81,10\rangle^2\) & \(\lbrack(0)^3;(1^2)^7,(21)^3\rbrack\) \\
\hline
\end{tabular}

}

\end{center}
\end{table}

\newpage

\renewcommand{\arraystretch}{1.1}

\begin{table}[ht]
\caption{Graph \(2\) of Category \(\mathrm{I}\) Continued}
\label{tbl:Cat1Gph2B}
\begin{center}

{\scriptsize

\begin{tabular}{|rrc|c|c|}
\hline
    No. &       \(c\) &  \(\lbrack q_1,q_2,q_3\rbrack_3\) &                                   \(\mathrm{G}_3^{(2)}{F_{c,\mu}}\) & \(\left\lbrack\mathrm{Cl}_3{S_i}\right\rbrack_{1\le i\le 13}\) \\
\hline
 \(31\) & \(54\,649\) & \(\lbrace 7\leftarrow 211\rightarrow 37\rbrace\) & \(\langle 81,14\rangle;\langle 81,8\rangle,\langle 81,10\rangle^2\) & \(\lbrack(0)^3;(1^2)^7,(21)^3\rbrack\) \\
 \(32\) & \(58\,387\) & \(\lbrace 19\leftarrow 7\rightarrow 439\rbrace\) &                    \(\langle 243,42\rangle;\langle 243,8\rangle^3\) &  \(\lbrack(0)^3;(21)^9,21^2\rbrack\) \\
 \(33\) & \(\mathbf{59\,031}\) & \(\lbrace 9\leftarrow 937\rightarrow 7\rbrace\) & \(\langle 2187,4670\rangle;\langle 3^6,52\rangle^3\) & \(\lbrack(0)^3;(21)^6,(2^2)^3,21^2\rbrack\) \\
 \(34\) & \(62\,109\) & \(\lbrace 67\leftarrow 9\rightarrow 103\rbrace\) & \(\langle 81,14\rangle;\langle 81,8\rangle,\langle 81,10\rangle^2\) & \(\lbrack(0)^3;(1^2)^7,(21)^3\rbrack\) \\
 \(35\) & \(63\,297\) & \(\lbrace 9\leftarrow 541\rightarrow 13\rbrace\) & \(\langle 81,14\rangle;\langle 81,8\rangle,\langle 81,10\rangle^2\) & \(\lbrack(0)^3;(1^2)^7,(21)^3\rbrack\) \\
 \(36\) & \(64\,519\) & \(\lbrace 7\leftarrow 13\rightarrow 709\rbrace\) & \(\langle 81,14\rangle;\langle 81,8\rangle,\langle 81,10\rangle^2\) & \(\lbrack(0)^3;(1^2)^7,(21)^3\rbrack\) \\
 \(37\) & \(65\,191\) & \(\lbrace 7\leftarrow 139\rightarrow 67\rbrace\) & \(\langle 81,14\rangle;\langle 81,8\rangle,\langle 81,10\rangle^2\) & \(\lbrack(0)^3;(1^2)^7,(21)^3\rbrack\) \\
 \(38\) & \(66\,969\) & \(\lbrace 7\leftarrow 1063\rightarrow 9\rbrace\) & \(\langle 81,14\rangle;\langle 81,8\rangle,\langle 81,10\rangle^2\) & \(\lbrack(0)^3;(1^2)^7,(21)^3\rbrack\) \\
 \(39\) & \(67\,249\) & \(\lbrace 7\leftarrow 13\rightarrow 739\rbrace\) & \(\langle 81,14\rangle;\langle 81,8\rangle,\langle 81,10\rangle^2\) & \(\lbrack(0)^3;(1^2)^7,(21)^3\rbrack\) \\
 \(40\) & \(68\,929\) & \(\lbrace 7\leftarrow 43\rightarrow 229\rbrace\) & \(\langle 81,14\rangle;\langle 81,8\rangle,\langle 81,10\rangle^2\) & \(\lbrack(0)^3;(1^2)^7,(21)^3\rbrack\) \\
 \(41\) & \(69\,939\) & \(\lbrace 9\leftarrow 19\rightarrow 409\rbrace\) & \(\langle 81,14\rangle;\langle 81,8\rangle,\langle 81,10\rangle^2\) & \(\lbrack(0)^3;(1^2)^7,(21)^3\rbrack\) \\
 \(42\) & \(71\,953\) & \(\lbrace 19\leftarrow 7\rightarrow 541\rbrace\) & \(\langle 81,14\rangle;\langle 81,8\rangle,\langle 81,10\rangle^2\) & \(\lbrack(0)^3;(1^2)^7,(21)^3\rbrack\) \\
 \(43\) & \(71\,991\) & \(\lbrace 9\leftarrow 19\rightarrow 421\rbrace\) &                    \(\langle 243,42\rangle;\langle 243,8\rangle^3\) &    \(\lbrack(0)^3;(21)^9,21^2\rbrack\) \\
 \(44\) & \(72\,541\) & \(\lbrace 7\leftarrow 43\rightarrow 241\rbrace\) & \(\langle 81,14\rangle;\langle 81,8\rangle,\langle 81,10\rangle^2\) & \(\lbrack(0)^3;(1^2)^7,(21)^3\rbrack\) \\
 \(45\) & \(77\,013\) & \(\lbrace 9\leftarrow 199\rightarrow 43\rbrace\) & \(\langle 81,14\rangle;\langle 81,8\rangle,\langle 81,10\rangle^2\) & \(\lbrack(0)^3;(1^2)^7,(21)^3\rbrack\) \\
 \(46\) & \(79\,513\) & \(\lbrace 7\leftarrow 307\rightarrow 37\rbrace\) & \(\langle 81,14\rangle;\langle 81,8\rangle,\langle 81,10\rangle^2\) & \(\lbrack(0)^3;(1^2)^7,(21)^3\rbrack\) \\
 \(47\) & \(80\,227\) & \(\lbrace 73\leftarrow 7\rightarrow 157\rbrace\) & \(\langle 81,14\rangle;\langle 81,8\rangle,\langle 81,10\rangle^2\) & \(\lbrack(0)^3;(1^2)^7,(21)^3\rbrack\) \\
 \(48\) & \(82\,899\) & \(\lbrace 61\leftarrow 9\rightarrow 151\rbrace\) & \(\langle 81,14\rangle;\langle 81,8\rangle,\langle 81,10\rangle^2\) & \(\lbrack(0)^3;(1^2)^7,(21)^3\rbrack\) \\
 \(49\) & \(83\,629\) & \(\lbrace 7\leftarrow 13\rightarrow 919\rbrace\) & \(\langle 81,14\rangle;\langle 81,8\rangle,\langle 81,10\rangle^2\) & \(\lbrack(0)^3;(1^2)^7,(21)^3\rbrack\) \\
 \(50\) & \(84\,409\) &\(\lbrace 13\leftarrow 151\rightarrow 43\rbrace\) & \(\langle 81,14\rangle;\langle 81,8\rangle,\langle 81,10\rangle^2\) & \(\lbrack(0)^3;(1^2)^7,(21)^3\rbrack\) \\
 \(51\) & \(85\,183\) & \(\lbrace 7\leftarrow 43\rightarrow 283\rbrace\) & \(\langle 81,14\rangle;\langle 81,8\rangle,\langle 81,10\rangle^2\) & \(\lbrack(0)^3;(1^2)^7,(21)^3\rbrack\) \\
 \(52\) & \(90\,097\) & \(\lbrace 7\leftarrow 211\rightarrow 61\rbrace\) & \(\langle 81,14\rangle;\langle 81,8\rangle,\langle 81,10\rangle^2\) & \(\lbrack(0)^3;(1^2)^7,(21)^3\rbrack\) \\
 \(53\) & \(91\,567\) &\(\lbrace 7\leftarrow 127\rightarrow 103\rbrace\) & \(\langle 81,14\rangle;\langle 81,8\rangle,\langle 81,10\rangle^2\) & \(\lbrack(0)^3;(1^2)^7,(21)^3\rbrack\) \\
 \(54\) & \(92\,241\) & \(\lbrace 9\leftarrow 37\rightarrow 277\rbrace\) & \(\langle 81,14\rangle;\langle 81,8\rangle,\langle 81,10\rangle^2\) & \(\lbrack(0)^3;(1^2)^7,(21)^3\rbrack\) \\
 \(55\) & \(93\,961\) & \(\lbrace 7\leftarrow 433\rightarrow 31\rbrace\) & \(\langle 81,14\rangle;\langle 81,8\rangle,\langle 81,10\rangle^2\) & \(\lbrack(0)^3;(1^2)^7,(21)^3\rbrack\) \\
 \(56\) & \(94\,877\) & \(\lbrace 9\leftarrow 811\rightarrow 13\rbrace\) &                                     \(\ast;\langle 243,8\rangle^3\) &    \(\lbrack(0)^3;(21)^9,21^2\rbrack\) \\
 \(57\) & \(94\,939\) &\(\lbrace 13\leftarrow 109\rightarrow 67\rbrace\) & \(\langle 81,14\rangle;\langle 81,8\rangle,\langle 81,10\rangle^2\) & \(\lbrack(0)^3;(1^2)^7,(21)^3\rbrack\) \\
 \(58\) & \(95\,157\) & \(\lbrace 9\leftarrow 109\rightarrow 97\rbrace\) & \(\langle 81,14\rangle;\langle 81,8\rangle,\langle 81,10\rangle^2\) & \(\lbrack(0)^3;(1^2)^7,(21)^3\rbrack\) \\
 \(59\) & \(98\,721\) & \(\lbrace 9\leftarrow 1567\rightarrow 7\rbrace\) & \(\langle 81,14\rangle;\langle 81,8\rangle,\langle 81,10\rangle^2\) & \(\lbrack(0)^3;(1^2)^7,(21)^3\rbrack\) \\
 \(60\) & \(99\,883\) & \(\lbrace 19\leftarrow 7\rightarrow 751\rbrace\) & \(\langle 81,14\rangle;\langle 81,8\rangle,\langle 81,10\rangle^2\) & \(\lbrack(0)^3;(1^2)^7,(21)^3\rbrack\) \\
\hline
\end{tabular}

}

\end{center}
\end{table}


\section{3-Towers of Length 3 over Quadratic Fields and Cyclic Cubic Fields}
\label{s:3ClassTowerLength3Comparison}

\noindent
It is very illuminating to compare three kinds of algebraic number fields
with distinct signatures, which share a common Artin pattern.
By a lucky coincidence, our computations yielded
a complex quadratic field \(F=\mathbb{Q}(\sqrt{d})\) with discriminant \(d=-17\,131\) and signature \((0,1)\),
two real quadratic fields \(F_i=\mathbb{Q}(\sqrt{d_i})\) with discriminants \(d_1=+8\,711\,453\), \(d_2=+9\,448\,265\) and signature \((2,0)\),
and four cyclic cubic fields \(F_{c,\mu}\) with conductor \(c=48\,393\) and signature \((3,0)\).
The conductor \(c\) belongs to Graph \(2\) of Category \(\mathrm{I}\),
since the combined cubic residue symbol is given by \(\lbrack q_1,q_2,q_3\rbrack_3=\lbrace 9\leftarrow 19\rightarrow 283\rbrace\).
All quadratic fields and three of the cyclic cubic fields,
having \(3\)-class groups \(\mathrm{Cl}_3{F}\) of type \((3,3)\), 
possess the Artin pattern
\(\mathrm{AP}=(\tau,\varkappa)\)
with \(\varkappa\sim (4231)\), of type \(\mathrm{G}.16\),
and \(\tau\sim\lbrack 32,21,21,21\rbrack\).

\noindent
\textbf{Notation.}
We are going to use
\textit{logarithmic type invariants} of abelian \(3\)-groups, for instance
\((321)\hat{=}(27,9,3)\), \((2^21)\hat{=}(9,9,3)\), and \((41^2)\hat{=}(81,3,3)\). 


\begin{figure}[ht]
\caption{Various Kinds of Fields with Tree Topology of Type \(\mathrm{G}.16\)}
\label{fig:TreeTopoG16}


{\tiny

\setlength{\unitlength}{0.8cm}

\begin{picture}(14,14)(0,-13)



\put(0,0){\line(0,-1){10}}
\multiput(-0.1,0)(0,-2){6}{\line(1,0){0.2}}

\put(-0.2,0){\makebox(0,0)[rc]{\(729\)}}
\put(0.2,0){\makebox(0,0)[lc]{\(3^6\)}}
\put(-0.2,-2){\makebox(0,0)[rc]{\(2\,187\)}}
\put(0.2,-2){\makebox(0,0)[lc]{\(3^7\)}}
\put(-0.2,-4){\makebox(0,0)[rc]{\(6\,561\)}}
\put(0.2,-4){\makebox(0,0)[lc]{\(3^8\)}}
\put(-0.2,-6){\makebox(0,0)[rc]{\(19\,683\)}}
\put(0.2,-6){\makebox(0,0)[lc]{\(3^9\)}}
\put(-0.2,-8){\makebox(0,0)[rc]{\(59\,049\)}}
\put(0.2,-8){\makebox(0,0)[lc]{\(3^{10}\)}}
\put(-0.2,-10){\makebox(0,0)[rc]{\(177\,147\)}}
\put(0.2,-10){\makebox(0,0)[lc]{\(3^{11}\)}}

\put(0,-10){\vector(0,-1){2}}
\put(0,-12.1){\makebox(0,0)[ct]{Order \(3^n\)}}



\put(7.8,0.4){\makebox(0,0)[rc]{\(\langle 54\rangle\)}}
\put(8,0){\circle*{0.2}}
\put(8.2,0.4){\makebox(0,0)[lc]{Fork, type c.21}}


\put(8,0){\line(-2,-1){4}}

\put(3.8,-2.4){\makebox(0,0)[rc]{\(\langle 301\rangle\)}}
\put(4,-2){\circle*{0.2}}

\put(5,-1.5){\makebox(0,0)[cc]{type G.16}}

\put(8,0){\line(-1,-1){2}}

\put(5.8,-2.4){\makebox(0,0)[rc]{\(\langle 305\rangle\)}}
\put(6,-2){\circle*{0.2}}

\put(4,-2){\line(0,-1){2}}

\put(3.8,-4.4){\makebox(0,0)[rc]{\(\langle 2048\rangle\)}}
\put(4,-4){\circle*{0.1}}

\put(6,-2){\line(0,-1){2}}

\put(5.8,-4.4){\makebox(0,0)[rc]{\(\langle 2058\rangle\)}}
\put(6,-4){\circle*{0.1}}

{\color{blue}
\put(2.8,-6.4){\makebox(0,0)[rc]{\(\#1;1\)}}
\put(2.9,-6.1){\framebox(0.2,0.2){}}
\put(3,-6){\vector(1,2){0.9}}
}

{\color{green}
\put(3.8,-6.4){\makebox(0,0)[rc]{\(2\)}}
\put(3.9,-6.1){\framebox(0.2,0.2){}}
\put(4,-6){\vector(0,1){1.9}}
}

\put(4.5,-7){\makebox(0,0)[cc]{strong \(\sigma\)-groups}}

{\color{blue}
\put(4.8,-6.4){\makebox(0,0)[rc]{\(1\)}}
\put(4.9,-6.1){\framebox(0.2,0.2){}}
\put(5,-6){\vector(1,2){0.9}}
}

{\color{green}
\put(5.8,-6.4){\makebox(0,0)[rc]{\(2\)}}
\put(5.9,-6.1){\framebox(0.2,0.2){}}
\put(6,-6){\vector(0,1){1.9}}
}


\put(8,0){\line(1,-2){2}}

{\color{orange}
\put(9.8,-4.4){\makebox(0,0)[rc]{\(\langle 619\rangle\)}}
\put(9.9,-4.1){\framebox(0.2,0.2){}}
\put(10,-4){\vector(-3,1){5.9}}
\put(11,-5){\makebox(0,0)[cc]{weak \(\sigma\)-groups}}
}

\put(11,-3.5){\makebox(0,0)[cc]{type G.16}}

\put(8,0){\line(1,-1){4}}

{\color{orange}
\put(11.8,-4.4){\makebox(0,0)[rc]{\(\langle 623\rangle\)}}
\put(11.9,-4.1){\framebox(0.2,0.2){}}
\put(12,-4){\vector(-3,1){5.9}}
}

\put(10,-4){\line(0,-1){2}}

\put(9.8,-6.4){\makebox(0,0)[rc]{\(\#1;4\)}}
\put(9.9,-6.1){\framebox(0.2,0.2){}}

\put(12,-4){\line(0,-1){2}}

\put(11.8,-6.4){\makebox(0,0)[rc]{\(4\)}}
\put(11.9,-6.1){\framebox(0.2,0.2){}}

\put(10,-6){\line(1,-2){2}}

{\color{purple}
\put(11.8,-10.4){\makebox(0,0)[rc]{\(\#2;1\)}}
\put(11.9,-10.1){\framebox(0.2,0.2){}}
\put(12,-10){\vector(-4,3){7.9}}
\put(13,-11){\makebox(0,0)[cc]{Schur \(\sigma\)-groups}}
}

\put(12,-6){\line(1,-2){2}}

{\color{purple}
\put(13.8,-10.4){\makebox(0,0)[rc]{\(1\)}}
\put(13.9,-10.1){\framebox(0.2,0.2){}}
\put(14,-10){\vector(-4,3){7.9}}
}


\put(5,-9){\makebox(0,0)[cc]{Topology Symbols:}}

{\color{orange}
\put(5,-10){\makebox(0,0)[cc]{\(\mathrm{G}.16\binom{1}{\rightarrow}\mathrm{c}.21\binom{2}{\leftarrow}\mathrm{G}.16\)}}
}

{\color{blue}
\put(3,-11){\makebox(0,0)[cc]{\(\mathrm{G}.16\binom{1}{\leftarrow}\mathrm{G}.16\)}}
}

{\color{green}
\put(7,-11){\makebox(0,0)[cc]{\(\mathrm{G}.16\binom{1}{\leftarrow}\mathrm{G}.16\)}}
}

{\color{purple}
\put(7,-12){\makebox(0,0)[cc]{\(\mathrm{G}.16\binom{1}{\rightarrow}\mathrm{G}.16\binom{1}{\rightarrow}\mathrm{c}.21\binom{2}{\leftarrow}\mathrm{G}.16\binom{1}{\leftarrow}\mathrm{G}.16\binom{2}{\leftarrow}\mathrm{G}.16\)}}
}


\end{picture}

}

\end{figure}
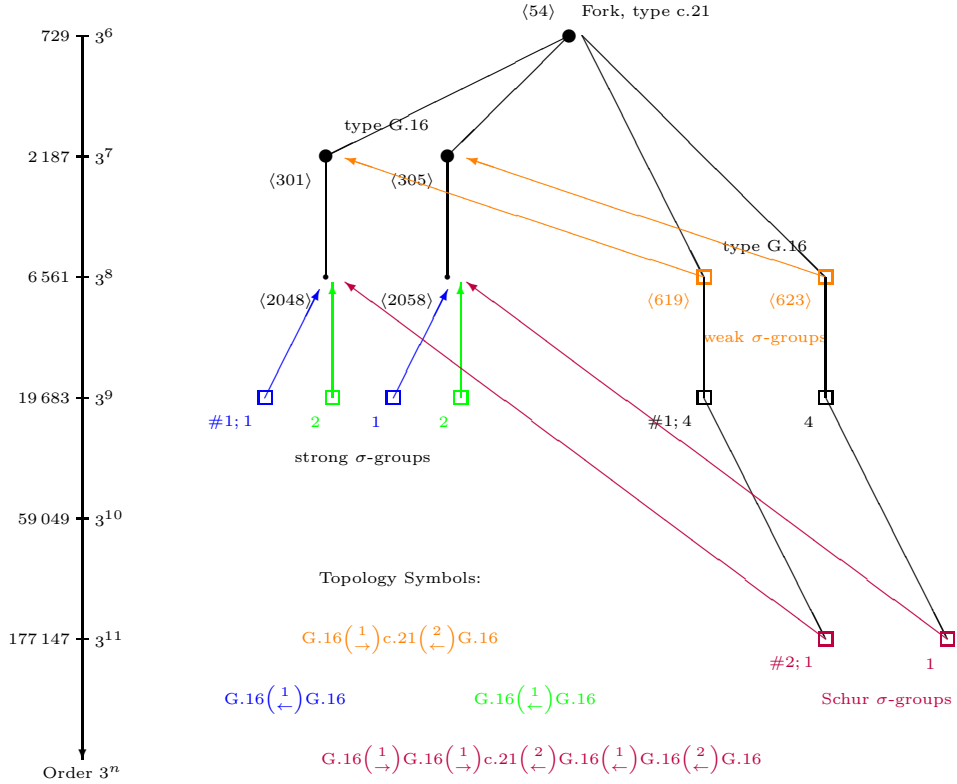


\noindent
Let \(F\) be a number field with \(\mathrm{Cl}_3\simeq (3,3)\)
and Artin pattern \(\mathrm{AP}(F)=(\tau(F),\varkappa(F))\),
where \(\varkappa(F)\sim (4231)\) is of type \(\mathrm{G}.16\)
and \(\tau(F)\sim\lbrack 32,21,21,21\rbrack\) indicates the ground state.
Denote by \(\mathfrak{M}=\mathrm{Gal}(F_3^{(2)}/F)\) the second \(3\)-class group
and by \(G=\mathrm{Gal}(F_3^{(\infty)}/F)\) the \(3\)-tower group of \(F\).


\begin{theorem}
\label{thm:ImgQdr}
\textbf{(Imaginary quadratic field with fork topology).} \\
If \(\tau^{(2)}{F}=\lbrack (32;321,(\mathbf{4}1^2)^3), (21;321,(\mathbf{3}1)^3)^3 \rbrack\)
for an imaginary quadratic field \(F=\mathbb{Q}(\sqrt{d})\), \(d<0\),
then \(\mathfrak{M}\simeq\langle 3^8,2048\vert 2058\rangle\),
and \(G\simeq\langle 3^8,619\vert 623\rangle-\#1;4-\#2;1\)
is a \textbf{Schur} \(\sigma\)-group.
\end{theorem}


\begin{theorem}
\label{thm:RealQdr}
\textbf{(Real quadratic field with child topology).} \\
If \(F=\mathbb{Q}(\sqrt{d})\), \(d>0\), is a real quadratic field, and \\
\(\tau^{(2)}{F}=\lbrack (32;321,(\mathbf{4}11)^3), (21;321,(\mathbf{2}1)^3)^3 \rbrack\),
respectively \\
\(\tau^{(2)}{F}=\lbrack (32;321,(\mathbf{3}11)^3), (21;321,(\mathbf{2}1)^3)^3 \rbrack\), \\
then \(\mathfrak{M}\simeq\langle 3^8,i\rangle\) with \(i=2048\), respectively \(i=2058\),
and
\(G\simeq\mathfrak{M}-\#1;j\) with \(j=1\),
respectively
\(j=2\),
is a \textbf{strong} \(\sigma\)-group.
\end{theorem}


\begin{theorem}
\label{thm:CyclCub}
\textbf{(Cyclic cubic field with fork topology).} \\
If \(\tau^{(2)}{F}=\lbrack (32;\mathbf{2^2}1,(\mathbf{3}1^2)^3), (21;\mathbf{2^2}1,(\mathbf{3}1)^3)^3 \rbrack\)
for a cyclic cubic field \(F\),
then \(\mathfrak{M}\simeq\langle 3^7,301\vert 305\rangle\),
and \(G\simeq\langle 3^8,619\vert 623\rangle\)
is a \textbf{weak} \(\sigma\)-group.
\end{theorem}

\begin{proof}
(Simultaneous proof of the Theorems
\ref{thm:ImgQdr},
\ref{thm:RealQdr},
and
\ref{thm:CyclCub}.)
Due to the various signatures of the different kinds of fields,
the torsionfree Dirichlet unit ranks are given by
\(r=0\), \(r=1\), \(r=2\), respectively.
According to the Shafarevich Theorem
\cite{Sh1964},
the relation rank of the \(3\)-tower group is therefore bounded by
\(d_2\le 2\), \(d_2\le 3\), \(d_2\le 4\), respectively.
Exclusively the Schur \(\sigma\)-groups
\(\langle 3^8,i\rangle-\#1;4-\#2;1\), \(i\in\lbrace 619,623\rbrace\),
have a balanced presentation with \(d_2=2\).
Thus they provide the unique possibility for the fastidious imaginary quadratic fields.
The real quadratic fields are happy with a strong \(\sigma\)-group, like
\(\langle 3^8,i\rangle-\#1;1\) or \(\langle 3^8,i\rangle-\#1;2\),
\(i\in\lbrace 2048,2058\rbrace\),
in dependence on the IPAD of second order \(\tau^{(2)}{F}\).
The cover of \(\langle 3^7,i\rangle\), \(i\in\lbrace 301,305\rbrace\),
does not contain any strong \(\sigma\)-group, let alone a Schur \(\sigma\)-group.
It is thus forbidden as second \(3\)-class group for both, real and imaginary quadratic fields.
Only the frugal cyclic cubic fields are satisfied
with a metabelianization of their \(3\)-tower group
\(\langle 3^8,619\vert 623\rangle\),
which is only a weak \(\sigma\)-group
(even a non-\(\sigma\) group would do it).
\end{proof}


\begin{example}
The smallest concrete realizations of the fields in the Theorems
\ref{thm:ImgQdr},
\ref{thm:RealQdr},
and
\ref{thm:CyclCub},
whose tree topologies are drawn in Figure
\ref{fig:TreeTopoG16},
have the following discriminants, conductors:

{\color{purple}
\(F=\mathbb{Q}(\sqrt{d})\) with discriminant \(d=-17\,131\),
}

{\color{blue}
\(F=\mathbb{Q}(\sqrt{d})\) with discriminant \(d=+8\,711\,453\),
}

{\color{green}
\(F=\mathbb{Q}(\sqrt{d})\) with discriminant \(d=+9\,448\,265\),
}

{\color{orange}
\(F\) cyclic cubic field with conductor \(c=48\,393\) (Graph \(2\) of Category \(\mathrm{I}\)).
}

\end{example}


\begin{theorem}
\label{thm:c21}
\textbf{(Three stage tower \(\mathbf{c.21}\)).}
Let \(F\) be a cyclic cubic field with \(3\)-class group \(\mathrm{Cl}_3{F}\simeq (3,3)\),
capitulation type \(\mathrm{c}.21\), \(\varkappa\sim (0231)\), and
(logarithmic) abelian type invariants of \textbf{second order}
\[\tau^{(2)}=\left(11;\lbrack 22;(211)^4\rbrack,\lbrack 21;211,(\mathbf{3}1)^3\rbrack,\lbrack 21;211,(21)^3\rbrack^2\right).\]
Then \(F\) has a \(3\)-class field tower of precise length \(\ell_3{F}=3\)
with group \(G=\mathrm{Gal}(F_3^{(\infty)}/F)\) isomorphic to
either \(\langle 2187,307\rangle\) or \(\langle 2187,308\rangle\),
the latter with double probability.
\end{theorem}

\begin{proof}
The isomorphism class of the metabelianization \(\mathfrak{M}=G/G^{\prime\prime}\) of \(G\)
is determined uniquely as \(\langle 729,54\rangle=U\) by the Artin pattern \(\mathrm{AP}=(\tau,\varkappa)\)
with \(\tau=\lbrack 22,21,21,21\rbrack\) and \(\varkappa=(0231)\).
The Shafarevich cover
\cite{Ma2016c}
with respect to \(F\),
\[\mathrm{cov}_F(\mathfrak{M})=\lbrace\langle 729,54\rangle,\langle 2187,307\rangle,\langle 2187,308\rangle\rbrace,\]
is non-trivial, finite and \textit{inhomogeneous}.
Since all members have admissible relation ranks \(3\le d_2\le 4\),
the decision between \(\ell_3{F}=2\) and \(\ell_3{F}=3\) requires abelian type invariants of second order \(\tau^{(2)}\),
because \(\langle 729,54\rangle\) has
\[\tau^{(2)}=\left(11;\lbrack 22;(211)^4\rbrack,\lbrack 21;211,(21)^3\rbrack^3\right).\]
Finally the order of the automorphism group of \(\langle 2187,307\rangle\) is bigger by a factor of \(2\)
than that of \(\langle 2187,308\rangle\).
\end{proof}

\begin{example}
\label{exm:c21}
All three fields of type \((3,3)\) in the quartets
for each of the conductors \(c=\mathbf{28\,791}\), \(\mathbf{32\,227}\), \(\mathbf{61\,087}\),
i.e. nine fields,
which belong to Graph \(1\) of Category \(\mathrm{I}\),
satisfy the conditions of Theorem
\ref{thm:c21}.
(The fourth field is of type \((9,3,3)\).)
\end{example}


\begin{theorem}
\label{thm:d19}
\textbf{(Three stage tower \(\mathbf{d.19}\)).}
Let \(F\) be a cyclic cubic field with \(3\)-class group \(\mathrm{Cl}_3{F}\simeq (3,3)\),
capitulation type \(\mathrm{d}.19\), \(\varkappa\sim (4043)\), and
abelian type invariants of \textbf{second order}
\[\tau^{(2)}=\left(11;\lbrack 22;(211)^4\rbrack,\lbrack 21;211,(21)^3\rbrack,\lbrack 111;\mathbf{(211)^4},(11)^9\rbrack,\lbrack 111;211,(111)^3,(11)^9\rbrack\right).\]
Then \(F\) has a \(3\)-class field tower of precise length \(\ell_3{F}=3\)
with group \(G=\mathrm{Gal}(F_3^{(\infty)}/F)\) isomorphic to \(\langle 2187,265\rangle\).
\end{theorem}

\begin{proof}
The isomorphism class of the metabelianization \(\mathfrak{M}=G/G^{\prime\prime}\) of \(G\)
is determined uniquely as \(\langle 729,41\rangle=D\) by the Artin pattern \(\mathrm{AP}=(\tau,\varkappa)\)
with \(\tau=\lbrack 22,21,111,111\rbrack\) and \(\varkappa=(4043)\).
The Shafarevich cover
\cite{Ma2016c}
with respect to \(F\),
\[\mathrm{cov}_F(\mathfrak{M})=\lbrace\langle 729,41\rangle,\langle 2187,263\rangle,\langle 2187,264\rangle,\langle 2187,265\rangle\rbrace,\]
is non-trivial, finite and \textit{inhomogeneous}.
Since all members have admissible relation ranks \(3\le d_2\le 4\),
the decision between \(\ell_3{F}=2\) and \(\ell_3{F}=3\) requires abelian type invariants of second order \(\tau^{(2)}\)
and is possible only for \(\langle 2187,265\rangle\),
since \(\langle 729,41\rangle\), \(\langle 2187,263\rangle\), \(\langle 2187,264\rangle\) have the same pattern
\[\tau^{(2)}=\left(11;\lbrack 22;(211)^4\rbrack,\lbrack 21;211,(21)^3\rbrack,\lbrack 111;211,(111)^3,(11)^9\rbrack^2\right).\qedhere\]
\end{proof}

\begin{example}
\label{exm:d19}
A single field for each of the conductors \(c=\mathbf{22\,581}\), \(\mathbf{34\,489}\), \(\mathbf{56\,547}\), \(\mathbf{79\,933}\),
i.e. four fields,
which belong to Graph \(1,1,2,2\) of Category \(\mathrm{II}\), respectively,
satisfy the conditions of Theorem
\ref{thm:d19}.
\end{example}


\begin{theorem}
\label{thm:333}
\textbf{(World record three stage tower of type \(\mathbf{(3,3,3)}\)).}
Let \(F\) be a cyclic cubic field with \(3\)-class group \(\mathrm{Cl}_3{F}\simeq (3,3,3)\),
capitulation type \(\varkappa\sim (O^3P^{10})\), and
abelian type invariants of first order \(\tau^{(1)}=\left\lbrack111;1111,(111)^3,(21)^9\right\rbrack\).
Then \(F\) has a \(3\)-class field tower of exact length \(\ell_3{F}=3\)
with group \(\mathrm{Gal}(F_3^{(\infty)}/F)\) isomorphic to one of the nine groups
\(\langle 6561,n\rangle\) with \(n\in\lbrace 261262,\ldots,261270\rbrace\).
\end{theorem}

\begin{proof}
There are four possible isomorphism classes of the metabelianization \(\mathfrak{M}=G/G^{\prime\prime}\) of \(G\),
namely \(\langle 2187,m\rangle\) with \(m\in\lbrace 5576,\ldots,5579\rbrace\)
sharing the common Artin pattern \(\mathrm{AP}=(\tau,\varkappa)\).
The cover
\cite{Ma2016c}
of the latter three groups is given by
\[\mathrm{cov}(\mathfrak{M})=\lbrace\langle 2187,5577\rangle,\langle 6561,261262\rangle,\langle 6561,261263\rangle,\langle 6561,261264\rangle\rbrace, \text{ resp.}\]
\[\mathrm{cov}(\mathfrak{M})=\lbrace\langle 2187,5578\rangle,\langle 6561,261265\rangle,\langle 6561,261266\rangle,\langle 6561,261267\rangle\rbrace, \text{ resp.}\]
\[\mathrm{cov}(\mathfrak{M})=\lbrace\langle 2187,5579\rangle,\langle 6561,261268\rangle,\langle 6561,261269\rangle,\langle 6561,261270\rangle\rbrace.\]
However, the Shafarevich cover \(\mathrm{cov}_F(\mathfrak{M})\) with respect to \(F\) consists of the non-metabelian groups with order \(6561\) only.
The reason for this fact is that the relation rank of the metabelian group \(\mathfrak{M}\) takes the forbidden value \(6\).
So, the Shafarevich cover is non-trivial, finite and \textit{homogeneous},
which dispenses us with calculating abelian type invariants of second order, \(\tau^{(2)}\).
\end{proof}

\begin{example}
\label{exm:333}
Two fields with conductors \(c\in\lbrace\mathbf{22\,581},\mathbf{34\,489}\rbrace\),
which belong to Graph \(1\) of Category \(\mathrm{II}\),
and two fields with conductors \(c\in\lbrace\mathbf{56\,547},\mathbf{79\,933}\rbrace\),
which belong to Graph \(2\) of Category \(\mathrm{II}\),
satisfy the conditions of Theorem
\ref{thm:333}.
\end{example}


\part{Applications}
\label{pt:Applications}

\chapter{Galois Action of Cyclic Fields}
\label{ch:GaloisAction}

\section{p-Capitulation Enforced by Galois Action}
\label{s:GaloisAction}

\noindent
The generating automorphism \(\sigma\)
of a cyclic number field \({F}/\mathbb{Q}\) of degree \(d\)
with Galois group \(\mathrm{Gal}({F}/\mathbb{Q})=\langle\sigma\rangle\)
acts on the class group \(\mathrm{Cl}_{F}\) of \({F}\)
and thus also on the higher \(p\)-class groups \(\mathrm{G}_p^{(n)}{{F}}\)
with \(n\in\mathbb{N}\cup\lbrace\infty\rbrace\)
\cite{Ma2013},
for a fixed prime number \(p\).
When \(d\) and \(p\) are coprime,
a remarkable restriction of the possibilities for the
metabelian second \(p\)-class group \(\mathfrak{M}=\mathrm{G}_p^{(2)}{{F}}\),
and consequently for the transfer kernel type \(\varkappa({F})\) of \({F}\),
is due to the fact that
the trace \(T_\sigma=\sum\limits_{\substack{i=0}}^{d-1}\sigma^i\) of \(\sigma\) annihilates
the commutator quotient of all the groups \(\mathrm{G}_p^{(n)}{{F}}\).


\begin{definition}
\label{dfn:AnnTrace}
Let \(p\) be a prime number and
\(G\) be a pro-\(p\) group with finite abelianization \(G/G^\prime\).
Suppose that \(d\ge 2\) is a fixed integer.
\(G\) is said to be a \(\sigma\)-\textbf{group of degree} \(d\),
if \(G\) possesses an automorphism \(\sigma\) of order \(d\)
whose trace
\(T_\sigma=\sum\limits_{\substack{j=0}}^{d-1}\,\sigma^{j}\in\mathbb{Z}\lbrack\mathrm{Aut}(G)\rbrack\)
annihilates \(G\) modulo \(G^\prime\),
that is, if there exists \(\sigma\in\mathrm{Aut}(G)\) such that \(\mathrm{ord}(\sigma)=d\) and
\(x^{T_\sigma}=\prod\limits_{\substack{j=0}}^{d-1}\,\sigma^{j}(x)\in G^\prime\) for all \(x\in G\).
\end{definition}


The following theorem is the cubic analogue of
\cite[Thm. 2.2]{AKMTT2020}.

\begin{theorem}
\label{thm:CycCub}
The \(p\)-class tower group \(G_p^{(\infty)}{{F}}\) and
all higher \(p\)-class groups \(G_p^{(n)}{{F}}\), \(n\ge 2\),
of a cyclic cubic number field \({F}\)
are \(\sigma\)-groups of degree \(3\), for any prime \(p\).
\end{theorem}

\begin{proof}
The generating automorphism \(\sigma\) of \(F/\mathbb{Q}\)
annihilates the class group \(\mathrm{Cl}_F\)
when it acts by its trace \(T_\sigma=\sum_{i=0}^2\,\sigma^i\in\mathbb{Z}\lbrack\langle\sigma\rangle\rbrack\),
since \(x^{T_\sigma}=\prod_{i=0}^2\,\sigma^i(x)=\mathrm{Norm}_{F/\mathbb{Q}}(x)\in\mathrm{Cl}_\mathbb{Q}=1\),
for all \(x\in\mathrm{Cl}_F\).
Of course, the same is true for all \(p\)-class groups \(\mathrm{Cl}_p{F}\) with primes \(p\).
Finally, we have isomorphisms
\(\mathrm{G}_p^{(n)}{F}/\left(\mathrm{G}_p^{(n)}{F}\right)^\prime\simeq\mathrm{Cl}_p{F}\),
for any \(n\in\mathbb{N}\cup\lbrace\infty\rbrace\)
\cite{Ma2013}.
\end{proof}


\renewcommand{\arraystretch}{1.1}

\begin{table}[ht]
\caption{Number of Automorphisms with Order \(3\) of \(p\)-Groups \(G\)}
\label{tbl:AutOrd3}
\begin{center}
{\scriptsize
\begin{tabular}{|rc|rr|cc|rl|rrr|}
\hline
    No. & \(p\) &    Order &    Id  & \(\nu\) & \(\mu\) &   \(\#\)Aut &       Factors &     \(o\) &     \(w\) &   \(s\) \\
\hline
  \(1\) & \(2\) &    \(4\) &  \(2\) &   \(3\) &   \(3\) &       \(6\) &       \(2,3\) &     \(2\) &     \(2\) &   \(2\) \\
  \(2\) & \(2\) &    \(8\) &  \(3\) &   \(1\) &   \(3\) &       \(8\) &       \(2^3\) &     \(0\) &     \(0\) &   \(0\) \\
  \(3\) & \(2\) &    \(8\) &  \(4\) &   \(0\) &   \(2\) &      \(24\) &     \(2^3,3\) &     \(8\) &     \(8\) &   \(2\) \\
\hline
  \(4\) & \(2\) &    \(8\) &  \(5\) &   \(1\) &   \(3\) &     \(168\) &   \(2^3,3,7\) &    \(56\) &     \(0\) &   \(0\) \\
\hline
\hline
  \(1\) & \(5\) &   \(25\) &  \(2\) &   \(3\) &   \(3\) &     \(480\) &   \(2^5,3,5\) &    \(20\) &    \(20\) &  \(20\) \\
  \(2\) & \(5\) &  \(125\) &  \(3\) &   \(2\) &   \(4\) &   \(12000\) & \(2^5,3,5^3\) &   \(500\) &   \(500\) &  \(20\) \\
  \(3\) & \(5\) &  \(125\) &  \(4\) &   \(0\) &   \(2\) &     \(500\) &   \(2^2,5^3\) &     \(0\) &     \(0\) &   \(0\) \\
  \(4\) & \(5\) &  \(625\) &  \(7\) &   \(1\) &   \(4\) &   \(50000\) &   \(2^4,5^5\) &     \(0\) &     \(0\) &   \(0\) \\
  \(5\) & \(5\) &  \(625\) &  \(8\) &   \(0\) &   \(3\) &   \(12500\) &   \(2^2,5^5\) &     \(0\) &     \(0\) &   \(0\) \\
  \(6\) & \(5\) &  \(625\) &  \(9\) &   \(0\) &   \(3\) &    \(5000\) &   \(2^3,5^4\) &     \(0\) &     \(0\) &   \(0\) \\
  \(7\) & \(5\) &  \(625\) & \(10\) &   \(0\) &   \(3\) &    \(5000\) &   \(2^3,5^4\) &     \(0\) &     \(0\) &   \(0\) \\
\hline
  \(8\) & \(5\) & \(3125\) &  \(3\) &   \(3\) &   \(5\) & \(7500000\) & \(2^5,3,5^7\) & \(12500\) & \(12500\) & \(500\) \\
  \(9\) & \(5\) & \(3125\) &  \(4\) &   \(1\) &   \(3\) &   \(62500\) &   \(2^2,5^6\) &     \(0\) &     \(0\) &   \(0\) \\
 \(10\) & \(5\) & \(3125\) &  \(5\) &   \(1\) &   \(3\) &  \(625000\) &   \(2^3,5^7\) &     \(0\) &     \(0\) &   \(0\) \\
 \(11\) & \(5\) & \(3125\) &  \(6\) &   \(1\) &   \(3\) &  \(625000\) &   \(2^3,5^7\) &     \(0\) &     \(0\) &   \(0\) \\
 \(12\) & \(5\) & \(3125\) &  \(7\) &   \(1\) &   \(3\) &  \(125000\) &   \(2^3,5^6\) &     \(0\) &     \(0\) &   \(0\) \\
 \(13\) & \(5\) & \(3125\) &  \(8\) &   \(0\) &   \(2\) &  \(156250\) &     \(2,5^7\) &     \(0\) &     \(0\) &   \(0\) \\
 \(14\) & \(5\) & \(3125\) &  \(9\) &   \(0\) &   \(2\) &   \(93750\) &   \(2,3,5^6\) &  \(1250\) &  \(1250\) &  \(50\) \\
 \(15\) & \(5\) & \(3125\) & \(10\) &   \(1\) &   \(3\) &  \(187500\) & \(2^2,3,5^6\) &  \(1250\) &  \(1250\) &  \(50\) \\
 \(16\) & \(5\) & \(3125\) & \(11\) &   \(0\) &   \(2\) &   \(62500\) &   \(2^2,5^6\) &     \(0\) &     \(0\) &   \(0\) \\
 \(17\) & \(5\) & \(3125\) & \(12\) &   \(0\) &   \(2\) &   \(93750\) &   \(2,3,5^6\) &  \(1250\) &  \(1250\) &  \(50\) \\
 \(18\) & \(5\) & \(3125\) & \(13\) &   \(0\) &   \(2\) &  \(156250\) &     \(2,5^7\) &     \(0\) &     \(0\) &   \(0\) \\
 \(19\) & \(5\) & \(3125\) & \(14\) &   \(0\) &   \(2\) & \(1875000\) & \(2^3,3,5^7\) & \(12500\) & \(12500\) & \(500\) \\
\hline
\end{tabular}
}
\end{center}
\end{table}


\noindent
In Table
\ref{tbl:AutOrd3},
we denote
the SmallGroups identifier
\cite{BEO2005}
by Id,
the nuclear rank by \(\nu\),
the \(p\)-multiplicator rank by \(\mu\)
\cite{Ma2015a}.
Further we give the order \(\#\)Aut of the automorphism group \(\mathrm{Aut}(G)\)
and its prime power factors.
\(o\) denotes the cardinality of the set \(\mathrm{Aut}_3(G):=\lbrace\sigma\in\mathrm{Aut}(G)\mid\mathrm{ord}(\sigma)=3\rbrace\),
\(w:=\#\lbrace\sigma\in\mathrm{Aut}_3(G)\mid (\forall x\in G)\,x^{1+\sigma+\sigma^2}\in G^\prime\rbrace\)
the number of \textit{weak} \(\sigma\)-automorphisms of degree \(3\), and
\(s:=\#\lbrace\sigma\in\mathrm{Aut}_3(G)\mid (\forall x\in G)\,x^{1+\sigma+\sigma^2}=1\rbrace\)
the number of \textit{strong} \(\sigma\)-automorphisms of degree \(3\).

\noindent
For the dominating part of the finite \(p\)-groups in Table
\ref{tbl:AutOrd3},
the failure of being a \(\sigma\)-group of degree \(3\)
is a consequence of \(\gcd(3,\#\mathrm{Aut}(G))=1\) already.
For the elementary abelian \(2\)-group \(\langle 8,5\rangle\) of rank \(3\), however,
the computation of \(w\) is required for the decision,
since \(v_3(\#\mathrm{Aut}(G))=1\).


\noindent
In view of our special situation with \(p\in\lbrace 2,5\rbrace\) and \(\mathrm{Cl}_p{F}=(p,p)\),
we tested finite metabelian \(p\)-groups \(G\) with \(G/G^\prime\simeq (p,p)\)
for the property of being a \(\sigma\)-group of degree \(3\).
The following theorem is the analogue of
\cite[Thm. 2.3]{AKMTT2020}
for degree \(3\).


\begin{theorem}
\label{thm:SigmaGroupsDegree3}
We characterize finite \(p\)-groups by their identifier in the SmallGroups Library
\cite{BEO2002,BEO2005},
which gives the order, \(\mathrm{ord}\), of the group and a counting number, \(\mathrm{id}\),
enclosed in angle brackets \(\langle\mathrm{ord},\mathrm{id}\rangle\).
\begin{enumerate}
\item
A finite \(2\)-group \(G\) with \(G/G^\prime\simeq (2,2)\) which is a \(\sigma\)-group of degree \(3\)
is isomorphic to either the \textbf{abelian group} \(\langle 4,2\rangle\) of type \((2,2)\)
or the \textbf{quaternion group} \(\langle 8,4\rangle\).
\item
There are no finite \(2\)-groups \(G\) with \(G/G^\prime\simeq (2,2,2)\) which are \(\sigma\)-groups of degree \(3\).
(Thus, there are no cyclic cubic fields of type \((2,2,2)\).)
\item
A finite \(5\)-group \(G\) with \(G/G^\prime\simeq (5,5)\) which is a \(\sigma\)-group of degree \(3\)
is isomorphic to either the \textbf{abelian group} \(\langle 25,2\rangle\) of type \((5,5)\)
or the \textbf{extra special group} \(\langle 125,3\rangle\)
or one of the three \textbf{Schur \(\sigma\)-groups} \(\langle 3125,i\rangle\) with \(i\in\lbrace 9,12,14\rbrace\)
(and TKT two \(3\)-cycles or one \(6\)-cycle or identity permutation)
or a descendant of one of the two groups
\(\langle 3125,i\rangle\) with \(i\in\lbrace 3,10\rbrace\).
\end{enumerate}
\end{theorem}

\begin{proof}
Using permutation representations,
we compiled a program script in Magma
\cite{MAGMA2022}
in order to test
whether an assigned \(p\)-group \(G\) with \(G/G^\prime\simeq (p,p)\) is a \(\sigma\)-group of degree \(3\),
for any prime number \(p\ne 3\).
See Table
\ref{tbl:AutOrd3}.
According to
\cite[Cor. 2.1]{AKMTT2020},
all descendants of a group \(R\) which is not a \(\sigma\)-group of degree \(3\)
share this property with their ancestor \(R\).
Consequently, the descendant tree \(\mathcal{T}^{(1)}\langle 8,3\rangle\)
\cite[Fig. 3.1, p. 419]{Ma2013}
is entirely forbidden for \(\mathrm{G}_2^{(2)}{F}\) of cyclic cubic fields \(F\),
and only the abelian root \(\langle 4,2\rangle\) and its terminal immediate descendant \(\langle 8,4\rangle\) are admissible.
Similarly, the complete coclass tree \(\mathcal{T}^{(1)}\langle 625,7\rangle\)
is forbidden for \(\mathrm{G}_5^{(2)}{F}\) of cyclic cubic fields \(F\).
Table
\ref{tbl:AutOrd3}
was computed by Algorithm
\ref{alg:TraceSigma},
rather than Algorithm
\ref{alg:ActionFrattiniQuotient}.
\end{proof}


\begin{algorithm}
\label{alg:TraceSigma}
(Action of the Trace of an Automorphism.) \\
\textbf{Input:}
degree \texttt{d}, compact presentation \texttt{sP} of group \texttt{G}. \\
\textbf{Code:}
implemented as an intrinsic procedure \texttt{TraceAut()}.
{\footnotesize
\texttt{
\begin{tabbing}
for \= for \= for \= for \= for \= \kill
intrinsic TraceAut(d::RngIntElt,sP::[RngIntElt]) \(\lbrace\,\rbrace\)\\
G := PCGroup(sP); // G as group with power-conjugate presentation\\
Id := ""; // identifier of G\\
if CanIdentifyGroup(Order(G)) then\+\\
   Id := IdentifyGroup(G);\-\\
end if;\\
n := Ngens(G); // number of generators of G\\
IG := Identity(G); // neutral element of G\\
DG := DerivedSubgroup(G); // commutator subgroup of G\\
A := AutomorphismGroup(G);\\
T,AP,R := PermutationRepresentation(A); // much CPU time (slow)\\
c := 0;\\
o := 0;\\
w := 0;\\
s := 0;\\
for a in AP do\+\\
   c := c+1; // counter of all automorphisms\\
   if (Order(a) eq d) then\+\\
      o := o+1; // sigma-aut of assigned degree \\
      bw := true;\\
      bs := true;\\
      for i in [1..n] do // let trace act on all generators\+\\
         tr := G.i\({}\,\hat{}\,{}\)0;\\
         for j in [0..d-1] do\+\\
            tr := tr\(\ast\)(G.i@((a@@T)\({}\,\hat{}\,{}\)j));\-\\
         end for; // j\\
         bw := bw and (tr in DG); // weak\\
         bs := bs and (tr eq IG); // strong\-\\
      end for; // i\\
      if (bw eq true) then\+\\
         w := w+1; // weak sigma-aut\-\\
      end if; // bw\\
      if (bs eq true) then\+\\
         s := s+1; // strong sigma-aut\-\\
      end if; // bs\-\\
   end if; // d\-\\
end for; // a\\
printf"\%o c=\%o, o=\%o, w=\%o, s=\%o\(\backslash\)n",Id,c,o,w,s;\\
end intrinsic; // TraceAut
\end{tabbing}
}
}
\noindent
\textbf{Output:}
Identifier \texttt{Id} and all counters \texttt{c}, \texttt{o}, \texttt{w}, \texttt{s} of various \(\sigma\)-automorphisms.
\end{algorithm}


\noindent
The following algorithm checks whether the Frattini quotient \(G/\Phi(G)\)
of an assigned \(p\)-group \(G\) has an action by the cyclic group \(C_3\) of order \(3\).
It can easily be modified by replacing \(C_3\simeq\langle 3,1\rangle\)
with another critical group, for instance,
the symmetric group \(S_{3}\simeq\langle 6,1\rangle\) of order \(6\) or
the metacyclic group \(M_{5}\simeq\langle 20,3\rangle\) of order \(20\).
The Magma function \texttt{pQuotient()} with parameters \(G,p,1\) returns
the Frattini quotient \(Q=G/\Phi(G)\) of \(G\)
together with the natural projection \(qp:\,G\to Q\).
The composition of mappings runs from the left to the right in Magma,
so \texttt{A.j\(\ast\)qp} means \(qp\circ A_j\),
and the inverse image \(qp^{-1}(Q_k)\) is denoted by \texttt{Q.k@@qp},
so \(B=(qp\circ A_j\circ qp^{-1})(Q_k)\) is the image of \(Q_k\)
under the composite mapping \(qp\circ A_j\circ qp^{-1}:\,Q\to Q\).
In the (rare) case that the subgroup \(U\le\mathrm{Aut}(Q)\)
cannot be identified, due to a huge order,
the group \(G\) is considered admissible,
in order to avoid an erroneous elimination.

\begin{algorithm}
\label{alg:ActionFrattiniQuotient}
(Action on the Frattini Quotient.) \\
\textbf{Input:}
group \texttt{G}, prime \texttt{p}. \\
\textbf{Code:}
implemented as a function \texttt{IsAdmissible()} with boolean return value.
{\footnotesize
\texttt{
\begin{tabbing}
for \= for \= for \= for \= for \= \kill
IsAdmissible := function(G,p);\+\\
   A := AutomorphismGroup(G);\\
   Q,qp := pQuotient(G,p,1);\\
   Lj := [];\\
   for j in [1..Ngens(A)] do\+\\
      Lk := [];\\
      for k in [1..Ngens(Q)] do\+\\
         B := (A.j\(\ast\)qp)(Q.k@@qp);\\
         Append(\(\sim\)Lk,<Q.k,B>);\-\\
      end for; // k\\
      h := hom<Q->Q|Lk>;\\ 
      Append(\(\sim\)Lj,h);\-\\
   end for; // j\\
   AQ := AutomorphismGroup(Q);\\
   U := sub<AQ|Lj>;\\
   o := Order(U);\\
   sN := [o,0];\\
   Boole := true;\\
   if CanIdentifyGroup(o) then\+\\
      sN := IdentifyGroup(U);\\
      PCU := SmallGroup(sN[1],sN[2]);\\
      Critical := SmallGroup(3,1); // cyclic group of order 3\\
      SU := Subgroups(PCU);\\
      Boole := false;\\
      for i in [1..\#SU] do\+\\
         if IsIsomorphic(Critical,SU[i]\({}\,\grave{}\,{}\)subgroup) then\+\\
            Boole := true;\\
            break;\-\\
         end if;\-\\
      end for; // i\-\\
   end if;\\
   return Boole;\-\\
end function; // IsAdmissible
\end{tabbing}
}
}
\noindent
\textbf{Output:}
\texttt{true}, if there is an action of \texttt{Critical} on the Frattini quotient of \texttt{G}.
\end{algorithm}


\section{Cyclic Cubic Fields of Type (2,2)}
\label{s:CycCub2x2}

We are now in the position to give a new proof of Derhem's result
on cyclic cubic fields \(F\) with \(\mathrm{Cl}_2{F}\simeq (2,2)\),
and to add a trivial supplement on the case \(\mathrm{Cl}_2{F}\simeq (2,2,2)\).

\begin{theorem}
\label{thm:Derhem}
(\textbf{Derhem, 1988},
\cite{Dh1988}.) \\
For a cyclic cubic field \(F\)
with \(2\)-class group \(\mathrm{Cl}_2{F}\simeq (2,2)\),
there are only two possibilities for the second \(2\)-class group \(\mathrm{G}_2^{(2)}{F}\),
either the \textbf{abelian group} \(\langle 4,2\rangle\)
with \(\varkappa(F)=(000)\), \(\tau(F)=\lbrack (1)^3\rbrack\), \(\ell_2{F}=1\)
or the \textbf{quaternion group} \(\langle 8,4\rangle\)
with \(\varkappa(F)=(123)\), \(\tau(F)=\lbrack (2)^3\rbrack\), \(\ell_2{F}=2\).

A cyclic cubic field \(F\) can never have a \(2\)-class group \(\mathrm{Cl}_2{F}\simeq (2,2,2)\).
\end{theorem}

\begin{proof}
These are immediate consequences of items 1. and 2. in Theorem
\ref{thm:SigmaGroupsDegree3}.
\end{proof}


\noindent
Table
\ref{tbl:CycCub2x2}
shows that cyclic cubic fields \(F\)
with elementary \(2\)-class group \(\mathrm{Cl}_2{F}\) of type \((2,2)\)
arise randomly without regard to the prime decomposition of the conductor \(c\) of \(F/\mathbb{Q}\).
There even occur cases where several, \(n=2\), members of a multiplet are simultaneously of type \((2,2)\).


\renewcommand{\arraystretch}{1.1}

\begin{table}[ht]
\caption{Various Multiplets of Cyclic Cubic Fields of Type \((2,2)\)}
\label{tbl:CycCub2x2}
\begin{center}
{\scriptsize
\begin{tabular}{|crccc|cc|}
\hline
 \(v_3{c}\) &      \(c\) &        Factors & \(t\) & \(m\) & \(n\) & \(\mathrm{G}_2^{(2)}{F_{c,\mu}}\) \\
\hline
      \(0\) &    \(163\) &          prime & \(1\) & \(1\) & \(1\) &            \(\langle 4,2\rangle\) \\
      \(0\) &    \(277\) &          prime & \(1\) & \(1\) & \(1\) &            \(\langle 8,4\rangle\) \\
\hline
      \(0\) &    \(679\) &       \(7,97\) & \(2\) & \(2\) & \(1\) &            \(\langle 4,2\rangle\) \\
      \(0\) &    \(703\) &      \(19,37\) & \(2\) & \(2\) & \(1\) &            \(\langle 8,4\rangle\) \\
      \(2\) & \(1\,413\) &    \(3^2,157\) & \(2\) & \(2\) & \(1\) &            \(\langle 4,2\rangle\) \\
      \(2\) &    \(711\) &     \(3^2,79\) & \(2\) & \(2\) & \(1\) &            \(\langle 8,4\rangle\) \\
      \(2\) & \(2\,169\) &    \(3^2,241\) & \(2\) & \(2\) & \(2\) & \(\langle 4,2\rangle,\langle 4,2\rangle\) \\
      \(0\) & \(6\,349\) &      \(7,907\) & \(2\) & \(2\) & \(2\) & \(\langle 4,2\rangle,\langle 4,2\rangle\) \\
      \(2\) & \(1\,899\) &    \(3^2,211\) & \(2\) & \(2\) & \(2\) & \(\langle 8,4\rangle,\langle 4,2\rangle\) \\
      \(0\) & \(3\,667\) &     \(19,193\) & \(2\) & \(2\) & \(2\) & \(\langle 8,4\rangle,\langle 4,2\rangle\) \\
\hline
      \(2\) & \(1\,197\) &   \(3^2,7,19\) & \(3\) & \(4\) & \(1\) &            \(\langle 4,2\rangle\) \\
      \(0\) & \(3\,913\) &    \(7,13,43\) & \(3\) & \(4\) & \(1\) &            \(\langle 4,2\rangle\) \\
      \(2\) & \(6\,489\) &  \(3^2,7,103\) & \(3\) & \(4\) & \(2\) & \(\langle 4,2\rangle,\langle 4,2\rangle\) \\
      \(0\) & \(6\,643\) &    \(7,13,73\) & \(3\) & \(4\) & \(2\) & \(\langle 4,2\rangle,\langle 4,2\rangle\) \\
\hline
\end{tabular}
}
\end{center}
\end{table}


\section{Cyclic Cubic Fields of Type (5,5)}
\label{s:CycCub5x5}


\begin{example}
\label{exm:CycCub5x5}
In the range \(1<c<1\,000\,000\) of conductors,
there are \(481\) occurrences of \(\mathrm{Cl}_5{F}\simeq (5,5)\).
The dominating part of \(463\) fields (\(96\%\)) has \(\mathrm{G}_5^{(2)}{F}\simeq\langle 125,3\rangle\)
with six total transfer kernels \(\varkappa(F)=(000000)\).
The leading, respectively trailing, example of the dominant part is 
\[c=6\,901=67\cdot 103, \text{ respectively } c=96\,733=7\cdot 13\cdot 1063.\]
Exceptions occur for the following \(18\) conductors only, confirming Corollary
\ref{cor:CycCub5x5}.
\end{example}

\renewcommand{\arraystretch}{1.1}

\begin{table}[ht]
\caption{Exceptional Singlets of Cyclic Cubic Fields of Type \((5,5)\)}
\label{tbl:CycCub5x5}
\begin{center}
{\scriptsize
\begin{tabular}{|r|rccc|c|l|c|}
\hline
    No. &        \(c\) &    Factors    & \(t\) & \(m\) & \(n\) & \(\mathrm{G}_5^{(2)}{F}\)  & \(\varkappa(F)\) \\
\hline
  \(1\) &  \(66\,313\) &   \(13,5101\) & \(2\) & \(2\) & \(1\) & \(\langle 3125,12\rangle\) & a \(6\)-cycle    \\
  \(2\) &  \(68\,791\) &         prime & \(1\) & \(1\) & \(1\) & \(\langle 3125,14\rangle\) & the identity     \\
  \(3\) &  \(77\,971\) &   \(103,757\) & \(2\) & \(2\) & \(1\) & \(\langle 3125,14\rangle\) & the identity     \\
  \(4\) &  \(87\,409\) &   \(7,12487\) & \(2\) & \(2\) & \(1\) & \(\langle 3125,12\rangle\) & a \(6\)-cycle    \\
  \(5\) & \(199\,621\) &         prime & \(1\) & \(1\) & \(1\) & \(\langle 3125,9\rangle\)  & two \(3\)-cycles \\
  \(6\) & \(317\,853\) &   \(9,35317\) & \(2\) & \(2\) & \(1\) & \(\langle 3125,12\rangle\) & a \(6\)-cycle    \\
  \(7\) & \(425\,257\) &  \(7,79,769\) & \(3\) & \(4\) & \(1\) & \(\langle 3125,14\rangle\) & the identity     \\
  \(8\) & \(464\,191\) & \(7,13,5101\) & \(3\) & \(4\) & \(1\) & \(\langle 3125,12\rangle\) & a \(6\)-cycle    \\
  \(9\) & \(481\,537\) &   \(7,68791\) & \(2\) & \(2\) & \(1\) & \(\langle 3125,14\rangle\) & the identity     \\
 \(10\) & \(545\,797\) & \(7,103,757\) & \(3\) & \(4\) & \(1\) & \(\langle 3125,14\rangle\) & the identity     \\
 \(11\) & \(596\,817\) & \(9,13,5101\) & \(3\) & \(4\) & \(1\) & \(\langle 3125,12\rangle\) & a \(6\)-cycle    \\
 \(12\) & \(619\,119\) &   \(9,68791\) & \(2\) & \(2\) & \(1\) & \(\langle 3125,14\rangle\) & the identity     \\
 \(13\) & \(678\,303\) &   \(9,75367\) & \(2\) & \(2\) & \(1\) & \(\langle 3125,12\rangle\) & a \(6\)-cycle    \\
 \(14\) & \(701\,739\) & \(9,103,757\) & \(3\) & \(4\) & \(1\) & \(\langle 3125,14\rangle\) & the identity     \\
 \(15\) & \(767\,623\) &         prime & \(1\) & \(1\) & \(1\) & \(\langle 3125,9\rangle\)  & two \(3\)-cycles \\
 \(16\) & \(786\,681\) & \(7,9,12487\) & \(3\) & \(4\) & \(1\) & \(\langle 3125,12\rangle\) & a \(6\)-cycle    \\
 \(17\) & \(894\,283\) &  \(13,68791\) & \(2\) & \(2\) & \(1\) & \(\langle 3125,14\rangle\) & the identity     \\
 \(18\) & \(909\,229\) &  \(487,1867\) & \(2\) & \(2\) & \(1\) & \(\langle 3125,14\rangle\) & the identity     \\
\hline
\end{tabular}
}
\end{center}
\end{table}

\bigskip
\noindent
All these \(481\) cyclic cubic fields have \(5\)-class towers of length \(\ell_5{F}=2\).
For the group \(\mathrm{G}_5^{(2)}{F}\simeq\langle 125,3\rangle\) of coclass \(1\) with cyclic commutator subgroup of order \(5\),
this follows from a theorem by Blackburn
\cite{Bl1957}.
In the exceptional cases with groups \(\mathrm{G}_5^{(2)}{F}\simeq\langle 3125,i\rangle\), \(i\in\lbrace 9,12,14\rbrace\),
it is a consequence of the fact that a metabelian Schur \(\sigma\)-group
cannot be the second derived quotient of a non-metabelian \(5\)-group.
Therefore, we always have \(\mathrm{G}_5^{(2)}{F}=\mathrm{G}_5^{(\infty)}{F}\) and \(\ell_5{F}=2\).

The \textit{fixed point capitulation problem} by \textbf{Olga Taussky} \(\mathbf{1970}\)
\cite[Rem. 1, p. 438]{Ta1970},
which has first been solved in 2011 by five imaginary quadratic fields
\cite[\S\ 3.5.2, p. 448]{Ma2013},
and later in 2015 by eleven cyclic quartic fields
\cite[Thm. 4.4, Tbl. 4--5]{AKMTT2020},
now also has \textit{nine solutions with cyclic cubic fields}
possessing the identity permutation type \(\varkappa=(123456)\)
and the group \(\mathrm{G}_5^{(2)}{F}\simeq\langle 3125,14\rangle\).

The numerical results are in perfect accordance with Theorem
\ref{thm:SigmaGroupsDegree3}.
However, it should be pointed out that
the abelian group \(\langle 25,2\rangle\) with  \(\ell_5{F}=1\)
and the unbalanced groups  \(\langle 3125,i\rangle\), \(i\in\lbrace 3,10\rbrace\),
did not occur as \(\mathrm{G}_5^{(2)}{F}\), up to now.


\begin{theorem}
\label{thm:CycCub5x5}
(\(\sigma\)-groups of degree \(3\) with type \((5,5)\)) \\
The transfer kernel type of a pro-\(5\) group \(\mathfrak{G}\)
with abelianization \(\mathfrak{G}/\mathfrak{G}^\prime\simeq (5,5)\)
which is a \(\sigma\)-group of degree \(3\)
is restricted to the following admissible types \\
\(\bullet\) \((000000)\), \\
\(\bullet\) two \(3\)-cycles, \\
\(\bullet\) a \(6\)-cycle, \\
\(\bullet\) the identity, \\
\(\bullet\) three \(2\)-cycles, \\
and descendant types of \((000000)\).
\end{theorem}

\begin{proof}
These are immediate consequences of item 3. in Theorem
\ref{thm:SigmaGroupsDegree3}.
\end{proof}


\begin{corollary}
\label{cor:CycCub5x5}
(Cyclic cubic fields of type \((5,5)\)) \\
The second \(5\)-class group \(\mathrm{G}_5^{(2)}{F}\)
of a cyclic cubic field \(F\)
with \(5\)-class group \(\mathrm{Cl}_5{F}\simeq (5,5)\)
is restricted to the isomorphism classes of the following groups\\
\(\bullet\) \(\langle 25,2\rangle\), \\
\(\bullet\) \(\langle 125,3\rangle\), \\
\(\bullet\) a descendant of \(\langle 3125,3\rangle\), \\
\(\bullet\) \(\langle 3125,9\rangle\), \\
\(\bullet\) \(\langle 3125,12\rangle\), \\
\(\bullet\) \(\langle 3125,14\rangle\), \\
\(\bullet\) a descendant of \(\langle 3125,10\rangle\). \\
All these finite \(5\)-groups are \(\sigma\)-groups of degree \(3\).
\end{corollary}


\chapter{Closed Andozhskii Groups}
\label{ch:Andozhskii}

\noindent
According to Koch and Venkov
\cite{KoVe1975},
\textit{Schur \(\sigma\)-groups} \(S\) are known to be
mandatory for realizations
\(S\simeq\mathrm{Gal}(\mathrm{F}_p^\infty{k}/k)\)
by \(p\)-class field towers of \textit{imaginary} quadratic fields \(k\),
with an odd prime \(p\).
They possess a balanced presentation \(d(S)=r(S)\) with coinciding
generator rank  \(d(S)=\dim_{\mathbb{F}_p}H^1(S,\mathbb{F}_p)\)
and relation rank \(r(S)=\dim_{\mathbb{F}_p}H^2(S,\mathbb{F}_p)\),
and an automorphism \(\sigma\in\mathrm{Aut}(S)\)
acting as inversion \(x\mapsto x^{-1}\)
on the commutator quotient \(S/\lbrack S,S\rbrack\).
However, in the older literature,
for instance Shafarevich
\cite[\S\ 6, pp. 88--91]{Sh1964},
there also appear Schur groups with balanced presentation,
but without a generator inverting \(\sigma\)-automorphism,
and they are called \textit{closed},
according to the original terminology by Schur.
In the present chapter,
we are interested in finite closed \(3\)-groups
given by Andozhskii
\cite{AnTs1974,An1975},
whose position in the descendant tree of 
\(3\)-groups \(G\) with elementary tricyclic
commutator quotient \(G/\lbrack G,G\rbrack\)
is illuminated in Figure
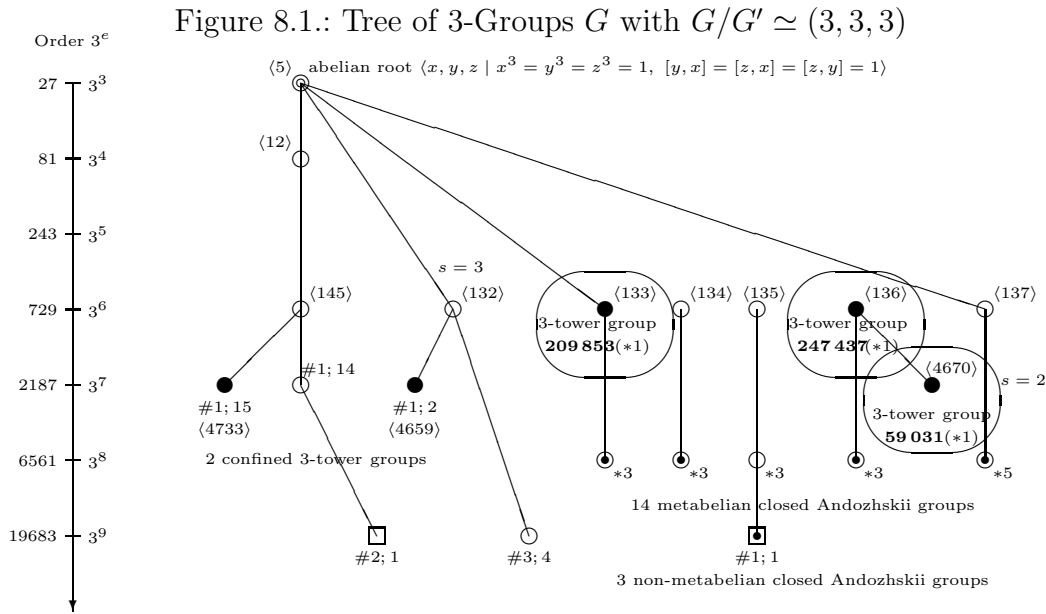
\ref{fig:Groups333}.


\begin{figure}[ht]
\caption{Tree of \(3\)-Groups \(G\) with \(G/G^\prime\simeq (3,3,3)\)}
\label{fig:Groups333}

{\tiny

\setlength{\unitlength}{1.0cm}
\begin{picture}(10,7.5)(-10.2,-7)


\put(-9,0.5){\makebox(0,0)[cb]{Order \(3^e\)}}
\put(-9,0){\line(0,-1){7}}
\multiput(-9.1,0)(0,-1){7}{\line(1,0){0.2}}
\put(-9.2,0){\makebox(0,0)[rc]{\(27\)}}
\put(-8.8,0){\makebox(0,0)[lc]{\(3^3\)}}
\put(-9.2,-1){\makebox(0,0)[rc]{\(81\)}}
\put(-8.8,-1){\makebox(0,0)[lc]{\(3^4\)}}
\put(-9.2,-2){\makebox(0,0)[rc]{\(243\)}}
\put(-8.8,-2){\makebox(0,0)[lc]{\(3^5\)}}
\put(-9.2,-3){\makebox(0,0)[rc]{\(729\)}}
\put(-8.8,-3){\makebox(0,0)[lc]{\(3^6\)}}
\put(-9.2,-4){\makebox(0,0)[rc]{\(2187\)}}
\put(-8.8,-4){\makebox(0,0)[lc]{\(3^7\)}}
\put(-9.2,-5){\makebox(0,0)[rc]{\(6561\)}}
\put(-8.8,-5){\makebox(0,0)[lc]{\(3^8\)}}
\put(-9.2,-6){\makebox(0,0)[rc]{\(19683\)}}
\put(-8.8,-6){\makebox(0,0)[lc]{\(3^9\)}}
\put(-9,-6){\vector(0,-1){1}}


\put(-6.1,0.1){\makebox(0,0)[rb]{\(\langle 5\rangle\)}}
\put(-5.9,0.1){\makebox(0,0)[lb]{abelian root \(\langle x,y,z\mid
x^3=y^3=z^3=1,\ \lbrack y,x\rbrack=\lbrack z,x\rbrack=\lbrack z,y\rbrack=1\rangle\)}}
\put(-6,0){\circle{0.2}}
\put(-6,0){\circle{0.1}}

\put(-6,0){\line(0,-1){1}}
\put(-6.1,-0.9){\makebox(0,0)[rb]{\(\langle 12\rangle\)}}
\put(-6,-1){\circle{0.2}}
\put(-6,-1){\line(0,-1){2}}
\put(-5.9,-2.9){\makebox(0,0)[lb]{\(\langle 145\rangle\)}}
\put(-6,-3){\circle{0.2}}
\put(-6,-3){\line(-1,-1){1}}
\put(-7,-4){\circle*{0.2}}
\put(-7,-4.3){\makebox(0,0)[cc]{\(\#1;15\)}}
\put(-7,-4.6){\makebox(0,0)[cc]{\(\langle 4733\rangle\)}}
\put(-6,-3){\line(0,-1){1}}
\put(-6,-3.9){\makebox(0,0)[lb]{\(\#1;14\)}}
\put(-6,-4){\circle{0.2}}
\put(-6,-4){\line(1,-2){1}}
\put(-5.1,-6.1){\framebox(0.2,0.2){}}
\put(-5,-6.3){\makebox(0,0)[cc]{\(\#2;1\)}}


\put(-5.8,-5){\makebox(0,0)[cc]{\(2\) confined \(3\)-tower groups}}

\put(-6,0){\line(2,-3){2}}
\put(-4.2,-2.5){\makebox(0,0)[lb]{\(s=3\)}}
\put(-3.9,-2.9){\makebox(0,0)[lb]{\(\langle 132\rangle\)}}
\put(-4,-3){\circle{0.2}}
\put(-4,-3){\line(-1,-2){0.5}}
\put(-4.5,-4){\circle*{0.2}}
\put(-4.5,-4.3){\makebox(0,0)[cc]{\(\#1;2\)}}
\put(-4.5,-4.6){\makebox(0,0)[cc]{\(\langle 4659\rangle\)}}
\put(-4,-3){\line(1,-3){1}}
\put(-3,-6){\circle{0.2}}
\put(-3,-6.3){\makebox(0,0)[cc]{\(\#3;4\)}}


\put(-6,0){\line(4,-3){4}}
\put(-1.9,-2.9){\makebox(0,0)[lb]{\(\langle 133\rangle\)}}
\put(-2,-3){\circle*{0.2}}
\put(-2,-3){\line(0,-1){2}}
\put(-2,-5){\circle{0.2}}
\put(-2,-5){\circle*{0.1}}
\put(-1.9,-5.1){\makebox(0,0)[lt]{\(\ast 3\)}}

\multiput(-2,-3.2)(3.3,0){2}{\oval(1.8,1.4)}

\put(-2.1,-3.2){\makebox(0,0)[cc]{\(3\)-tower group}}
\put(-2.1,-3.5){\makebox(0,0)[cc]{\(\mathbf{209\,853}(\ast 1)\)}}

\put(-0.9,-2.9){\makebox(0,0)[lb]{\(\langle 134\rangle\)}}
\put(-1,-3){\circle{0.2}}
\put(-1,-3){\line(0,-1){2}}
\put(-1,-5){\circle{0.2}}
\put(-1,-5){\circle*{0.1}}
\put(-0.9,-5.1){\makebox(0,0)[lt]{\(\ast 3\)}}

\put(-0.2,-2.9){\makebox(0,0)[lb]{\(\langle 135\rangle\)}}
\put(0,-3){\circle{0.2}}
\put(0,-3){\line(0,-1){2}}
\put(0.6,-5.6){\makebox(0,0)[cc]{\(14\) metabelian closed Andozhskii groups}}
\put(0,-5){\circle{0.2}}
\put(0.1,-5.1){\makebox(0,0)[lt]{\(\ast 3\)}}
\put(0,-5){\line(0,-1){1}}
\put(-0.1,-6.1){\framebox(0.2,0.2){}}
\put(0,-6){\circle*{0.1}}
\put(0,-6.3){\makebox(0,0)[cc]{\(\#1;1\)}}
\put(0.6,-6.6){\makebox(0,0)[cc]{\(3\) non-metabelian closed Andozhskii groups}}

\put(1.4,-2.9){\makebox(0,0)[lb]{\(\langle 136\rangle\)}}
\put(1.3,-3){\circle*{0.2}}
\put(1.3,-3){\line(0,-1){2}}
\put(1.3,-5){\circle{0.2}}
\put(1.3,-5){\circle*{0.1}}
\put(1.4,-5.1){\makebox(0,0)[lt]{\(\ast 3\)}}

\put(1.2,-3.2){\makebox(0,0)[cc]{\(3\)-tower group}}
\put(1.2,-3.5){\makebox(0,0)[cc]{\(\mathbf{247\,437}\!(\ast 1)\)}}

\put(1.3,-3){\line(1,-1){1}}
\put(2.3,-4){\circle*{0.2}}
\put(2.2,-3.9){\makebox(0,0)[lb]{\(\langle 4670\rangle\)}}

\put(2.3,-4.2){\oval(1.8,1.4)}

\put(2.3,-4.4){\makebox(0,0)[cc]{\(3\)-tower group}}
\put(2.3,-4.7){\makebox(0,0)[cc]{\(\mathbf{59\,031}(\ast 1)\)}}

\put(-6,0){\line(3,-1){9}}
\put(3.1,-2.9){\makebox(0,0)[lb]{\(\langle 137\rangle\)}}
\put(3,-3){\circle{0.2}}
\put(3,-3){\line(0,-1){2}}
\put(3.2,-4){\makebox(0,0)[lb]{\(s=2\)}}
\put(3,-5){\circle{0.2}}
\put(3,-5){\circle*{0.1}}
\put(3.1,-5.1){\makebox(0,0)[lt]{\(\ast 5\)}}


\end{picture}
}
\end{figure}


In \S\
\ref{s:AndozhskiiIdentified},
we identify  the \(17\) \textit{closed Andozhskii groups}. In \S\
\ref{s:RealizationAndozhskii},
we show that two of them can be realized as Galois groups
of \(3\)-class field towers of \textit{cyclic cubic fields}.
Incidentally,
we point out that related but not closed \(3\)-groups
(see \(j=2\) in Table
\ref{tbl:HBC})
are realized by lots of \(3\)-class field towers over 
\textit{totally complex \(S_3\)-fields} \(K\),
which are unramified extensions of imaginary quadratic fields
\(k=\mathrm{Q}(\sqrt{d})\)
with \(3\)-class group \(\mathrm{Cl}_3{k}\simeq (3,3)\),
capitulation type \(\mathrm{H}.4\), \(\varkappa(k)\sim (4111)\),
and three abelian type invariants \(\tau(k)\sim (111,111,111,21)\) of rank \(3\).


\section{Identification of Closed Andozhskii Groups}
\label{s:AndozhskiiIdentified}

\noindent
\begin{theorem}
\label{thm:Andozhskii}
Among the finite \(3\)-groups \(G\)
with commutator quotient \(G/G^\prime\simeq (3,3,3)\),
there exist precisely \(14\) metabelian closed
groups \(S\) of order \(\#(S)=3^8\)
with identifiers
\begin{equation}
\label{eqn:Metab}
S\simeq\langle 6561,217700+i\rangle \text{ where } 1\le i\le 6 \text{ or } 10\le i\le 17,
\end{equation}
and \(3\) non-metabelian closed
groups \(S\) of order \(\#(S)=3^9\)
with identifiers
\begin{equation}
\label{eqn:NonMetab}
S\simeq\langle 6561,217700+i\rangle-\#1;1 \text{ where } 7\le i\le 9.
\end{equation}
They possess a trivial Schur multiplier \(M(S)=H_2(S,\mathbb{Q}/\mathbb{Z})=0\)
and a balanced presentation \(d(S)=r(S)\)
with coinciding generator rank
\(d(S)=\dim_{\mathbb{F}_p}H^1(S,\mathbb{F}_p)\)
and relation rank
\(r(S)=\dim_{\mathbb{F}_p}H^2(S,\mathbb{F}_p)\).
The class is \(\mathrm{cl}(S)=3\) for soluble length \(\mathrm{sl}(S)=2\)
and \(\mathrm{cl}(S)=4\) for \(\mathrm{sl}(S)=3\).
\end{theorem}

\begin{proof}
By a search in the SmallGroups database
\cite{BEO2005},
extended to order \(3^9\) by the \(p\)-group generation algorithm
\cite{Nm1977,Ob1990},
the finite closed Andozhskii \(3\)-groups \(S\)
are identified.
There are \(14\) hits of order \(\#(S)=3^8\)
and only three hits of order \(\#(S)=3^9\).
The non-metabelian groups are characterized by
their identifiers defined by the ANUPQ package
\cite{GNO2006}.
\end{proof}


\begin{corollary}
\label{cor:Andozhskii}
Each of the \(17\) closed Andozhskii groups \(S\) in Theorem
\ref{thm:Andozhskii}
shares a common Artin pattern \((\varkappa,\alpha)\)
with its unique ancestor \(A\simeq\langle 729,130+j\rangle\),
as shown in Table
\ref{tbl:HBC},
where
\begin{equation}
\label{eqn:Connection}
j=
\begin{cases}
3 & \text{ for } 1\le i \le 3, \\
4 & \text{ for } 4\le i \le 6, \\
5 & \text{ for } 7\le i \le 9, \\
6 & \text{ for } 10\le i \le 12, \\
7 & \text{ for } 13\le i \le 17. \\
\end{cases}
\end{equation}
\end{corollary}

\begin{proof}
According to the theorem on the antitony of the Artin pattern
\cite[\S\S\ 5.1--5.4, pp. 78--87]{Ma2016},
it suffices to calculate the \textit{stable} transfer kernels
of the five ancestors \(A\) of the \(17\) closed groups in Theorem
\ref{thm:Andozhskii}.
They are of order \(\#(A)=3^6\)
and have much simpler presentations.
It turns out that the transfer kernels are \textit{harmonically balanced},
that is, permutations in the symmetric group \(S_{13}\).
\end{proof}

\renewcommand{\arraystretch}{1.1}

\begin{table}[ht]
\caption{TKT \(\varkappa\) and AQI \(\tau\) of Closed Andozhskii \(3\)-Groups \(S\), \(S/S^{\prime}\simeq (3,3,3)\)}
\label{tbl:HBC}
\begin{center}
{\scriptsize
\begin{tabular}{|c||c|c|c|c|c|c|c|c|c|c|c|c|c||c|}
\hline
 \(j\) & \multicolumn{13}{|c||}{\(\varkappa\) respectively \(\tau\)} & Action \\
\hline
 \(2\) &  1 &  2 &  3 &   6 & 11 &  9 & 10 &   4 & 13 &   5 & 12 &  8 &   7 & \(\langle 24,12\rangle\) \\
       & 22 & 22 & 22 & 211 & 22 & 22 & 22 & 211 & 22 & 211 & 22 & 22 & 211 & \\
\hline
 \(3\) &  9 &  2 &  3 &   6 & 10 &   8 &   4 &  11 & 12 &  13 &   5 &  1 &   7 & \(\langle 6,2\rangle\) \\
       & 22 & 22 & 22 & 211 & 22 & 211 & 211 & 211 & 22 & 211 & 211 & 22 & 211 & \\
\hline
 \(4\) &  1 &  7 &  3 &  2 & 10 &  9 &   4 &  11 &   5 &  12 & 13 &  8 &  6 & \(\langle 4,1\rangle\) \\
       & 22 & 22 & 22 & 22 & 22 & 22 & 211 & 211 & 211 & 211 & 22 & 22 & 22 & \\
\hline
 \(5\) &  9 &  7 &  3 &  2 &   4 &   8 & 11 &  10 & 13 &   5 & 12 &  1 &  6 & \(\langle 3,1\rangle\) \\
       & 22 & 22 & 22 & 22 & 211 & 211 & 22 & 211 & 22 & 211 & 22 & 22 & 22 & \\
\hline
 \(6\) & 12 &  7 &  3 &  2 &  9 &  5 &  8 &  1 & 10 & 11 &  4 &  13 &  6 & \(\langle 6,2\rangle\) \\
       & 22 & 22 & 22 & 22 & 22 & 22 & 22 & 22 & 22 & 22 & 22 & 211 & 22 & \\
\hline
 \(7\) & 10 &  7 &  3 &   6 &   8 &  4 &  1 &  9 & 13 &   5 & 12 &  11 &  2 & \(\langle 24,3\rangle\) \\
       & 22 & 22 & 22 & 211 & 211 & 22 & 22 & 22 & 22 & 211 & 22 & 211 & 22 & \\
\hline
\end{tabular}
}
\end{center}
\end{table}


\noindent
We recall that the transfer kernel (capitulation kernel) \(\ker(T_i)\)
of an Artin transfer homomorphism \(T_i:\,G/G^\prime\to M_i/M_i^\prime\)
\cite{Ma2016}
from a \(3\)-group \(G\) with \(G/G^\prime\simeq (3,3,3)\)
to one of its \(13\) maximal subgroups \(M_i\), \(1\le i\le 13\),
is called of \textit{Taussky type} A, if the meet
\(\ker(T_i)\cap M_i>1\) is non-trivial,
and of \textit{Taussky type} B, if \(\ker(T_i)\cap M_i=1\).

\begin{corollary}
\label{cor:Taussky}
The Artin patterns \((\varkappa(A),\tau(A))\) of the six groups
\(A=\langle 729,130+j\rangle\) with \(2\le j\le 7\)
share the common property that the Taussky type of their transfer kernels \(\ker(T_i)\)
is determined uniquely by the AQI \(M_i/M_i^\prime\) of the corresponding maximal subgroup \(M_i\):
\begin{equation}
\label{eqn:Taussky}
\begin{aligned}
\tau(A)_i=M_i/M_i^\prime\simeq (211) & \Longleftrightarrow \varkappa(A)_i=\ker(T_i)\cap M_i>1, \text{ Taussky type }\mathrm{A}, \\
\tau(A)_i=M_i/M_i^\prime\simeq (22)  & \Longleftrightarrow \varkappa(A)_i=\ker(T_i)\cap M_i=1, \text{ Taussky type }\mathrm{B},
\end{aligned}
\end{equation}
for all \(1\le i\le 13\). Here the abelian quotient invariants are logarithmic.
\end{corollary}

\begin{proof}
This follows by comparing the \(1\)-dimensional transfer kernels in Table
\ref{tbl:HBC}
to the planes in Table
\ref{tbl:Planes}.

\end{proof}


\section{Realization as 3-Class Field Tower Groups}
\label{s:RealizationAndozhskii}

\noindent
Since the groups in Theorem
\ref{thm:Andozhskii}
are non-\(\sigma\) groups,
they cannot be realized by any quadratic field,
neither imaginary nor real. Therefore,
we investigated the possible Galois actions (Table
\ref{tbl:HBC})
on the five ancestors \(A=\mathrm{SmallGroup}(729,130+j)\).
It turned out that the unique non-metabelian case
\(j=5\) can only be realized by \textit{cyclic cubic} fields,
\(j=4\) by cyclic quartic fields, and
\(j\in\lbrace 3,6,7\rbrace\) by \textit{cyclic cubic or sextic} fields.
We show that certain metabelian descendants \(S\) for \(j\in\lbrace 3,6\rbrace\)
can actually be realized as Galois groups
\(\mathrm{Gal}(\mathrm{F}_3^\infty{K}/K)\simeq S\)
of maximal unramified pro-\(3\) extensions
of cyclic cubic fields \(F\) with conductors
\(c=\mathbf{59\,031},\mathbf{209\,853}\),
\(\mathbf{247\,437},\mathbf{263\,017},\mathbf{329\,841},\mathbf{377\,923}\)
and \(3\)-class group \(\mathrm{Cl}_3{F}\simeq (3,3,3)\).


\begin{theorem}
\label{thm:Main1}
If a number field \(K/\mathbb{Q}\)
with elementary tricyclic \(3\)-class group \(\mathrm{Cl}_3{K}\simeq (3,3,3)\)
possesses the Artin pattern \((\varkappa(K),\tau(K))\)
with harmonically balanced capitulation type 
\(\varkappa(K)\sim (9,2,3,6,10,8,4,11,12,13,5,1,7)\)
and  abelian type invariants
\(\tau(K)\sim ((22)^{3},211,22,(211)^3,22,(211)^2,22,211)\),
then \(K/\mathbb{Q}\) must be cyclic cubic or sextic,
and has a metabelian \(3\)-class field tower with automorphism group
\(\mathrm{Gal}(\mathrm{F}_3^\infty{K}/K)\simeq\)
\begin{equation}
\label{eqn:Main1}
\langle 3^8,217700+i\rangle,\ 1\le i\le 3,
\text{ or } \langle 3^7,4660+k\rangle,\ 1\le k\le 4,
\text{ or } \langle 3^6,133\rangle.
\end{equation}
\end{theorem}


\begin{theorem}
\label{thm:Main2}
If a number field \(K/\mathbb{Q}\)
with elementary tricyclic \(3\)-class group \(\mathrm{Cl}_3{K}\simeq (3,3,3)\)
possesses the Artin pattern \((\varkappa(K),\tau(K))\)
with harmonically balanced capitulation type 
\(\varkappa(K)\sim (12,7,3,2,9,5,8,1,10,11,4,13,6)\)
and abelian type invariants
\(\tau(K)\sim ((22)^{11},211,22)\),
then \(K/\mathbb{Q}\) must be cyclic cubic or sextic,
and has a metabelian \(3\)-class field tower with automorphism group
\(\mathrm{Gal}(\mathrm{F}_3^\infty{K}/K)\simeq\)
\begin{equation}
\label{eqn:Main2}
\langle 3^8,217700+i\rangle,\ 10\le i\le 12,
\text{ or } \langle 3^7,4669+k\rangle,\ 0\le k\le 3,
\text{ or } \langle 3^6,136\rangle.
\end{equation}
\end{theorem}

\begin{proof}
Theorems
\ref{thm:Main1}
and
\ref{thm:Main2}
are immediate consequences of Table
\ref{tbl:HBC}.
\end{proof}


\noindent
In April 2002,
we used the Voronoi algorithm and the Euler product method
in order to compute the \(15851\) cyclic cubic fields \(F\)
with conductors \(c_{F/\mathbb{Q}}<10^5\)
and their class numbers \(1\le h_F\le 1953\).
Among the fields,
\(4785\) occur as singlets,
\(7726\) in doublets,
\(3132\) in quartets,
and \(208\) in octets.
Twenty years later, in July 2022,
we have confirmed these results,
extended by the class group structures,
as reported in detail in \S\
\ref{s:ClassGroups}.
The cyclic cubic fields \(F\) were constructed 
as ray class fields over the rational number field,
using the class field theoretic routines of Magma
\cite{MAGMA2022}.
Additionally, we constructed the \(13\)
unramified cyclic cubic relative extensions \(E_i/F\),
whenever the class group of \(F\) was \(\mathrm{Cl}_3{F}\simeq (3,3,3)\),
which was the primary goal for the reconstruction
in view of the intended realization of Andozhskii groups.

\begin{example}
\label{exm:Andozhskii}
Indeed, there are \(15851\) fields in \(9457\) multiplets:
\(4785\) singlets,
\(3863\) doublets,
\(783\) quartets,
and \(26\) octets,
as shown in Table
\ref{tbl:StatMult3}.
One of the four fields with \(c=\mathbf{59\,031}=3^2\cdot 7\cdot 937\)
and \(c\in\lbrace\mathbf{247\,437},\mathbf{329\,841},\mathbf{377\,923}\rbrace\),
respectively \(c=\mathbf{209\,853}=3^2\cdot 7\cdot 3331\) and \(c=\mathbf{263\,017}\),
i.e. quartets of Category I and Graph 2,
were hits
of the Artin pattern in the Main Theorem
\ref{thm:Main2}, respectively
\ref{thm:Main1}.
The density of such fields is horrificly sparse.
The conductors \(c>10^5\) are due to an extended systematic investigation
of harmonically balanced transfer kernels \(\varkappa\).
We point out that \(c=\mathbf{271\,737}\) leads to some of the three groups
\(\langle 729,130+j\rangle\) with \(j\in\lbrace 2,5,7\rbrace\),
and \(c=\mathbf{59\,031}\) is the single hit of a group with order \(2187\).
\end{example}


\part{Future Research}
\label{pt:Future}

\noindent
Among the \textit{remaining challenges} for future activities of the young researchers
we mention the following open problems.

\begin{itemize}
\item
Theorems
\ref{thm:Cat3Gph679} --
\ref{thm:Cat1Gph12}
have been obtained by thorough observations and analysis
of computational results.
It would be desirable to prove these theorems
by exploiting results in the later chapters of Ayadi's Thesis
\cite{Ay1995}
on principal factors (Parry invariants) and
bicyclic bicubic subfields of the \(3\)-genus field.
\item
An interesting supplement would be the determination
of the \(3\)-class field tower group of cyclic cubic fields
in quartets of Category IV, Graphs 1--3, and in octets with \(t=4\).
\item
It should still be within the reach of Magma
to investigate the \(6\), respectively \(31\),
unramified cyclic quintic extensions \(E/F\)
of \textbf{cyclic quintic fields} \(F\) with elementary \(5\)-class group of
bicyclic type \(\mathrm{Cl}_5{F}\simeq (5,5)\), respectively
tricyclic type \(\mathrm{Cl}_5{F}\simeq (5,5,5)\),
and to establish theoretical foundations
for the capitulation \(\ker(T_{E/F})\)
and \(5\)-class field tower \(F_5^\infty\).
\end{itemize}


\backmatter

\chapter{Conclusion}
\label{ch:Conclusion}


\section{Construction Process and Statistics of Cyclic Cubic Fields}
\label{s:Construction}


\subsection{Computational Techniques}
\label{ss:Techniques}

\noindent
In \(2002\), we computed the regulators \(\mathrm{Reg}(F)\) and class numbers \(h(F)\) 
of all \(15851\) cyclic cubic fields \(F\) with conductor \(c<10^5\),
based on generating polynomials in \(\mathbb{Z}\lbrack X\rbrack\) by Marie Nicole Gras
\cite[p. 90]{MNGr1975},
that is,
\(P(X)=X^3+X^2+\frac{1-c}{3}X-\frac{c(3+a)-1}{27}\) for \(\gcd(3,c)=1\),
simplified by a Tschirnhausen transformation, respectively
\(P(X)=X^3-\frac{c}{3}X-\frac{ac}{27}\) for \(9\mid c\),
where \(c=\frac{a^2+27b^2}{4}\).
For the \textit{regulator}, we used the two-dimensional Voronoi algorithm
\cite{Vo1896},
implemented in Delphi (i.e. object oriented Pascal),
for the \textit{class number}, we used the analytic Euler product formula.

In \(2013\), we computed a list of 
the first \(251\) conductors \(c\) of cyclic cubic number fields \(F=F_c\), 
having a \(3\)-class group \(\mathrm{Cl}_3{F}\) of type \((3,3)\). 
These conductors occur in the range \(657 \le c\le 26\,523\).
They are divisible by two or three primes, \(2\le t\le 3\).
Inspired by Amandine Leriche and Henri Cohen
\cite[pp. 337--339]{Co2000}, 
we used a very convenient parametrization of generating polynomials 
\(P(X) = X^3 - 3cX - cu\in\mathbb{Z}\lbrack X\rbrack\)
for the cyclic cubic fields \(F\).
These results were the basis of preliminary statistics in \S\
\ref{ss:Statistics}.

In \(2016\), we extended the list of conductors up to \(c\le 10^5\)
without strict book keeping of exact numbers.
The construction process was completely different:
Instead of using generating polynomials,
we employed the class field package of Magma
\cite{BCP1997,BCFS2022,MAGMA2022},
as implemented by Claus Fieker
\cite{Fi2001},
and obtained cyclic cubic fields \(F=F_c\) as subfields of \textit{ray class fields} \(\mathbb{Q}_c\)
modulo admissible conductors \(c\) over the rational number field \(\mathbb{Q}\).
Then, we constructed the \textit{unramified cyclic cubic extensions} \(E/F\) and
determined the \(3\)-capitulation kernel of \(F\) in each superfield \(E\) 
(the transfer kernel type, TKT, \(\varkappa\)),
and the structure of the \(3\)-class group \(\mathrm{Cl}_3{E}\) of every extension \(E\) 
(the transfer target type, TTT, \(\tau\)).
Finally, we applied our new algorithm,
the strategy of \textit{pattern recognition} via Artin transfers
\cite{Ma2020},
in order to determine the second \(3\)-class group \(\mathfrak{M}=\mathrm{G}_3^2{F}\) 
of the number field \(F\) 
by seeking the \textit{Artin pattern} \(\mathrm{AP}=(\tau,\varkappa)\),
which consists of targets \(\tau\) and kernels \(\varkappa\) of transfer homomorphisms,
in the descendant tree \(\mathcal{T}(R)\) with root \(R=\langle 9,2\rangle\).
The tree was constructed by successive applications of
the \(p\)-group generation algorithm by Mike F. Newman
\cite{Nm1977}
and Eamonn A. O'Brien
\cite{Ob1990},
which is also described in
\cite{HEO2005}.

In \(2022\), we repeated and double checked the computations of \(2002\) and \(2016\)
with modern techniques, very thorough book keeping and extensive statistics.
Some metabelian \(3\)-groups in Table
\ref{tbl:Metabelian33}
could only be identified unambiguously by means of the commutator subgroup.
This required the computation of the \(3\)-class group \(\mathrm{Cl}_3{F_3^{(1)}}\)
of the first Hilbert \(3\)-class field of many cyclic cubic fields \(F\),
which is of absolute degree \(\lbrack F_3^{(1)}:\mathbb{Q}\rbrack=3\cdot 9=27\).
In this manner, we had to distinguish between
\(\mathfrak{M}=\langle 729,34|35|36\rangle\) with \(\mathfrak{M}^\prime=(3,3,3,3)\) and
\(\mathfrak{M}=\langle 729,37|38|39\rangle\) with \(\mathfrak{M}^\prime=(3,3,9)\).
On the other hand, for extensive classes of cyclic cubic fields \(F\),
initially we did not know whether
\(\mathfrak{M}=\langle 81,9\rangle\),
which was known to arise for certain bicyclic biquadratic fields
\cite{ATTDM2016},
or \(\mathfrak{M}=\langle 243,28|29|30\rangle\),
known for real quadratic fields
\cite{Ma2012a}.
Since \(\mathrm{Cl}_3{F_3^{(1)}}\) turned out to be \((3,9)\),
rather than \((3,3)\),
the group \(\langle 81,9\rangle\) could be eliminated (cannot occur for cyclic cubic fields). 


\subsection{Statistics}
\label{ss:Statistics}

\noindent
For the calculation of relative frequencies (percentages),
we need some absolute frequencies as references.
We restrict to \(3\)-class group \(\mathrm{Cl}_3{F}\) of type \((3,3)\).
Due to the explanations given in Remark
\ref{rmk:3ClassRank}
and Definition
\ref{dfn:Categories},
the \(251\) conductors in the range \(657\le c\le 26\,523\)
of our computations in 2013 consist of \\
\(\bullet\) \(19\) with \(t = 2\) and \(9\mid c\), corresponding to \(38\) fields, and \\
\(69\) with \(t = 2\) and \(\gcd(c,3) = 1\), corresponding to \(138\) fields, \\
since fields with \(t = 2\) arise in doublets, and \\
\(\bullet\) \(85\) with \(t = 3\) and \(9\mid c\), corresponding to \(10\cdot 3 + 7\cdot 2 + 68\cdot 4 = 316\) fields, and \\
\(78\) with \(t = 3\) and \(\gcd(c,3) = 1\), corresponding to \(14\cdot 3 + 11\cdot 2 + 53\cdot 4 = 276\) fields, \\
since fields with \(t = 3\) may arise in triplets (Category I),
doublets (Category II) or quartets (Category III) --- always with respect to type \((3,3)\) alone.

There is a single instance with \textit{dominating two-stage towers}:
For \(t = 2\), \(9\mid c\), 
the extraspecial group \(\langle 27,4\rangle\) is populated most densely with \(20/38=53\%\). 

However, \textit{abelian towers} are most frequent in all other situations:
For \(t = 2\), \(\gcd(c,3) = 1\), 
the single stage \(3\)-class towers are dominating with \(96/138=70\%\). 
For \(t = 3\), \(9\mid c\), 
the single stage \(3\)-class towers are dominating with \(220/316=70\%\). 
For \(t = 3\), \(\gcd(c,3) = 1\), 
the single stage \(3\)-class towers are dominating with \(184/276=67\%\).
Note that the \(3\)-class tower of \(F\) consists of a single stage 
if and only if \(\mathrm{G}_3^2{F}\simeq\langle 9,2\rangle\) is abelian.


This \textit{preliminary legacy statistics} is superseded by our most recent computations in 2022,
presented in the \textit{final statistics} of the Tables
\ref{tbl:FinalStatA}
and
\ref{tbl:FinalStatB}.

\renewcommand{\arraystretch}{1.1}

\begin{table}[ht]
\caption{Statistics of Second \(3\)-Class Group \(\mathfrak{M}\) and Tower Length \(\ell_3{F}\)}
\label{tbl:FinalStatA}
\begin{center}
{\tiny
\begin{tabular}{|r|r|r||c|c|c||c|c|c||c|}
\hline
 Abs.     & Rel.       & Ref.      & \(\mathfrak{M}\)     & cc & \(\ell_3{F}\) & \(t\) & \(m\) & \(v\) & Remark     \\
\hline
 \(4785\) & \(30.2\%\) & \(15851\) & \(\langle 1,1\rangle\)    & \(0\) & \(0\) & \(1\) & \(1\) & \(0\) & Singlets   \\
\hline
 \(7726\) & \(48.7\%\) & \(15851\) &                           &       &       & \(2\) & \(2\) &       & Doublets   \\
\hline
          &            &           &                           &       &       &       &       &       & regular    \\
 \(6910\) & \(\mathbf{89.4\%}\) &  \(7726\) & \(\mathbf{\langle 3,1\rangle}\)    & \(0\) & \(1\) &       &       & \(0\) & \(h_3=3\)  \\
  \(704\) &  \(9.1\%\) &  \(7726\) &                           &       &       &       &       &       & \(h_3=9\)  \\
  \(420\) & \(\mathbf{59.7\%}\) &   \(704\) & \(\mathbf{\langle 9,2\rangle}\)    & \(1\) & \(1\) &       &       & \(1\) & \(h_3=9\)  \\
  \(284\) & \(40.3\%\) &   \(704\) & \(\langle 27,4\rangle\)   & \(1\) & \(2\) &       &       & \(2\) & \(h_3=9\)  \\
  \(112\) &  \(1.4\%\) &  \(7726\) &                           &       &       &       &       &       & \(h_3>9\)  \\
\hline
          &            &           &                           &       &       &       &       &       & singular   \\
   \(74\) & \(\mathbf{66.1\%}\) & \(112\) & \(\mathbf{\langle 81,3\rangle}\) & \(2\) & \(2\) & & & \(3\) & \(h_3=27\) \\
\hline
          &            &           &                           &       &       &       &       &   & super-singular \\
   \(18\) & \(\mathbf{16.1\%}\) & \(112\) & \(\mathbf{\langle 243,14\rangle}\) & \(2\) & \(2\) & & & \(\ge 4\) & \(h_3=27\) \\
    \(5\) &  \(4.5\%\) &   \(112\) & \(\langle 243,13\rangle\) & \(2\) & \(\ge 2\) &   &       & \(4\) & \(h_3=27\) \\
    \(3\) &  \(2.7\%\) &   \(112\) & \(\langle 729,12\rangle\) & \(3\) & \(\ge 2\) &   &       & \(4\) & \(h_3=27\) \\
    \(3\) &  \(2.7\%\) &   \(112\) & \(\langle 729,17|20\rangle\) & \(3\) & \(\ge 2\) & &      & \(4\) & \(h_3=27\) \\
    \(2\) &  \(1.8\%\) &   \(112\) & \(\langle 243,15\rangle\) & \(2\) & \(2\) &       &       & \(4\) & \(h_3=27\) \\
    \(7\) &  \(6.3\%\) &   \(112\) &                           &       &       &   & \(\ge 4\) & \(h_3>27\) \\
\hline
 \(3132\) & \(19.8\%\) & \(15851\) &                           &       &       & \(3\) & \(4\) &       & Quartets   \\
 \(2316\) & \(73.9\%\) &  \(3132\) &                           &       &       &       &       &       & Cat. III   \\
\hline
 \(1820\) & \(\mathbf{78.6\%}\) & \(2316\) & \(\mathbf{\langle 9,2\rangle}\) & \(1\) & \(1\) & & &     & III/1--4   \\
\hline
  \(148\) &  \(6.4\%\) &  \(2316\) &                           &       &       &       &       &       & III/5      \\
   \(82\) & \(\mathbf{55.4\%}\) & \(148\) & \(\mathbf{\langle 243,28|29|30\rangle}\) & \(1\) & \(2\) & & & & III/5  \\
   \(32\) & \(21.6\%\) &   \(148\) & \(\langle 81,7\rangle\)   & \(1\) & \(2\) &       &       &       & III/5      \\
   \(16\) & \(10.8\%\) &   \(148\) & \(\langle 243,27\rangle\) & \(1\) & \(2\) &       &       &       & III/5      \\
   \(10\) &  \(6.8\%\) &   \(148\) & \(\langle 243,25\rangle\) & \(1\) & \(2\) &       &       &       & III/5      \\
    \(8\) &  \(5.4\%\) &   \(148\) & \(\langle 729,37|38|39\rangle\) & \(2\) & \(\ge 2\) & &   &       & III/5      \\
\hline
  \(124\) &  \(5.4\%\) &  \(2316\) &                           &       &       &       &       &       & III/6      \\
  \(120\) & \(\mathbf{96.8\%}\) & \(124\) & \(\mathbf{\langle 81,7\rangle}\) & \(1\) & \(2\) & & &     & III/6      \\
    \(2\) &  \(1.6\%\) &   \(124\) &\(\langle 2187,250\rangle\)& \(2\) & \(\ge 2\) &   &       &       & III/6      \\
    \(2\) &  \(1.6\%\) &   \(124\) &\(\langle 2187,251|252\rangle\)& \(2\) & \(\ge 2\) & &     &       & III/6      \\
\hline
  \(136\) &  \(5.9\%\) &  \(2316\) &                           &       &       &       &       &       & III/7      \\
  \(108\) & \(\mathbf{79.4\%}\) & \(136\) & \(\mathbf{\langle 81,7\rangle}\) & \(1\) & \(2\) & & &     & III/7      \\
   \(12\) &  \(8.8\%\) &   \(136\) & \(\langle 729,41\rangle\) & \(2\) & \(\ge 2\) &   &       &       & III/7      \\
    \(8\) &  \(5.9\%\) &   \(136\) & \(\langle 2187,65|67\rangle\)& \(3\) & \(\ge 3\) & &      &       & III/7      \\
    \(4\) &  \(2.9\%\) &   \(136\) & \(\langle 729,37|38|39\rangle\) & \(2\) & \(\ge 2\) & &   &       & III/7      \\
    \(4\) &  \(2.9\%\) &   \(136\) &\(\langle 2187,253|254\rangle\)& \(2\) & \(\ge 2\) & &     &       & III/7      \\
\hline
   \(28\) &  \(1.2\%\) &  \(2316\) &                           &       &       &       &       &       & III/8      \\
   \(20\) & \(\mathbf{71.4\%}\) & \(28\) & \(\mathbf{\langle 729,34|35|36\rangle}\) & \(2\) & \(\ge 2\) & & & & III/8      \\
\hline
   \(60\) &  \(2.6\%\) &  \(2316\) &                           &       &       &       &       &       & III/9      \\
   \(56\) & \(\mathbf{93.3\%}\) & \(60\) & \(\mathbf{\langle 81,7\rangle}\) & \(1\) & \(2\) & & &      & III/9      \\
\hline
\end{tabular}
}
\end{center}
\end{table}


\renewcommand{\arraystretch}{1.1}

\begin{table}[ht]
\caption{Second \(3\)-Class Group \(\mathfrak{M}\) and Tower Length \(\ell_3{F}\) Continued}
\label{tbl:FinalStatB}
\begin{center}
{\tiny
\begin{tabular}{|r|r|r||c|c|c||c|c|c||c|}
\hline
 Abs.     & Rel.       & Ref.      & \(\mathfrak{M}\)     & cc & \(\ell_3{F}\) & \(t\) & \(m\) & \(v\) & Remark     \\
\hline
  \(392\) & \(12.5\%\) &  \(3132\) &                           &       &       &       &       &       & Cat. I     \\
  \(152\) & \(38.8\%\) &   \(392\) &                           &       &       &       &       &       & I/1        \\
   \(36\) & \(\mathbf{23.7\%}\) & \(152\) & \(\mathbf{\langle 81,10\rangle}\)  & \(1\) & \(2\) & & &   & I/1        \\
   \(21\) & \(13.8\%\) &   \(152\) & \(\langle 243,8\rangle\)  & \(2\) & \(2\) &       &       &       & I/1        \\
   \(18\) & \(11.8\%\) &   \(152\) & \(\langle 81,8\rangle\)   & \(1\) & \(2\) &       &       &       & I/1        \\
   \(18\) & \(11.8\%\) &   \(152\) & \(\langle 81,14\rangle\)  & \(2\) & \(2\) &       &       &       & I/1        \\
   \(18\) & \(11.8\%\) &   \(152\) & \(\langle 243,28|29|30\rangle\) & \(1\) & \(2\) & &       &       & I/1        \\
    \(9\) &  \(5.9\%\) &   \(152\) & \(\langle 243,25\rangle\) & \(1\) & \(2\) &       &       &       & I/1        \\
    \(9\) &  \(5.9\%\) &   \(152\) & \(\langle 729,54\rangle\) & \(2\) & \(3\) &       &       &       & I/1        \\
    \(5\) &  \(3.3\%\) &   \(152\) & \(\langle 243,46\rangle\) & \(3\) & \(2\) &       &       &       & I/1        \\
    \(4\) &  \(2.6\%\) &   \(152\) & \(\langle 243,47\rangle\) & \(3\) & \(2\) &       &       &       & I/1        \\
    \(3\) &  \(2.0\%\) &   \(152\) & \(\langle 2187,303\rangle\) & \(2\) & \(\ge 2\) & &       &       & I/1        \\
\hline 
  \(240\) & \(61.2\%\) &   \(392\) &                           &       &       &       &       &       & I/2        \\
   \(96\) & \(\mathbf{40.0\%}\) & \(240\) & \(\mathbf{\langle 81,10\rangle}\)  & \(1\) & \(2\) & & &   & I/2        \\
   \(48\) & \(20.0\%\) &   \(240\) & \(\langle 81,8\rangle\)   & \(1\) & \(2\) &       &       &       & I/2        \\
   \(48\) & \(20.0\%\) &   \(240\) & \(\langle 81,14\rangle\)  & \(2\) & \(2\) &       &       &       & I/2        \\
   \(30\) & \(12.5\%\) &   \(240\) & \(\langle 243,8\rangle\)  & \(2\) & \(2\) &       &       &       & I/2        \\
    \(3\) &  \(1.3\%\) &   \(240\) & \(\langle 243,42\rangle\) & \(3\) & \(2\) &       &       &       & I/2        \\
    \(3\) &  \(1.3\%\) &   \(240\) & \(\langle 729,52\rangle\) & \(2\) & \(\ge 2\) &   &       &       & I/2        \\
    \(3\) &  \(1.3\%\) &   \(240\) & \(\langle 2187,301|305\rangle\) & \(2\) & \(3\) & &       &       & I/2        \\
    \(1\) &  \(0.4\%\) &   \(240\) & \(\langle 2187,4670\rangle\) & \(5\) & \(2\) &  &       &       & I/2        \\
\hline
  \(368\) & \(11.7\%\) &  \(3132\) &                           &       &       &       &       &       & Cat. II    \\
  \(188\) & \(51.1\%\) &   \(368\) &                           &       &       &       &       &       & II/1       \\
   \(76\) & \(\mathbf{40.4\%}\) & \(188\) & \(\mathbf{\langle 81,7\rangle}\)  & \(1\) & \(2\) & & &    & II/1       \\
   \(76\) & \(\mathbf{40.4\%}\) & \(188\) & \(\mathbf{\langle 81,13\rangle}\) & \(2\) & \(2\) & & &    & II/1       \\
    \(6\) &   \(3.2\%\) &  \(188\) & \(\langle 729,37|38|39\rangle\) & \(2\) & \(\ge 2\) & &   &       & II/1       \\
    \(4\) &   \(2.1\%\) &  \(188\) & \(\langle 729,41\rangle\) & \(2\) & \(3\) &       &       &       & II/1       \\
    \(2\) &   \(1.1\%\) &  \(188\) & \(\langle 729,372\rangle\) & \(3\) & \(\ge 2\) &  &       &       & II/1       \\
    \(2\) &   \(1.1\%\) &  \(188\) & \(\langle 2187,248|249\rangle\) & \(2\) & \(\ge 2\) & &   &       & II/1       \\
    \(2\) &   \(1.1\%\) &  \(188\) & \(\langle 2187,65|67\rangle\) & \(3\) & \(\ge 3\) & &     &       & II/1       \\
    \(2\) &   \(1.1\%\) &  \(188\) & \(\langle 2187,4595\rangle\) & \(4\) & \(\ge 2\) & &      &       & II/1       \\
    \(2\) &   \(1.1\%\) &  \(188\) & \(\langle 2187,4606\rangle\) & \(4\) & \(\ge 2\) & &      &       & II/1       \\
    \(2\) &   \(1.1\%\) &  \(188\) & \(\langle 2187,5577\rangle\) & \(3\) & \(3\) &    &       &       & II/1       \\
\hline
  \(180\) &  \(48.9\%\) &  \(368\) &                           &       &       &       &       &       & II/2       \\
   \(66\) & \(\mathbf{36.7\%}\) & \(180\) & \(\mathbf{\langle 81,7\rangle}\)  & \(1\) & \(2\) & & &    & II/2       \\
   \(66\) & \(\mathbf{36.7\%}\) & \(180\) & \(\mathbf{\langle 81,13\rangle}\) & \(2\) & \(2\) & & &    & II/2       \\
   \(16\) &   \(8.9\%\) &  \(180\) & \(\langle 729,37|38|39\rangle\) & \(2\) & \(\ge 2\) & &   &       & II/2       \\
    \(4\) &   \(2.2\%\) &  \(180\) & \(\langle 729,41\rangle\) & \(2\) & \(3\) &       &       &       & II/2       \\
    \(2\) &   \(1.1\%\) &  \(180\) & \(\langle 2187,253|254\rangle\) & \(2\) & \(\ge 2\) & &   &       & II/2       \\
    \(2\) &   \(1.1\%\) &  \(180\) & \(\langle 2187,4595\rangle\) & \(4\) & \(\ge 2\) & &      &       & II/2       \\
    \(2\) &   \(1.1\%\) &  \(180\) & \(\langle 2187,5577\rangle\) & \(3\) & \(3\) &    &       &       & II/2       \\
\hline
   \(56\) &  \(1.8\%\) &  \(3132\) &                           &       &       &       &       &       & Cat. IV    \\
\hline
  \(208\) &  \(1.3\%\) & \(15851\) &                           &       &       & \(4\) & \(8\) &       & Octets     \\
\hline 
\end{tabular}
}
\end{center}
\end{table}

\newpage

\section{Main results}
\label{s:MainResults}

\noindent
The \(2\)-class tower of cyclic cubic fields \(F\)
with \(2\)-class group \(\mathrm{Cl}_2{F}\simeq (2,2)\)
is completely settled by our Theorems
\ref{thm:CycCub},
\ref{thm:SigmaGroupsDegree3},
\ref{thm:Derhem},
and by the detailed criteria distinguishing the two cases \(\ell_2{F}\in\lbrace 1,2\rbrace\) in
\cite{Dh1988,CtDh1992}.
The second fruitful application of restrictions enforced by Galois action in Theorem
\ref{thm:CycCub5x5}
and Corollary
\ref{cor:CycCub5x5}
justifies the concrete numerical results for cyclic cubic fields \(F\)
with \(5\)-class group \(\mathrm{Cl}_5{F}\simeq (5,5)\) in Table
\ref{tbl:CycCub5x5}.
In particular, Olga Taussky's famous \(1970\) \textit{fixed point capitulation} problem
\cite[\S\ 3.5.2, p. 448]{Ma2013}
which has first been solved by five \textit{imaginary quadratic} fields
\cite{Ma2013},
and later by eleven \textit{cyclic quartic} fields
\cite[Thm. 4.4, Tbl. 4--5]{AKMTT2020},
now also has nine solutions with \textit{cyclic cubic} fields
having the identity permutation type \(\varkappa=(123456)\).
There only remains the open question whether
the abelian group \(\langle 25,2\rangle\) or
descendants of the groups
\(\langle 3125,i\rangle\) with \(i\in\lbrace 3,10\rbrace\)
will be realized as \(\mathrm{G}_5^{(2)}{F}\)
by cyclic cubic fields \(F\) with bigger conductors \(c>10^6\).

Our investigation of the \(3\)-class tower of cyclic cubic fields \(F\)
with \textit{elementary tricyclic} \(3\)-class group \(\mathrm{Cl}_3{F}\simeq (3,3,3)\)
is a striking novelty.
Similar attempts with imaginary quadratic fields of type \((3,3,3)\),
where all capitulation kernels are of order \(\#L=3\) (lines),
successfully yielded the Artin pattern \(\mathrm{AP}=(\tau,\varkappa)\) by means of arithmetic computations
\cite[\S\ 7.2, Tbl. 2--4, pp. 308--311]{Ma2015b}
but were doomed to failure group theoretically,
since the complexity of descendant trees became unmanageable
\cite[\S\ 7.4, p. 312]{Ma2015b},
\cite[\S\ 10, p. 54]{Ma2015c},
\cite[Thm. 8.2, p. 174]{Ma2017},
\cite[\S\ 8, pp. 98--99]{Ma2016c},
\cite[\S\ 2, Example, p. 6]{Ma2016d}.
Therefore, we were delighted that cyclic cubic fields of type \((3,3,3)\)
impose much less severe requirements on the second \(3\)-class group
\(\mathfrak{M}=\mathrm{G}_3^{(2)}{F}\),
since capitulation kernels of order \(\#P=9\) (planes) and even \(\#O=27\) (full space) are admissible.





\end{document}